\documentclass[12pt]{amsart}
\font\emailfont=cmtt10

\usepackage{amsmath,amsthm,amsfonts,amscd,flafter,epsf}

\newcommand\divis{\mathfrak d}

\newcommand\commentable[1]{#1}
\newcommand\Id{\mathrm{Id}}

\newcommand{\Tors}{\mathrm{Tors}}

\newcommand{\HF}{HF}

\newtheorem{theorem}{Theorem}[section]
\newtheorem{prop}[theorem]{Proposition}

\newtheorem{lemma}[theorem]{Lemma}

\newtheorem{example}[theorem]{Example}
\newtheorem{defn}[theorem]{Definition}

\newtheorem{remark}[theorem]{Remark}

\def\endproof{\relax\ifmmode\expandafter\endproofmath\else
  \unskip\nobreak\hfil\penalty50\hskip.75em\hbox{}\nobreak\hfil\bull
  {\parfillskip=0pt \finalhyphendemerits=0 \bigbreak}\fi}
\def\endproofmath$${\eqno\bull$$\bigbreak}
\def\bull{\vbox{\hrule\hbox{\vrule\kern3pt\vbox{\kern6pt}\kern3pt\vrule}\hrule}}

\newcounter{bean}
\newcommand{\Q}{\mathbb{Q}}
\newcommand{\R}{\mathbb{R}}

\newcommand{\C}{\mathbb{C}}

\newcommand{\Z}{\mathbb{Z}}

\newcommand{\OneHalf}{\frac{1}{2}}

\newcommand{\Zmod}[1]{\Z/{#1}\Z}

\newcommand{\Ker}{\mathrm{Ker}}
\newcommand{\CoKer}{\mathrm{Coker}}
\newcommand{\Coker}{\mathrm{Coker}}
\newcommand{\ind}{\mathrm{ind}}
\newcommand{\Image}{\mathrm{Im}}

\newcommand{\grad}{\vec\nabla}

\newcommand{\cm}{\cdot}
\newcommand\Sections{\mbox{$\Gamma$}}

\newcommand{\CDisk}{D}

\newcommand{\ModSWfour}{\mathcal{M}}
\newcommand{\ModFlow}{\ModSWfour}

\newcommand{\Sobol}[2]{L^{#1}_{#2}}

\newcommand{\Maps}{\mathrm{Map}}

\newcommand{\SpinC}{{\mathrm{Spin}}^c}

\newcommand{\Proj}{\Pi}

\newcommand{\goesto}{\mapsto}

\newcommand{\DBar}{\overline{\partial}}

\newcommand\Wedge{\Lambda}

\newcommand\loc{\mathrm{loc}}
\newcommand\ext{\mathrm{ext}}
\newcommand\Pic[1]{\mathrm{Pic}^{#1}}
\newcommand\Hom{\mathrm{Hom}}

\newcommand\abuts\Rightarrow
\newcommand\Sym{\mathrm{Sym}}

\newcommand\End{\mathrm{End}}

\newcommand\TRGras{G{\R}}
\newcommand{\Diag}{D}
\newcommand{\AbelJacobi}{\Theta}
\newcommand{\AlmostCxStruct}{\mathcal J}
\newcommand{\AlmostCx}{\AlmostCxStruct}
\newcommand{\Gl}{\mathrm{Gl}}    
\newcommand{\ModDeg}{\mathcal N}
\newcommand{\UnparModDeg}{{\widehat\ModDeg}}
\newcommand\ModRect{\Mod(\Square)}
\def\sqr#1#2{{\vcenter{\vbox{\hrule height.#2pt
	\hbox{\vrule width.#2pt height#1pt \kern#1pt
	\vrule width.#2pt}
	\hrule height.#2pt}}}}
\def\Square{\mathchoice\sqr64\sqr64\sqr{2.1}3\sqr{1.5}3}
\newcommand\Rect{\Square}

\newcommand\SigOne{\Sigma}
\newcommand\Lu{D_{\widetilde u}}
\newcommand\Lv{D_{\widetilde v}}
\newcommand\Dinf{D_{\infty}}
\newcommand\met{\mathfrak g}

\newcommand{\CInftyLoc}{C^\infty_\loc}

\newcommand{\sj}{\mathfrak j}

\newcommand\interior{\mathrm{int}}
\newcommand\SpinCz{\mathfrak S}

\newcommand\ThreeCurveComp{\Sigma-\alpha_{1}-\ldots-\alpha_{g}-\beta_{1}-\ldots-\beta_{g}-\gamma_1-...-\gamma_g}
\newcommand\SpanA{{\mathrm{Span}}([\alpha_i]_{i=1}^g)}
\newcommand\SpanB{{\mathrm{Span}}([\beta_i]_{i=1}^g)}
\newcommand\SpanC{{\mathrm{Span}}([\gamma_i]_{i=1}^g)}
\newcommand\Filt{\mathcal F}
\newcommand\HFinfty{\HFinf}
\newcommand\CFinfty{\CFinf}

\newcommand\Area{\mathcal A}

\newcommand\x{\mathbf x}
\newcommand\w{\mathbf w}
\newcommand\z{\mathbf z}
\newcommand\p{\mathbf p}
\newcommand\q{\mathbf q}
\newcommand\y{\mathbf y}
\newcommand\Harm{\mathcal H}

\newcommand\sW{\mathcal W}
\newcommand\sX{\mathcal X}

\newcommand\sZ{\mathcal Z}

\newcommand\BigO{O}

\newcommand\ModCent{\ModFlow^{\mathrm{cent}}
\left({\mathbb S}\longrightarrow \Sym^{g-1}(\Sigma_{1})\times 
\Sym^2(\Sigma_{2})\right)}
\newcommand\ModSphere{\ModFlow\left({\mathbb S}\longrightarrow 
\Sym^{g-1}(\Sigma_{1})\times \Sym^2(\Sigma_{2})\right)}
\newcommand\ModSpheres\ModSphere
\newcommand\CF{CF}
\newcommand\cyl{\mathrm{cyl}}
\newcommand\CFa{\widehat{CF}}
\newcommand\CFp{\CFb}
\newcommand\CFm{\CF^-}
\newcommand\CFleq{\CF^{\leq 0}}
\newcommand\HFleq{\HF^{\leq 0}}

\newcommand\CFme{\CFleq}

\newcommand\HFpred{\HFp_{\rm red}}

\newcommand\HFmred{\HFm_{\rm red}}
\newcommand\HFred{\HF_{\rm red}}

\newcommand\HFp{\HFb}
\newcommand\HFpm{HF^{\pm}}
\newcommand\HFm{\HF^-}
\newcommand\CFinf{CF^\infty}
\newcommand\HFinf{HF^\infty}
\newcommand\CFb{CF^+}
\newcommand\HFa{\widehat{HF}}
\newcommand\HFb{HF^+}
\newcommand\gr{\mathrm{gr}}
\newcommand\Mas{\mu}
\newcommand\UnparModSp{\widehat \ModSp}
\newcommand\UnparModFlow\UnparModSp
\newcommand\Mod\ModSp

\newcommand{\cald}{{\mathcal D}}

\newcommand\PD{\mathrm{PD}}

\newcommand\Paths{\mathcal{B}}

\newcommand{\spinc}{\mathfrak s}

\newcommand{\spinct}{\mathfrak t}
\newcommand{\Exp}{\mathrm{Exp}}

\newcommand\Real{\mathrm Re}
\newcommand\brD{F}
\newcommand\brDisk{\brD}

\newcommand\ModMaps{\mathcal M}
\newcommand\ModSp\ModMaps

\newcommand\ProdSig[1]{\Sigma^{\times{#1}}}
\newcommand\ProductForm{\omega_0}
\newcommand\ProdForm{\ProductForm}
\newcommand\Energy{E}

\newcommand\Cinfty{C^{\infty}}
\newcommand\Ta{{\mathbb T}_{\alpha}}
\newcommand\Tb{{\mathbb T}_{\beta}}
\newcommand\Tc{{\mathbb T}_{\gamma}}
\newcommand\Td{{\mathbb T}_{\delta}}

\newcommand\Diff{\mathrm{Diff}}

\newcommand\dbar{\overline\partial}
\newcommand\Map{\mathrm{Map}}
\newcommand\Strip{\mathbb{D}}

\newcommand\alphas{\mbox{\boldmath$\alpha$}}

\newcommand\etas{\mbox{\boldmath$\eta$}}
\newcommand\betas{\mbox{\boldmath$\beta$}}
\newcommand\gammas{\mbox{\boldmath$\gamma$}}
\newcommand\deltas{\mbox{\boldmath$\delta$}}

\newcommand\Dom{\mathcal D}
\newcommand\PerClass[1]{{\PerDom}_{#1}}
\newcommand\PerClasses[1]{{\Pi}_{#1}}
\newcommand\PerDom{\mathcal P}
\newcommand\RenPerDom{\mathcal Q}

\newcommand\csum{*}
\newcommand\CurveComp{\Sigma-\alpha_{1}-\ldots-\alpha_{g}-\beta_{1}-\ldots-\beta_{g}}

\newcommand\Mult{m}

\newcommand\NumDoms{m}

\newcommand\Fstar[1]{F^{*}_{#1}}

\newcommand\Fm[1]{F^{-}_{#1}}
\newcommand\Fp[1]{F^{+}_{#1}}
\newcommand\Fa[1]{{\widehat F}_{#1}}
\newcommand\Finf[1]{F^{\infty}_{#1}}
\newcommand\Fleq[1]{F^{\leq 0}_{#1}}

\newcommand\Ring{\mathbb A}
\newcommand\HFstar{\HF^*}

\newcommand\orient{\mathfrak o}

\commentable{
\title[{Holomorphic disks and three-manifold invariants}]
{Holomorphic disks and topological invariants for closed three-manifolds}
\author[Peter Ozsv{\'a}th]{Peter Ozsv\'ath}
\thanks{PSO was supported by NSF grant number DMS 9971950 and a Sloan 
Research Fellowship}
\address{Department of
Mathematics,  Princeton University, New Jersey 08540 \newline
\indent{\emailfont{petero@math.princeton.edu}}}

\author[Zolt{\'a}n Szab{\'o}]{Zolt{\'a}n Szab{\'o}} 
\thanks{ZSz was supported by NSF grant number DMS 9704359, a Sloan 
Research Fellowship, and a Packard Fellowship}
\address{Department of
Mathematics,  Princeton University, New Jersey 08540 \newline
\indent{\emailfont{szabo@math.princeton.edu}}}}
\date{Oct 2000. Revised: Dec 2002}

\begin{document}

\begin{abstract} 
The aim of this article is to introduce certain topological invariants
for closed, oriented three-manifolds $Y$, equipped with a $\SpinC$
structure.  Given a Heegaard splitting of $Y=U_{0}\cup_{\Sigma}U_{1}$,
these theories are variants of the Lagrangian Floer homology for the
$g$-fold symmetric product of $\Sigma$ relative to certain totally
real subspaces associated to $U_{0}$ and $U_{1}$.
\end{abstract}

\maketitle
\section{Introduction}

Let $Y$ be a connected, closed, oriented three-manifold, equipped with
a $\SpinC$ structure $\spinc$. Our aim in this paper is to define
certain Floer homology groups $\HFa(Y,\spinc)$, $\HFp(Y,\spinc)$,
$\HFm(Y,\spinc)$, $\HFinf(Y,\spinc)$, and $\HFred(Y,\spinc)$ using
Heegaard splittings of $Y$. For calculations and applications of these
invariants, we refer the reader to the sequel, ~\cite{HolDiskTwo}.

Recall that a Heegaard splitting of $Y$ is a decomposition $Y
=U_{0}\cup_{\Sigma}U_{1}$, where $U_{0}$ and $U_{1}$ are handlebodies
joined along their boundary $\Sigma$. The splitting is determined by
specifying a connected, closed, oriented two-manifold $\Sigma$ of
genus $g$ and two collections $\{\alpha_1,...,\alpha_g\}$ and
$\{\beta_1,...,\beta_g\}$ of simple, closed curves in $\Sigma$.

The invariants are defined by studying the $g$-fold symmetric product
of the Riemann surface $\Sigma$, a space which we denote by
$\Sym^g(\Sigma)$: i.e. this is the quotient of the $g$-fold product of
$\Sigma$, which we denote by $\Sigma^{\times g}$, by the action of the
symmetric group on $g$ letters. There is a quotient map
$$\pi\colon \Sigma^{\times g}\longrightarrow
\Sym^g(\Sigma).$$ $\Sym^g(\Sigma)$ is a smooth manifold; in fact,
a complex structure on $\Sigma$ naturally gives rise to a complex
structure on $\Sym^g(\Sigma)$, for which $\pi$ is a holomorphic map.

In~\cite{FloerLag}, Floer considers a homology theory defined for a
symplectic manifold and a pair of Lagrangian submanifolds, whose
generators correspond to intersection points of the Lagrangian
submanifolds (when the Lagrangians are in sufficiently general
position), and whose boundary maps count pseudo-holomorphic disks with
appropriate boundary conditions. We spell out a similar theory, where
the ambient manifold is $\Sym^g(\Sigma)$ and the submanifolds playing
the role of the Lagrangians are tori $\Ta=\alpha_{1}\times \ldots
\times\alpha_{g}$ and $\Tb=\beta_{1}\times\ldots\times\beta_{g}$.
These tori are half-dimensional totally real submanifolds with respect
to any complex structure on the symmetric product induced from a
complex structure on $\Sigma$. These tori are transverse to one
another when all the $\alpha_{i}$ are transverse to the $\beta_{j}$.
To bring $\SpinC$ structures into the picture, we  fix a point $z\in 
\Sigma-\alpha_{1}-\ldots-\alpha_{g}-\beta_{1}-\ldots-\beta_{g}$.
We show in Section~\ref{subsec:SpinCStructures} that the choice of $z$ 
induces a natural map from the intersection points $\Ta\cap 
\Tb$ to the set of $\SpinC$ structures over $Y$.

While the submanifolds $\Ta$ and $\Tb$ in $\Sym^{g}(\Sigma)$ are not
{\em a priori} Lagrangian, we show that certain constructions from
Floer's theory can still be applied, to define a chain complex
$\CFinf(Y,\spinc)$. This complex is freely generated by pairs
consisting of an intersection point of the tori (which represents
the given $\SpinC$ structure) and an integer which keeps track of the
intersection number of the holomorphic disks with the subvariety
$\{z\}\times \Sym^{g-1}(\Sigma)$; and its differential counts
pseudo-holomorphic disks in $\Sym^g(\Sigma)$ satisfying appropriate
boundary conditions. Indeed, a natural filtration on the complex gives
rise to an auxiliary collection of complexes $\CFm(Y,\spinc)$,
$\CFp(Y,\spinc)$, and $\CFa(Y,\spinc)$.  We let $\HFm$, $\HFinf$,
$\HFp$, and $\HFa$ denote the homology groups of the corresponding
complexes.

These homology groups are relative $\Zmod{\divis(\spinc)}$-graded
Abelian groups, where $\divis(\spinc)$ is the integer given by
$$\divis(\spinc)=\gcd_{\xi\in H_{2}(Y;\Z)}\langle
c_{1}(\spinc),\xi\rangle,$$ where $c_1(\spinc)$ denotes the first
Chern class of the $\SpinC$ structure.  In particular, when
$c_1(\spinc)$ is a torsion class (which is guaranteed, for example,
if $b_1(Y)=0$), then the groups are relatively $\Z$-graded.

Moreover, we define actions $$U\colon \HFinf(Y,\spinc)
\longrightarrow \HFinf(Y,\spinc)$$ and $$(H_1(Y,\Z)/\Tors)\otimes
\HFinf(Y,\spinc)\longrightarrow \HFinf(Y,\spinc),$$ which decrease
the relative degree in $\HFinf(Y,\spinc)$ by two and one
respectively. These induce actions on $\HFa$, $\HFp$, and $\HFm$
(although the induced $U$-action on $\HFa$ is trivial), endowing
the homology groups with the structure of a module over
$\Z[U]\otimes_\Z\Wedge^*(H_1(Y;\Z)/\Tors)$. We show
in Section~\ref{sec:DefHF} that the quotient
$\HFp(Y,\spinc)/U^d\HFp(Y,\spinc)$
stabilizes for all sufficiently
large exponent $d$, and we let $\HFred(Y,\spinc)$ denote the group so
obtained.  After defining the groups, we turn to their topological invariance:

\begin{theorem}
\label{intro:Invariance}
The invariants $\HFa(Y,\spinc)$, 
$\HFm(Y,\spinc)$, $\HFinf(Y,\spinc)$, $\HFp(Y,\spinc)$,
and $\HFred(Y,\spinc)$, thought of as modules over
$\Z[U]\otimes_\Z \Wedge^*(H_1(Y;\Z)/\Tors)$, are topological
invariants of $Y$ and $\spinc$, in the sense that they are
independent of the Heegaard splitting, the choice of attaching
circles, the basepoint $z$, and the complex structures used in their
definition.
\end{theorem}

See also Theorem~\ref{thm:Invariance} for a more precise, technical
statement.  The proof of the above theorem consists of many steps, and
indeed, they take up the rest of the present
paper.

In Section~\ref{sec:TopPrelim}, we recall the topological
preliminaries on Heegaard splittings and symmetric products used
throughout the paper. In Section~\ref{sec:Analysis}, we describe the
modifications to the usual Lagrangian set-up which are necessary to
define the totally real Floer homologies for the Heegaard splittings.
In Subsection~\ref{subsec:Transversality}, we address the issue of
smoothness for the moduli spaces of disks.  In
Subsection~\ref{subsec:EnergyBounds}, we prove {\em a priori} energy
estimates for pseudo-holomorphic disks which are essential for proving
compactness results for the moduli spaces. 

With these pieces in place, we define the Floer homology groups in
Section~\ref{sec:DefHF}. We begin with the technically simpler case of
three-manifolds with $b_1(Y)=0$, in Subsection~\ref{subsec:QHS}. We
then turn to the case where $b_1(Y)>0$ in
Section~\ref{subsec:DefHFBOneBig}. In this case, we must work with a
special class of Heegaard diagrams (so-called {\em admissible}
diagrams) to obtain groups which are independent of the isotopy class
of Heegaard diagram. The precise type of Heegaard diagram needed
depends on the $\SpinC$ structure in question, and the variant of
$\HF(Y,\spinc)$ which one wishes to consider.  We define the types of
Heegaard diagrams in Section~\ref{subsec:Admissibility}, and discuss
some of the additional algebraic structure on the homology theories
when $b_1(Y)>0$ in Subsection~\ref{subsec:DefAct}.  With these
definitions in hand, we turn to the construction of admissible
Heegaard diagrams required when $b_1(Y)>0$ in
Section~\ref{sec:Special}.

After defining the groups, we show that they are independent of
initial analytical choices (complex structures) which go into their
definition. This is established in Section~\ref{section:CxStructure},
using chain homotopies which follow familiar constructions in
Lagrangian
Floer homology.  Thus, the groups now depend on the Heegaard diagram.

In Section~\ref{sec:Isotopies}, we turn to the question of their
topological invariance. To show that we have a topological invariant
for three-manifolds, we must show that the groups are invariant under
the three basic Heegaard moves: isotopies of the attaching circles,
handleslides among the attaching circles, and stabilizations of the
Heegaard diagram. Isotopy invariance is established in
Subsection~\ref{subsec:Isotopies}, and its proof is closely modeled on
the invariance of Lagrangian Floer homology under exact Hamiltonian
isotopies.

To establish handleslide invariance, we show that a handleslide
induces a natural chain homotopy between the corresponding chain
complexes. With a view towards this application, we describe in
Section~\ref{sec:HolTriangles} the chain maps induced by counting
holomorphic triangles, which are associated to three $g$-tuples of
attaching circles. Indeed, we start with the four-dimensional
topological preliminaries of this construction in
Subsection~\ref{subsec:TriangleTopPrelim}, and turn to the Floer
homological construction in later subsections. In fact, we set up this
theory in more generality than is needed for handleslide invariance,
to make our job easier in the sequel~\cite{HolDiskTwo}.

With the requisite naturality in hand, we turn to the proof of
handleslide invariance in Section~\ref{sec:HandleSlides}. This starts
with a model calculation in $\#^g(S^1\times S^2)$
(c.f. Subsection~\ref{subsec:STwoTimesSOne}), which we transfer to an
arbitrary three-manifold in Subsection~\ref{subsec:Naturality}.

In Section~\ref{sec:Stabilization}, we prove stabilization
invariance. In the case of $\HFa$, the result is quite
straightforward, while for the others, we must establish certain
gluing results for holomorphic disks.

In Section~\ref{sec:Conclusion} we assemble the various components 
of the proof of Theorem~\ref{intro:Invariance}.

\subsection{On the Floer homology package}

Before delving into the constructions, we pause for a moment to
justify the profusion of Floer homology groups. Suppose for simplicity
that $b_1(Y)=0$.

Given a Heegaard diagram for $Y$, the complex underlying
$\CFinf(Y,\spinc)$ can be thought of as a variant of Lagrangian Floer
homology in $\Sym^g(\Sigma)$ relative to the subsets $\Ta$ and $\Tb$,
and with coefficients in the ring of Laurent polynomials
$\Z[U,U^{-1}]$ to keep track of the homotopy classes of connecting
disks.  This complex in itself is independent of the choice of
basepoint in the Heegaard diagram (and hence gives a homology theory
which is independent of the choice of $\SpinC$ structure on
$Y$). Indeed (especially when $b_1(Y)=0$) the homology groups of this
complex turn out to be uninteresting
(c.f. Section~\ref{HolDiskTwo:sec:HFinfty} of~\cite{HolDiskTwo}).

However, the choice of basepoint $z$ gives rise to a $\Z$-filtration
on $\CFinf(Y,\spinc)$ which respects the action of the polynomial
subalgebra $\Z[U]\subset
\Z[U,U^{-1}]$. Indeed, the filtration has the following form: there
is a $\Z[U]$-subcomplex $\CFm(Y,\spinc)\subset \CFinf(Y,\spinc)$,
and for $k\in\Z$, the $k^{th}$ term in the filtration is given by
$U^{k}\CFm(Y,\spinc)\subset\CFinf(Y,\spinc)$.  It is now the 
chain homotopy type of $\CFinf$ as a {\em filtered} complex which gives 
an interesting three-manifold invariant. To detect this object, we
consider the invariants $\HFm$, $\HFp$, $\HFa$, and $\HFinf$ which
are the homology groups of 
$$
\CFm,
{\mbox{\hskip.4in}}
\frac{\CFinf}{\CFm}, 
{\mbox{\hskip.4in}}
\frac{U^{-1}\cm \CFm(Y,\spinc)}{\CFm(Y,\spinc)},
{\mbox{\hskip.4in}}
{\text{and}}
{\mbox{\hskip.4in}}
\CFinf
$$ respectively.  From their construction, it is clear that there are
relationships between these various homology groups including, in
particular, a long exact sequence relating $\HFm$, $\HFinf$, and
$\HFp$. So, although $\HFinf$ in itself contains no interesting
information, we claim that its subcomplex, quotient complex, and
indeed the connecting maps all do. 

\subsection{Further developments.}
We give more motivation for these invariants, and their relationship
with gauge theory, in the introduction to the
sequel,~\cite{HolDiskTwo}. Indeed, first computations and applications
of these Floer homology groups are given in that paper. See
also~\cite{HolDiskThree} where a corresponding smooth four-manifold invariant
constructed, and~\cite{HolDiskGraded} where we endow the Floer
homology groups with an absolute grading, and give topological
applications of this extra structure.

\subsection{Acknowledgements.}
We would like to thank Stefan Bauer, John Morgan, Tom Mrowka, Rob
Kirby, and Andr{\'a}s Stipsicz for helpful discussions during the
course of the writing of this paper.

\section{Topological preliminaries}
\label{sec:TopPrelim}

In this section, we recall some of the topological ingredients
used in the definitions of the Floer homology theories: Heegaard
diagrams, symmetric products, homotopy classes of connecting disks,
$\SpinC$ structures and their relationship with Heegaard diagrams.

\subsection{Heegaard diagrams}
\label{subsec:HeegaardSplittings}

A genus $g$ Heegaard splitting of a connected, closed, oriented
three-manifold $Y$ is a decomposition of $Y=U_0\cup_\Sigma U_1$ where
$\Sigma$ is an oriented, connected, closed $2$-manifold with genus
$g$, and $U_0$ and $U_1$ are handlebodies with $\partial
U_0=\Sigma=-\partial U_1$.  Every closed, oriented three-manifold
admits a Heegaard decomposition. For modern surveys on the theory of
Heegaard splittings, see~\cite{Scharlemann} and~\cite{Zieschang}.

\begin{figure}
\mbox{\vbox{\epsfbox{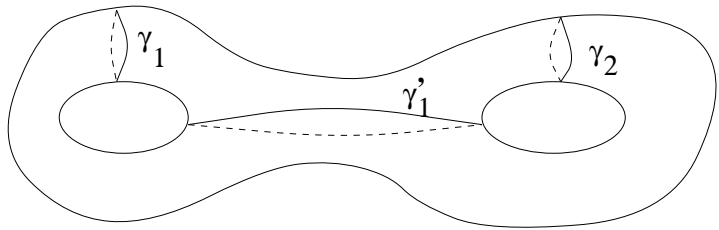}}}
\caption{\label{fig:hs} Handlesliding $\gamma_1$ over $\gamma_2$}
\end{figure}

A handlebody $U$ bounding $\Sigma$ can be described using Kirby 
calculus. $U$ is obtained from $\Sigma$ by first attaching $g$ 
two-handles along $g$ disjoint, simple closed curves 
$\{\gamma_{1},\ldots,\gamma_{g}\}$ which are linearly independent in 
$H_{1}(\Sigma;\Z)$, and then one three-handle. The curves
$\gamma_{1},\ldots,\gamma_{g}$ are called {\em attaching circles 
for $U$}. Since the three-handle is unique, $U$ is determined by the 
attaching circles.
Note that the attaching circles are not uniquely determined by $U$. 
For example, they can be moved by isotopies. But more importantly, 
if $\gamma_{1},\ldots,\gamma_{g}$ are attaching circles for $U$, then 
so are $\gamma_{1}',\gamma_{2},\ldots,\gamma_{g}$, where $\gamma_{1}$ 
is obtained by ``sliding'' the handle of $\gamma_{1}$ over another 
handle, say, 
$\gamma_{2}$; i.e. $\gamma_{1}'$ is any simple, closed 
curve 
which is disjoint from the $\gamma_{1},\ldots,\gamma_{g}$ with the 
property that $\gamma_{1}',\gamma_{1}$ and $\gamma_{2}$ bound an 
embedded 
pair of pants in $\Sigma-\gamma_{3}-\ldots-\gamma_{g}$ (see Figure~\ref{fig:hs} 
for an illustration in the $g=2$ case). 

In view of these remarks, one can concretely think of a genus $g$
Heegaard splitting of a closed three-manifold $Y=U_{0}\cup_{\Sigma}
U_{1}$ as specified by a genus $g$ surface $\Sigma$, and a pair of
$g$-tuples of curves in $\Sigma$, $\alphas=\{\alpha_1,...,\alpha_g\}$
and $\betas=\{\beta_1,...,\beta_g\}$, which are $g$-tuples of
attaching circles for the  $U_0$- and $U_1$-handlebodies
respectively. The triple $(\Sigma,\alphas,\betas)$ is called a {\em
Heegaard diagram}.

Note that Heegaard diagrams have a 
Morse-theoretic interpretation as follows (see for instance~\cite{GompfStipsicz}).
If $f\colon Y\longrightarrow [0,3]$ is a self-indexing 
Morse function on $Y$ with
one minimum and one maximum, then $f$ induces a Heegaard decomposition
with surface $\Sigma=f^{-1}(3/2)$, $U_0=f^{-1}[0,3/2]$,
$U_1=f^{-1}[3/2,3]$. The attaching circles $\alphas$ and $\betas$ are
the intersections of $\Sigma$ with the ascending and descending
manifolds for the index one and two critical points respectively (with
respect to some choice of Riemannian metric over $Y$).
We will call such a Morse function on $Y$ {\em compatible} with the 
Heegaard diagram $(\Sigma,\alphas,\betas)$.

\begin{defn} Let $(\Sigma,\alphas,\betas)$ and $(\Sigma',\alphas',\betas')$
be a pair of Heegaard diagrams. We say that the Heegaard diagrams are
{\em isotopic} if $\Sigma=\Sigma'$ and there are two one-parameter
families $\alphas_t$ and $\betas_t$ of $g$-tuples of curves, moving by
isotopies so that for each $t$, both the $\alphas_t$ and the
$\betas_t$ are $g$-tuples of smoothly embedded, pairwise disjoint
curves. We say that $(\Sigma',\alphas',\betas')$ is obtained from
$(\Sigma,\alphas,\betas)$ by {\em handleslides} if $\Sigma=\Sigma'$
and $\alphas'$ are obtained by handleslides amongst the $\alphas$, and
$\betas'$ is obtained by handleslides amongst the $\betas$. Finally,
we say that $(\Sigma',\alphas',\betas')$ is obtained from
$(\Sigma,\alphas,\betas)$ by {\em stabilization}, if $\Sigma'\cong
\Sigma\#E$, and
$\alphas'=\{\alpha_{1},\ldots,\alpha_{g},\alpha_{g+1}\}$,
$\betas'=\{\beta_{1},\ldots,\beta_{g},\beta_{g+1}\}$, where $E$ is a
two-torus, and $\alpha_{g+1}$, $\beta_{g+1}$ are a pair of curves in
$E$ which meet transversally in a single point. Conversely, in this
case, we say that $(\Sigma,\alphas,\betas)$ is obtained from
$(\Sigma',\alphas',\betas')$ by {\em destabilization}. Collectively,
we will call isotopies, handleslides, stabilizations, and
destabilizations of Heegaard diagrams {\em Heegaard moves}.
\end{defn}

Recall the following basic result (compare~\cite{Reidemeister} and~\cite{Singer}):

\begin{prop}
\label{prop:HeegaardMoves}
Any two Heegaard diagrams $(\Sigma,\alphas,\betas)$ and
$(\Sigma',\alphas',\betas')$ which specify the same three-manifold
are diffeomorphic after a finite sequence of Heegaard moves.
\end{prop}

For the above statement, two Heegaard diagrams
$(\Sigma,\alphas,\betas)$ and $(\Sigma',\alphas',\betas')$ are said
to be diffeomorphic if there is an orientation-preserving
diffeomorphism of $\Sigma$ to $\Sigma'$ which carries $\alphas$ to
$\alphas'$ and $\betas$ to $\betas'$.

Most of Proposition~\ref{prop:HeegaardMoves} follows from the usual
handle calculus (as described, for example, in~\cite{GompfStipsicz}).
Introducing a canceling pair of index one and two critical points
increases the genus of the Heegaard surface by one. After possible
isotopies and handleslides, this corresponds to the stablization
procedure described above.  {\em A priori}, we might have to introduce
canceling pairs of critical points with index zero and one, or two and
three. (The two and three case is dual to the index zero and one case,
so we consider only the latter.) To consider new index zero critical
points, we have to relax the notion of attaching circles: any set
$\{\alpha_1,...\alpha_d\}$ of pairwise disjoint, embedded circles in
$\Sigma$ which bound disjoint, embedded disks in $U$ and span the
image of the boundary homomorphism $\partial\colon
H_2(U,\Sigma;\Z)\longrightarrow H_1(\Sigma,\Z)$ is called an {\em
extended set of attaching circles for $U$} (i.e., here we have $d\geq
g$).  Introducing a canceling zero and one pair corresponds to
preserving $\Sigma$, but introducing a new attaching circle (which
cancels with the index zero critical point). Pair cancellations
correspond to deleting an attaching circle which can be homologically
expressed in terms of the other attaching
circles. Proposition~\ref{prop:HeegaardMoves} is established once we
see that handleslides using these additional attaching circles can be
expressed in terms of handleslides amongst a minimal set of attaching
circles. To this end, we have the following lemmas:

\begin{lemma}
\label{lemma:NewAttachingCircle}
Let $\{\alpha_1,...,\alpha_g\}$ be a set of attaching circles in
$\Sigma$ for $U$. Suppose that $\gamma$ is a simple, closed
curve which is disjoint from $\{\alpha_1,...,\alpha_g\}$. Then,
either $\gamma$ is null-homologous or there is some $\alpha_i$ with
the property that $\gamma$ is isotopic to a curve obtained by
handlesliding $\alpha_i$ across some collection of the $\alpha_j$ for
$j\neq i$.
\end{lemma}

\begin{proof}
If we surger out the $\alpha_1,...,\alpha_g$, we replace $\Sigma$ by
the two-sphere $S^2$, with $2g$ marked points
$\{p_1,q_1,...,p_g,q_g\}$ (i.e. the pair $\{p_i,q_i\}$ corresponds to
the zero-sphere which replaced the circle $\alpha_i$ in $\Sigma$).
Now, $\gamma$ induces a Jordan curve $\gamma'$  in
this two-sphere. If $\gamma'$ does not separate any of the $p_i$ from
the corresponding $q_i$, then it is easy to see that 
the original curve $\gamma$ had to be null-homologous. On the other hand, if
$p_i$ is separated from $q_i$, then it is easy to see that $\gamma$ is
obtained by handlesliding $\alpha_i$ across some collection of the
$\alpha_j$ for $j\neq i$. 
\end{proof}

\begin{lemma}
\label{lemma:AnyTwoSubsets}
Let $\{\alpha_1,...,\alpha_d\}$ be an extended set of attaching circles in
$\Sigma$ for $U$.
Then, any two $g$-tuples of these circles which form a set of
attaching circles for $U$ are related by a series of isotopies and handleslides.
\end{lemma}

\begin{proof}
This is proved by induction on $g$. The case $g=1$ is obvious: if two
embedded curves in the torus represent the same generator in homology,
they are isotopic.

Next, if the two subsets have some element, say $\alpha_1$, in common,
then we can reduce the genus, by surgering out $\alpha_1$. This gives
a new Riemann surface $\Sigma'$ of genus $g-1$ with two marked
points. Each isotopy of a curve in $\Sigma'$ which
crosses one of the marked points corresponds to a handleslide in
$\Sigma$ across $\alpha_1$. Thus, by the inductive hypothesis, the two
subsets are related by isotopies and handleslides.

Consider then the case where the two subsets are
disjoint, labeled $\{\alpha_1,...,\alpha_g\}$ and
$\{\alpha_1',...,\alpha_g'\}$. Obviously, $\alpha_1'$ is not
null-homologous, so, according to
Lemma~\ref{lemma:NewAttachingCircle}, after renumbering, we can obtain
$\alpha_1'$ by handlesliding $\alpha_1$ across some collection of the
$\alpha_i$ ($i=2,...,g$). Thus, we have reduced to the case where the
two subsets are not disjoint.
\end{proof}

\vskip.2cm
\noindent{\bf{Proof of Proposition~\ref{prop:HeegaardMoves}.}}
Given any two Heegaard diagrams of $Y$, we connect corresponding
compatible 
Morse
functions through a generic family $f_t$ of functions, and equip $Y$
with a generic metric. The genericity ensures that the gradient
flow-lines for each of the $f_t$ never flow from higher- to
lower-index critical points. In particular, at all but finitely many
$t$ (where there is cancellation of index one and two critical
points), we get induced Heegaard diagrams for $Y$, whose extended
sets of attaching circles undergo only handleslides and pair
creations and cancellations.

Suppose, now that two sets of attaching circles
$\{\alpha_1,...,\alpha_g\}$ and $\{\alpha_1',...,\alpha_g'\}$ for $U$
can be extended to sets of attaching circles
$\{\alpha_1,...,\alpha_d\}$ and $\{\alpha_1',...,\alpha_d'\}$ for $U$,
which are related by isotopies and handleslides. We claim that the
original sets $\{\alpha_1,...,\alpha_g\}$ and
$\{\alpha_1',...,\alpha_g'\}$ are related by isotopies and
handleslides, as well. To see this, suppose that $\alpha_i'$ (for some
fixed $i\in\{1,...,d\}$) is obtained by handle-sliding $\alpha_i$ over
some $\alpha_j$ (for $j=1,...,d$), then since $\alpha_i'$ can be made
disjoint from all the other $\alpha$-curves, we can view the extended
subset $\{\alpha_1,...,\alpha_d,\alpha_i'\}$ as a set of attaching
circles for $U$. Thus, Lemma~\ref{lemma:AnyTwoSubsets} applies,
proving the claim for a single handleslide amongst the
$\{\alpha_1,...,\alpha_d\}$, and hence also for arbitrary many
handleslides.  The proposition then follows.
\qed
\vskip.2cm

In light of Proposition~\ref{prop:HeegaardMoves}, we see that any
quantity associated to Heegaard diagrams which is unchanged by
isotopies, handleslides, and stabilization is actually a topological
invariant of the underlying three-manifold. Indeed, we will need a
slight refinement of Proposition~\ref{prop:HeegaardMoves}. To this
end, we will find it convenient to fix an additional reference point
$z\in\CurveComp$.

\begin{defn}
The collection data $(\Sigma,\alphas,\betas,z)$ is called a {\em
pointed Heegaard diagram}. Heegaard moves which are supported in a
complement of $z$ -- i.e. during the isotopies, the curves never cross
the basepoint $z$, and for handleslides, the pair of pants does not
contain $z$ -- are called {\em pointed Heegaard moves}.
\end{defn}

\subsection{Symmetric products}
\label{subsec:SymmetricProducts}

In this section, we review the topology of symmetric products. For
more details, see~\cite{MacDonald}.

The {\em diagonal} $\Diag$ in $\Sym^{g}(\Sigma)$ consists of those $g$-tuples 
of points in $\Sigma$, where at least two entries coincide.

\begin{lemma}
    \label{lemma:IdentifyPiOneSym}
Let $\Sigma$ be a genus $g$ surface. Then
$\pi_1(\Sym^g(\Sigma))\cong H_1(\Sym^g(\Sigma))\cong H_1(\Sigma)$.
\end{lemma}

\begin{proof}
We begin by proving the isomorphism on the level of homology.
There is an obvious map 
$$H_1(\Sigma)\rightarrow H_1(\Sym^g(\Sigma))$$
induced from the inclusion $\Sigma\times \{x\}\times ...\times \{x\}\subset
\Sym^g(\Sigma)$. To invert this, note that a curve (in general
position) in $\Sym^g(\Sigma)$ corresponds to a  map of a 
$g$-fold cover of $S^1$ to $\Sigma$,
giving us a homology class in $H_{1}(\Sigma)$. This gives a well-defined map
$H_1(\Sym^g(\Sigma))\longrightarrow H_1(\Sigma)$, since a cobordism
$Z$ in
$\Sym^g(\Sigma)$, which meets the diagonal transversally 
gives rise to a branched cover ${\widetilde Z}$ which maps to
$\Sigma$. It is easy to see that these two maps are inverses of each
other.

To see that $\pi_1(\Sym^{g}(\Sigma))$ is Abelian, consider a
null-homologous curve $\gamma\colon S^1\longrightarrow
\Sym^g(\Sigma)$, which misses the diagonal.
As above, this corresponds to a map ${\widehat \gamma}$
of a $g$-fold cover of the circle  into $\Sigma$,
which is null-homologous; i.e. there is a map of a
two-manifold-with-boundary $F$ into $\Sigma$, $i\colon
F\longrightarrow \Sigma$, with 
$i|\partial F={\widehat \gamma}$. By increasing the genus of $F$ if
necessary, we can extend the $g$-fold covering of the circle to
a branched $g$-fold covering of the disk
$\pi\colon F\longrightarrow \CDisk$. Then, the map sending $z\in
\CDisk$ to the image of $\pi^{-1}(z)$ under $i$ induces the requisite
null-homotopy of $\gamma$. 
\end{proof}

The isomorphism above is Poincar\'e dual to the one
induced from the Abel-Jacobi map $$\AbelJacobi\colon
\Sym^g(\Sigma)\rightarrow\Pic{g}(\Sigma)$$ which
associates to each divisor the corresponding (isomorphism class of)
line bundle. Here, $\Pic{g}(\Sigma)$ is the set of isomorphism 
classes of degree $g$ line bundles over $\Sigma$, which in turn is 
isomorphic to the torus
$$\frac{H^1(\Sigma,\R)}{H^1(\Sigma,\Z)}\cong T^{2g}.$$
Since, $H_1(\Pic{g}(\Sigma))=H^1(\Sigma,\Z)$, we obtain an isomorphism
$$\mu \colon
H_1(\Sigma;\Z)\longrightarrow H^1(\Sym^g(\Sigma);\Z).$$

The cohomology of $\Sym^g(\Sigma)$ was studied in~\cite{MacDonald}. It
is proved there that the cohomology ring is generated by the image of
the above map $\mu$, and one additional two-dimensional cohomology
class, which we denote by $U$, which is Poincar\'e dual to the
submanifold $$\{x\}\times \Sym^{g-1}(\Sigma)\subset \Sym^g(\Sigma),$$
where $x$ is any fixed point in $\Sigma$.

As is implicit in the above discussion, a holomorphic structure $\sj$
on $\Sigma$ naturally endows the symmetric product $\Sym^{g}(\Sigma)$
with a holomorphic structure, denoted $\Sym^g(\sj)$. This structure
$\Sym^g(\sj)$ is specified by the property that the natural quotient
map $$\pi\colon\Sigma^{\times g}\longrightarrow \Sym^{g}(\Sigma)$$ is
holomorphic (where the product space is endowed with a product
holomorphic structure). Indeed, this complex structure can be
K\"ahler, since first any Riemann surface has a projective embedding,
inducing naturally a projective embedding on the $g$-fold product
$\Sigma^{\times g}$, so that elementary geometric invariant theory (as
explained in Chapter 10 of~\cite{Harris}) endows $\Sym^g(\Sigma)$,
its quotient by the symmetric group on $g$ letters (a finite group acting
holomorphically) with the structure of a projective algebraic variety.

As is usual in the study of Gromov invariants and Lagrangian Floer
theory, we must understand the holomorphic spheres in our manifold
$\Sym^g(\Sigma)$. To this end, we study how the first Chern class
$c_1$ (of the tangent bundle $T\Sym^g(\Sigma)$) evaluates on homology
classes which are representable by spheres. First, we identify these
homology classes. To this end, we introduce a little notation. If $X$
is a connected space endowed with a basepoint $x\in X$, let
$\pi_2'(X)$ denote the quotient of $\pi_2(X,x)$ by the action of
$\pi_1(X,x)$. Note that this group is independent of the choice of
basepoint $x$, and also that the natural Hurewicz homomorphism
from $\pi_2(X,x)$ to $H_2(X;\Z)$ factors through $\pi_2'(X)$.

\begin{prop}
\label{prop:PiTwoSym}
Let $\Sigma$ be a Riemann surface of genus $g>1$, then 
$$\pi_2'(\Sym^g(\Sigma))\cong \Z.$$
Furthermore, if $\{A_i,B_i\}$ is a symplectic basis
for $H_1(\Sigma)$, then there is a generator of $\pi_2'(\Sym^g(\Sigma))$, denoted $S$, 
whose image under the Hurewicz homomorphism is Poincar\'e dual to
$$(1-g)U^{g-1} + \sum_{i=1}^g\mu(A_i)\mu(B_i)U^{g-2}.$$
In the case where $g>2$, $\pi_1(\Sym^g(\Sigma))$ acts trivially on
$\pi_2(\Sym^g(\Sigma))$ and thus
$$\pi_2(\Sym^g(\Sigma))\cong \Z.$$
\end{prop}

\begin{proof}   
The isomorphism 
$\pi_2'(\Sym^g(\Sigma))\cong \Z$
is given by the intersection number with the submanifold $x\times
\Sym^{g-1}(\Sigma)$, for generic $x$. Specifically, if we take a
hyperelliptic structure on $\Sigma$, the hyperelliptic involution
gives rise to a sphere
$S_0\subset \Sym^2(\Sigma)$, which we can then use to construct a sphere
$S=S_0\times x_3\times ... \times x_g\subset \Sym^g(\Sigma)$. 
Clearly, $S$ maps to $1$ under this isomorphism.

Consider a sphere $Z$ in the kernel of this map. By moving $Z$ into
general position, we can arrange that $Z$ meets $x\times
\Sym^{g-1}(\Sigma)$ transversally in finitely many points. By splicing
in homotopic translates of $S$ (with appropriate signs) at the
intersection points, we can find a new sphere $Z'$ homotopic to $Z$
which misses $x\times \Sym^{g-1}(\Sigma)$; i.e. we can think of $Z'$
as a sphere in $\Sym^{g}(\Sigma-x)$. Note that since this splicing
construction makes no reference to a basepoint, this operation is
taking place in $\pi_2'(\Sym^g(\Sigma))$.  We claim that
$\pi_2(\Sym^{g}(\Sigma-x))=0$, for $g>2$.

One way to see that $\pi_2(\Sym^{g}(\Sigma-x))=0$ is to observe that 
$\Sigma-x$ is homotopy equivalent to the wedge of $2g$ circles or, 
equivalently, the complement in $\C$ of $2g$ points 
$\{z_{1},\ldots,z_{2g}\}$. Now, 
$\Sym^{g}(\C-\{z_{1},\ldots,z_{2g}\})$ can be thought of as the 
space of monic degree $g$ polynomials $p$ in one variable,
with $p(z_{i})\neq 
0$ for $i=1,\ldots,2g$. Considering the coefficients of $p$, this is 
equivalent to considering $\C^{g}$ minus $2g$ generic hyperplanes.
A theorem of Hattori~\cite{Hattori} states that the 
homology groups of the universal covering space of this complement are 
trivial except in dimension zero or $g$. This proves that $\pi _2(\Sym ^g(\Sigma -x))=0$ and hence completes the proof that
$\pi _2'(\Sym ^g(\Sigma))=\Z$ for $g>2$.

In the case where $g=2$ it is easy to see
that $\Sym ^2(\Sigma)$ is diffeomorphic to the blowup of $T^4$
(indeed, the Abel-Jacobi map gives the map to the torus, and the 
exceptional sphere is the sphere $S_{0}\subset \Sym^{2}(\Sigma)$ 
induced from the hyperelliptic involution on the genus two Riemann 
surface). In this case, the calculation of $\pi_2'$ is straightforward.

To verify the second claim, note that
the Poincar\'e dual of $S$ is characterized by the fact that:
\begin{eqnarray*}
\PD[S]\cup U=\PD[1]&{\text{and}}&
\PD[S]\cup\mu(A_{i})\cup \mu(B_{j})=0,
\end{eqnarray*}
where the latter equation holds for all $i,j=1,\ldots,g$.
It is easy to see that $(1-g)U^{g-1} + \sum_{i=1}^g\mu(A_i)\mu(B_i)U^{g-2}$
satisfies these properties, as claimed.

Finally, in the case where $g>2$, we verify that the action of
$\pi_1(\Sym^g(\Sigma),\{x\times...\times x\})$ is trivial. Fix maps
$$\gamma\colon S^1\longrightarrow \Sym^g(\Sigma)$$ and $$\sigma\colon
S^2\longrightarrow \Sym^g(\Sigma).$$ According to
Lemma~\ref{lemma:IdentifyPiOneSym}, we can arrange after a homotopy
that $\gamma=\gamma_1\times \{x,..,x\}$ where $\gamma_1\colon
S^1\longrightarrow \Sigma$, and according to our calculation of
$\pi_2'$, we can arrange that $\sigma$ has the form
$\sigma=\{x\}\times \sigma_1$, where $\sigma_1\colon
S^2\longrightarrow \Sym^{g-1}(\Sigma)$. Now, the map
$$\gamma\vee\sigma\colon S^1\vee S^2\longrightarrow
\Sym^g(\Sigma)$$ admits an obvious extension $$\gamma_1\times
\sigma_1\colon S^1\times S^2\longrightarrow \Sym^g(\Sigma).$$ Since
the action of $\pi_1(S^1\times S^2)$ on $\pi_2(S^1\times S^2)$ is
trivial, the claim now follows immediately.
\end{proof}

The evaluation of the first Chern class on the generator $S$ is given in the following:

\begin{lemma}
\label{lemma:ChernClass}
The first Chern class of $\Sym^{g}(\Sigma_{g})$ is given by 
$$c_1=U - \sum_{i=1}^g \mu(A_i)\mu(B_i).$$
In particular, $\langle c_1,[S]\rangle = 1.$
\end{lemma}

\begin{proof}
See~\cite{MacDonald} for the calculation of $c_1$. The rest follows from this, together with
Proposition~\ref{prop:PiTwoSym}. 
\end{proof}

\subsection{Totally real tori}

Fix a Heegaard diagram $(\Sigma,\alphas,\betas)$. 
There is a naturally induced 
pair of smoothly embedded, $g$-dimensional tori 
\begin{eqnarray*}
\Ta=\alpha_1\times...\times\alpha_g
&{\text{and}}&
\Tb=\beta_1\times...\times\beta_g
\end{eqnarray*}
in $\Sym^g(\Sigma)$. More precisely $\Ta$ consists of those $g$-tuples
of points $\{x_1,...,x_g\}$ for which $x_i\in \alpha_i$ for
$i=1,...,g$. 

These tori enjoy a certain
compatibility with any complex structure on $\Sym^g(\Sigma)$ induced
(as in Section~\ref{subsec:SymmetricProducts}) from $\Sigma$. 

\begin{defn}
Let $(Z,J)$ be a complex manifold, and $L\subset Z$ be a
submanifold. Then, $L$ is called {\em totally real} if none of its
tangent spaces contains a $J$-complex line, i.e.
$T_\lambda L \cap J T_\lambda L =(0)$
for each $\lambda\in L$.
\end{defn}

\begin{lemma}
Let $\Ta\subset \Sym^g(\Sigma)$ be the torus induced from a 
set of attaching circles $\alphas$. Then, $\Ta$ is a
totally real submanifold of $\Sym^g(\Sigma)$ (for any complex
structure induced from $\Sigma$).
\end{lemma}

\begin{proof}
Note that the projection map
$\pi\colon\Sigma^{\times g}\longrightarrow \Sym^{g}(\Sigma)$
is a holomorphic local diffeomorphism away from the diagonal
subspaces (consisting of those $g$-tuples for which at least two of
the points coincide).
Since $\Ta\subset \Sym^{g}(\Sigma)$ misses the diagonal, the claims 
about $\Ta$ follow immediately from the fact that
$\alpha_{1}\times\ldots\times\alpha_{g}\subset \Sigma^{\times g}$ 
is a 
totally real submanifold (for the product complex structure), 
which follows easily from the definitions.
\end{proof}

Note also that if all the $\alpha_{i}$ curves meet all the 
$\beta_{j}$ curves transversally, then the tori $\Ta$ and $\Tb$ meet 
transversally. We will make these transversality 
assumptions as needed.

\subsection{Intersection points and disks}
\label{subsec:Disks}

Let $\x,\y \in \Ta \cap \Tb$ be a pair of intersection points.
Choose a pair of paths $a\colon [0,1]\longrightarrow \Ta$,
$b\colon [0,1]\longrightarrow \Tb$ from $\x$ to $\y$ in $\Ta$
and $\Tb$ respectively. The difference $a-b$, then, gives a loop in
$\Sym^g(\Sigma)$. 

\begin{defn}
Let $\epsilon(\x,\y)$ denote the image of $a-b$ under the map
$$
\frac{H_1(\Sym^g(\Sigma))}{H_1(\Ta)\oplus H_1(\Tb)}\cong
\frac{H_1(\Sigma)}{[\alpha_1],...,[\alpha_g], [\beta_1],
...,[\beta_g]}
\cong H_1(Y;\Z).
$$ 
Of course, $\epsilon(\x,\y)$ is independent of the choice of the
paths $a$ and $b$.
\end{defn}

It is worth emphasizing that $\epsilon$ can be calculated in $\Sigma$, 
using the identification between $\pi_{1}(\Sym^{g}(\Sigma))$ and 
$H_{1}(\Sigma)$ described in Lemma~\ref{lemma:IdentifyPiOneSym}. 
Specifically, writing $\x=\{x_{1},\ldots,x_{g}\}$ and 
$\y=\{y_{1},\ldots,y_{g}\}$, 
we can think of the path $a\colon [0,1]\longrightarrow \Ta$ as a 
collection of arcs in 
$\alpha_{1}\cup\ldots\cup\alpha_{g}\subset\Sigma$, whose boundary 
(thought of as a zero-chain in $\Sigma$) is given by
$\partial a = y_{1}+\ldots+y_{g}-x_{1}-\ldots-x_{g}$;
similarly, 
we think of the path $b\colon [0,1]\longrightarrow \Tb$ as a 
collection of arcs in 
$\beta_{1}\cup\ldots\cup\beta_{g}\subset\Sigma$, whose boundary is given by
$\partial b = y_{1}+\ldots+y_{g}-x_{1}-\ldots-x_{g}$. Thus, the difference
$a-b$ is a closed one-cycle in $\Sigma$, whose image
in $H_{1}(Y;\Z)$ is the difference $\epsilon(\x,\y)$ defined above.

Clearly $\epsilon$ is additive, in the sense that
$$\epsilon(\x,\y)+\epsilon(\y,\z)=\epsilon(\x,\z),$$ so $\epsilon$
allows us to partition the intersection points of $\Ta\cap \Tb$ into
equivalence classes, where $\x\sim \y$ if $\epsilon(\x,\y)=0$.

We will study holomorphic disks connecting $\x$ and $\y$. These can be
naturally partitioned into homotopy classes of disks with certain boundary
conditions.  To describe this, we consider the unit disk 
$\Strip$ in $\C$, and let $e_1\subset \partial \Strip$
denote the arc where $\Real(z)\geq 0$, and $e_2\subset \partial \Strip$
denote the arc where $\Real(z)\leq 0$. 
When $g>2$, let $\pi_{2}(\x,\y)$ 
denote the set of homotopy classes of maps
$$\left\{u\colon \Strip\longrightarrow
\Sym^g(\Sigma)\Bigg|
\begin{array}{l}
u(-i)=\x, u(i)=\y \\ 
u(e_1)\subset \Ta, u(e_2)\subset \Tb
\end{array}\right\}.$$
In the case where $g=2$, we let $\pi_2(\x,\y)$ denote the quotient
of this set by the additional equivalence relation that $\phi_1\sim \phi_2$
if they become homotopic after splicing sime fixed sphere $\psi\in\pi_2(\Sym^g(\Sigma))$.
In general,
$\pi_2(\x,\y)$ is empty if $\epsilon(\x,\y)\neq 0$.

The set $\pi_{2}(\x,\y)$ is equipped with certain algebraic
structure. Note that there is a natural splicing action
$$\pi_{2}'(\Sym^{g}(\Sigma)) *
\pi_{2}(\x,\y)\longrightarrow
\pi_{2}(\x,\y).$$ 
Also, if we take a Whitney disks connecting $\x$ to $\y$, and one 
connecting $\y$ to $\z$, we can ``splice'' them, to get a Whitney disk 
connecting $\x$ to $\z$. This operation gives rise to
a generalized multiplication
$$\csum\colon \pi_{2}(\x,\y) \times \pi_{2}(\y,\z)\longrightarrow \pi_{2}(\x,\z),$$
which is easily seen to be associative.
As a special case, when $\x=\y$, we see that $\pi_{2}(\x,\x)$ is a 
group. 

\begin{defn}
    \label{def:Additive}
    Let $A$ be a collection of functions
    $\{A_{\x,\y}\colon \pi_{2}(\x,\y)\longrightarrow 
    \Z\}_{\x,\y\in\Ta\cap\Tb}$, satisfying the property that
    $$A_{\x,\y}(\phi)+A_{\y,\z}(\psi)=A_{\x,\z}(\phi\csum\psi),$$
    for each $\phi\in\pi_{2}(\x,\y)$, $\psi\in\pi_{2}(\y,\z)$. Such a 
    collection $A$ is called an {\em additive assignment}.
\end{defn}

For example, for each fixed basepoint $z\in\CurveComp$, the map which
sends a Whitney disk $u$ to the algebraic intersection number
$$n_{z}(u)=\# u^{-1}(\{z\}\times\Sym^{g-1}(\Sigma))$$ descends to
homotopy classes, to give an additive assignment $$n_{z}\colon
\pi_{2}(\x,\y)\longrightarrow \Z.$$ This assignment can be used to
define the domain belonging to a Whitney disk:

\begin{defn} Let
$\Dom_{1},\ldots,\Dom_{\NumDoms}$ denote the closures of the
components of
$\Sigma-\alpha_{1}-\ldots-\alpha_{g}-\beta_{1}-\ldots-\beta_{g}$. Given
a Whitney disk $u\colon \Strip\longrightarrow
\Sym^g(\Sigma)$, the {\em domain associated to $u$} is the formal linear combination of the 
domains $\{\Dom_{i}\}_{i=1}^{\NumDoms}$:
$$\Dom(u)=\sum_{i=1}^{\NumDoms} n_{z_{i}}(u)\Dom_{i},$$
where $z_{i}\in\Dom_{i}$ are points in the interior of $\Dom_{i}$. 
If all the coefficients $n_{z_i}(u)\geq 0$, then we write
$\cald(u)\geq 0$.
\end{defn}

This quantity is obviously independent of the  choice of $z_i$, 
and indeed, $\Dom(u)$ depends only on the homotopy class of $u$. 

\begin{defn}
For a pointed Heegaard diagram $(\Sigma,\alphas,\betas,z)$,
a {\em periodic domain} is a
two-chain $\PerDom{}=\sum_{i=1}^\NumDoms a_i\cald_i$ whose boundary is
a sum of $\alpha$- and $\beta$-curves, and whose $n_z(\PerDom)=0$.
For each $\x\in\Ta\cap\Tb$, a class $\phi\in\pi_{2}(\x,\x)$ with
$n_{z}(\phi)=0$ is called a {\em periodic class}. The set
$\PerClasses{\x}(z)$ of periodic classes is naturally a subgroup of
$\pi_{2}(\x,\x)$. The domain belonging to a periodic class is, of
course, a periodic domain.
\end{defn}

The algebraic topology of the $\pi_{2}(\x,\y)$ is described in the 
following:

\begin{prop}
    \label{prop:WhitneyDisks}
    Suppose $g>1$. 
    For all $\x\in\Ta\cap\Tb$, there is an isomorphism
    $$\pi_{2}(\x,\x)\cong \Z\oplus H^{1}(Y;\Z);$$
    which identifies the subgroup of periodic classes 
    $$\PerClasses{\x}(z)\cong H^{1}(Y;\Z).$$
    In general, for each $\x,\y\in\Ta\cap\Tb$, if 
    $\epsilon(\x,\y)\neq 0$, then $\pi_{2}(\x,\y)$ is empty; 
    otherwise,
    $$\pi_{2}(\x,\y)\cong \Z \oplus H^{1}(Y;\Z)$$
    as principal $\pi_{2}'(\Sym^{g}(\Sigma))\times \PerClasses{\x}(z)$ 
    spaces. 
\end{prop}

For each $\x\in\Ta\cap\Tb$, the above proposition shows that the
natural map which associates to a periodic class in
$\PerClasses{\x}(z)$ its periodic domain is an isomorphism of groups
(when $g>1$).

\begin{proof}
    Suppose that $g>2$.
    The space $\pi_{2}(\x,\x)$ is naturally identified with the 
    fundamental group of the space $\Omega(\Ta,\Tb)$
    of paths in $\Sym^{g}(\Sigma)$ 
    joining $\Ta$ to $\Tb$, based at the constant ($\x$) path. 
    Evaluation maps (at the two endpoints of the paths)
    induce a Serre fibration (with fiber the path-space of 
    $\Sym^{g}(\Sigma)$):
    $$
    \begin{CD}
        \Omega\Sym^{g}(\Sigma)@>>>
        \Omega(\Ta,\Tb)
        @>>>
        \Ta\times \Tb,
    \end{CD}
    $$
    whose associated homotopy long exact sequence gives:
    $$
    \begin{CD}
        0  @>>>\Z\cong\pi_{2}(\Sym^{g}(\Sigma))
        @>>> \pi_{1}(\Omega(\Ta,\Tb))
        @>>> \pi_{1}(\Ta\times \Tb)
        @>>>\pi_{1}(\Sym^{g}(\Sigma)).
    \end{CD}
    $$
    But under the identification $\pi_{1}(\Sym^{g}(\Sigma))\cong 
    H^{1}(\Sigma;\Z)$, the images of $\pi_{1}(\Ta)$ and $\pi_{1}(\Tb)$
    correspond to 
    $H^{1}(U_{0};\Z)$ and $H^{1}(U_{1};\Z)$ respectively. Hence, after 
    comparing with the cohomology long exact sequence for $Y$, we can 
    reinterpret the above as a short exact sequence:
    $$\begin{CD}
    0@>>> \Z @>>> \pi_{2}(\x,\x) @>>> H^{1}(Y;\Z) @>>> 0.
    \end{CD}$$
    The homomorphism $n_{z}\colon \pi_{2}(\x,\x)\longrightarrow \Z$ 
    provides a splitting for the sequence. The proposition in the case where
    $g>2$ follows. The case where $g=2$ follows similarly, only now one must
    divide by the action of $\pi_1(\Sym^g(\Sigma))$.

    In the case where $\x\neq \y$ and $\epsilon(\x,\y)=$, then
    $\pi_2(\x,\y)$ is non-empty, so the above reasoning applies.
\end{proof}

\begin{remark}
        The above result, of course, fails when $g=1$. However, it is 
        still clear that
        $\pi_{2}(\x,\y)\longrightarrow \Z\oplus H^{1}(Y;\Z)$ is  
        injective,
        and that is  the only part of this result which is required
        for the Floer homology constructions described below
        to work. (Note also that the only three-manifolds which admit genus one
        Heegaard diagrams are lens spaces and  $S^{2}\times S^{1}$.)
\end{remark}

\subsection{Periodic domains and surfaces in $Y$}
\label{subsec:PerDom}

Given a periodic domain $\PerDom$, there is a map from a surface-with-boundary 
$\Phi\colon F \longrightarrow \Sigma$ representing
$\PerDom$, in the sense that $\Phi_{*}[F] = \PerDom$ as chains
(where here $[F]$ is a fundamental cycle of $F$).   Typically,
such representatives can be ``inefficient'': $\Phi$ need not be
orientation preserving, so $F$ can be quite complicated.  However, for
chains of the form $\PerDom+\ell[\Sigma]$ with no negative
coefficients, we can choose $F$ in a special manner, according to the
following.

\begin{lemma}
    \label{lemma:RepPerDom}
    Consider a chain $\PerDom+\ell[\Sigma]$ with $\ell$ sufficiently large
    that $n_{z'}(\PerDom+\ell[\Sigma])\geq 0$ for all 
    $z'\in\CurveComp$.
    Then there is an oriented two-manifold with 
    boundary $F$ and a map $\Phi\colon F \longrightarrow \Sigma$ 
    with $\Phi_{*}[F]=\PerDom+\ell [\Sigma]$ with the property that
    $\Phi$ is nowhere orientation reversing and the restriction of 
    $\Phi$ to each boundary component of $F$ is a diffeomorphism onto 
    its image.
\end{lemma}

\begin{proof}
    Write 
    $$
         \PerDom+\ell[\Sigma]= \sum_{i=1}^{\NumDoms}n_{i}\cald_{i},
    $$
    (where, by assumption, $n_{i}\geq 0$). If $\cald$ is the domain 
    $\cald_{i}$, then we let $\Mult(\cald)$ denote the coefficient 
    $n_{i}$. 
    The surface $F$ is 
    constructed as an identification space from
    $$
    X=\coprod_{i=1}^{\NumDoms}\coprod_{j=1}^{n_{i}}\cald_{i}^{(j)},$$
    where $\cald_{i}^{(j)}$ is a diffeomorphic copy of the domain 
    $\cald_{i}$. 
    
    The $\alpha$-curves are divided up by the $\beta$-curves into 
    subsets, which we call $\alpha$-arcs; and similarly, the 
    $\beta$-curves are divided up by the $\alpha$-curves
    into $\beta$-arcs. Each $\alpha$ or 
    $\beta$-arc 
    $c$ is contained in two (not necessarily distinct) domains, 
    $\cald_{1}(c)$ and $\cald_{2}(c)$. We order the domains so that
    $$\Mult(\cald_{1}(c))\leq \Mult(\cald_{2}(c)).$$

    $F$ is obtained from $X$ by the following identifications.  For
    each $\alpha$-arc $a$, if $x\in a$, then for
    $j=1,\ldots,\Mult(\cald_{1}(a))$, we identify 
    $$\Big(x^{(j)}\in\cald^{(j)}_{1}(a)\Big)\sim \Big(x^{(j+\delta_{a})}\in
    \cald^{(j+\delta_a)}_{2}(a)\Big),$$ where
    $\delta_{a}=\Mult(\cald_{2}(a))-\Mult(\cald_{1}(a))$.  Similarly,
    for each $\beta$-arc $b$, if $x\in b$, then for
    $j=1,\ldots,\Mult(\cald_{1}(b))$, we identify
    $$\Big(x^{(j)}\in\cald^{(j)}_{1}(b)\Big)\sim \Big(x^{(j)}\in
    \cald^{(j)}_{2}(b)\Big).$$ The map $\Phi$, then, is induced from the
    natural projection map from $X$ to $\Sigma$.

    It is easy to verify that the space $F$ is actually a manifold-with-boundary as claimed. 
\end{proof}

Let $\Phi\colon F \longrightarrow \Sigma$ be a representative for a 
periodic domain $\PerDom+\ell[\Sigma]$ as constructed in 
Lemma~\ref{lemma:RepPerDom}. $\Phi$ can be extended to a map 
into the three-manifold:
$${\widehat \Phi}\colon {\widehat F}\longrightarrow Y$$
by gluing copies of the attaching disks for the index one and two 
critical points (with appropriate multiplicity) along the boundary of
$F$. This gives us a concrete correspondence between periodic domains
and homology classes in $Y$ which, in the case where $\Ta$ meets
$\Tb$, is Poincar\'e dual to the isomorphism of
Proposition~\ref{prop:WhitneyDisks}.

One can also think of the intersection numbers $n_{z}$ as taking place 
in $Y$. To set this up, note that each 
(oriented) attaching circle
$\alpha_{i}$ naturally gives rise to a cohomology class 
$\alpha_{i}^{*}\in H^{2}(Y;\Z)$. This class is, by definition, 
Poincar\'e dual to the closed curve $\gamma\subset U_0\subset Y$ which is the
difference between the two flow-lines connecting 
the corresponding index one critical point $a_{i}\in U_{0}\subset Y$ with 
the index zero critical point. The sign of $\alpha_i^*$ is specified
by requiring
that the linking number of $\gamma$ with $\alpha_i$ in $U_0$ is  $+1$.  

\begin{lemma}
    \label{lemma:MovePoint} Let $z_{1},z_{2}\in\CurveComp$ be a pair
    of points which are separated by $\alpha_{1}$, in the sense that
    there is a curve $z_t$ from $z_{1}$ to $z_{2}$ which is
    disjoint from $\alpha_{2},\ldots,\alpha_{g}$, and $\#(\alpha_1\cap
    z_t)=+1$. Then, if $\PerDom$ is a periodic domain (with respect to
    some possibly different base-point), then
    $$n_{z_{1}}(\PerDom)-n_{z_{2}}(\PerDom)=\langle H(\PerDom),
    \alpha_{1}^{*}\rangle,$$ where $H(\PerDom)\in H_{2}(Y;\Z)$ is the
    homology class belonging to the periodic domain.
\end{lemma}

\begin{proof}
    For $i=1,2$, let $\gamma_{i}$ be the gradient flow line passing
    through $z_{i}$ (connecting the index zero to the index three
    critical point). Clearly, $n_{z_{i}}(\PerDom)=\# \gamma_{i}\cap
    \PerDom$. Now the difference $\gamma_{1}-\gamma_{2}$ is a closed
    loop in $Y$, which is clearly homologous to a loop in $U_{0}$
    which meets the attaching disk for $\alpha_{1}$ in a single
    transverse point (and is disjoint from the attaching disks for
    $\alpha_i$ for $i\neq 1$). The formula then follows.
\end{proof}
    
\subsection{$\SpinC$ structures}
\label{subsec:SpinCStructures}

Fix a point $z\in
\Sigma-\alpha_1-...-\alpha_g-\beta_1-...-\beta_g$. In this section we
define a natural map $$s_z\colon \Ta\cap \Tb\longrightarrow
\SpinC(Y).$$

To construct this, it is convenient to use Turaev's formulation of
$\SpinC$ structures in terms of homology classes of vector fields
(see~\cite{Turaev}; see also~\cite{KMcontact}). Fix a Riemannian
metric $g$ over a closed, oriented three-manifold $Y$.
Following~\cite{Turaev}, two unit vector fields $v_{1}, v_{2}$ are
said to be {\em homologous} if they are homotopic in the complement of
a three-ball in $Y$ (or, equivalently, in the complement of finitely
many disjoint three-balls in $Y$).  Denote the space of homology
classes of unit vector fields over $Y$ by $\SpinC(Y)$.  Fixing an
ortho-normal trivialization $\tau$ of the tangent bundle $TY$, there
is a natural one-to-one correspondence between vector fields over $Y$
and maps from $Y$ to $S^{2}$, which descends to homology classes
(where we say that two maps $f_{0},f_{1}\colon Y\longrightarrow S^{2}$
are {\em homologous} if they are homotopic in the complement of a
three-ball). Fixing a generator $\mu$ for $H^{2}(S^2;\Z)$, it follows
from elementary obstruction theory that the assignment which
associates to a map from $Y$ to $S^{2}$ the pull-back of $\mu$ induces
an identification between the space of homology classes of maps from
$Y$ to $S^{2}$ and the cohomology group $H^{2}(Y;\Z)$. Hence, we
obtain a one-to-one correspondence, depending on the trivialization
$\tau$: $$\delta^{\tau}\colon \SpinC(Y)\longrightarrow H^{2}(Y;\Z).$$
More canonically, if $v_{1}$ and $v_{2}$ are a pair of nowhere
vanishing vector fields over $Y$, then the difference
$$\delta(v_{1},v_{2})=\delta^{\tau}(v_{1})-\delta^{\tau}(v_{2})\in
H^{2}(Y;\Z)$$ is independent of the trivialization $\tau$, since any
two trivializations $\tau$ and $\tau'$ differ by the action of a map
$g\colon Y\longrightarrow SO(3)$, and, as is elementary to check,
$$\delta^{g\cm \tau}(v)-\delta^{\tau}(v)=g^{*}(w),$$ where $w$ is the
generator of $H^{2}(SO(3);\Z)\cong \Zmod{2}$.  Moreover, since (for
any fixed $v\in\SpinC(Y)$) the map $\delta(v,\cm)$ defines a
one-to-one correspondence between $\SpinC(Y)$ and $H^{2}(Y;\Z)$, and
$\delta(v_{1},v_{2})+\delta(v_{2},v_{3})=\delta(v_{1},v_{3})$, the
space $\SpinC(Y)$ is naturally an affine space for $H^{2}(Y;\Z)$. It
is convenient to write the action additively, so that if $a\in
H^{2}(Y;\Z)$ and $v\in \SpinC(Y)$, then $a+v\in\SpinC(Y)$ is
characterized by the property that $\delta(a+v,v)=a$. Moreover, given
$v_{1},v_{2}\in\SpinC(Y)$, we let $v_{1}-v_{2}$ denote
$\delta(v_{1},v_{2})$.

Thus, one could simply define the space of $\SpinC$ structures over $Y$ to be 
the space of homology classes of vector fields. The 
correspondence with the more traditional definition of $\SpinC$ 
structures is given by associating to the vector $v$ the 
``canonical'' $\SpinC$ structure associated to the reduction of the 
structure group of $TY$ to $SO(2)$ (for this, and other equivalent 
formulations, see~\cite{Turaev}). 

The natural map $s_{z}$ is defined as follows.
Let $f$ be a Morse function on $Y$ compatible with the attaching circles
$\alphas$, $\betas$, see Section~\ref{subsec:HeegaardSplittings}.
Then each $\x\in \Ta\cap \Tb$ determines a $g$-tuple
of trajectories for the 
gradient flow of $f$ connecting the
index one critical
points to index two critical points. Similarly $z$ gives a trajectory 
connecting the index zero critical point with the index three
critical point.
Deleting tubular neighborhoods of these $g+1$ trajectories,
we obtain a subset of $Y$ where the gradient vector field $\grad f$
does not vanish.
Since each trajectory connects critical points of different
parities, the gradient vector field has index 0 on all the 
boundary spheres of the subset, so it
can be extended as a nowhere vanishing vector field over $Y$.
The homology class of the nowhere vanishing vector field obtained in 
this manner (after renormalizing, to make it a unit vector field)
gives the $\SpinC$ structure $s_{z}(\x)$. Clearly
$s_z(\x)$ does not depend on the  choice of the  compatible Morse function $f$
or the extension of the vector field  $\grad f$ to the balls.

Now  we investigate how $s_z(\x)\in \SpinC(Y)$ depends on $\x$ and $z$. 

\begin{lemma}
\label{lemma:VarySpinC}
Let $\x,\y\in \Ta\cap \Tb$. Then we have 
\begin{equation}
        \label{eq:VaryIntPoint}
        s_z(\y)-s_z(\x)=\PD[\epsilon(\x,\y)].
\end{equation}
Furthermore if $z_1, z_2 \in \Sigma
-\alpha_1-...-\alpha_g-\beta_1-...-\beta_g$ can be connected in
$\Sigma$ by an arc $z_t$ from $z_1$ to $z_2$ which is disjoint from
the $\betas$, whose intersection number $\#(\alpha_i\cap z_t)=1$, and
$\#(\alpha_j\cap z_t)$ for $j\neq i$ vanishes, then for all $\x\in \Ta\cap \Tb$ we
have
\begin{equation}
\label{eq:VaryBasePoint}
s_{z_2}(\x)-s_{z_1}(\x)=\alpha_i^*,
\end{equation}
where $\alpha_i^*\in H^2(Y,\Z)$ is Poincar\'e dual to the homology
class in $Y$ induced from a curve $\gamma$ in $\Sigma$ with
$\alpha_i\cdot \gamma=1$, and whose intersection number with all other
$\alpha_j$ for $j\neq i$ vanishes.
\end{lemma}

\begin{proof}
    Given $\x\in\Ta\cap\Tb$, let $\gamma_{\x}$ denote the 
    $g$ trajectories for $\grad f$ connecting the index one
    to the index 
    two critical points which contains the $g$-tuple $\x$; 
    similarly, given $z\in \CurveComp$, let $\gamma_{z}$ denote the 
    corresponding trajectory from the index zero to the index
    three critical point.
    
    Thus, if $\x,\y\in\Ta\cap\Tb$, 
    $\gamma_{\x}-\gamma_{\y}$ is a closed loop in $Y$. 
    A 
    representative for $s_{z}(\x)$ is obtained by modifying the 
    vector field $\grad f$ in a neighborhood of 
    $\gamma_{\x}\cup \gamma_{z}$. It follows then that 
    $s_{z}(\x)-s_{z}(\y)$ can be represented by a cohomology class 
    which is compactly supported in a neighborhood of 
    $\gamma_{\x}-\gamma_{\y}$ (we can use the same vector field to 
    represent both $\SpinC$ structures outside this neighborhood). 
    
    It follows that the difference $s_{z}(\x)-s_{z}(\y)$ is some 
    multiple of the Poincar\'e dual of $\gamma_{\x}-\gamma_{\y}$ (at 
    least if the curve is connected; though the following argument is 
    easily seen to apply in the disconnected case as well). To find 
    out which multiple, we fix a disk $D_{0}$ transverse to 
    $\gamma_\x-\gamma_\y$; to find such a disk take some $x_{i}\in 
    \x$ so that $x_{i}\not\in \y$ (if no such $x_i$ can be found,
    then $\x=\y$, and Equation~\eqref{eq:VaryIntPoint}
    is trivial), and let
    $D_{0}$ be a small neighborhood of $x_{i}$ in $\Sigma$. Our 
    representative $v_{\x}$ of $s_{z}(\x)$ can be chosen to agree with $\grad 
    f$ near $\partial D_{0}$; and the representative $v_{\y}$
    for $s_{z}(\y)$ 
    can be chosen to agree with $\grad f$ over $D_{0}$. 
    With respect to any fixed trivialization of $TY$,
    the two maps from $Y$ to $S^{2}$ corresponding to $v_{\x}$ and 
    $v_{\y}$ agree on $\partial D_{0}$. It makes sense, then, to compare the 
    difference between the degrees
    $\deg_{D_{0}}(v_{\x})$ and $\deg_{D_{0}}(v_{\y})$ (maps from the 
    disk to the sphere, relative 
    to their boundary). Indeed, 
    $$s_{z}(\x)-s_{z}(\y)=
    \left(\deg_{D_{0}}(v_{\x})-\deg_{D_{0}}(v_{\y})\right)\PD(\gamma_{\x}-\gamma_{\y})
    $$
    
    To calculate this difference, take 
    another disk $D_{1}$ with the same boundary as $D_{0}$, so that 
    $D_{0}\cup D_{1}$ bounds a three-ball in $Y$ containing 
    the index one critical point 
    corresponding to $x_{i}$ (and no other critical point); thus we 
    can assume that $v_{\x}\equiv \grad f$ over $D_{1}$. Now, since 
    $v_{\x}$ does not vanish inside this three-ball, we have:
    $$0 = \deg_{D_{0}}(v_{\x})+\deg_{D_{1}}(v_{\x})
    = \deg_{D_{0}}(v_{\x})+\deg_{D_{1}}(\grad f).$$
    Thus,
    $$\deg_{D_{0}}(v_{\x})-\deg_{D_{0}}(v_{\y})
    =-\deg_{D_{1}}(\grad f)-\deg_{D_{0}}(\grad f)
    =1,$$
    since $\grad f$ vanishes with winding number $-1$ around the 
    index $1$ critical points of $f$.
    It follows from this calculation that 
    $v_{\x}-v_{\y}=\PD(\gamma_{\x}-\gamma_{\y})$.
    Letting $a\subset \alpha_{1}\cup \ldots 
    \cup \alpha_{g}$ be a collection of arcs with 
    $\partial a = 
    \y-\x$,
    and $b\subset \beta_{1}\cup \ldots \cup \beta_{g}$ be 
    such a collection with $\partial b = \y-\x$, we know that
    $a-b$ 
    represents $\epsilon(\x,\y)$. On the other hand, if $a_{i}\subset 
    a$ is one of the arcs which connects $x_{i}$ to $y_{i}$, then it 
    is easy to see that $a_{i}$ is homotopic relative to its boundary 
    to the segment in $U_{0}$ formed by joining the two gradient trajectories 
    connecting $x_{i}$ and $y_{i}$ to the index one critical point.
    It follows from this (and the analogous statement in $U_{1}$) that
    $a-b$ is homologous to $\gamma_{\y}-\gamma_{\x}$. 
    Equation~\eqref{eq:VaryIntPoint} follows.
    
    Equation~\eqref{eq:VaryBasePoint} follows from similar
    considerations.  Note first that $s_{z_{1}}(\x)$ agrees with
    $s_{z_{2}}(\y)$ away from $\gamma_{z_{1}}-\gamma_{z_{2}}$.
    Letting now $D_{0}$ be a disk which meets $\gamma_{z_{1}}$
    transversally in a single positive point (and is disjoint from
    $\gamma_{z_{2}}$), and $D_{1}$ be a disk with the same boundary as
    $D_{0}$ so that $D_{0}\cup D_{1}$ contains the index zero critical
    point, we have that
    $$\deg_{D_{0}}(v_{z_{1}})-\deg_{D_{0}}(v_{z_{2}})
    =-\deg_{D_{1}}(v_{z_{1}})-\deg_{D_{0}}(v_{z_{2}})
    =-\deg_{D_{1}}(\grad f)-\deg_{D_{0}}(\grad f)=-1$$ (note now that
    $\grad f$ vanishes with winding number $+1$ around the index zero
    critical point).  It follows that
    $s_{z_{1}}(\x)-s_{z_{2}}(\x)=-\PD(\gamma_{z_{1}}-\gamma_{z_{2}})$.
    Now, $\gamma_{z_1}-\gamma_{z_2}$ is easily seen to be Poincar\'e
    dual to $\alpha_i^*$.
\end{proof}

It is not difficult to generalize the above discussion to give a
one-to-one correspondence between $\SpinC$ structures and homotopy
classes of paths of $\Ta$ to $\Tb$ (having fixed the base point
$z$). This is closely related to Turaev's notion of ``Euler systems''
(see~\cite{Turaev}).

There is a natural involution on the space of $\SpinC$ structures
which carries the homology class of the vector field $v$ to the
homology class of $-v$. We denote this involution by the map
$\spinc\mapsto {\overline \spinc}$. Sometimes, ${\overline \spinc}$
is called the {\em conjugate} $\SpinC$ structure to $\spinc$.

There is also a natural map
$$c_1\colon \SpinC(Y)\longrightarrow H^2(Y;\Z),$$
the first Chern class. This is defined by
$c_1(\spinc)=\spinc-{\overline \spinc}$. Equivalently, if $\spinc$ is
represented by the vector field $v$, then 
$c_1(\spinc)$ is the first Chern class of the orthogonal complement of 
$v$, thought of as an oriented real two-plane (hence complex line)
bundle over $Y$. It is clear that
$c_1({\overline\spinc})=-c_1(\spinc)$.

\section{Analytical Aspects}
\label{sec:Analysis}

Lagrangian Floer homology (see~\cite{FloerLag}) is a homology theory
associated to a pair $L_0$ and $L_1$ of Lagrangian submanifolds in a
symplectic manifold. Its boundary map counts certain
pseudo-holomorphic disks whose boundary is mapped into the union of
$L_0$ and $L_1$.  Our set-up here differs slightly from Floer's: we
are considering a pair of totally real submanifolds, $\Ta$ and $\Tb$,
in the symmetric product. It is the aim of this section to show that
the essential analytical aspects -- the Fredholm theory,
transversality, and compactness -- carry over to this context.  We
then turn our attention to orientations. In the final subsection, we
discuss certain disks, whose boundary lies entirely in either $\Ta$ or
$\Tb$. 

\subsection{Nearly symmetric almost-complex structures}

We will be counting pseudo-holomorphic disks in $\Sym^g(\Sigma)$,
using a restricted class of almost-complex structures over
$\Sym^g(\Sigma)$ (which can be thought of as a suitable elaboration of
the taming condition from symplectic geometry).

Recall that an almost-complex structure $J$ over a symplectic manifold
$(M,\omega)$ is said to {\em tame} $\omega$ if $\omega(\xi,J\xi)>0$
for every non-zero tangent vector $\xi$ to $M$. This is an open
condition on $J$. 

The quotient map $$\pi\colon \ProdSig{g} \longrightarrow
\Sym^g(\Sigma)$$ induces a covering space of
$\Sym^g(\Sigma)-{\Diag}$, where $\Diag\subset \ProdSig{g}$ is the
diagonal, see Subsection~\ref{subsec:SymmetricProducts}. Let $\eta$ be
a K\"ahler form over $\Sigma$, and $\omega_0 =\eta^{\times
g}$. Clearly, $\omega_0$ is invariant under the covering action, so it
induces a K\"ahler form $\pi_*(\omega_0)$ over $\Sym^g(\Sigma)-\Diag$.

\begin{defn}
\label{def:NearlyStandard}
Fix a K\"ahler structure $(\sj,\eta)$ over $\Sigma$, 
a finite collection of points
$$\{z_i\}_{i=1}^m\subset \CurveComp,$$
and an open set $V$ with
$$\left(\{z_i\}_{i=1}^m\times \Sym^{g-1}(\Sigma) \bigcup \Diag\right)\subset V \subset
\Sym^g(\Sigma)$$ 
and 
$${\overline V}\cap (\Ta\cup\Tb)=\emptyset.$$
An almost-complex structure $J$ on
$\Sym^g(\Sigma)$ is called $(\sj,\eta,V)$-{\em nearly symmetric}
if \begin{itemize}
\item $J$ tames $\pi_*(\omega_0)$ over $\Sym^g(\Sigma)-{\overline V}$
\item $J$ agrees with $\Sym^g(\sj)$ over $V$
\end{itemize}
The space of $(\sj,\eta,V)$-nearly symmetric almost-complex
structures will be denoted $\AlmostCx(\sj,\eta,V)$.
\end{defn}

Note that since $\Ta$ and $\Tb$ are Lagrangian with respect to
$\pi_*(\omega_0)$, and $J$ tames $\pi_*(\omega_0)$, the tori $\Ta$ and
$\Tb$ are totally real for $J$.

The space $\AlmostCx(\sj,\eta,V)$ is a subset of the set of all
almost-complex structures, and as such it can be endowed with Banach
space topologies ${\mathcal C}^\ell$ for any $\ell$. In fact,
$\Sym^g(\sj)$ is $(\sj,\eta,V)$-nearly symmetric for any choice of $\eta$
and $V$; and the space $\AlmostCx(\sj,\eta,V)$ is an open neighborhood
of $\Sym^g(\sj)$ in the space of almost-complex structures which agree
with $\Sym^g(\sj)$ over $V$. 

Unless otherwise specified, we choose the points $\{z_i\}_{i=1}^m$
so that there is some
$z_i$ in each connected component of $\CurveComp$.

\subsection{Fredholm Theory}
\label{subsec:FredholmTheory}

We recall the Fredholm theory for pseudo-holomorphic disks, with appropriate 
boundary conditions. For more details, we refer the reader 
to~\cite{FloerUnregularized}, see also
\cite{OhFloer},
\cite{FloerHoferSalamon}, and~\cite{FOOO}.

To set this up we assume that $\Ta$ and $\Tb$ meet transversally,
i.e. that each $\alpha_{i}$ meets each $\beta_{j}$ transversally.

We consider the moduli space of holomorphic strips connecting $\x$
to $\y$, suitably generalized as follows. 
Let $\Strip=[0,1]\times i\R \subset \C$ be the
strip in the complex plane. Fix a path $J_s$ 
of almost-complex structures over $\Sym^g(\Sigma)$. 
Let $\ModFlow_{J_s}(\x,\y)$ be the set of
maps satisfying the following conditions:
$$\ModFlow_{J_s}(\x,\y)=\left\{u\colon \Strip\longrightarrow
\Sym^g(\Sigma)\Bigg|
\begin{array}{ll}
u(\{1\}\times\R)\subset \Ta \\
u(\{0\}\times\R)\subset \Tb \\
\lim_{t\goesto -\infty} u(s+it)=\x \\
\lim_{t\goesto +\infty} u(s+it)=\y \\
\frac{du}{ds}+J(s) \frac{du}{dt}=0 &
\end{array}\right\}. $$
For $\phi\in\pi_2(\x,\y)$, the space $\ModFlow_{J_s}(\phi)$ denotes
the subset consisting of maps as above which represent the given 
homotopy class  (or equivalence class, when $g=2$) $\phi$.
The translation action on $\Strip$ endows this moduli space with an
$\R$ action.  The space of {\em unparameterized $J_s$-holomorphic disks} is
the quotient
$$\UnparModSp_{J_s}(\phi)=\frac{\ModSp_{J_s}(\phi)}{\R}.$$ The word
``disk'' is used, in view of the holomorphic identification of the
strip with the unit disk in the complex plane with two
boundary points removed (and maps in the moduli space
extend across these points, in view of the asymptotic conditions).

We will be considering moduli space $\ModFlow_{J_s}(\x,\y)$, where
$J_s$ is a one-parameter family of nearly symmetric almost-complex structures: i.e. where we
have some fixed $(\sj,\eta,V)$ for which each $J_s$ is 
$(\sj,\eta,V)$-nearly symmetric (see
Definition~\ref{def:NearlyStandard}) for each $s\in[0,1]$.

In the definition of nearly-symmetric almost-complex structure, the
almost-complex structure in a neighborhood of $\Diag$ is fixed to help
prove the required energy bound, c.f. Subsection~\ref{subsec:EnergyBounds}.
Moreover, the complex structure in a neighborhood of the
$\{z_i\}_{i=1}^m\times\Sym^{g-1}(\Sigma)$ is fixed to establish the following:

\begin{lemma}
\label{lemma:NonNegativity}
If $u\in\ModSp_{J_s}(\phi)$ is any $J_s$-holomorphic disk, then 
$\cald(u)\geq 0$.
\end{lemma}

\begin{proof}
In a neighborhood of $\{z_i\}_{i=1}^m\times \Sym^{g-1}(\Sigma)$, we are using an
integrable complex structure, so the disk $u$ must either be contained
in the subvariety (which is excluded by the boundary conditions) or it
must meet it non-negatively. 
\end{proof}

Let $E$ be a vector bundle over $[0,1]\times \R$ equipped with a 
metric and compatible connection, 
$p$, $\delta$ be positive real numbers, 
and $k$ be a non-negative integer. The 
$\delta$-weighted Sobolev space of sections of $E$, written
$\Sobol{p}{k,\delta}([0,1]\times \R,E)$, is the space of sections 
$\sigma$ for which the norm 
$$\|\sigma\|_{\Sobol{p}{k,\delta}(E)}=
\sum_{\ell=0}^{k}
\int_{[0,1]\times \R}|\nabla^{(\ell)} 
\sigma(s+it)|^{p}e^{\delta \tau(t)}ds\wedge dt$$
is finite. Here, $\tau\colon \R\longrightarrow \R$ is a 
smooth function with $\tau(t)=|t|$ provided that $|t|\geq 1$. 

Fix some $p>2$.
Let $\Paths_{\delta}(\x,\y)$ denote the space of maps
$$u\colon [0,1]\times \R\longrightarrow\Sym^{g}(\Sigma)$$
in $\Sobol{p}{1,\loc}$, 
satisfying the boundary conditions
\begin{eqnarray*}
u(\{1\}\times \R)\subset \Ta, &{\text{and}}& 
u(\{0\}\times \R)\subset \Tb,
\end{eqnarray*}
which are asymptotic to $\x$ and $\y$ as $t\goesto -\infty$ and 
$+\infty$, in the following sense.
There is a real number $T>0$ and sections 
\begin{eqnarray*}
    \xi_{-}\in 
    \Sobol{p}{1,\delta}\Big([0,1]\times (-\infty,-T], T_{\x}\Sym^{g}(\Sigma)\Big)
    &{\text{and}}& 
    \xi_{+}\in 
    \Sobol{p}{1,\delta}\Big([0,1]\times [T,\infty), T_{\y}\Sym^{g}(\Sigma)\Big)
\end{eqnarray*}
with the property that 
\begin{eqnarray*}
    u(s+it)= \exp_{\x}(\xi_{-}(s+it))
    &{\text{and}}&
    u(s+it)= \exp_{\y}(\xi_{+}(s+it),
\end{eqnarray*}
for all $t<-T$ and $t>T$ respectively. Here,
$\exp$ denotes the usual exponential map for some Riemannian 
metric on $\Sym^{g}(\Sigma)$. 
Note that $\Paths_{\delta}(\x,\y)$ can be naturally given the 
structure of a Banach manifold, whose tangent space at any 
$u\in\Paths_{\delta}(\x,\y)$ is given by
$$\Sobol{p}{1,\delta}(u)
:= \left\{\xi\in \Sobol{p}{1,\delta}([0,1]\times \R, 
u^{*}(T\Sym^{g}(\Sigma)))\Big| 
\begin{array}{l}
\xi(1,t)\in T_{u(1+it)}(\Ta), \forall t\in\R\\
\xi(0,t)\in T_{u(0+it)}(\Tb), \forall t\in\R
\end{array}\right\}.$$
Moreover, at each $u\in\Paths_{\delta}(\x,\y)$, we denote the
space of sections 
$$\Sobol{p}{\delta}(\Lambda^{0,1}u)
:= \Sobol{p}{\delta}\Big([0,1]\times \R, u^{*}(T\Sym^{g}(\Sigma))\Big)
$$
These Banach spaces fit together to form a bundle ${\mathcal
L}^{p}_{\delta}$ over $\Paths_{\delta}(\x,\y)$.  
At each
$u\in\Paths_{\delta}(\x,\y)$, 
$\DBar_{J_s}u=\frac{d}{ds}+J(s)\frac{d}{dt}$ lies in the space $\Sobol{p}{\delta}(\Lambda^{0,1}(u))$ 
and is  zero exactly
when $u$ is a $J_s$-holomorphic map. (Note that our definition of
$\DBar_{J_s}$ implicitly uses the natural trivialization of the
the bundle $\Lambda^{0,1}$ over $\Strip$, which is why the bundle does
not appear in the definition of $\Sobol{p}{\delta}(\Lambda^{0,1}u)$, but does appear in its notation.)
This assignment fits together over $\Paths_{\delta}(\x,\y)$ to induce
a Fredholm section of ${\mathcal L}^{p}_{\delta}$. The linearization
of this section is denoted $$D_{u}
\colon \Sobol{p}{1,\delta}(u)\longrightarrow
\Sobol{p}{\delta}(\Lambda^{0,1}u),$$
and it is given by the formula
$$D_u(\nu)=\frac{d\nu}{ds} + J(s)\frac{d\nu}{dt} + (\nabla_\nu J(s))\frac{du}{dt}.$$
Since the intersection of $\Ta$ and $\Tb$ is transverse,
this linear map is Fredholm for all sufficiently small non-negative 
$\delta$. Indeed, there is some $\delta_{0}>0$ with the property 
that any map $u\in\ModFlow(\x,\y)$ lies in 
$\Paths_{\delta}(\x,\y)$, for all $0\leq \delta < \delta_{0}$.

The components of $\Paths_{\delta}(\x,\y)$ can be partitioned
according to homotopy classes $\phi\in\pi_{2}(\x,\y)$.  The index of
$D_u$, acting on the unweighted space ($\delta=0$) descends to a
function on $\pi_{2}(\x,\y)$. Indeed, the index is calculated by the
Maslov index $\Mas$ of the map $u$ (see~\cite{FloerMaslov},
\cite{SalamonZehnder}, \cite{SpecFlow},
\cite{Viterbo}).  We conclude the subsection with a result about the
Maslov index which will be of relevance to us later:

\begin{lemma}
\label{lemma:MasClass}
Let $S\in\pi_2'(\Sym^g(\Sigma))$ be the positive generator. Then for
any $\phi\in\pi_2(\x,\y)$, we have that
$$\Mas(\phi+k[S])=\Mas(\phi)+2k.$$
In particular, if $O_\x\in \pi_2(\x,\x)$ denotes the class of the
constant map, then
$$\Mas(O_\x+kS)=2k.$$
\end{lemma}

\begin{proof} It follows from the excision principle for the index that attaching a 
topological sphere $Z$ to a disk changes the Maslov index by $2\langle c_1, [Z]\rangle$
(see \cite{FloerMaslov},  \cite{McDuffSalamon}).
On the other hand for the positive 
generator we have $\langle c_1, [S]\rangle=1$  according
to Lemma~\ref{lemma:ChernClass}.
\end{proof}

\subsection{Transversality}
\label{subsec:Transversality}

Given a Heegaard diagram $(\Sigma,\alphas,\betas)$ 
for which all the $\alpha_{i}$ meet the $\beta_{j}$ transversally, 
the tori $\Ta$ and $\Tb$ meet transversally, so
the holomorphic disks connecting
$\Ta$ with $\Tb$ are naturally endowed with a Fredholm deformation
theory. 

Indeed, the usual arguments from Floer theory
(see~\cite{FloerUnregularized}, \cite{OhFloer} and
\cite{FloerHoferSalamon}) can be modified to prove the following result:

\begin{theorem}
\label{thm:Transversality}
	Fix a Heegaard diagram $(\Sigma,\alphas,\betas)$ with the
	property that each $\alpha_i$ meets $\beta_j$ transversally,
	and fix $(\sj,\eta,V)$ as in
	Definition~\ref{def:NearlyStandard}. Then, for a dense set of
	paths $J_s$ of $(\sj,\eta,V)$-nearly symmetric almost-complex
	structures, the moduli
	spaces $\ModFlow_{J_s}(\x,\y)$ are all smoothly cut out by the
	defining equations.
\end{theorem}

In the above statement, ``dense'' is meant in the $C^{\infty}$
topology on
the path-space of $\AlmostCx(\sj,\eta,V)$.

\vskip.5cm
\noindent{\bf{Proof of Theorem~\ref{thm:Transversality}}.}
This is a modification of the usual proof of transversality, 
see~\cite{FloerUnregularized}, \cite{OhFloer} and
\cite{FloerHoferSalamon}.

Recall (see for instance Theorem~5.1 of~\cite{OhFloer}) that if $u$ is
any non-constant holomorphic disk, then there is a dense set of points
$(s,t)\in [0,1]\times \R$ satisfying the two conditions that
$du_{(s,t)}\neq 0$ and $u(s,t)\cap u(s,\R-\{t\})=\emptyset$. By
restricting to an open neighborhood of the boundary of $\Strip$ (note
that we have assumed that ${\overline V}$ is disjoint from $\Ta$ and
$\Tb$), it follows that we can find such an $(s,t)$ with
$u(s,t)\not\in {\overline V}$. By varying the path $J_s$ in a
neighborhood of $u(s,t)$, the usual arguments show that $u$ is a
smooth point for the parameterized moduli space ${\mathfrak M}$,
consisting of pairs $(J_s,u)$ for which $\DBar_{J_s}u=0$. The result
then follows from the Sard-Smale theorem, applied to the Fredholm
projection from ${\mathfrak M}$ to the space of paths of
nearly-symmetric almost-complex structures.
\qed
\vskip.3cm

Under certain topological hypotheses, one can achieve
transversality by placing the curves $\alphas$ and $\betas$ in general
position, but leaving the almost-complex structure fixed: indeed,
letting $J_s$ be the constant path $\Sym^g(\sj)$. We return to this in
Proposition~\ref{prop:MoveToriTransversality}, after setting up more of the
theory of holomorphic disks in $\Sym^g(\Sigma)$.

\subsection{Energy bounds}
\label{subsec:EnergyBounds}

Let $\Omega$ be a domain in $\C$. 
Recall that the energy of a
map $u\colon \Omega\rightarrow X$  to a Riemannian manifold $(X,g)$ is 
given by $$\Energy(u)=\OneHalf\int_{\Omega}|d u|^2.$$ 

Fix 
$\phi\in \pi_2(\x,\y).$
In order to get the usual compactness results for holomorphic disks 
representing $\phi$, we need an 
{\em a priori} energy bound for any holomorphic representative $u$ for
$\phi$.

Such a bound exists in the symplectic context.  Suppose that
$(X,\omega)$ is a compact symplectic manifold, with a tame
almost-complex structure $J$, then there a constant $C$ for which
$$\Energy(u)\leq C \int_{\Omega}u^*(\omega),$$ for each
$J$-holomorphic map $u$. When the $u$ has Lagrangian boundary
conditions, the integral on the right-hand-side depends only on the
homotopy class of the map. This principle holds in our context as
well, according to the following lemma.

\begin{lemma}
\label{lemma:EnergyBound}
Fix a path $J_s$ in the space of nearly-symmetric almost-complex
structures. Then, for each pair of intersection points
$\x,\y\in\Sym^g(\Sigma)$, and  $\phi\in
\pi_2(\x,\y)$, there is an upper bound on the energy of any holomorphic
representative of $\phi$.
\end{lemma}

\begin{proof}
Given 
$$u\colon (\Strip,\partial \Strip)\longrightarrow
(\Sym^g(\Sigma),\Ta\cup \Tb),$$ we consider the 
lift 
$${\widetilde u}\colon (\brDisk,\partial
\brDisk)\longrightarrow (\Sigma^{\times g},\pi^{-1}(\Ta\cup \Tb))$$
obtained by pulling back the branched covering space $\pi\colon
\Sigma^{\times g}\longrightarrow \Sym^g(\Sigma)$. (That is to say, $F$
is defined to be the covering space of the image of
$u$ away from the diagonal
$\Diag\subset \Sym^g(\Sigma)$,
and in a
neighborhood of $\Diag$, $F$ is defined as a subvariety of
$\Sigma^{\times g}$ -- it is here that we are using the fact that each of the
$J_s$ agree with the standard complex structure near $\Diag$.)

We break the energy integral into two regions:
$$\Energy(u)=\int_{u^{-1}(\Sym^g(\Sigma)-V)}|du|^2 +
\int_{u^{-1}(V)}|du|^2.$$
To estimate the integral on $\Sym^g(\Sigma)-V$, we use the fact that
each $J_s$ tames $\pi_*(\omega_0)$, from which it follows that there is a
constant $C_1$ for which 
\begin{eqnarray}
\Energy(u|_{\Sym^g(\Sigma)-V})&\leq& C_1
\int_{u^{-1}(\Sym^g(\Sigma)-V)}u^*(\pi_*(\omega_0))
=
\frac{C_1}{g!}\int_{{\widetilde u}^{-1}(\ProdSig{g}-{\widetilde V})}{\widetilde
u}^{*}(\omega_0),
\label{ineq:AwayFromDiag}
\end{eqnarray}
where ${\widetilde V}=\pi^{-1}({\overline V})$. 

To estimate the other integrand, choose a K\"ahler form $\omega$ over
$\Sym^g(\Sigma)$.  Over $V$ all the $J_s$ agree with
$\Sym^g(\sj)$, so
$u$ is $\Sym^g(\sj)$-holomorphic in that region,
so there is some constant $C_2$ with
the property that
\begin{equation}
\label{eq:KahlerDominates}
\Energy(u|_V) \leq C_2 \int_{u^{-1}(V)} u^*(\omega)
\end{equation}
(the constant $C_2$ depends on the Riemannian metric used
over $\Sym^g(\Sigma)$ and the choice of K\"ahler form $\omega$). 
Moreover, the right hand side can be calculated using ${\widetilde u}$
according to the
following formula:
\begin{equation}
\label{eq:BranchFormula}
\int_{u^{-1}(V)}u^*(\omega)=\frac{1}{g!}\int_{{\widetilde
u}^{-1}({\widetilde V})}{\widetilde u}^*(\pi^*(\omega)).
\end{equation}

Now, 
fix any two-form $\omega_1$ over $\ProdSig{g}$. 
Then there is a constant $C_3$ with the
following property. Let $${\widetilde u}\colon F
\longrightarrow \ProdSig{g}$$ be any map which is
$\sj^{\times g}$-holomorphic on ${\widetilde u}^{-1}({\widetilde V})$. 
Then
we have the inequality
\begin{equation}
\label{ineq:BoundEnergyDiag}
\int_{{\widetilde u}^{-1}({\widetilde V})}{\widetilde u}^*(\omega_1) \leq C_3
\int_{{\widetilde u}^{-1}({\widetilde V})}
{\widetilde u}^*(\ProductForm).
\end{equation}
This holds for the constant with the property that for each tangent 
vector $\xi$ to $\ProdSig{g}$ and 
$$
\omega_1(\xi,J\xi)\leq C_3\ProdForm(\xi,J\xi),
$$
where $J=\sj^{\times g}$.
Such a constant can be found since $\ProdSig{g}$ is compact and 
$\ProdForm(\cdot ,J \cdot)$ determines a non-degenerate quadratic form on
each tangent space $T\ProdSig{g}$.

Applying Inequality~\eqref{ineq:BoundEnergyDiag} for the form
$\omega_1=\pi^*(\omega)$, and combining with
Inequality~\eqref{ineq:AwayFromDiag}, we find a constant $C_0$ with the property that
\begin{equation}
\label{eq:TopologicalBound}
\Energy(u)\leq C_0\int_{\brDisk}{\widetilde u}^*(\omega_0).
\end{equation}

Moreover, with respect to the symplectic form $\ProductForm$, the preimage
under $\pi$ of $\Ta$ and $\Tb$ are both Lagrangian. This
gives a topological interpretation to the right-hand-side of
Equation~\eqref{eq:TopologicalBound}:
\begin{equation}
\label{eq:CohomInterp}
\int_{\brDisk}{\widetilde u}^*(\ProductForm)=\langle \ProductForm,
[\brDisk,\partial\brDisk]\rangle,
\end{equation}
which makes sense since $\ProductForm$ defines a relative cohomology
class in $H^2(\ProdSig{g},\pi^{-1}(\Ta\cup \Tb))$.  Note that the correspondence
$u\mapsto{\widetilde u}$ induces a right inverse,
up to a multiplicative constant, to the map on
homology
$$\pi_*\colon H_2(\ProdSig{g},\pi^{-1}(\Ta\cup\Tb))\longrightarrow H_2(\Sym^g(\Sigma),\Ta\cup\Tb);$$
thus,
the homology class $[\brDisk,\partial\brDisk]$ depends only on
the relative homology class of $u$, thought of as a class in
$H_2(\Sym^g(\Sigma),\Ta\cup\Tb)$  -- 
in particular, it depends only on the 
equivalence class $\phi\in\pi_{2}(\x,\y)$ of $u$). 

Thus, given a class $\phi\in\pi_2(\x,\y)$,
this gives us an {\em a priori}
bound on the $\ProductForm$-energy of the (branched) lift of any
holomorphic disk $u\in \ModFlow_{J_s}({\phi})$, combining
Inequality~\eqref{eq:KahlerDominates},
Equation~\eqref{eq:BranchFormula}, Inequality~\eqref{eq:TopologicalBound},
and Equation~\eqref{eq:CohomInterp}, we get that
\begin{equation}
\label{eq:EnergyBound}
\Energy(u)\leq C_0 \langle
\ProductForm, [\brDisk,\partial\brDisk]\rangle, 
\end{equation}
(for some constant $C_0$ independent of the class $\phi\in\pi_2(\x,\y)$).
\end{proof}

\subsection{Holomorphic disks in the symmetric product}
  
Suppose that the path $J_s$ is constant, and it is given by
$\Sym^g(\sj)$ for some complex structure $\sj$ over $\Sigma$. Then, the
space of holomorphic disks connecting $\x,\y$ can be given an alternate
description, using only maps between one-dimensional complex manifolds.

\begin{lemma}
    \label{lemma:Correspondence}
    Given any holomorphic disk $u\in\ModFlow(\x,\y)$, there is a 
    branched $g$-fold covering space $p\colon {\widehat 
    \Strip}\longrightarrow \Strip$ and a holomorphic map ${\widehat 
    u}\colon {\widehat \Strip}\longrightarrow \Sigma $, with the 
    property that for each $z\in\Strip$, $u(z)$ is the image under 
    ${\widehat u}$ of the pre-image $p^{-1}(z)$.
\end{lemma}

\begin{proof}
Given a holomorphic map $u\colon \Strip\longrightarrow
\Sym^g(\Sigma)$ which does not lie in the diagonal, 
we can find a branched $g!$-fold cover 
$p\colon {\widetilde \Strip}\longrightarrow \Strip$ pulling back the canonical $g!$-fold
cover $\pi\colon \Sigma^{\times g}\longrightarrow \Sym^g(\Sigma)$,
i.e. making the following diagram commutative:
$$\begin{CD}
{\widetilde \Strip}@>{\widetilde u}>> \Sigma^{\times g} \\
@V{p}VV @V{\pi}VV \\
\Strip @>{u}>> \Sym^g(\Sigma).
\end{CD}
$$
Indeed, ${\widetilde \Strip}$ inherits an action by the symmetric
group on $g$ letters $S_g$, and ${\widetilde u}$ is equivariant for
the action (and its quotient is $u$). Let $\Pi_1\colon \Sigma^{\times
g}\longrightarrow \Sigma$ denote projection onto the first
factor. Then, the composite map 
$$\Pi_1\circ {\widetilde u}\colon {\widetilde \Strip}\longrightarrow
\Sigma$$ is invariant under the action of $S_{g-1}\subset S_g$
consisting of permutations which fix the first letter. Then, we let
${\widehat \Strip}={\widetilde \Strip}/S_{g-1}$, and ${\widehat u}$ be
the induced map from ${\widehat \Strip}$ to $\Sigma$. It is easy to 
verify that ${\widehat u}$ has the desired properties.
\end{proof}

\begin{remark}
    It is straightforward to find appropriate topological conditions
    on ${\widehat u}|{\partial{\widehat \Strip}}$ to give a one-to-one 
    correspondence between flows in $\ModFlow(\x,\y)$ and certain 
    pairs $(p\colon {\widehat \Strip}\longrightarrow \Strip, {\widehat 
    u}\colon {\widehat \Strip}\longrightarrow \Sigma)$.
\end{remark}
    
Let $\cald _1,...,\cald _m$ denote the connected components of $\Sigma
-\alpha _1,-...-\alpha _g-\beta_1,...-\beta _g$. Fix a basepoint $z_i$
inside each $\cald _i$. Then for any $\phi \in \pi_2(\x,\y)$ we define
the {\it domain} associated to $\phi$, as a formal linear combination
of components:
$$\cald(\phi)=\sum _{i=1} ^m n_{z_i}(\phi)\cdot \cald _i.$$ 
Similarly the {\it area} of $\phi$ is given by
$${\mathcal A}(\phi)= \sum _{i=1} ^m n_{z_i}(\phi)\cdot {\rm
Area}_{\eta}({\cald _i}),$$
where $\eta$ is the K\"ahler form on $\Sigma$.
This area gives us a concrete way to understand the energy bound from
the 
previous section since, as is easy to verify,
$$\int_{F}{\widetilde u}^*(\omega_0)=(g!){\mathcal A}(\phi).$$

As an application of Lemma~\ref{lemma:Correspondence}, we observe that
for certain special homotopy classes of maps in $\pi_2(\x,\y)$
transversality can also be achieved by moving the curves $\alphas$ and
$\betas$, following the approach of Oh~\cite{OhTransversality}. (This
observation will prove helpful in the explicit calculations of
Section~\ref{sec:HandleSlides} and also Section~\ref{HolDiskTwo:sec:Examples}
of~\cite{HolDiskTwo}.)

To state it, we need the following:

\begin{defn}
A domain $\cald(\phi)$ is called {\em $\alpha$-injective} if all of its
multiplicities are $0$ or $1$, if its interior (i.e. the interior of
the region with multiplicity $1$) is
disjoint from each $\alpha_i$ for $i=1,...,g$, and its boundary
contains intervals in each $\alpha_i$. 
\end{defn}

\begin{prop}
\label{prop:MoveToriTransversality}
Let $\phi\in\pi_2(\x,\y)$ be an $\alpha$-injective homotopy class, and
fix a complex structure $\sj$ over $\Sigma$. Then, for generic
perturbations of the $\alphas$, the moduli space $\ModFlow(\phi)$ of
$\Sym^g(\sj)$-holomorphic disks is smoothly cut out by its defining
equation.
\end{prop}

\begin{proof}
The hypotheses
ensure that for all  $t\in \R$, we have that
$u(1+it)=\{a_1,...,a_g\}\in\Ta$ where $a_i\not\in u(1+it')$ for any
$t'\neq t$. This is true because the $\alpha$-injectivity hypothesis
ensures that the corresponding map ${\widehat u}\colon
\brDisk\longrightarrow \Sigma$, coming from
Lemma~\ref{lemma:Correspondence}, is injective (with injective
linearization, by elementary complex analysis) on the region mapping
to the $\alpha$-curves $p^{-1}(\{1\}\times \R)$. Thus,
following~\cite{OhTransversality}, by varying the $\alpha_i$ in a
neighborhood of the $a_i$, one can see that the map $u$ is a smooth
point in a parameterized moduli space (parameterized now by variations
in the curves). Thus, according to the Sard-Smale theorem, for generic
small variations in the $\alphas$, the corresponding moduli spaces are
smooth.
\end{proof}

\subsection{Orientability}
\label{subsec:Orientability}

In this subsection, we show that the moduli spaces of flows 
$\ModFlow(\x,\y)$ are orientable. As is usual 
in the gauge-theoretic set-up, 
this is done by proving triviality of the determinant line bundle of the 
linearization of the equations (the $\DBar$-equation) 
which cut out the moduli spaces. A thorough treatment of orientability
in general can be found in~\cite{FOOO}.

For some fixed $p>2$ and some real $\delta>0$ (both of which we 
suppress from the notation),
consider the space ${\mathcal B}(\x,\y)={\mathcal B}_{\delta}(\x,\y)$
of  maps discussed in Subsection~\ref{subsec:FredholmTheory}.
The moduli spaces of holomorphic disks 
are finite-dimensional subspaces of this Banach 
manifold.

Recall that for a family $F_{x}$ of Fredholm operators parameterized by an auxiliary space 
$X$, the virtual vector spaces
$\ker F_{x}-\CoKer F_{x}$ naturally fit together
to give rise to an element in the $K$-theory of $X$ 
(see~\cite{KTheory}), the {\em virtual index bundle}. 
The determinant of this is a real line bundle over $X$, the {\em 
determinant line bundle} of the family $F_{x}$.

\begin{prop}
    \label{prop:Orientable} There is a trivial line bundle over
    ${\mathcal B}(\x,\y)$ whose restriction to the moduli space
    $\ModFlow_{J_s}(\x,\y)\subset {\mathcal B}(\x,\y)$ is naturally
    identified with the determinant line bundle for the linearization
    $\det(D_{u})$, where $J_s$ is any path of $(\sj,\eta,V)$-nearly
    symmetric almost-complex structures.
\end{prop}

As we shall see, the main ingredient in the above proposition is the fact 
that the totally real subspaces $\Ta$ and $\Tb$ have trivial tangent 
bundles. We shall give the proof after a
preliminary discussion.

Let  $L_{0}(t)$ and $L_{1}(t)$ be a pair of paths of totally real 
subspaces of $\C^{n}$, indexed by $t\in \R$ which are asymptotically 
constant as $t\goesto\pm \infty$, i.e. there are totally real subspaces 
$L_{0}^{-}$, $L_{1}^{-}$, $L_{0}^{+}$, and $L_{1}^{+}$ with the 
property that
\begin{eqnarray*}
    \lim_{t\goesto \pm \infty}L_{0}(t)=L_{0}^{\pm} &&
    \lim_{t\goesto \pm \infty}L_{1}(t)=L_{1}^{\pm}.
\end{eqnarray*}
Suppose moreover that $L_{0}^{-}$ and $L_{1}^{-}$ are transverse, and 
similarly $L_{0}^{+}$ and $L_{1}^{+}$ are transverse, too. 
Then, the $\DBar$ on $\C^{g}$-valued functions on 
the strip, satisfying  boundary conditions specified by 
the paths $L_{0}(t)$ and $L_{1}(t)$ 
$$\DBar\colon 
\left\{f \in \Sobol{p}{1}([0,1]\times \R; \C^{g})
\Bigg|
\begin{array}{ll}
    f(1+it)\in L_{0}(t), \\
    f(0+it)\in L_{1}(t), \\
    \DBar f = 0 
\end{array}\right\} 
\longrightarrow \Sobol{p}{}([0,1]\times \R;\C^{g})$$
is Fredholm. 
Thus, the $\DBar$ operator induces a family of Fredholm operators 
indexed by the space
$${\mathcal P}=
\left\{L_{0},L_{1}\colon [0,1] \longrightarrow \TRGras(g)
\Bigg| \begin{array}{l}
L_{0}(0)=L_{0}^{-}, \\
L_{1}(0)=L_{1}^{-}, \\
L_{0}(1)=L_0^+, \\
L_1(1)=L_1^+
\end{array}
\right\}
$$
(after reparameterizing the paths in ${\mathcal P}$ to be indexed by 
$\R\cup\{\pm\infty\}$ rather than $[0,1]$),
where $\TRGras(g)$ denotes the Grassmannian of totally real 
$g$-dimensional subspaces of $\C^{g}$. 

The index of the linearization $D_{u}$ of the $\DBar_{J_s}$ operator on maps
of the disk into $\Sym^{g}(\Sigma)$ can be related to index of the
$\DBar$-operators over ${\mathcal P}$, as follows.  First observe that
$D_u$ depends on a path $J_s$ of almost-complex structures.  
However, we can
connect the family to the constant path $\Sym^g(\sj)$, 
without changing
the index bundle. Next, fix a contraction of the unit disk to
$-i$.
Together with a connection over $T\Sym^g(\Sigma)$,
this induces a trivialization for any $u\in
{\mathcal B}(\x,\y)$ of the pull-back of the complex tangent bundle of
$\Sym^{g}(\Sigma)$ (induced from $\Sym^g(\sj)$). Via
these trivializations, the one-parameter family of totally real
subspaces
\begin{eqnarray*}
    \left\{t\mapsto T_{u(1+it)}\Ta\subset 
    T_{u(1+it)}\Sym^{g}(\Sigma)\right\}, &&
    \{t\mapsto T_{u(0+it)}\Tb\subset T_{u(0+it)}\Sym^{g}(\Sigma)\}
\end{eqnarray*}
induce one-parameter families $L_{0}(t)$ and $L_{1}(t)$ of totally
real subspaces of $\C^{g}$. Indeed, if we use a connection over
$T\Sym^g(\Sigma)$ which trivial along $\Tb$,
and we choose the
contraction of our disk to preserve the left arc, both $t=0$ and $t=1$
endpoints of the families can be viewed as a fixed (i.e. independent of
the particular choice of $u$). Thus, we have a map
$$\Psi\colon {\mathcal B}(\x,\y)\longrightarrow
{\mathcal P},$$ together with an identification between the pull-back
of the (virtual) index bundle for $\DBar$ and the (virtual) index
bundle for $D_{u}$ (over the moduli space $\ModFlow(\x,\y)\subset
{\mathcal B}(\x,\y)$).

We wish to study the index bundle over ${\mathcal P}$. 
There is a ``difference'' map
$$\delta\colon \TRGras(g)\times \TRGras(g) \longrightarrow 
\frac{\Gl_{g}(\C)}{\Gl_{g}(\R)},$$
where $\delta(L_{0},L_{1})$ is the equivalence class of any matrix 
$A\in\Gl_{g}(\C)$ with the property that $A L_{0}=L_{1}$.  (By 
taking the difference with $\R^{g}\subset \C^{g}$, we obtain a 
diffeomorphism between $\TRGras(g)$ and the homogeneous space 
$\frac{\Gl_{g}(\C)}{\Gl_{g}(\R)}$.) In this space, we have a
Maslov cycle 
$$Z_{\mu}\subset \frac{\Gl_{g}(\C)}{\Gl_{g}(\R)}= \{[A]\big| 
\R^{g}+A(\R^{g})\neq \C^{n}\}.$$
Of course, $L_{0}$ and $L_{1}$ meet transversally if and only if their 
difference $\delta(L_{0},L_{1})$ does not lie in the Maslov cycle.

Let $[a_{0}]=\delta(L_{0}^{-}, L_{1}^{-})$, $[a_1]=\delta(L_0^+,L_1^+)$.
The difference map gives us a map 
$$\Phi\colon {\mathcal P} \longrightarrow {\mathcal Q} = \{A \colon [0,1]\longrightarrow 
\frac{\Gl_{g}(\C)}{\Gl_{g}(\R)} \big| A(0)=[a_{0}],  A(1)=[a_1]\}.$$
In this notation,
then, the numerical index of the $\DBar$ operator associated to a pair
of paths $L_{0}(t)$ and $L_{1}(t)$ in ${\mathcal P}$ is calculated by
the intersection number of the difference with the Maslov cycle:
$$\ind(\DBar(L_{0}(t),L_{1}(t)))=\delta(L_{0}(t),L_{1}(t))\cap
Z_{\mu}.$$ Moreover, we could work entirely over ${\mathcal Q}$:
${\mathcal Q}$ is identified with the subspace of ${\mathcal P}$ where
$L_{0}(t)\equiv
\R^{g}$, so there is an index bundle over ${\mathcal Q}$, and the 
index bundle for $\DBar$ over ${\mathcal P}$ is easily seen to be the 
pull-back of this index bundle over ${\mathcal Q}$ (since the index 
bundle over ${\mathcal P}$ is trivial over the 
fiber of $\Phi$). 

Clearly, 
$$\pi_{2}\left(\frac{\Gl_{g}(\C)}{\Gl_{g}(\R)}\right)\cong
\pi_{1}({\mathcal Q}).$$
Now, if $g\geq 2$, it is easy to see that 
$$\pi_{2}\left(\frac{\Gl_{g}(\C)}{\Gl_{g}(\R)}\right)\cong\Zmod{2}.$$
Thus, there is no {\em a priori} reason for the determinant line bundle 
$\det(\DBar)\longrightarrow {\mathcal Q}$ to be trivial: its first 
Stiefel-Whitney class may evaluate non-trivially on the non-trivial 
homotopy class $\Zmod{2}$.  Proposition~\ref{prop:Orientable}
is established by giving a suitable lift of the composite map 
$\Phi\circ\Psi$.

\vskip0.3cm
\noindent{\bf{Proof of Proposition~\ref{prop:Orientable}.}}
Continuing the above notation, 
fix matrices $a_{0}, a_1 \in\Gl_{g}(\C)$, and
consider the space
$${\widetilde {\mathcal Q}}=\{A\colon [0,1]\longrightarrow \Gl_{g}(\C)\big|
A(0)=a_{0}, A(1)=a_1 \}.$$
Since $\pi_{2}(\Gl_{g}(\C))=0$, we 
see that the index bundle of the $\DBar$ operator over ${\widetilde Q}$ is 
orientable. Thus, to establish orientability of the determinant line 
bundle over the moduli spaces of flows, we lift 
$\Phi\circ \Psi$ to a map
$${\widetilde \Phi}\colon 
{\mathcal B}(\x,\y)\longrightarrow {\widetilde Q}.$$
To define this lift, note that the tangent spaces to $\Ta$ and $\Tb$ 
respectively can be 
trivialized by ordering and orienting the attaching circles $\alphas$ 
and $\betas$. This in turn gives rise to a complex trivialization of 
the restrictions of $T\Sym^{g}(\Sigma)$ to $\Ta$ and $\Tb$ respectively 
(induced from the identifications
$T \Sym^{g}(\Sigma)|_{\Ta}\cong \Ta\otimes_{\R}\C$,
$T \Sym^{g}(\Sigma)|_{\Tb}\cong \Tb\otimes_{\R}\C$ arising from the 
corresponding totally 
real structures).
Given a holomorphic disk $u$, then, we let $A(t)$ 
denote the matrix corresponding to the linear transformation from 
$\C^{g}$ to itself given by parallel transporting the vector space
$\C^{g}\cong T_{u(1+it)}\Ta\otimes_{\R}\C$ to 
$T_{u(1)}\Ta\otimes_{\R}\C\cong \C^{g}$, 
using the arc which is the image under $u$
of the path prescribed by the fixed contraction of
$\Strip$.
Now, the composite $\Phi\circ\Psi$ factors through
the projection from ${\widetilde {\mathcal Q}}$ 
to ${\mathcal Q}$, so the pull-back of the determinant of the 
index bundle
is trivial since it is trivial over ${\widetilde{\mathcal Q}}$.
\qed
\vskip0.3cm

We would like to choose orientations for all moduli spaces in a
consistent manner. To this end, we construct ``coherent orientations''
closely following~\cite{FloerHofer}.  Note that splicing gives an
identification $$\det(u_1)\wedge\det(u_2)\cong \det(u_1 * u_2),$$
where $u_1\in\pi_2(\x,\y)$ and $u_2\in\pi_2(\y,\w)$ are a pair of
maps.

\begin{defn}
\label{defn:CoherentSystem}
A {\em coherent system of orientations for $\spinc$}, ${\mathfrak o}$, 
is a choice of non-vanishing sections ${\mathfrak o}(\phi)$ of the determinant line
bundle over each $\phi\in\pi_2(\x,\y)$ 
for each $\x,\y\in{\mathcal S}$ and each
$\phi\in\pi_2(\x,\y)$, which are compatible with gluing in the sense that
$${\mathfrak o}({\phi_1})\wedge {\mathfrak o}({\phi_2})={\mathfrak o}({\phi_1 *\phi_2}),$$
under the identification coming from splicing, and
$${\mathfrak o}({u*S})={\mathfrak o}(u),$$ 
under the identification coming from the canonical orientation for
the moduli space of holomorphic spheres.
\end{defn}

To construct these it is useful to have the following:

\begin{defn}
\label{def:CompSetPath}
Let $(\Sigma,\alphas,\betas,z)$ be a Heegaard diagram representing
$Y$, and let $\spinc$ be a $\SpinC$ structure for $Y$. A {\em
complete set of paths for $\spinc$} is an enumeration
$\{\x_0,...,\x_m\}={\mathcal S}$ of all the intersection points of
$\Ta$ with $\Tb$ representing $\spinc$, and a collection of homotopy classes
$\theta_i\in\pi_2(\x_0,\x_i)$ for $i=1,...,m$ with $n_z(\theta_i)=0$. 
\end{defn}

Fix periodic classes $\phi_1,...,\phi_b\in\pi_2(\x,\x)$ representing a
basis for $H^1(Y;\Z)$, and non-vanishing
sections of the determinant line bundle for bundle for the
homotopy classes $\theta_1,...,\theta_m$ and $\phi_1,...,\phi_b$.
These data uniquely determine coherent system of orientations by splicing, 
since any homotopy class $\phi\in\phi_2(\x_i,\x_j)$ can be uniquely
written as $$\phi=a_1\phi_1+...+a_b\phi_b-\theta_i+\theta_j.$$

\subsection{Degenerate disks}

Fix a nearly symmetric almost-complex structure $J$ over
$\Sym^g(\Sigma)$. For each $\x\in\Ta\cap \Tb$, the moduli space of
{\em $\alpha$-degenerate disks} is the set of maps
$$\ModDeg_J(\x)=\left\{ u\colon [0,\infty)\times \R\longrightarrow
\Sym^g(\Sigma)\Bigg|
\begin{array}{l}
u(\{0\}\times \R)\subset \Ta \\
\lim_{z\goesto \infty}u(z)=\x \\
\frac{du}{ds}+J\frac{du}{dt}=0
\end{array}
\right\}.$$
Equivalently, we can think of $\ModDeg_J(\x)$ as the moduli space of
$J$-holomorphic maps of the unit disk $\CDisk$ in $\C$ to
$\Sym^g(\Sigma)$, which carry $\partial\CDisk$ into $\Ta$, and $i$ to
$\x$.  This also gives rise to a finite-dimensional moduli space,
partitioned according to the homotopy classes of maps satisfying these
boundary conditions, a set which we denote by $\pi_2(\x)$ (again, when
$g=2$, we divide out by the larger equivalence relation as in the
definition of $\pi_2(\x,\y)$). Suppose that $g>1$.  Since the map
$\pi_1(\Ta)\longrightarrow
\pi_1(\Sym^g(\Sigma))$ is injective, it follows that $\pi_2(\x)\cong\Z$
under the map $n_z(u)$; equivalently, if $O_\x\in\pi_2(\x)$ is the
homotopy class of the constant, then any other is given by $O_\x+k[S]$
for $k\in\Z$. Of course, when $g=1$, there is only one homotopy class, and
that is $O_\x$.

Note that there is a two-dimensional automorphism group acting on
$\ModDeg_J(\x)$ (pre-composing $u$ by either a purely imaginary
translation or a real dilation), and we denote the quotient space by
$\UnparModDeg_J(\x)$. If $\phi\in\pi_2(\x)$ is a homotopy class,
then we let $\ModDeg_J(\phi)$, resp. $\UnparModDeg_J(\phi)$, denote
its corresponding component in $\ModDeg_J(\x)$,
resp. $\UnparModDeg_J(\x)$. 

In studying smoothness properties of $\ModDeg_J(\x)$, it is useful to
have the following result concerning the complex structures $\Sym^g(\sj)$:

\begin{lemma}
\label{lemma:Genericj}
Given a finite collection of points $\{\x_i\}_{i=1}^n$ in
$\Sym^g(\Sigma)$, the set of complex structures $\sj$ over $\Sigma$ for
which there is a $\Sym^g(\sj)$-holomorphic sphere containing at least
one of the $\x_i$ has real codimension two.
\end{lemma}

\begin{proof}
The spheres in $\Sym^g(\Sigma)$ for the complex structure $\Sym^g(\sj)$
are all contained in the set of critical points for the Abel-Jacobi
map $$\AbelJacobi\colon \Sym^g(\Sigma)\longrightarrow
\Pic{g}(\Sigma)\cong H^1(\Sigma;S^1),$$ which is a degree one holomorphic map.
Thus, the set of spheres is contained in a subset of real
codimension two.
\end{proof}

\begin{prop}
\label{prop:CompactDegDisk}
Suppose $\x\in\Ta$ is not contained in any $\Sym^g(\sj)$-holomorphic sphere in $\Sym^g(\Sigma)$.
Then, there is a contractible neighborhood ${\mathcal U}$ of $\Sym^g(\sj)$ in
$\AlmostCx(\sj,\eta,V)$ with the property that for generic
$J\in{\mathcal U}$, the moduli space
$\UnparModDeg_J(O_\x+[S])$ is a compact, formally zero-dimensional
space which is smoothly cut out by its defining equations.
\end{prop}

\begin{proof}
To investigate compactness, note first that a sequence of elements in
$\UnparModDeg_J(\x)$ has a subsequence which either bubbles off spheres, or
additional disks with boundaries lying in $\Ta$. However, 
it is impossible for a sequence to bubble off a
null-homotopic disk with boundary lying in $\Ta$, since such disks must
be constant, as they have no energy (according to the
proof of Lemma~\ref{lemma:EnergyBound}, see
Equation~\eqref{eq:EnergyBound}). Moreover, sequences in $O_\x+[S]$ cannot bubble off 
homotopically non-trivial disks, because then one of the components in the
decomposition would have negative $\omega$-integral, and such homotopy
classes have no holomorphic representatives. 

This argument also rules out bubbling off spheres, except in the
special case where the subsequence converges to a single sphere (more
precisely, the constant disk mapping to $\x$, attached to some
sphere). But this is ruled out by our hypothesis on $\sj$, which ensures that
for any $J$ sufficiently close to $\Sym^g(\sj)$, the $J$-holomorphic
spheres are disjoint from $\x$.

To prove smoothness, note first that any holomorphic disk in
$\UnparModDeg(\x)$ for $\Sym^g(\sj)$ has a dense set of injective
points. To see this, fix any point $z'\in
\Sigma-\alpha_1-...-\alpha_g$. The intersection number of
$\{z'\}\times\Sym^{g-1}(\Sigma)$ with $u$ is $+1$, and both are
varieties; it follows that there is a single point of intersection,
i.e. there is only one $(s,t)$ for which $u(s,t)\in \{z'\}\times
\Sym^{g-1}(\Sigma)$.
Thus, $u$ is injective in a neighborhood of
$(s,t)$.  It follows that for any $J$ sufficiently close to
$\Sym^g(\sj)$, all the $J$-holomorphic degenerate disks are injective in
a neighborhood of $u(s,t)$. Thus, according to the usual proof of
transversality, these pairs are all smooth points in the parameterized
moduli space.  Thus, the result follows from the Sard-Smale theorem.
\end{proof}

Proceeding as earlier, we can orient the moduli spaces 
$\UnparModDeg_J$.
Our aim now is to prove the following:

\begin{theorem}
\label{thm:GromovInvariant}
Fix a finite collection $\{\x_i\}$ of points in $\Sym^{g}(\Sigma)$,
and an almost-complex structure $\sj$ over $\Sigma$
for which each 
$\Sym^{g}(\sj)$-holomorphic sphere misses $\{\x_i\}$.
Then, there is a contractible open neighborhood of $\Sym^g(\sj)$,
${\mathcal U}$, in the space
of nearly-symmetric almost-complex structures, with the property that
for generic $J\in {\mathcal U}$, the total signed number of points in
$\UnparModDeg_J(O_{\x_i}+[S])$ is zero.
\end{theorem}

Thinking of $\Sigma$ as the connected sum of $g$ tori, each of which
contains exactly one $\alpha_i$, we can endow $\Sigma$ with a complex
structure with long connected sum tubes.

\begin{prop}
\label{prop:ModDegEmpty}
If $\sj$ is sufficiently stretched out along the connected sum tubes,
then the moduli space ${\widehat \ModDeg}_\sj(O_\x+[S])$ is empty for any $\x\in\Ta$.
\end{prop}

\begin{proof}
Fix a genus one Riemann surface $E$. Let $\sj_t$
denote the complex structure on $\Sigma$, thought of as the connected
sum of $g$ copies of $E$, connected along cylinders isometric to
$S^1\times [-t,t]$. As $t\goesto\infty$, the Riemann surface degenerates to the wedge
product of $g$
copies of $E$, $\bigvee_{i=1}^g E_i$.

If for each $\sj_t$, the moduli space were non-empty, we could take
the Gromov limit of a sequence $u_t$ in ${\widehat
\ModDeg}_{\sj_t}(O_\x+[S])$
to obtain a holomorphic map $u_\infty$ into $\Sym^g(E_1\vee...\vee
E_g)$ (a linear chain of $g$ tori meeting in $g-1$ nodes). (In this
argument, we have a one-parameter family of symmetric products, which
we can embed into a fixed K\"ahler manifold, where we can apply the
usual Gromov compactness theorem, see also
Section~\ref{sec:Stabilization}.)  The latter symmetric product
decomposes into irreducible components $$\bigcup_{\{k_1,...,k_g \in \Z
\big| 0\leq k_i \leq g,~~k_1+...+k_g=g\}}
\Sym^{k_1}(E_1)\times...\times\Sym^{k_g}(E_g).$$
These components meet along loci containing the connected sum points
for the various $E_i$.  Moreover, the torus $\Ta$ can be viewed as a
subset of the irreducible component $E_1\times...\times E_g$
(corresponding to all $k_i=1$).  The Gromov limit $u_\infty$ then
consists of a holomorphic disk $v$ with boundary mapping into $\Ta$,
and a possible collection of spheres bubbling off into the other
irreducible components.  But $\pi_2(E_1\times...\times
E_g,\alpha_1\times...\times\alpha_g)=0$, so it follows that $v$ is
constant, mapping to $\x\in\Ta$ (which is disjoint from the connected
sum points). Since $v$ misses the other components of the symmetric
product, it cannot meet any of the spheres, so $v$ is the Gromov limit
of the $u_t$. But, $n_z(v)=0$, while we have assumed that
$n_z(u_t)=1$.
\end{proof}

\begin{lemma} 
\label{lemma:CompactCobordism}
Let $\x$, $\sj$, and ${\mathcal U}$
be as in Proposition~\ref{prop:CompactDegDisk}.
Suppose that $J_1,J_2\in {\mathcal U}$ are a pair of generic
almost-complex structures, in the sense that $\UnparModDeg_{J_s}(\x)$
is smooth for $s=0$ and $1$. Then, these moduli
spaces are compactly cobordant.
\end{lemma}

\begin{proof} 
We connect $J_1$ and $J_2$ by a generic path $\{J_s\}$ in ${\mathcal
U}$.  As in the proof of Proposition~\ref{prop:CompactDegDisk}, this
gives rise to the required compact cobordism. Note that the
possibility of bubbling off a sphere is ruled out, choosing ${\mathcal
U}$ small enough to ensure that $\x$ is disjoint from all
$J$-holomorphic spheres with $J\in{\mathcal U}$.
\end{proof}

\vskip.3cm
\noindent{\bf Proof of Theorem~\ref{thm:GromovInvariant}.}
Let $\sj$ be any complex structure over $\Sigma$ for which the
$\Sym^{g}(\sj)$-holomorphic spheres miss the $\{\x_i\}$.  Let
${\mathcal U}\subset \bigcap_{i=1}^n {\mathcal U}_i$ be a
contractible, open subset of the the open subsets ${\mathcal U}_i$
given to us by Proposition~\ref{prop:CompactDegDisk} for the points
$\x_i\in\Ta$.  According to Lemma~\ref{lemma:CompactCobordism}, the
number of points $\#\UnparModDeg_J(O_{\x_i}+[S])$ is independent of
$J$, i.e. it depends only on the complex structure $\sj$ over
$\Sigma$.  In fact, if $J$ is a generic $\sj$-nearly-symmetric
almost-complex structure, and $J'$ is a sufficiently close
$\sj'$-nearly-symmetric almost-complex structure, then the moduli
spaces are identified. It follows that
$\#\UnparModDeg_J(O_{\x_i}+[S])$ is a locally-constant function of the
complex structure $\sj$. Since the space of complex structures for
which the $\Sym^{g}(\sj)$-holomorphic spheres miss $\{\x_i\}$ is
connected, the theorem follows from
Proposition~\ref{prop:ModDegEmpty}.
\qed
\vskip.3cm

\subsection{Structure of moduli spaces}

\begin{theorem}
\label{thm:StructMod}
Let $(\Sigma,\alphas,\betas)$ be a Heegaard diagram with curves in
general position. For a generic path $J_s$ of nearly-symmetric
almost-complex structures, we have the following. There is no
non-constant $J_s$-holomorphic disk $u$ with $\Mas(u)\leq 0$. Moreover
for each $\phi\in\pi_2(\x,\y)$ with 
$\Mas(\phi)=1$,  the quotient space
$$\UnparModSp(\phi)=\frac{\ModSp(\phi)}{\R}$$ is a compact,
zero-dimensional manifold.
\end{theorem}

\begin{proof}
The first part follows directly from the  Theorem~\ref{thm:Transversality}.

Compactness follows from the usual compactification theorem for
holomorphic curves
(see~\cite{FloerUnregularized} and also~\cite{Gromov}, 
\cite{ParkerWolfson}, \cite{Ye}), which holds thanks to the energy bound
(Lemma~\ref{lemma:EnergyBound}).

Specifically, the compactness theorem says that a sequence of points
in the moduli spaces converges to an ideal disk, with possible broken
flowlines, boundary degenerations, and bubblings of
spheres. Broken flowlines are excluded by the additivity of the Maslov
index, and the transversality result
Theorem~\ref{thm:Transversality}. Spheres and boundary
degenerations both carry Maslov index at
least two, so these kinds of degenerations are excluded as well. 
\end{proof}

\section{Definition of the Floer homology groups}
\label{sec:DefHF}

We are now ready to define the Floer homology groups for
three-manifolds.  In Subsection~\ref{subsec:QHS}, we consider the
technically simpler case of three-manifolds with $b_1(Y)=0$.  We then
turn our attention to the issues which arise when pushing this
definition to the case of three-manifolds with $b_1(Y)>0$: the cyclic
gradings in Subsection~\ref{subsec:Grading}, the ``admissibility
hypotheses'' on the Heegaard diagrams required for topological
invariance of the constructions in
Subsection~\ref{subsec:Admissibility}. (We will return to the
construction of such Heegaard diagrams in
Subsection~\ref{sec:Special}.)  With these technical pieces in place,
we proceed as before to define the Floer homology groups when
$b_1(Y)>0$, in Subsection~\ref{subsec:BOneBig}. These groups can be
endowed the additional structure of the action by $H_1(Y;\Z)$, which
is constructed in Subsection~\ref{subsec:DefAct}. An additional
modification -- Floer homology with ``twisted coefficients'' -- is
introduced in Section~\ref{HolDiskTwo:sec:TwistedCoeffs}
of~\cite{HolDiskTwo}.

\subsection{Floer homologies when $b_1(Y)=0$}
\label{subsec:QHS}

Let $(\Sigma,\alphas,\betas,z)$ be a pointed
Heegaard diagram with genus $g>0$ for a rational homology three-sphere $Y$,
where the $\alphas$ and $\betas$ are in  general position, and choose a 
$\SpinC$ structure $\spinc \in \SpinC(Y)$. We  let $\SpinCz$ be the 
set of intersection points $\x\in\Ta\cap\Tb$ with $s_z(\x)=\spinc$. 
We fix also the following auxiliary data:
\begin{itemize}
\item a coherent orientation system, in the sense of Definition~\ref{defn:CoherentSystem}
(note that this is not necessary when defining Floer homology groups
with $\Zmod{2}$ coefficents)
\item a generic complex structure $\sj$ over $\Sigma$ (generic in the sense
that each intersection point $\x\in\Ta\cap\Tb$ is disjoint from the
$\Sym^g(\sj)$-holomorphic spheres in $\Sym^g(\Sigma)$ -- see
Lemma~\ref{lemma:Genericj}),
\item a generic path of
nearly-symmetric almost-complex structure $J_s$ over $\Sym^g(\Sigma)$,
contained in the open subset ${\mathcal U}$ of
Theorem~\ref{thm:GromovInvariant} (associated to the subset
$\Ta\cap\Tb$ and also $V$),
\end{itemize} 

Let $\CFa(\alphas,\betas,\spinc)$ denote the free
Abelian group generated by the points in $\SpinCz\subset\Ta\cap\Tb$
This group 
can be endowed with a relative grading\footnote{A relatively graded Abelian group is one which is
generated by elements partitioned into equivalence classes ${\mathfrak
S}$, with a relative grading function $\gr\colon {\mathfrak S}\times
{\mathfrak S}\longrightarrow \Z$, satisfying
$\gr(\x,\y)+\gr(\y,\w)=\gr(\x,\w)$ for each $\x,\y,\w\in{\mathfrak
S}$. When the corresponding theory for four-manifolds is developed,
this relative $\Z$-grading can be lifted to an absolute $\Q$-grading,
see~\cite{HolDiskThree}.}, defined by
\begin{equation}
\label{eq:RelativeGrading}
\gr(\x,\y)=\Mas(\phi)-2 n_z(\phi),
\end{equation} 
where $\phi$ is any element
$\phi\in\pi_2(\x,\y)$, and 
$\Mas$ is the Maslov index. 
In view of Proposition~\ref{prop:WhitneyDisks} and Lemma~\ref{lemma:MasClass},
this integer is independent of the
choice of Whitney disk $\phi\in\pi_2(\x,\y)$.

Let $$\partial\colon \CFa(\alphas,\betas,\spinc)\longrightarrow
\CFa(\alphas,\betas,\spinc)$$ be the map defined by: $$\partial \x =
\sum_{\{\y\in\SpinCz \big| \gr(\x,\y)=1\}}
\#\left(\UnparModSp_0(\x,\y)\right)\y,$$
where $\UnparModSp_0(\x,\y)=\UnparModSp(\phi)$ for the element
$\phi\in\pi_2(\x,\y)$ with $n_z(\phi)=0$ and $\Mas(\phi)=1$.  Note
that by Proposition~\ref{prop:WhitneyDisks} and
Lemma~\ref{lemma:MasClass}, there is at most one such homotopy class.
Also, counting these 
holomorphic disks in $\Sym^g(\Sigma)$ is equivalent to counting
holomorphic disks in $\Sym^g(\Sigma-z)$, in view of
Lemma~\ref{lemma:NonNegativity}.
(We have suppressed the path $J_s$ from the notation, but one should
bear in mind that $\partial$ does depend on the path $J_s$. When it is
important to call attention to this dependence, we write
$\partial_{J_s}$, see the proof of Theorem~\ref{thm:IndepCxStruct}
below.) 

The count appearing in the above boundary operator is, as usual,
meant to signify a signed (oriented) count of points in the compact,
zero-dimensional moduli spaces (see Theorem~\ref{thm:StructMod}, and
Subsection~\ref{subsec:Orientability}), and as such, it depends on a
coherent orientation system as defined in
Definition~\ref{defn:CoherentSystem}. As we shall see in
Lemma~\ref{lemma:OrientationSystems} in the present case,
different such choices give rise to isomorphic chain complexes,
so we shall usually drop them from the notation.

\begin{theorem}
When $b_1(Y)=0$, the pair $(\CFa(\alphas,\betas,\spinc),\partial)$ is
a chain complex; i.e.  $\partial^2 = 0$.
\end{theorem}

\begin{proof}
This follows in the usual manner from the compactifications of the
one-dimensional moduli spaces $\UnparModFlow(\phi)$ with
$\mu(\phi)=2$ (together with the gluing descriptions of the
neighborhoods of the ends).  Note that if $\x'\in  \Ta\cap \Tb-\SpinCz$, then
$\epsilon(\x,\x')\neq 0$, so there are no flows connecting $\x$ to $\x'$.
We note also that there are
no spheres in $\Sym^g(\Sigma-z)$ or
degenerate holomorphic disks (whose boundary lies entirely in $\Ta$ or
$\Tb$),
so the only boundary
components in the compactification consist of broken flow-lines.
\end{proof}

\begin{defn}
The Floer homology groups $\HFa(\alphas,\betas,\spinc)$ are
the homology groups of the
complex $(\CFa(\alphas,\betas,\spinc),\partial)$.
\end{defn}

Next, let $\CFinf(\alphas,\betas,\spinc)$ be the free Abelian group generated
by pairs $[\x,i]$ where $\x\in\SpinCz $, and $i\in \Z$ is an integer.
We give the generators a relative grading defined by
$$\gr([\x,i],[\y,j])=\gr(\x,\y)+2i-2j.$$

Let 
$$\partial^\infty\colon \CFinf(\alphas,\betas,\spinc)\longrightarrow
\CFinf(\alphas,\betas,\spinc)$$ be the map defined by: 
\begin{equation}
\label{eq:DefBoundaryInfty}
\partial^\infty
[\x,i] = \sum_{\y\in\pi_2(\x,\y)}\sum_{\{\phi\in\pi_2(\x,\y)\big|\Mas(\phi)=1\}}
\#\left(\UnparModSp(\phi)\right)[\y,i-n_z(\phi)].
\end{equation}
Although we have written the above expression as a double-sum, 
Proposition~\ref{prop:WhitneyDisks} and Lemma~\ref{lemma:MasClass} ensure that
for given $\x$ and $\y$, there is at most one homotopy class $\phi\in\pi_2(\x,\y)$
with $\Mas(\phi)=1$.

\begin{theorem}
\label{thm:DSquaredZero}
When $b_1(Y)=0$, the pair
$(\CFinf(\alphas,\betas,\spinc),\partial^\infty)$ is a chain complex;
i.e.  $(\partial^\infty)^2 = 0$.
\end{theorem}

\begin{proof}
As is usual in Floer's theory, one considers the ends of the moduli
spaces $\UnparModFlow(\phi)$, where $\phi\in\pi_2(\x,\w)$ satisfies
$\Mas(\phi)=2$. This space has {\em a priori} three
kinds of ends:
    \begin{list}
	{(\arabic{bean})}{\usecounter{bean}\setlength{\rightmargin}{\leftmargin}}
\item those corresponding to ``broken flow-lines'', i.e. a pair 
$u\in\ModFlow(\x,\y)$ and $v\in\ModFlow(\y,\w)$ with
$\Mas(u)=\Mas(v)=1$
\item
\label{item:Spheres} those which correspond to a sphere bubbling off, i.e. another 
$v\in\ModFlow(\x,\w)$ and a holomorphic sphere
$S\in\Sym^g(\Sigma)$ 
which meets $v$
\item 
\label{item:BoundaryDeg}
those which correspond to ``boundary bubbling'', i.e. we have 
a $v\in\ModFlow(\x,\w)$, and a holomorphic map $u$ from the disk,
whose boundary is mapped into $\Ta$ or $\Tb$, which meet in  a point on the boundary.
\end{list}

(In principle, several of the above degenerations could happen at once
-- multiple broken flows, spheres, and boundary degenerations, but
these multiple degenerations  are easily ruled out by dimension counts
and the transversality theorem, Theorem~\ref{thm:Transversality}.)

In the Cases~(\ref{item:Spheres}) and (\ref{item:BoundaryDeg}), we
argue that $[v]=\phi-\ell[S]$ (note that a disk whose boundary lies
entirely inside $\Ta$ or $\Tb$ also has a corresponding domain
$\cald(u)$, which, in this case, must be a multiple of $\Sigma$; if
$u$ is pseudo-holomorphic, then $\cald(u)=\ell[\Sigma]$ for $\ell\geq
0$ according to Lemma~\ref{lemma:NonNegativity}, and if $\cald(u)=0$,
the disk must be constant).  Thus it follows from
Lemma~\ref{lemma:MasClass} that $\Mas([v])=\Mas(\phi)-2\ell$.  From
transversality (Theorem~\ref{thm:Transversality}), it follows that $\ell=1$ and 
$v$ must be constant; so that, in particular, $\x=\w$.  Now for generic
$\sj$, we know that the holomorphic spheres miss the intersection
points $\Ta\cap\Tb$; hence, the case of spheres bubbling off is
excluded.

Thus, when $\x\neq \w$, boundary bubbles are excluded, so, counting
the ends of the moduli space $\ModFlow(\phi)$, we get that $$0=\sum_\y
\sum_{\{\psi\in\pi_2(\x,\y),
\zeta\in\pi_2(\y,\w)\big|\psi\csum\zeta=\phi\}}
\Big(\#\ModFlow(\psi)\Big)\cm\Big(\#\ModFlow(\zeta)\Big).$$
When $\x=\w$, there are additional terms, corresponding to the
boundary bubbles (and the gluing descriptions of the ends), giving a relation
$$0=\#\UnparModDeg^\alpha(\x)+\#\UnparModDeg^\beta(\x) +\sum_\y \sum_{\{\psi\in\pi_2(\x,\y),
\zeta\in\pi_2(\y,\x)\big|\psi\csum\zeta=\phi\}}
\Big(\#\ModFlow(\psi)\Big)\cm \Big(\#\ModFlow(\zeta)\Big),$$
see for example~\cite{FOOO}.
But the terms 
$\#\UnparModDeg^\alpha(\x)$ and $\#\UnparModDeg^\beta(\x)$ both vanish, according to
Theorem~\ref{thm:GromovInvariant}.

From the additivity of $n_z$ under juxtapositions of flow-lines,
it follows that the double sums considered above are coefficients of
$\partial^2[\x,i]$.
\end{proof}

There is a chain map
$$U\colon \CFinf(\alphas,\betas,\spinc)\longrightarrow \CFinf(\alphas,\betas,\spinc),$$
which lowers degree by two, defined by
$$U[\x,i]=[\x,i-1].$$

Let $\CFm(\alphas,\betas,\spinc)$ denote the subgroup of $\CFinf(\alphas,\betas,\spinc)$
which is freely generated by pairs $[\x,i]$, where
$i<0$. 
Let $\CFp(\alphas,\betas,\spinc)$ denote the quotient group
$\CFinf(\alphas,\betas,\spinc)/\CFm(\alphas,\betas,\spinc)$.

\begin{lemma}
The group $\CFm(\alphas,\betas,\spinc)$ is a subcomplex of
$\CFinf(\alphas,\betas,\spinc)$, so we have a short exact sequence of chain
complexes: 
$$
\label{eq:ExactSequence}
\begin{CD}
0 @>>> CF^-(\alphas,\betas,\spinc)@>{i}>>
\CFinf(\alphas,\betas,\spinc)@>{\pi}>>\CFp(\alphas,\betas,\spinc) @>>>0\end{CD}.  
$$
\end{lemma}

\begin{proof}
The fact that $\CFm(\alphas,\betas,\spinc)$ is a subcomplex is an easy consequence
of Lemma~\ref{lemma:NonNegativity}. 
\end{proof}

Clearly, $U$ is restricts to an endomorphism of $\CFm(\alphas,\betas,\spinc)$
(which lowers degree by $2$), and hence it also induces an
endomorphism of the quotient $\CFp(\alphas,\betas,\spinc)$. Sometimes, for
clarity we denote the induced actions on these complexes (or their
homologies) by $U^-$ or $U^+$. It is also 
easy to see that there is a short
exact sequence 
$$
\begin{CD} 0@>>> \CFa(\alphas,\betas,\spinc) @>{\iota}>> 
\CFp(\alphas,\betas,\spinc)
@>{U^+}>>\CFp(\alphas,\betas,\spinc)@>>>0,
\end{CD}
$$
where $\iota(\x)=[\x,0]$. In view of this, we declare the $U$ action
on $\CFa(\alphas,\betas,\spinc)$ to be trivial.

\begin{defn}
Let $\HFinf(\alphas,\betas,\spinc)$, $\HFm(\alphas,\betas,\spinc)$, and
$\HFp(\alphas,\betas,\spinc)$ denote the homologies of the complexes
$\CFm(\alphas,\betas,\spinc)$, $\CFinf(\alphas,\betas,\spinc)$, and
$\CFp(\alphas,\betas,\spinc)$ respectively, thought of as $\Z[U]$ modules.
\end{defn}

Note of course, that these constructions can be carried out with
coefficients in any ring $\Lambda$: one defines the corresponding
chain complexes as free $\Lambda$-modules, and the homology groups
obtained in this way are modules over the polynomial algebra
$\Lambda[U]$. We have no particular use for this construction, as $\Z$
is a ``universal'' case, though we do point out that if we had chosen
to use $\Lambda=\Zmod{2}$, then the issues of orientation (and choices of
orientation system) would become unnecessary.

In the interest of conciseness, we have suppressed additional data --
notably, complex structures (and their perturbations) and orientation
systems -- from the notation of these homology groups. In fact, 
we will show that the homology groups are
independent of these choices.

There is another algebraic construction which gives rise to a
finitely-generated variant of Floer homology, which we define with the help
of the following:

\begin{lemma}
\label{lemma:Hom}
If $k$ is sufficiently large, then $$\Image (U^+)^k = \Image
(\pi_\ast),$$ where $\pi_*\colon \HFinf(\alphas,\betas,\spinc)
\longrightarrow \HFp(\alphas,\betas,\spinc)$ is the naturally induced map on homology. 
Similarly $$\Ker (U^-)^k = \Ker (i_*),$$ for $i_*:
\HFm(\alphas,\betas,\spinc)\longrightarrow
\HFinf(\alphas,\betas,\spinc)$.
\end{lemma}

\begin{proof}
We begin with the claim about $\HFm$. It is clear that the kernel of
$i_*$ consists of those homology classes $\xi\in
\HFm(\alphas,\betas,\spinc)$ for which $(U^-)^k\xi = 0$ for sufficiently
large integers $k$. To see that $k$ can be chosen independent of $\xi$,
observe that $\CFm$ is a finitely generated chain complex over the
ring $\Z[U^{-}]$, and hence its homology is finitely generated over
$\Z[U^{-}]$, as well. Thus, the sequence of submodules of
$\HFm(\alphas,\betas,\spinc)$
$$\Ker(U^-) \subseteq \Ker(U^-)^2
\subseteq
\Ker(U^-)^3 \subseteq ... $$ must stabilize. 

Now, we claim that for $k$ large enough that $\Ker(U^-)^k=\Ker (i_*)$,
it is also the case that $\Image(U^+)^k=\Image(\pi_*)$. First, since
$U$ is an automorphism of $\HFinf(\alphas,\betas,\spinc)$, it is
clear that $\Image(U^+)^k\supseteq\Image(\pi_*)$. Conversely, if 
$\xi\in\Image(U^+)^k$, then writing
$\xi=(U^+)^{k}\eta$ , we have that
$\delta \xi = (U^-)^{k}\delta \eta=0$ (since
$i_*\circ \delta=0$); thus $\xi\in\Image(\pi_*)$.
\end{proof}

\begin{defn}
Let 
$$\HFpred(\alphas,\betas,\spinc)={\HFp(\alphas,\betas,\spinc)}/{\Image (U^+)^k},$$
for sufficiently large $k$. Similarly, let 
$$\HFmred(\alphas,\betas,\spinc)=\Ker (U^-)^k\subset \HFm(\alphas,\betas,\spinc).$$
\end{defn}

\begin{prop}
The boundary homomorphism of the long exact sequence induces a
$U$-equivariant isomorphism $$\HFpred(\alphas,\betas,\spinc)\cong
\HFmred(\alphas,\betas,\spinc).$$
Moreover, both are finitely generated $\Z$ modules.
\end{prop}

\begin{proof}
The isomorphism follows readily from Lemma \ref{lemma:Hom} and the
long exact sequence. Since $\HFmred(\alphas,\betas,\spinc)$ is a
finitely-generated $\Z[U^-]$-module which is annihilated by $(U^-)^k$,
it follows immediately that this module is also finitely generated over $\Z$.
\end{proof}

In view of the above result, we will denote 
$\HFpred(\alphas,\betas,\spinc)\cong \HFmred(\alphas,\betas,\spinc)$ simply
by $\HFred(\alphas,\betas,\spinc)$.

\subsection{Constructions when $b_1(Y)>0$}
\label{subsec:DefHFBOneBig}

\subsubsection{Grading}
\label{subsec:Grading}
    
As before, we fix a $\SpinC$ structure $\spinc \in \SpinC(Y)$, and let
$\SpinCz$ be the set of intersection points $\x\in\Ta\cap\Tb$ with
$s_z(\x)=\spinc$. 
The expression used to define relative grading on $\SpinCz$ as in
Equation~\eqref{eq:RelativeGrading} is now
well-defined only modulo an indeterminacy, given by
the Maslov indices of periodic classes. Indeed, this indeterminacy is given 
by the following more familiar quantity
\begin{equation}
\label{eq:Indeterminacy}
\divis(\spinc)=\gcd_{\xi\in H_{2}(Y;\Z)}\langle
c_{1}(\spinc),\xi\rangle,
\end{equation} in view of the following result, which is
proved in Subsection~\ref{subsec:MasIndex}:

\begin{theorem}
\label{thm:Grading}
Fix a $\SpinC$ structure $\spinc\in 
\SpinC(Y)$. Then for each  $\x\in \Ta\cap \Tb$ with
$s_z(\x)=\spinc$,  and for each periodic class
$\psi \in \PerClasses{\x}$ we have 
$$\mu(\psi)=\langle c_1(\spinc), H(\psi)\rangle ,$$
where $H(\psi)\in H_2(Y;\Z)$ is the homology class corresponding to
the periodic class $\psi$.
\end{theorem}

\subsubsection{Admissibility}
\label{subsec:Admissibility}

To ensure compactness of the index one moduli spaces connecting
intersection points, we will need to use only certain special kinds of
Heegaard diagrams. It turns out that these conditions are somewhat
different for the various theories.

\begin{defn} 
\label{def:NonTorsionAdmissible}
A pointed Heegaard diagram is called {\em strongly
admissible for  the $\SpinC$
structure $\spinc$} if  for every non-trivial
periodic domain $\cald$ with $$\langle
c_1(\spinc),H(\cald) \rangle = 2n \geq 0, $$ $\cald$ has some
coefficient $>n$. A pointed Heegaard diagram is called {\em weakly
admissible for $\spinc$} if for each non-trivial periodic domain
$\cald$ with $$\langle c_1(\spinc),H(\cald) \rangle = 0, $$ $\cald$
has both positive and negative coefficients.  
\end{defn}

\begin{remark}
\label{rmk:AdmissibleAtOnce}
Note that for a $\SpinC$ structure with $c_1(\spinc)$ torsion, the
weak and strong admissibility conditions coincide. Also note that if a
Heegaard diagram is strongly admissible for any torsion $\SpinC$
structure then in fact it is weakly admissible for all $\SpinC$
structures.
\end{remark}

We have the following geometric reformulation of the weak
admissibility condition (for all $\SpinC$ structures):

\begin{lemma}
\label{lemma:EnergyZero}
A Heegaard diagram is weakly admissible for all $\SpinC$ structures
if and only if $\Sigma$ can be endowed
with a volume form for which each periodic domain has total signed
area equal to zero.
\end{lemma}

\begin{proof}
The existence of such a volume form obviously implies weak
admissibility, since each non-trivial domain has positive area.

Assume, conversely, that each non-trivial periodic domain has both
positive and negative coefficients.  By changing the volume form, we
are free to make each domain in $\CurveComp$ have arbitrary positive
area. Thus, the claim now reduces to some linear algebra. We say that
a vector subspace $V\subset \R^m$ is {\em balanced} if each of its
non-zero vectors has both positive and negative components.  The
claim, then, follows form the fact that a vector subspace of $\R^m$
which is balanced admits an orthogonal vector each of whose
coefficients is positive.

This fact is true by induction on the dimension of the ambient vector
space (and it is vacuously true for $m=1$). Now, suppose $V$ is a
balanced subspace of $\R^m$, and let $\Pi_i\colon \R^m
\longrightarrow
\R^{m-1}$ denote the projection map
$\Pi_i(x_1,...,x_m)=(x_1,...,\widehat{x_i},...,x_m)$. Either
$\Pi_i(V)$ is also balanced, or $V$ contains a vector $v$ whose
$i^{th}$ component is $+1$, all other components are non-positive, and
at least one of them is negative. In this latter case, we construct the
required positive orthogonal vector as follows.  Apply the induction
hypothesis to find a vector
$\xi=(\xi_1,...,\xi_{i-1},0,\xi_{i+1},...,\xi_m)$ with $\xi_j>0$ for
$i\neq j$, which is orthogonal to $V\cap \R^{m-1}$. The required
vector, then, is $\xi-\langle v,\xi\rangle e_i$. 

If, on the other hand, all $i$ of the vector spaces $\Pi_i(V)$ are
balanced, then by induction we can find vectors
$\xi=(0,\xi_2,...,\xi_m)$ and $\eta=(\eta_1,0,\eta_3,...,\eta_m)$ with
$\xi_i>0$ for $i\neq 1$, and $\eta_i>0$ for $i\neq 2$. Then,
$\xi+\eta$ is our required vector.
\end{proof}

The following two lemmas are, ultimately, the reasons for introducing the
admissibility hypotheses.

\begin{lemma}
\label{lemma:WeakFiniteness}
Suppose that $(\Sigma,\alphas,\betas,z)$ is weakly admissible for the
$\SpinC$ structure $\spinc$, and fix
integers $j,k \in\Z$. Then, for
each $\x,\y\in{\mathfrak S}$, there are only finitely many
$\phi\in\pi_2(\x,\y)$ for which $\Mas(\phi)=j$, $n_z(\phi)=k$, and
$\cald(\phi)\geq 0$.
\end{lemma}

\begin{proof}
Fix some initial
$\psi\in\pi_2(\x,\y)$ with $\Mas(\psi)=j$. 
Then, in view of Theorem~\ref{thm:Grading},
any
other $\phi\in\pi_2(\x,\y)$ with $\Mas(\psi)=j$ has the form
$$\phi=\psi+\PerClass{\x} -
\frac{\langle c_1(\spinc),H(\PerDom) \rangle}{2} S,$$ 
where $\PerClass{x}$ is some periodic class, 
$\PerDom$ its associated
periodic domain, and $S$ is the positive generator of $\pi_2(\Sym^g(\Sigma))$.
If $n_z(\psi)=n_z(\phi)$, this forces $\cald(\phi)=\cald(\psi)+\PerDom$ for some
periodic domain whose associated homology class is annihilated by
$c_1(\spinc)$. If $\cald(\phi)\geq 0$, then
$\PerDom\geq -\cald(\psi)$.

Thus, the lemma follows from the observation that for any fixed
$\psi\in\pi_2(\x,\y)$, there are only finitely many periodic domains
$\PerDom$ in the set $$Q=\left\{ \PerDom\in\PerClasses{\x} |
\langle c_1(\spinc),H(\PerDom)\rangle = 0, \PerDom\geq -\cald(\psi)\right\}.$$
We see this as follows. Let $m$ denote the total number of domains
(components in $\CurveComp$). We can think of $Q$ as lattice points in
the $m$-dimensional vector space generated by the domains
$\cald_i$. Given $p\in Q$, written as $p=\sum a_i\cald_i$, we let
$\|p\|$ denote its naturally induced Euclidean norm
$$\|p\|=\sqrt{\sum_{i=1}^m |a_i|^2}.$$ If $Q$ had infinitely many
elements, we could find a sequence of $\{p_j\}_{j=1}^\infty \subset Q$
with $\|p_j\|\goesto \infty$. In particular, the sequence
$\frac{p_j}{\|p_j\|}$ has a subsequence which converges to a unit
vector in the vector space of periodic domains with real coefficients
which annihilate $c_1(\spinc)$. We write the vector as $p=\sum b_i
\cald_i$. Since the coefficients of $p_j$ are bounded below, but the
lengths of the $p_j$ diverge, it follows that all the coefficients of
$p$ are non-negative. Of course, if the polytope the subspace of
$H_2(Y;\Z)$ annihilated by $c_1(\spinc)$, corresponding to periodic
domains with only non-negative multiplicities has a non-trivial real
vector, then it must also have a non-trivial rational vector. After
clearing denominators, we obtain a non-zero periodic domain (with integer
coefficients) annihilating $c_1(\spinc)$, with only non-negative
coefficients. This contradicts the hypothesis of weak admissibility.
\end{proof}

\begin{lemma}
\label{lemma:StrongFiniteness}
For a strongly admissible pointed Heegaard diagram, and an integer
$j$,  there are only
finitely many $\phi\in\pi_2(\x,\y)$ with $\Mas(\phi)=j$ and 
$\cald(\phi)\geq 0$.
\end{lemma}

\begin{proof} 
Fix a reference $\psi\in\pi_2(\x,\y)$ with $\Mas(\psi)=j$.  Then, as in the
previous lemma, 
any other class $\phi\in\pi_2(\x,\y)$ with $\Mas(\psi)=j$ can be
written as $$\phi=\psi-\PerClass{\x}+
\frac{\langle c_1(\spinc),H(\PerDom) \rangle}{2} S.$$  Thus, each class $\phi$
with $\cald(\phi)\geq 0$ corresponds to a periodic domain $\PerDom$ with
$$-\PerDom+
\frac{\langle c_1(\spinc),H(\PerDom)\rangle}{2} [\Sigma] \geq 
-\cald(\psi). $$ 

The lemma follows from the fact that (for fixed $\psi$) there are only finitely many
periodic domains satisfying this inequality.  This follows as in the
proof of Lemma~\ref{lemma:WeakFiniteness}: an infinite number of such
periodic domains would give rise to a a real periodic domain $\PerDom$
for which $$-\PerDom+
\frac{\langle c_1(\spinc),H(\PerDom)\rangle}{2} [\Sigma] \geq 0,$$ 
from which it is easy to see that there must be an integral periodic
domain with the same property. But such a periodic domain would
violate the strong admissibility hypothesis.
\end{proof}

We will establish the existence of admissible Heegaard diagrams in
Section~\ref{sec:Special}.

\subsubsection{The chain complex}
\label{subsec:BOneBig}

Let $(\Sigma,\alphas,\betas,z)$ be a pointed Heegaard diagram for a
three-manifold with $b_1(Y)>0$, and fix a $\SpinC$ structure $\spinc$.
Suppose moreover that the Heegaard diagram is strongly
$\spinc$-admissibile.  In this case, we define the groups
$\CFa(\alphas,\betas,\spinc,\orient)$,
$\CFinf(\alphas,\betas,\spinc,\orient)$ as in
Subsection~\ref{subsec:QHS}. Note that when $b_1(Y)>0$, we 
do  include the coherent orientation system in the notation, 
since now the Floer homologies do depend on this choice -- 
we shall see in Lemma~\ref{lemma:OrientationSystems}
that there are in principle $2^{b_1(Y)}$ different possible chain
complexes corresponding to variations in this choice.
Equation~\eqref{eq:RelativeGrading} now endows
$\CFa(\alphas,\betas,\spinc,\orient)$ with a relative
$\Z/\divis(\spinc)$-grading, where $\divis(\spinc)$ is given in
Equation~\eqref{eq:Indeterminacy}. We can define the subgroup
$\CFm(\alphas,\betas,\spinc,\orient)$ and quotient group
$\CFp(\alphas,\betas,\spinc,\orient)$ as before.  We endow
$\CFinf(\alphas,\betas,\spinc,\orient)$ with the differential from
Equation~\eqref{eq:DefBoundaryInfty}, endowing
$\CFp(\alphas,\betas,\spinc,\orient)$ with the induced differential $$
\partial^+[\x,i] = \sum_{\y\in\pi_2(\x,\y)}\sum_{\{\phi\in\pi_2(\x,\y)\big|\Mas(\phi)=1,
i\geq n_z(\phi)\}}
\#\left(\UnparModSp(\phi)\right)[\y,i-n_z(\phi)].
$$

\begin{theorem}
\label{thm:DSquaredZeroBOneBig}
Let $Y$ be a three-manifold equipped with a $\SpinC$ structure $\spinc$.
Then, 
\begin{itemize}
\item if $(\Sigma,\alphas,\betas,z)$ is strongly $\spinc$-admissible,
then $\CFinf(\alphas,\betas,\spinc,\orient,\partial^\infty)$ is a chain complex,
with subcomplex $\CFm$ and quotient complex $\CFp$.
\item if $(\Sigma,\alphas,\betas,z)$ is weakly $\spinc$-admissible,
then $\CFp(\alphas,\betas,\spinc,\orient,\partial^+)$ is a chain complex
with subcomplex $\CFa(\alphas,\betas,\spinc,\orient,{\widehat \partial})$.
\end{itemize}
\end{theorem}

\begin{proof}
When $(\Sigma,\alphas,\betas,\spinc,\orient)$ is strongly
$\spinc$-admissible, the key point is to observe that the boundary
operators Equation~\eqref{eq:DefBoundaryInfty} is actually a finite
sum. This follows from
Lemma~\ref{lemma:StronglyAdmissible}. Similarly, when the diagram is
only weakly admissible, Lemma~\ref{lemma:WeakFiniteness} ensures
that the differentials for $\CFa$ and $\CFp$ are finite sums.
With these remarks, the proof proceeds
exactly as in the proof of Theorem~\ref{thm:DSquaredZero}.
\end{proof}

\subsubsection{Coherent orientation systems}

Although the above chain complexes depend on the choice of orientation system,
its isomorphism type depends on the orientation system only through its
equivalence class, in the following sense.
Let ${\mathfrak o}$ and ${\mathfrak o}'$ be a pair of systems of
coherent orientations for $\spinc$.  Then, we define their difference
$\delta=\delta({\mathfrak o},{\mathfrak o}')\in
\Hom(H^1(Y;\Z),\Zmod{2})$ as follows.  Let $\phi\in\pi_2(\x,\x)$ be
the periodic class representing some cohomology class $H\in
H^1(Y;\Z)$. Then, the section $\orient$ of the determinant line bundle
over the component specified by $\phi$ is either a positive multiple
of $\orient'$, in which case we let $\delta(H)=0$, or it is a negative
multiple of $\orient'$, in which case we let $\delta(H)=1$. We say
that two systems of coherent orientations are {\em equivalent} if
their difference $\delta$ vanishes.
Clearly, there are $2^{b_1(Y)}$ inequivalent choices of orientation conventions.

\begin{lemma}
\label{lemma:OrientationSystems}
If $\orient$ and $\orient'$ are equivalent orientation systems in the
above sense, then the chain complexes
$\CFinf(\alphas,\betas,\spinc,{\mathfrak o})$ and
$\CFinf(\alphas,\betas,\spinc,{\mathfrak o}')$ (and the corresponding
$\CFm$, $\CFp$, and $\CFa$)
are isomorphic.
\end{lemma}

\begin{proof}
Let $\orient$ and $\orient'$ be a pair of isomorphic
coherent orientation systems,
Fix a reference point $\x_0\in\SpinCz$. 
Given any other $\x\in\SpinCz$ and path $\phi\in\pi_2(\x_0,\x)$,
there is a sign $\sigma(\x)\in\{\pm 1\}$ 
with the property that for $\orient(\phi)=\sigma(\x)\cm \orient'(\phi)$.
As the notation suggests, if $\orient$ and $\orient'$ are isomorphic
orientation systems, the number 
$\sigma(\x)$ is independent of the choice of $\phi$, so we
can define a map 
$$f\colon\CFinf(\alphas,\betas,\spinc,{\mathfrak o})\longrightarrow
\CFinf(\alphas,\betas,\spinc,{\mathfrak o})$$
by
$f([\x,i])=\sigma(\x)\cm[\x,i]$. It is straightforward to verify
that $f$ induces an isomorphism of chain complexes. 
\end{proof}

\subsubsection{Additional algebra: the $H_1(Y;\Z)/\Tors$ and $U$-actions}
\label{subsec:DefAct}

As in the case where $b_1(Y)=0$, we define $\HFm$, $\HFinf$, $\HFp$,
and $\HFa$ to be the homologies of the corresponding chain complexes;
and, as before, all of these homology groups come with the structure
of a $\Z[U]$ module, where $U$ lowers the relative grading by
two. Similarly, we can form a reduced group $\HFred$ using the
procedure outlined earlier.  Moreover, when $b_1(Y)>0$, there is a new
algebraic object: an action of $H_1(Y;\Z)/\Tors$ the Floer homology
groups.  Recall from the proof of Proposition~\ref{prop:WhitneyDisks}
that the choice of basepoint gives an isomorphism
\begin{equation}
\label{eq:CohomConfig}
H^1(\Omega(\Ta,\Tb);\Z)\cong
\Hom(\pi_1(\Omega(\Ta,\Tb)),\Z)\cong \pi_2(\Sym^g(\Sigma))
\oplus \Hom(H^1(Y,\Z),\Z).
\end{equation}

\begin{prop}
\label{prop:Action}
There is a natural action of $H^1(\Omega(\Ta,\Tb))$ lowering degree by
one on $\HFinf(Y,\spinc)$, $\HFp(Y,\spinc)$, $\HFm(Y,\spinc)$ and
$\HFa(Y,\spinc)$. Furthermore, this induces actions 
of the exterior algebra
$\Wedge^* (H_1(Y;\Z)/\Tors)\subset \Wedge^*(H^1(\Omega(\Ta,\Tb));\Z)$ on each group.
\end{prop}

To define this action, let $\zeta\in Z^1(\Omega(\Ta,\Tb);\Z)$ be a
one-cocycle in the space of paths connecting $\Ta$ to $\Tb$. We define
a map $$A_\zeta\colon \CFinf(Y,\spinc)\longrightarrow
\CFinf(Y,\spinc)$$ which lowers degree by one, by the formula
$$A_{\zeta}([\x,i])=\sum_{\y\in {\mathfrak S}}
\sum_{\{\phi\in\pi_2(\x,\y)|\Mas(\phi)=1\}}
\zeta(\phi)\cdot\left(\#\UnparModFlow(\phi)\right)[\y,i-n_z(\phi)].$$
By $\zeta(\phi)$, we mean the following. Choose any representative $u$
for the homotopy class $\phi$, and view it as an arc in
$\Omega(\Ta,\Tb)$ which connects the constant paths $\x$ and $\y$. 
If we choose a different representative for the same homotopy class,
then the corresponding paths will be homotopic (as arcs in
$\Omega(\Ta,\Tb)$ connecting $\x$ to $\y$), so the evaluation of
$\zeta$ is independent of the particular choice (since $\zeta$ is a cocycle). 

We turn to the proof of Proposition~\ref{prop:Action}, which we break into several lemmas.

\begin{lemma}
\label{lemma:ChainMap}
$A_\zeta$ is a chain map. 
\end{lemma}

\begin{proof}
This is a variant on the usual proof that $\partial^2=0$. Suppose that
$\phi\in\pi_2(\x,\w)$ satisfies
$\Mas(\phi)=2$, and let $k=n_z(\phi)$. Then, since $\zeta(\phi_1\csum\phi_2)=\zeta(\phi_1)+\zeta(\phi_2)$
(since $\zeta$ is a cocycle), we get that
\begin{eqnarray*}
0 &=& \zeta(\phi) \cdot \left(\#(\text{ends of $\UnparModFlow(\phi)$})\right) \\
&=& \sum_{\{\phi_1,\phi_2 | \phi=\phi_1 \csum \phi_1,
\Mas(\phi_1)=\Mas(\phi_2)=1\}} \left(\zeta(\phi_1)+\zeta(\phi_2)\right) \left(\#\UnparModFlow(\phi_1)\right)\cm
\left(\#\UnparModFlow(\phi_2)\right).
\end{eqnarray*}
(Note that boundary degenerations do not contribute to the above sum,
as in the proof that $\partial^2=0$.) 
Summing over all $\phi\in\pi_2(\x,\w)$ with
$n_z(\phi)=k$ and $\Mas(\phi)=2$, we get the
$[\w,i-k]$-coefficient of 
$\left(\partial \circ A_\zeta +
A_\zeta \circ \partial\right) [\x,i]$. 
\end{proof}

\begin{lemma}
\label{lemma:NullHomotopic}
If $\zeta$ is a coboundary, then $A_\zeta$ is chain homotopic to zero.
\end{lemma}

\begin{proof}
If $\zeta$ is a coboundary, then there is a zero-cochain $B$ (a possibly discontinuous map from
$\Omega(\Ta,\Tb)$ to $\Z$) with the property that if $\gamma$ is an
arc in $\Omega(\Ta,\Tb)$ (a one-simplex), then 
$\zeta(\gamma)=B(\gamma(0))-B(\gamma(1))$. Let
$$H([\x,i])= B(\x)[\x,i],$$
where the evaluation of $B$ on $\x$ is performed by viewing the latter
as a constant path from $\Ta$ to $\Tb$.
Then, it follows from the definitions that
$$A_\zeta=\partial\circ H - H\circ\partial.$$
\end{proof}

{\noindent{\bf Proof of Proposition~\ref{prop:Action}.}}
Together, Lemmas~\ref{lemma:ChainMap} and \ref{lemma:NullHomotopic}
show that the $A_\zeta$ descends to a well-defined action of
$H^1(\Omega(\Ta,\Tb))$ on $\HFinfty$. To see that the action descends
to the exterior algebra, we must verify that the composite
$A_\zeta\circ A_\zeta=0$ in homology.

To see this,
we think of $A_\zeta$ using codimension one constraints. 
Specifically,
we begin with a map $f\colon \Omega(\Ta,\Tb)\longrightarrow S^1$
representing $\zeta$. Given a
generic point $p\in S^1$, and we let $V=f^{-1}(p)$, so that the
action of $\zeta$ is given by
$$A_\zeta([\x,i])=\sum_\y
\sum_{\{\phi\in\pi_2(\x,\y)\big|\Mas(\phi)=1\}}
a(\zeta,\phi)[\y,i-n_z(\phi)],$$
where
$$
a(\zeta,\phi)=
\# \{u\in\Mod(\phi)\big| u([0,1]\times \{0\})\in V\}.
$$

Fix a homotopy class $\phi\in\pi_2(\x,\w)$ with $\Mas(\phi)=2$. We
consider the one-manifold 
$$
\Xi=\Big\{s\in [0,\infty), u\in \ModFlow(\phi)\Big| u([0,1]\times \{s\})\in
V,~~~~~ u([0,1]\times \{-s\})\in V'\Big\}.$$ where $V$, $V'$ are the
preimages of $p$ and $p'$ under $f$. Choosing generic $p$ and $p'$,
the
one-manifold $\Xi$ has no boundary at $s=0$. The ends as $s\goesto
\infty$ (disregarding boundary degenerations, which do not contribute
algebraically), are modeled on $$\Big\{u_1\in\ModFlow(\phi_1)\Big|
u_1([0,1]\times \{0\})\in V\Big\}\times
\Big\{u_2\in\ModFlow(\phi_2)\Big| u_2([0,1]\times \{0\})\in V'\Big\},$$
where $\phi=\phi_1\csum\phi_2$. 
On the one hand, the number of points, counted with sign, must vanish;
on the other hand, it is the $[\w,i-n_z(\phi)]$ coefficient of
$A_\zeta\circ A_\zeta$. It follows that the action of
$H_1(Y;\Z)/\Tors$ on $\HFinf(Y,\spinc)$ descends to an action of the exterior algebra

The chain map $A_\zeta$, and the chain homotopy from Lemma~\ref{lemma:NullHomotopic}
preserve $\CFm(Y,\spinc)$, so it is follows
that $A_\zeta$ induces actions on $\HFp$ and $\HFm$. The action on
$\HFa$ is defined in an analogous manner, as well.
\qed
\vskip.2cm

Although this gives an action of all of $H^1(\Omega(\Ta,\Tb);\Z)$ on
the Floer homologies, the interesting part of the action is induced
from $H_1(Y;\Z)$ (c.f. the isomorphism from
Equation~\eqref{eq:CohomConfig}): it is a straightforward verification
that when $g>2$, the additional $\Z$ summand acts trivially on all the
Floer homology groups.

%
%

\begin{remark}
\label{rmk:GeoRep}
A geometric realization of the action of
$H_1(Y;\Z)/\Tors$ can be given as follows. 
Let $\gamma\in \Sigma$ be a curve which misses the
intersection points between the $\alpha_i$ and $\beta_j$, and let 
$[\gamma]$ be its induced homology class in $H_1(Y;\Z)$. Then,
$$A_{[\gamma]}([\x,i])=\sum_{\y}\sum_{\{\phi\in\pi_2(\x,\y)|\Mas(\phi)=1\}}
a(\gamma,\phi)\cm [\y,j-n_z(\phi)],$$ where $$ a(\gamma,\phi)=\#\{u\in
\ModFlow(\phi)\big| u(1\times 0)\in (\gamma\times
\Sym^{g-1}(\Sigma))\cap \Ta\} $$ or, equivalently, $a(\gamma,\phi)$ is
the product of $\#\UnparModFlow(\phi)$ with the intersection number in
$\Ta$
between the codimension one submanifold
$\left(\gamma\cap \Sym^{g-1}(\Sigma)\right)\cap \Ta$ and 
the curve in $\Ta$ obtained by restricting $u$ to
$u(\{1\}\times \R)$, where $u$ is any representative of $\phi$.
\end{remark}

\section{Special Heegaard moves}
\label{sec:Special}

In Section~\ref{sec:DefHF}, for three-manifolds with $b_1(Y)>0$, we
required that the Heegaard diagram satisfy additional admissibility
hypotheses.  It is the purpose of this section is to construct such
Heegaard diagrams, and indeed to show that any admissible diagrams are
isotopic through such diagrams.

To this end, we will be considering certain special isotopies. Let
$\gamma$ be an oriented simple closed curve in $\Sigma$. By {\em
winding along $\gamma$} we mean the diffeomorphism of $\Sigma$
obtained by integrating a vector field $X$ supported in a tubular
neighborhood of $\gamma$, where it satisfies the property that
$d\theta(X)>0$, with respect to a coordinate system $(t,\theta)\in
(-\epsilon,\epsilon)\times S^1$ in the tubular neighborhood of
$\gamma=\{0\}\times S^1$.

\begin{figure}
	\mbox{\vbox{\epsfbox{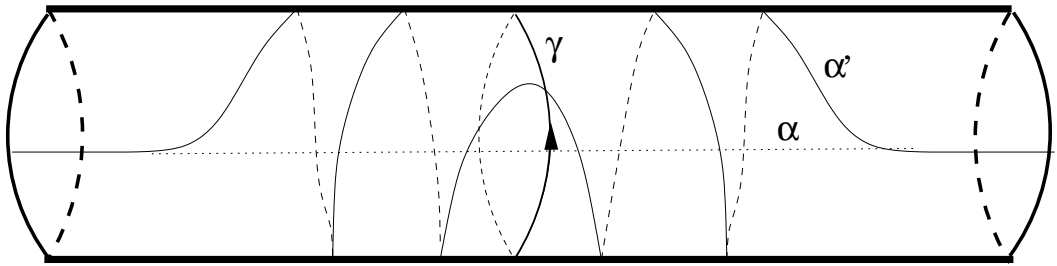}}}
	\caption{\label{fig:Winding}
	{\bf Winding transverse to $\alpha$.}
	We have pictured a cylindrical subregion of $\Sigma$, where
        $\alpha$ is the horizontal curve, which we wind twice along the
	vertical circle $\gamma$ (in the direction indicated) to obtain $\alpha'$.}
\end{figure}

Choose a curve $\gamma$ transverse to $\alpha_1$, meeting it in a
single transverse point, and which is disjoint from the other
$\alpha_i$ for $i\neq 1$, and suppose that $\phi$ is some
diffeomorphism which winds along $\gamma$. Suppose, moreover, that
$\phi(\alpha_1)$ meets $\alpha_1$ transversally in the neighborhood of
$\gamma$, meeting it there in $2k$ points. Then, we say that $\phi$
winds $\alpha_1$ along $\gamma$ $k$ times.  See Figure~\ref{fig:Winding}. 

We have the following notion:

\begin{defn}
Fix a $\SpinC$ structure $\spinc$ over $Y$. A pointed Heegaard diagram
$$(\Sigma,\alphas,\betas,z)$$ is called $\spinc$-realized if there is a
point $\x\in\Ta\cap\Tb$ with the property that $s_z(\x)=\spinc$.
\end{defn}

\begin{lemma}
\label{lemma:spincRealized}
Fix $Y$ and a $\SpinC$ structure $\spinc$. Then, $Y$ admits an 
$\spinc$-realized pointed Heegaard diagram.
\end{lemma}

\begin{proof}
Begin with any Heegaard diagram $(\Sigma,\alphas,\betas)$ for $Y$ and
let $\gammas$ be a collection of pairwise disjoint curves which are
dual to the $\alphas$, in the sense that for all $i$ and $j$,
$$\#(\alpha_i\cap\gamma_j)=\delta_{i,j}$$ (the right hand side is
Kronecker delta, and the left hand side denotes both the geometric and
algebraic intersection numbers of the curves).  By isotoping the
$\betas$ if necessary, we can arrange that $\Tb\cap\Tc\neq \emptyset$.
Choose a basepoint $z$ distinct from $\alphas$, $\betas$, and
$\gammas$ (indeed, choose $z$ to be disjoint from the neighborhood of
the $\gammas$ where the winding is performed).

Let $\x=\{x_1,...,x_g\}\in\Tb\cap \Tc$, labelled so that $x_i\in
\beta_i\cap\gamma_i$ for $i=1,...,g$. 
Each time we wind $\alpha_i$ along
$\gamma_i$. we create a new pair of intersection points near $x_i$
between
$\beta_i$
and the
new copy of $\alpha_i$. Winding along each $\gamma_i$ $k$ times, then, we can label
these intersection points $x_i^\pm(1), x_i^\pm(2),...,x_i^\pm(k)$
(ordered in decreasing order of their distance to $x_i$, and with sign
distinguishing which side of $\gamma_i$ -- in its tubular neighborhood
-- they lie in). Thus, we have induced intersection points
$$\x(i_1,...,i_g)=\{x_1^+(i_1),...,x_g^+(i_g)\}\in \Ta'\cap \Tb$$
labeled by $i_1,...,i_g \in 1,...,k$. 
Note that with our conventions, the short arc in $\alpha_i$
connecting $x_i(k)$ to $x_i(k+1)$, followed by the short arc in $\beta_i$
with the same endpoints, is homologous to $\gamma_i$ in $\Sigma$.

No matter how many times we wind $\alpha_i$ along
$\gamma_i$, the $\SpinC$ structure of the 
farthest intersection point $\x(1,...,1)$ remains fixed (this
is clear from the definition of $s_z(\x)$: the winding isotopy induces
an isotopy between the induced non-vanishing vector fields induced
over $Y$). Moreover, by the definition of the difference map $\epsilon$ 
introduced in Subsection~\ref{subsec:Disks}, together with Lemma~\ref{lemma:VarySpinC},
we have that
$$\spinc_z(\x(i_1,...,i_g))-\spinc_z(\x(j_1,...,j_g))=
\Big((i_1-j_1)\PD[\gamma_1]+...+(i_g-j_g)\PD[\gamma_g]\Big). $$

Thus, we can find Heegaard diagrams which realize the $\SpinC$
structures which differ from some fixed $\SpinC$ structure $\spinc_0$
by non-positive multiples of the $[\gamma_1]$,...,$[\gamma_g]$.
Moreover,  if we choose parallel copies
$\{\gamma^-_1,...,\gamma^-_g\}$ of the $\gammas$, only with the
opposite orientations, and wind along those in addition, we can
realize all $\SpinC$ structures which differ from $\spinc_0$ by
arbitrary multiples the $[\gamma_1]$,...,$[\gamma_g]$.
Now, it is easy to see that the group $H^2(Y;\Z)$
is generated by the Poincar\'e duals of the $\gammas$. Hence, we can
realize all $\SpinC$ structures.
\end{proof}

Winding can be used also to arrange for strong admissibility. For this,
it is useful to have the following:

\begin{defn}
\label{def:RenormalizedPeriodicDomain}
An {\em $\spinc$-renormalized periodic domain} is a two-chain
$\RenPerDom=\sum a_i\cald_i$ in
$\Sigma$ whose boundary is a sum of the curves $\alphas$ and
$\betas$ (with multiplicities),  satisfying the additional property that
$$n_z(\RenPerDom)=-\frac{\langle c_1(\spinc),H(\RenPerDom)\rangle}{2}.$$
\end{defn}

Of course, the group of  $\spinc$-renormalized periodic domains is isomorphic to
the group of periodic domains. (The periodic domain $\PerDom$ gives
rise to the renormalized periodic domain 
$\PerDom-\frac{\langle c_1(\spinc),H(\PerDom)\rangle}{2} [\Sigma]$.)

\begin{lemma}
\label{lemma:StronglyAdmissible}
Fix $Y$ and a $\SpinC$ structure $\spinc$. Then, $Y$ admits a strongly 
$\spinc$-admissible pointed Heegaard diagram.
\end{lemma}

\begin{proof}
In view of Lemma~\ref{lemma:spincRealized}, we can start with an
$\spinc$-realized Heegaard diagram.  We will show that after winding
the $\alphas$ sufficiently many times along curves $\gammas$ as in the
proof of the previous lemma, we obtain a pointed Heegaard diagram for
which each renormalized $\spinc$-periodic domain has both positive and
negative coefficients. Such a Heegaard diagram is strongly
$\spinc$-admissible.

We find it convenient to use rational coefficients for our periodic
domains.  Write $b=b_1(Y)$, and choose a basis
$\{\RenPerDom_1,...,\RenPerDom_b\}$ for the $\Q$-vector space of
renormalized periodic domains. Note that a renormalized periodic
domain $\RenPerDom$ is uniquely determined by the part of its boundary
which is spanned by $\{[\alpha_1],...,[\alpha_g]\}$; i.e. the map
which associates to a periodic domain the $\alphas$-coefficients of
its boundary gives an injection of vector spaces.  After a change of
basis of the $\{\RenPerDom_i\}$ and reordering the $\alphas$, we can
assume that for all $i=1,...,b$, $$ \partial \RenPerDom_i = \alpha_i +
\sum_{j=b+1}^g a_{i,j} \alpha_i + b_{i,j} \beta_j,$$ where $a_{i,j}, b_{i,j}\in\Q$.

For each $i=1,...,b$ 
choose points $w_i\in \gamma_i$ which are not contained in any of the
$\alphas$ or $\betas$.
Let 
$$c_i=\max_{j=1,...,b} |n_{w_i}(\RenPerDom_j)|,$$
and then choose some integer $N$ so that 
$$N> b\cdot \left(\max_{i=1,...,b} c_i\right).$$

Choose parallel copies $\gamma_i^-$ of the $\gamma_i$ for $i=1,...,b$, and let 
$\{\RenPerDom_1',...,\RenPerDom_b'\}$ be the new periodic domains,
obtained after winding the curves  $\{\alpha_1,...,\alpha_b\}$ $N$ times along the
$\{\gamma_1,...,\gamma_b\}$ and $N$ times in the opposite direction
along the 
$\{\gamma_1^-,...,\gamma_b^-\}$. Note that
\begin{eqnarray*}
n_{w_i}(\RenPerDom_i') &=&
n_{w_i}(\RenPerDom_i)+N \\
&>& n_{w_i}(\RenPerDom_i)+ b \cm c_i \\
&\geq& (b-1)\max_{j=1,...,i-1,i+1,...,b} |n_{w_i}(\RenPerDom_j)| \\
&=& (b-1)\max_{j=1,...,i-1,i+1,...,b} |n_{w_i}(\RenPerDom_j')|.
\end{eqnarray*} 
In a similar manner, we see that
$$n_{w_i^-}(\RenPerDom_i')<-(b-1) \max_{j=1,...,i-1,i+1,...,b} |n_{w_i^-}(\RenPerDom_j')|.$$
It is a straightforward matter, then, to verify that for any linear
combination of the $\RenPerDom_i'$, one can find some point $w$ for which
$n_{w}$ is positive, and another $w'$ for which $n_{w'}$ is negative.
\end{proof}

Indeed,  an elaboration of this argument gives the following refinement.
But first, we give a definition.

Recall that a (generic) pointed isotopy between two Heegaard diagrams
can subdivided into a sequence of isotopies where, at seach stage,
there is one curve $\alpha_i\in\alphas$ and one curve
$\beta_j\in\betas$ whose number of intersection points either
increases by two (pair creation) or drops by two (pair annihilation),
while all other curves $\alpha_k$ and $\beta_\ell$ when $(i,j)\neq
(k,\ell)$ remain transverse throughout the isotopy.

\begin{defn}
Two strongly $\spinc$-isotopic pointed Heegaard diagrams are said to
be {\em strongly $\spinc$-isotopic} if all the intermediate Heegaard
diagrams in the isotopy are also strongly $\spinc$-isotopic.
\end{defn}

\begin{lemma}
\label{lemma:StronglyAdmissibleIsotopy}
Suppose that two strongly $\spinc$-admissible pointed Heegaard
diagrams are isotopic (via an isotopy which does not cross the
basepoint $z$), then they are strongly $\spinc$-isotopic.
\end{lemma}

\begin{proof}
First, note that if $(\Sigma,\alphas,\betas,z)$ is strongly
$\spinc$-admissible, then if we choose curves along which to wind the
$\alphas$ (disjoint from the basepoint $z$), then the winding gives
an isotopy through strongly $\spinc$-admissible pointed Heegaard
diagrams. The reason for this is that, in the complement of a small
neighborhood of the winding region, the various renormalized periodic domains
remain unchanged; thus, if some renormalized periodic domain 
has positive coefficients, then it  retains this property as it undergoes
winding.

Thus, it suffices to show that if two Heegaard diagrams are isotopic
(via an isotopy which we can assume without loss of generality takes
place only among the $\betas$ -- taking $\betas$ to $\betas'$), then
if we wind their $\alpha$-curves simultaneously along some collection
of $\gammas$ to obtain $\alphas'$, then the pointed Heegaard diagrams
$(\Sigma,\alphas',\betas,z)$ and $(\Sigma,\alphas',\betas',z)$ are
isotopic through strongly $\spinc$-admissible Heegaard diagrams. To
see this, we choose $\gamma_i$ curves and their translates
$\gamma_i^-$ as in the proof of
Lemma~\ref{lemma:StronglyAdmissible}. Now, we choose constants $$c_i =
\sup_{t\in[0,1]} \max_{i=1,...,b} |n_{w_i}(\RenPerDom_i(t))|,$$
where we think of $t\in[0,1]$ as the parameter in some isotopy taking
$\betas$ to $\betas'$, and $\RenPerDom_i(t)$ is the corresponding
one-parameter family of renormalized periodic domains. (Strictly
speaking, the point $w_i$ generically lies on the
translates of the $\beta_i$ for finitely many $t$, so that for those
values of $t$, the multiplicity $n_{w_i}(\RenPerDom_i(t))$ does not
make sense as we have defined it; for those values of $t$, we use a
small perturbation $w_i'\in\gamma_i$ of the basepoint $w_i$.)
Using these constants $c_i$ as in the proof of
Lemma~\ref{lemma:StronglyAdmissible}, the present lemma follows.
\end{proof}

\begin{remark}
\label{rmk:SpinCRealized}
Note that this lemma also proves that any two isotopic $\spinc$-realized
pointed Heegaard diagrams are isotopic through $\spinc$-realized
Heegaard diagrams.
\end{remark}

Suppose that two weakly $\spinc$-admissible pointed Heegaard diagrams
are isotopic (via an isotopy which does not cross the basepoint $z$),
then they are said to be {\em weakly $\spinc$-isotopic}.

\begin{lemma}
\label{lemma:WeaklyAdmissibleIsotopy}
Any weakly $\spinc$-admissible Heegaard diagram is weakly
$\spinc$-isotopic to a strongly ${\spinc}$-admissible Heegaard diagram.
\end{lemma}

\begin{proof}
It is easy to see that if we take a weakly $\spinc$-admissible
Heegaard diagram, then Lemma~\ref{lemma:StronglyAdmissible} provides
an isotopy to a strongly $\spinc$-admissible Heegaard diagram, and
that the given isotopy is an isotopy through weakly $\spinc$-admissible
Heegaard diagrams.
\end{proof}

\section{Independence of complex structures}
\label{section:CxStructure}

In Section~\ref{sec:DefHF}, we defined various chain complexes, whose
definition required a Heegaard diagram (satisfying appropriate
admissibility hypotheses), and analytical choices -- a complex
structure on $\Sigma$, and a one-parameter perturbation $J$. Our aim
here is to prove that the homology groups of the chain complexes are
independent of the latter choices, and hence define an invariant
of pointed Heegaard diagrams. More precisely, we prove the following:

\begin{theorem}
\label{thm:IndepCxStruct}
Let $Y$ be a closed, oriented three-manifold equipped with a $\SpinC$
structure $\spinc$, and let $(\Sigma,\alphas,\betas,z)$ be a strongly
$\spinc$-admissible Heegaard diagram, endowed with an equivalence
class $\orient$ of coherent orientation system. Then, the homology
groups $\HFa(\alphas,\betas,\spinc,\orient)$, $HF^\pm
(\alphas,\betas,\spinc,\orient)$ and $\HFinf
(\alphas,\betas,\spinc,\orient)$, thought of as modules over
$\Z[U]\otimes_\Z\Wedge^* H_1(Y;\Z)/\Tors$, are independent of the
choice of complex structure $\sj$ on $\Sigma$ and the path $J_s$.
\end{theorem}

We prove first Theorem~\ref{thm:IndepCxStruct} in the case where $Y$ is a rational homology
three-sphere:

\vskip.3cm
\noindent{\bf{Proof of Theorem~\ref{thm:IndepCxStruct} when $b_1(Y)=0$}.}
First, we argue that if we fix $\sj$ over $\Sigma$, the Floer homology
groups are independent of the choice of $(\eta,\sj,V)$-nearly symmetric
path $J_s$ (in ${\mathcal U}$). Suppose we have two paths $J_s(0)$ and
$J_s(1)$ in ${\mathcal U}$. Since ${\mathcal U}$ is simply-connected,
we can connect them by a two-parameter family $J\colon [0,1]\times
[0,1]\longrightarrow {\mathcal U}$, thought of as a one-parameter
family of paths indexed by $t\in [0,1]$, writing $J_s(t)$ for the path
obtained by fixing $t$. In fact, we can arrange that $J_s(t)$ is
independent of $t$ for $t$ near $0$ and $1$, so that $J_s(t)$ can be
naturally extended to all $t\in\R$.  Then (as is familiar in Floer
theory) we can define an associated chain map 
$$
\Phi^\infty_{J_{s,t}}\colon
(\CFinf(\alphas,\betas,\spinc),\partial^\infty_{J_s(0)})\longrightarrow
(\CFinf(\alphas,\betas,\spinc),\partial^\infty_{J_s(1)}),
$$ by 
\begin{equation}
\label{eq:DefPhi}
\Phi^\infty_{J_{s,t}}[\x,i]=\sum_{\y}\sum_{\{\phi\in\pi_2(\x,\y)|\Mas(\phi)=0\}}
\#\left(\Mod_{J_{s,t}}(\phi)\right)\cdot [\y,i-n_z(\phi)], 
\end{equation}
where $\Mod_{J_{s,t}}(\phi)$ denotes holomorphic disks with a time-dependent
complex structure on the target,
i.e. $\Mod_{J_{s,t}}(\phi)$ consists of maps
$$\left\{u\colon \Strip\cong [0,1]\times \R\longrightarrow \Sym^{g}(\Sigma)
\Bigg|
\begin{array}{l}
	\frac{du}{ds}+J(s,t)\frac{du}{dt}=0, \\
    u(\{1\}\times \R)\subset \Ta, u(\{0\}\times \R)\subset \Tb \\
    \lim_{t\goesto -\infty}u(s+it)=\x, 
    \lim_{t\goesto +\infty}u(s+it)=\y 
    \end{array}
    \right\},$$
which represent the homotopy class $\phi$.
Note that the argument of Lemma~\ref{lemma:EnergyBound} still
goes over to give an energy bound on this moduli space.

The usual arguments from Floer theory then apply to
show that $\Phi^\infty_{J_{s,t}}$ is a chain map which induces an isomorphism in
homology. We outline these briefly.

The transversality theorem (Theorem~\ref{thm:Transversality}) shows
that for a generic path $J_s(t)$,
the zero-dimensional components of the moduli spaces $\Mod_{J_{s,t}}(\x,\y)$ are
smoothly cut out and compact, as in the proof of
Theorem~\ref{thm:StructMod}. Thus, the map $\Phi^\infty_{J_{s,t}}$ is
well-defined. To show that it is a chain map, we consider the ends of
the one-dimensional moduli spaces $\Mod_{J_{s,t}}(\psi)$ with
$\Mas(\psi)=1$.  We claim that the only ends of these moduli spaces
correspond to products $$\left(\coprod_{\phi\csum\phi'=\psi}
\Mod_{J_{s,t}}(\phi)\times
\UnparModFlow_{J_s(1)}(\phi')\right)
\coprod
\left(\coprod_{\phi'\csum\phi=\psi}\UnparModFlow_{J_s(0)}(\phi') \times \Mod_{J_{s,t}}(\phi)\right),$$ 
where the $\phi$ homotopy classes all have 
$\Mas(\phi)=0$ and $\phi'$ have $\Mas(\phi')=1$.
Counted with sign, these represent the coefficients of
$\partial^\infty_{J_s(1)} \circ \Phi^\infty_{J_{s,t}}-
\Phi^\infty_{J_{s,t}}\circ
\partial^\infty_{J_s(0)}$. This follows from Gromov's compactness, together with the
observation that there can be no spheres bubbling off 
or boundary degenerations, as both carry
Maslov index at least two. Hence, $\Phi^\infty_{J_{s,t}}$ is a chain map.

To see that $\Phi^\infty_{J_{s,t}}$ induces an isomorphism in
homology, we show that the composite $\Phi^\infty_{J_{s,t}}\circ
\Phi^\infty_{J_{s,1-t}}$ is chain homotopic to the identity map. The
chain homotopy is constructed using a homotopy $J_{s,t,\tau}$ between
two two-parameter families of complex structures (once again, we let
$J_{s,t}(\tau)$ denote the two-parameter family in $s$ and $t$, with
$\tau\in[0,1]$ fixed); i.e. $J_{s,t}(0)$ is the family of complex
structures obtained by juxtaposing $J_{s,t}$ with $J_{s,1-t}$, while
$J_{s,t}(1)=J_s(0)$ is independent of $t$. We can define a moduli space
$$\ModFlow_{J_{s,t,\tau}}(\phi)=
\bigcup_{\tau \in[0,1]} \ModFlow_{J_{s,t}(\tau)}(\phi),$$
for each fixed homotopy class $\phi\in\pi_2(\x,\y)$.
For generic $J_{s,t,\tau}$, this is a manifold of dimension $\Mas(\phi)+1$.
We define a map
$$H^\infty_{J_{s,t,\tau}}([\x,i]) =
\sum_\y\sum_{\{\phi\in\pi_2(\x,\y)|\Mas(\phi)=-1\}}
\# \left(\ModFlow_{J_{s,t,\tau}}(\phi)\right)\cdot[\y,i-n_z(\phi)].$$
Note that $H^\infty_{J_{s,t,\tau}}$ has degree $-1$. To see, then that $H^\infty_{J_{s,t,\tau}}$ is the chain
homotopy between $\Phi^\infty_{J_t}\circ \Phi^\infty_{J_{1-t}}$ and the identity
map, we consider ends of the moduli spaces $\ModFlow_{J_{s,t,\tau}}(\psi)$
with $\Mas(\psi)=0$.
These have three kinds of ends: 
the $\tau=0$ end, which
corresponds to the composite 
$\Phi^\infty_{J_t}\circ \Phi^\infty_{J_{1-t}}$, the
$\tau=1$ end, corresponding to the identity map, and those ends which
correspond to splittings
$$\left(\coprod_{\phi\csum\phi'=\psi} \Mod_{J_{s,t,\tau}}(\phi)\times
\UnparModFlow_{J_s(0)}(\phi')\right)
\coprod
\left(\coprod_{\phi'\csum\phi=\psi}\UnparModFlow_{J_s(0)}(\phi') \times \Mod_{J_{s,t,\tau}}(\phi)\right),$$ 
where $\Mas(\phi)=-1$, and $\Mas(\phi')=1$,
i.e. this corresponds to $\partial^\infty _{J_s(0)}\circ  
H^\infty_{J_{s,t,\tau}} -
H^\infty_{J_{s,t,\tau}} \circ \partial^\infty_{J_s(0)}$. Once
again, there are no other possible bubbles generically.

Similar chain maps can be defined on $\HFa$, $\HFp$ and $\HFm$ as
well. For $\HFa$, the corresponding map ${\widehat \Phi}_{J_{s,t}}$
counts $\phi$ only if $n_z(\phi)=0$.  For $\HFp$ and $\HFm$, we can
let $\Phi^+_{J_{s,t}}$ and $\Phi^-_{J_{s,t}}$ be the induced maps from
$\Phi^\infty_{J_{s,t}}$, because $\phi$ admits no holomorphic
representative if $n_z(\phi)<0$. It is clear from their definition that
$\Phi^\infty$, $\Phi^+$, and $\Phi^-$ commute with the corresponding
$U$-actions.

As a consequence, we see that the Floer homologies for fixed
$\sj$ are independent of the choice of generic 
$(\eta,\sj,V)$-nearly symmetric path. Thus, it follows also that
the homology groups do not depend on the $\eta$ and $V$.

Next, we see that the chain complex is independent of the complex
structure $\sj$ on the Riemann surface. To this end, we observe that
the chain complexes remain unchanged under small perturbations of the
path of almost-complex structures $J_s$, provided that we still have
{\em a priori} energy bounds after the perturbation. Furthermore, we
can approximate a $\sj$-nearly-symmetric path $J_s$ by
$\sj'$-nearly-symmetric paths $J_s'$, with $\sj'$ close to $\sj$. This
shows that the Floer homology is also independent of the choice of
$\sj$, since the space of allowed complex structures over $\Sigma$ is
connected (see~Lemma~\ref{lemma:Genericj}), as it is obtained from the
space of all complex structures by removing a codimension two subset.
\qed
\vskip.3cm

We turn attention to Theorem~\ref{thm:IndepCxStruct} in the case where $b_1(Y)>0$.
We assume that the pointed Heegaard diagram $(\Sigma,\alphas,\betas,z)$
is strongly $\spinc$-admissible.

As in the proof when $b_1(Y)=0$, consider a fixed complex structure
$\sj$ over $\Sigma$, and let $J_{s,t}$ be a one-parameter family
$(\eta,\sj,V)$-nearly symmetric paths in ${\mathcal U}$. Note that the
non-negativity result of Lemma~\ref{lemma:NonNegativity} applied to
the parameterized moduli spaces $\ModFlow_{J_{s,t}}(\phi)$, together
with admissibility and Lemma~\ref{lemma:StrongFiniteness}, ensure that
the sum defining $\Phi^{\infty}_{J_{s,t}}$ from
Equation~\eqref{eq:DefPhi} is a finite sum. 
To see that the map respects the module
structure, we have the following:

\begin{lemma}
\label{lemma:ActionCommutes}
For any $\zeta\in H_1(Y,\Z)/\Tors$, 
$$A_\zeta\circ (\Phi^\infty_{J_{s,t}})=(\Phi^\infty_{J_{s,t}})\circ A_\zeta $$
as a map from $H_*(\CFinfty(\Ta,\Tb),\partial^\infty_{J_s(0)}) \longrightarrow 
H_*(\CFinfty(\Ta,\Tb),\partial^\infty_{J_s(1)}) $
\end{lemma}

\begin{proof}
Let $V$ be a codimension one constraint
in $\Omega(\Ta,\Tb)$ representing
the class $\zeta\in H^1(\Omega(\Ta,\Tb);\Z)$, chosen to miss all the
constant paths (corresponding to the intersection
points $\Ta\cap\Tb$). 

Consider the map $$h\colon \CFinfty(\Ta,\Tb)\longrightarrow
\CFinfty(\Ta,\Tb),$$ 
defined by $$h([\x,i])=\sum_\y \sum_{\{\phi\in
\pi_2(\x,\y)| \Mas(\phi)=0\}} \# \{(r,u)\in \Mod_{J_{s,t}}(\phi)|
u([0,1]\times r)\in V\} [\y,i-n_z(\phi)].$$ We claim that 
\begin{equation}
\label{eq:HomotopyFormula}
A_\zeta\circ \Phi^\infty_{J_{s,t}} - \Phi^\infty_{J_{s,t}}\circ
A_\zeta = \partial_{J_s(1)} \circ h - h \circ \partial_{J_s(0)}.
\end{equation}
This follows by considering the ends of the one-dimensional moduli
spaces $$\Xi=\{(r,u)\in \R\times \ModFlow_{J_{s,t}}(\psi)\big|
u([0,1]\times \{r\})\in V\}$$ where $\Mas(\psi)=1$. The ends where
$r\goesto \pm \infty$ correspond to the commutator of $A_\zeta$ and
$\Phi^\infty$, while the ends where the maps $u\in\ModFlow_{J_{s,t}}$ 
bubble off correspond to the commutator of $h$ with the
corresponding boundary maps. 

Equation~\eqref{eq:HomotopyFormula}, of course, says that $A_\zeta$
commutes with $\Phi^\infty_{J_{s,t}}$, on the level of homology.
\end{proof}

\vskip.2cm
\noindent{\bf{Proof of Theorem~\ref{thm:IndepCxStruct} when $b_1(Y)>0$}.}
The proof proceeds exactly as in the case where
$b_1(Y)=0$. Lemma~\ref{lemma:ActionCommutes} is used to prove that the
induced isomorphisms of Floer homologies are $\Wedge^*
H_1(Y;\Z)/\Tors$-module isomorphisms.
\qed
\vskip.2cm

\section{First steps towards topological invariance}
\label{sec:Isotopies}

\subsection{Overview of topological invariance}

According to Theorem~\ref{thm:IndepCxStruct}, the Floer homology
groups as defined in Section~\ref{sec:DefHF} are invariants of
strongly $\spinc$-admissible pointed Heegaard diagram.

Loosely speaking, to show they are actually invariants of the
underlying three-manifold, we must show that they are invariant under
the three basic moves: isotopies, handleslides, and stabilizations. In
fact, the invariants we have are invariants of {\em pointed} Heegaard
diagrams, and we must restrict ourselves to moves which never cross
the basepoint $z$ (i.e. the isotopies do not cross $z$, and the pairs
of pants involved in the handleslides do not contain $z$ either). As
we shall see in Proposition~\ref{prop:PointedHeegaardMoves}, an
invariant of Heegaard diagrams which is invariant under these
restricted moves still gives a topological invariant, since we can
trade an isotopy which crosses the basepoint $z$ for a sequence of
handleslides which do not. 

In fact, when $b_1(Y)>0$, we allow ourselves only strongly
$\spinc$-admissible Heegaard diagrams for $Y$ and allow ourselves an
even more restricted set of Heegaard moves which connect such
diagrams. By adapting the arguments from Section~\ref{sec:Special}
(see Proposition~\ref{prop:VerySpecialMoves} below), we see that a
quantity associated to admissible Heegaard diagrams which is invariant
under only these kinds of Heegaard moves still gives a topological
invariant of the underlying three-manifold.

After establishing these topological preliminaries in
Section~\ref{subsec:PointedMoves}, we establish isotopy invariance
of the Floer homologies in Subsection~\ref{subsec:Isotopies} (and a
corresponding version for the weakly $\spinc$-admissible required to
define $\HFa$ and $\HFp$ in Subsection~\ref{subsec:Weakly}), returning
to handleslide invariance in Section~\ref{sec:HandleSlides} and
stabilization invariance in Section~\ref{sec:Stabilization}.

\subsection{Topological invariants and special Heegaard moves}
\label{subsec:PointedMoves}

A quantity associated to pointed Heegaard diagrams which is invariant
under pointed Heegaard moves is a three-manifold invariant, according
to the following:

\begin{prop}
\label{prop:PointedHeegaardMoves}
Any two Heegaard diagrams $(\Sigma,\alphas,\betas,z)$ and
$(\Sigma',\alphas',\betas',z')$ which specify the same three-manifold
are diffeomorphic after a finite sequence of pointed Heegaard moves
(i.e. Heegaard moves supported in the complement of the basepoint).
\end{prop}

\begin{proof}
Given a sequence of Heegaard moves, we can clearly introduce isotopies
as needed to arrange that no handleslides, only isotopies (of the
$\alphas$ or the $\betas$), cross the basepoint $z$.

Since the roles of $\alphas$, $\betas$ are symmetric, it suffices to
consider the case where the isotopy of $\beta_i$, say $\beta_1$,
crosses $z$. We denote the new isotopic curve by ${\overline\beta}_1$.
We claim that $\beta_1$ can be moved by a series of handle-slides and
isotopies to ${\overline\beta}_1$ all of which are supported in
$\Sigma-z$. This can be seen by first surgering out
$\beta_2,...,\beta_g$ to get a torus $T^2$ with $2g-2$ marked
points. Clearly, in $T^2-z$ the curves induced by $\beta_1$ and
${\overline\beta}_1$ are isotopic.  We can follow the isotopy by moves
$\Sigma-z$ where isotopies across the marked points $T^2-z$ are
replaced by handle-slides in $\Sigma-z$.  See
Figure~\ref{fig:largeisotopy} for the $g=3$ case.

\begin{figure}
\mbox{\vbox{\epsfbox{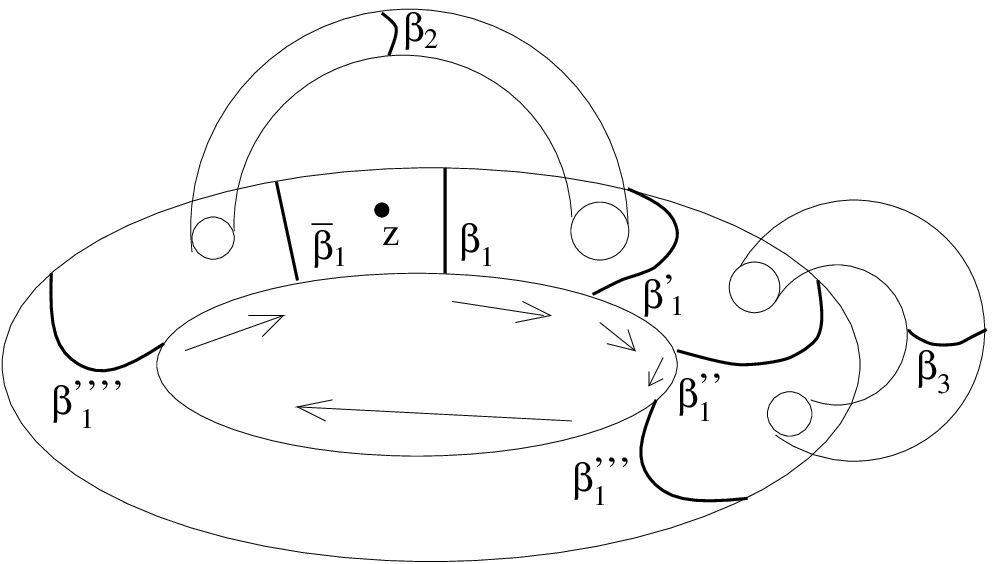}}}
\caption{\label{fig:largeisotopy} Moving $\beta_1$ to ${\overline \beta }_1$
in $\Sigma -z$}
\end{figure} 
\end{proof}

When $b_1(Y)>0$, in view of the additional admissibility hypotheses on
the Heegaard diagrams, we need the following refinement
Proposition~\ref{prop:PointedHeegaardMoves}:

\begin{prop}
\label{prop:VerySpecialMoves}
Let $Y$ be a three-manifold, equipped with a $\SpinC$ structure
$\spinc\in\SpinC(Y)$. Then, there is strongly $\spinc$-admissible
pointed Heegaard diagram for $Y$, and any two strongly
$\spinc$-admissible pointed Heegaard diagrams can be connected by a
finite sequence of pointed Heegaard moves where, at each intermediate
stage, the Heegaard diagrams are all strongly $\spinc$-admissible
pointed Heegaard diagrams.
\end{prop}

\begin{proof}
Existence is ensured by Lemma~\ref{lemma:StronglyAdmissible}.  Given
two strongly $\spinc$-admissible Heegaard diagrams,
Proposition~\ref{prop:PointedHeegaardMoves} gives a sequence of
pointed Heegaard moves which connect them.  Now, by introducing
additional isotopies as in
Lemma~\ref{lemma:StronglyAdmissibleIsotopy}, we can arrange for all
the isotopies to go through strongly $\spinc$-isotopic Heegaard
diagrams. 
\end{proof}

\subsection{Isotopy invariance}
\label{subsec:Isotopies}

Theorem~\ref{thm:IndepCxStruct} implies that the Floer homologies
remain unchanged under isotopies of the attaching circles which
preserve the transversality of the $\alphas$ and $\betas$. To show
isotopy invariance of the homology groups in general, we must allow a
larger class of isotopies which allows us to introduce new
intersections between the attaching circles, in a controlled manner.

\begin{figure}
\mbox{\vbox{\epsfbox{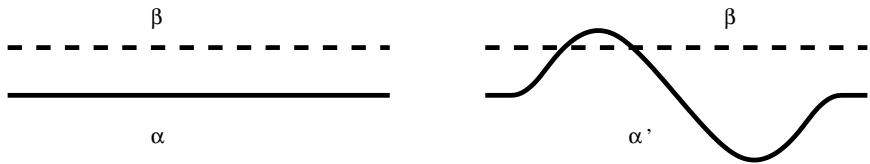}}}
\caption{\label{fig:HamIsot} Moving the curve $\alpha$ to $\alpha'$ 
by a 
Hamiltonian isotopy. In this 
manner, we can introduce a pair of canceling intersection points 
between $\alpha'$ and the curve $\beta$.}
\end{figure}

Such moves are provided by exact Hamiltonian motions of the attaching
circles in $\Sigma$. Recall that on a symplectic manifold, a
one-parameter family of real-valued functions $H_{t}$ naturally gives
rise to a unique one-parameter family of Hamiltonian vector fields
$X_{t}$, specified by $$\omega(X_{t},\cdot)=dH_{t}$$ where the
left-hand-side denotes the contraction of the symplectic form $\omega$
with the vector field.  A one-parameter family of diffeomorphisms
$\Psi_{t}$ is said to be an exact Hamiltonian isotopy if it is
obtained by integrating a Hamiltonian vector field, i.e. if
$\Psi_{0}=\Id$, and $$\frac{d\Psi_{t}}{dt}=X_{t}.$$ By taking a
positive bump function $h$ supported in a neighborhood of a point
which lies on $\alpha_{1}$, and letting $f\colon \R \longrightarrow
[0,1]$ be a non-negative smooth function whose support is $(0,1)$ we
can consider Hamiltonian $H_t=f(t)h$. The corresponding diffeomorphism
moves the curve $\alpha_{1}$ slightly (without moving any of the other
$\alpha$-curves). See Figure~\ref{fig:HamIsot} for an
illustration. (The picture takes place in a small Euclidean patch of
$\Sigma$, which meets two curves $\alpha$ and $\beta$; the isotopy is
used to displace $\alpha$, to give a new curve, $\alpha'$).

With these preliminaries in place we state the main result of this subsection:

\begin{theorem}
\label{thm:Isotopies}
If $Y$ is a three-manifold equipped with a $\SpinC$ structure
$\spinc$, and endowed with two isotopic $\spinc$-strongly admissible
Heegaard diagrams, then there is an identification between equivalence
classes of orientation systems for the two diagrams and corresponding
identification between the homology groups 
(thought of as $\Z[U]\otimes_\Z\Wedge^* H_1(Y;\Z)/\Tors$-modules)
\begin{eqnarray*}
\HFa(\Sigma,\alphas,\betas,\orient)&{\cong}&
\HFa(\Sigma',\alphas',\betas',\orient') \\
\HFpm(\Sigma,\alphas,\betas,\orient)&{\cong}&
\HFpm(\Sigma',\alphas',\betas',\orient') \\
\HFinf(\Sigma,\alphas,\betas,\orient)&{\cong}&
\HFinf(\Sigma',\alphas',\betas',\orient')
\end{eqnarray*}
\end{theorem}

\vskip.3cm
\noindent{\bf{Proof of Theorem~\ref{thm:Isotopies} when $b_1(Y)=0$.}}
First of all, isotopies which preserve the condition that the
$\alphas$ are transverse to the $\betas$ can be thought of as
variations in the metric (or equivalently, the complex structure $\sj$) 
on $\Sigma$. Hence, the invariance of the
Floer homology under such isotopies follows from
Theorem~\ref{thm:IndepCxStruct}. 

It suffices then to show that the homology remains unchanged when a
pair of canceling intersection points between, say, $\alpha_1$ and
$\beta_1$ is introduced. Such an isotopy can be realized by moving the
$\alphas$ by an exact Hamiltonian diffeomorphism of $\Sigma$. 
We assume that the exact Hamiltonian is supported on $\Sigma -
\alpha_2 -...- \alpha_g - z $
(i.e. that the corresponding vector fields $X_t$ are supported on the
subset), and that the Hamiltonian is supported in $t\in [0,1]$.  The
isotopy $\{\Psi_{t}\}$ induces an isotopy of $\Ta$.  We must show that
this isotopy of $\Ta$ induces a map on Floer homology, by imitating
the usual constructions from Lagrangian Floer theory.

The map is induced by counting points in the zero-dimensional
components of the moduli spaces with a time-dependent constraint.  To
be precise, here, we will be using a fixed one-parameter family $J$ of
nearly-symmetric complex structure. Recall that for this notion, we
needed to fix a collection of points $\{z_i\}_{i=1}^m\subset
\CurveComp$. For the present argument, it suffices to choose
$\{z_i\}_{i=1}^m=\{z\}$.  

Now, we have homotopy classes of disks $\pi_2^{\Psi_t}(\x,\y)$, where
$\x\in\Ta\cap\Tb$ and $\y\in
\Psi_1(\Ta)\cap\Tb$, which
denote homotopy classes of maps $u\colon \Strip\longrightarrow \Sym^g(\Sigma)$
satisfying
\begin{equation}
\label{eq:DefTimeDependentBoundary}
\left\{u\colon \Strip\longrightarrow
\Sym^g(\Sigma)\Bigg|
\begin{array}{l}
    u(1+it)\in \Psi_{t}(\Ta), \forall t\in\R, \\ 
    u(0+it)\in \Tb, \forall t\in\R, \\
    \lim_{t\goesto -\infty}u(s+it)=\x, \\
    \lim_{t\goesto +\infty}u(s+it)=\y
\end{array}\right\}
\end{equation}  
(dividing out by the action of $\pi_1$ on $\pi_2$ in the case where $g=2$).
We think of this set as a set of homotopy class of Whitney disks
with ``dynamic boundary conditions.''  Fixing such a class
$\phi\in\pi_2^{\Psi_t}(\x,\y)$, we have the moduli space of maps
$\ModFlow^{\Psi_{t}}(\phi)$ satisfying the above boundary conditions
and also $\DBar_{J_s} u = 0$.  It is easy to see that if
$\pi_2^{\Psi_t}(\x,\y)\neq \emptyset$, then $\x\in \Ta\cap \Tb$ and
$\y\in\Psi_{1}(\Ta)\cap \Tb$ lie in equivalence classes $\SpinCz$ and
$\SpinCz'$ corresponding to the same $\SpinC$ structure $\spinc$, for
the fixed base-point $z$.  Moreover, in this case
$\pi_2^{\Psi_t}(\x,\y)\cong \Z$.  Note that by juxtaposition, a single
$\phi_0\in\pi_2^{\Psi_t}(\x_0,\y_0)$ naturally gives rise to an
identification between homotopy classes
$\pi_2(\x_0,\x_0)\cong\pi_2(\y_0,\y_0)$ and a corresponding
identification between isomorphism classes of coherent orientation
systems for $(\Sigma,\alphas,\betas,\z)$ and $(\Sigma,\alphas',\betas,z)$.
(This point is more important, of course, in the case where $b_1(Y)>0$).

Now, we define the map associated to the isotopy 
$$
\Gamma^{\infty}_{\Psi_{t}}\colon \CFinf(\alphas,\betas,z)\longrightarrow
\CFinf(\alphas',\betas,z)
$$
by the formula
\begin{equation}
\label{eq:DefGamma}
\Gamma^{\infty}_{\Psi_{t}}([\x,i])
=
\sum_{\{\y\in\SpinCz'\}}\sum_{\{\phi\in\pi_2^{\Psi_t}(\x,\y)|\Mas(\phi)=0\}}
\#\left(\ModFlow^{\Psi_{t}}(\x,\y)\right)\cm [\y,i-n_z(\phi)],
\end{equation}
where $\Mas(\phi)$ is the expected dimension of the moduli space
$\ModFlow_{\Psi_t}(\phi)$, and the number
$\#\ModFlow^{\Psi_{t}}(\phi)$ is the signed number of points in this
zero-dimensional moduli space. Note that this is a finite sum for each
given $[\x,i]$, since $\pi_2^{\Psi_t}(\x,\y)\cong \Z$.  The map is
well-defined up to one overall sign. 

The important observation is that the moduli 
spaces considered have Gromov compactifications.
This follows from the energy bounds on the moduli spaces of
disks: given a class 
$\phi\in\pi_{2}^{\Psi_t}(\x,\y)$ with $\x\in\Ta\cap\Tb$ and 
$\y\in\Psi_1(\Ta)\cap\Tb$, we must show that there is a bound 
(depending only on $\phi$) on the energy of any holomorphic disk
$u\in\ModFlow^{\Psi_{t}}(\phi)$. 

Specifically, suppose that 
$u_0$ is some fixed representative of the homotopy class $\phi$. 
We will that there is some fixed constant
which bounds the energy of another such representative
$u_1$. Specifically, since 
$u_{0}$ and $u_{1}$ are homotopic
holomorphic disks which 
connect $\x$ to $\y$, i.e. there is a homotopy
$$U\colon \Strip\times [0,1]\longrightarrow \Sym^{g}(\Sigma),$$
with $U(s+it,0)=u_{0}(s+it), U(s+it,1)=u_{1}(s+it)$, and 
$U(1+it,\tau)\in\Psi_{t}(\Ta)$, $U(0+it,\tau)\in\Tb$. 
Then, there are corresponding
lifts 
\begin{eqnarray*}
    {\widetilde u}_{0} \colon F_{0}\longrightarrow 
    \Sigma^{\times g} &{\text{and}}&
    {\widetilde u}_{1} \colon F_{1}\longrightarrow 
    \Sigma^{\times g}
\end{eqnarray*}
coming from pulling back the branched covering $\pi\colon \ProdSig{g}
\longrightarrow \Sym^{g}(\Sigma)$  as in 
Subsection~\ref{subsec:EnergyBounds}. The homotopy of $u_{0}$ and 
$u_{1}$ lifts to give a homology between ${\widetilde u}_{0}$ and 
${\widetilde u}_{1}$, so that if $\omega_{0}$ denotes the product form 
on $\Sigma^{\times g}$, the difference
$$\int_{F_{1}}{\widetilde u}_{1}^{*}\omega_{0} - 
\int_{F_{0}}{\widetilde u}_{0}^{*}\omega_{0}$$
is calculated by integrating the pull-back of $\omega_{0}$ over the 
cover (in $\ProdSig{g}$) of  
the restriction $U|(\partial \Strip)\times[0,1]$. Indeed, the form $\omega_{0}$
vanishes over the image of $\{0\}\times\R\subset \partial \Strip$, 
so we need to bound the integral of the pull-back of $\omega_{0}$ via 
the $g!$ maps
$$f\colon \R\times [0,1]\longrightarrow \Sigma^{\times g},$$
given by
$$f(t,\tau)={\widetilde U}(1+it,\tau)$$ 
(where $\R\times [0,1]$ is one of the $g!$ components of 
the covering space of $\{1\}\times \R \times [0,1]\subset
\Strip \times [0,1]$ induced from $\ProdSig{g}$).
Thus, the map $f$ satisfies
$f(t,\tau)\in \Psi_{t}({\widetilde \Ta})$ 
(where ${\widetilde Ta}$ is the preimage of $\Ta$ in $\ProdSig{g}$)
for all $t,\tau$.
Taking derivatives of this condition, we get that
$$\frac{df(t,\tau)}{dt}\equiv X_{t}f(t,\tau) 
\pmod{(\Psi_{t})_{*}(T{\widetilde \Ta})}$$
and
$$\frac{df(t,\tau)}{d\tau}\equiv 0 
\pmod{(\Psi_{t})_{*}(T{\widetilde \Ta})}.$$
In view of the fact that $\omega_{0}$ vanishes on all tangent spaces 
$(\Psi_{t})_{*}(T{\widetilde \Ta})$, we get that
\begin{eqnarray*}
    \int f^{*}\omega_{0}  &=&
    \int \omega_{0}(\frac{df}{dt}, 
    \frac{df}{d\tau}) dt\wedge  d\tau \\
    &=& \int \omega_{0}(X_{t}, \frac{df}{d\tau}) \\
    &=& \int \langle dH_{t}, \frac{df}{d\tau}\rangle dt \wedge d\tau 
    \\
    &=& \int \big(H_{t}(f(t,1))-H_{t}(f(t,0))\big) dt.
\end{eqnarray*}
Since $H_{t}$ is identically zero outside 
$t\in [0,1]$, the above integral is bounded by a constant 
$$ k=\sup_{(t,\w)\in \R\times\Sym^g(\Sigma)} H_t(\w) - 
\inf_{(t,\w)\in \R\times\Sym^g(\Sigma)} H_t(\w)$$
which is independent of our choice of $u_0$ and $u_1$.
Thus, we see that
$$\int_{F_{1}}{\widetilde u}_{1}^{*}(\omega_{0})-
\int_{F_{0}}{\widetilde u}_{0}^{*}(\omega_{0})\leq k.$$
But now, in view of the proof of Lemma~\ref{lemma:EnergyBound}, we see that
if $\omega_{1}$ denotes the symplectic form on $\Sym^{g}(\Sigma)$ 
(which tames the product complex structure), then 
\begin{eqnarray*}
    \int_{\Strip} u_{1}^{*}(\omega_{1})
    &\leq & C \int_{F_{0}} {\widetilde u}_1(\omega_0) \\
    &\leq & C \Bigg(\int_{F_{0}} {\widetilde u}_0(\omega_0) +
			\left(\int_{F_{1}} {\widetilde u}_1(\omega_0) -
				\int_{F_{0}} {\widetilde u}_0(\omega_0)\right)\Bigg) \\
    &\leq & C \left(\int_{F_{0}} {\widetilde u}_0(\omega_0) + k\right)
\end{eqnarray*}
which in turn is independent of our choice of $u_1$. (Indeed, for the case
where $g=2$, it is easy to see that the key property used for this bound
is not that $u_1$ is homotopic to $u_0$, but merely that the two maps
induce the same relative homology class.)

With the energy bounds in place, now, we can show that 
$\Gamma^{\infty}$
induces an isomorphism in $\HFinfty$. 

As in the proof of Theorem~\ref{thm:IndepCxStruct} above, we verify
that $\Gamma^\infty_{\Psi_t}$ is a chain map. Consider
$\psi\in\pi_2^{\Psi_t}(\x,\y)$ with $\Mas(\psi)=1$. Then,
$\ModFlow_{\Psi_t}(\psi)$ is a non-compact space, which can be
compactified by: 
$$
	\partial \ModFlow_{\Psi_t}(\psi) = 
	\left(\coprod_{\phi\csum\phi'=\psi}\ModFlow_{\Psi_t}(\phi)\times
	\UnparModFlow(\phi')\right)
	\coprod 
	\left(\coprod_{\phi'\csum\phi=\psi} \UnparModFlow(\phi') \times 
	\ModFlow_{\Psi_t}(\phi)\right).
$$
In the first decomposition, $\phi$ ranges over those elements 
of $\pi_2^{\Psi_t}(\x,\y')$ (where $\y'$ is any point in 
$\Psi_1(\Ta)\cap\Tb$)
with $\Mas(\phi)=0$, while $\phi'$ ranges over those homotopy classes 
$\phi'\in\pi_2(\y',\y)$ for the tori $\Psi_1(\Ta)$ and $\Tb$, which
satisfy $\Mas(\phi')=1$, and $\phi\csum\phi'=\psi$ (using the obvious 
juxtaposition operation).
In the second, $\phi'$ ranges over those elements of $\pi_2(\x,\x')$
(where $\x'$ is any point in $\Ta\cap\Tb$) 
where $\Mas(\phi')=1$, while $\phi$ ranges over those elements of
$\pi_2^{\Psi_t}(\x',\y)$ with $\Mas(\phi)=0$, and $\phi'\csum\phi=\psi$.
By counting points with sign, it
follows that $\Gamma^{\infty}_{\Psi_t}$ is a chain map. 

To see that $\Gamma^\infty_{\Psi_t}$ induces an isomorphism in homology,
observe that the composite $\Gamma^\infty_{\Psi_t}\circ
\Gamma^\infty_{\Psi_{1-t}}$ is chain homotopic to the identity map. The
chain homotopy is constructed using a homotopy $\Phi_{t,\tau}$ between
two one-parameter families of isotopies -- thought of as $\tau\mapsto
\Psi_{t}(\tau)$ -- which connects the juxtaposition of $\Psi_t$ with
$\Psi_{1-t}$ at $\tau=0$, to the stationary identity isotopy at
$\tau=1$.  Letting $$\ModFlow^{\Phi_{t,\tau}}(\phi)=
\bigcup_{\tau\in[0,1]}\ModFlow^{\Phi_{t}(\tau)}(\phi),$$
we define
$$H^\infty([\x,i]) = 
\sum_{\y}\sum_{\{\phi\in\pi_2(\x,\y)|\Mas(\phi)=-1\}}
\#\left(\ModFlow^{\Phi_{t,\tau}}(\phi)\right)
[\y,i-n_z(\phi)],$$
where we are implicitly using the identification $\pi_2(\x,\y)\cong
\pi_2^{\Phi_t(\tau)}(\x,\y)$ for each $\tau\in[0,1]$.
Note that $H^\infty$ has degree $-1$. To see, then that $H^\infty$ is the chain
homotopy between $\Gamma^\infty_{\Psi_t}\circ \Gamma^\infty_{\Psi_{1-t}}$ and the identity
map, we consider ends of the moduli spaces 
$\ModFlow^{\Phi_{t,\tau}}(\psi)$,
where $\psi$ has $\Mas(\psi)=0$.
These spaces have three kinds of ends: those at to $\tau=0$, which
correspond to the composite 
$\Gamma^\infty_{\Psi_t}\circ \Gamma^\infty_{\Psi_{1-t}}$, those
at
to $\tau=1$, corresponding to the identity map, and those
at the  splittings
$$
\left(\coprod_{\phi\csum\phi'=\psi}\ModFlow^{\Phi_{t,\tau}}(\phi)\times
\UnparModFlow(\phi')\right)
\coprod 
\left(\coprod_{\phi'\csum\phi=\psi} \UnparModFlow(\phi') \times 
\ModFlow^{\Phi_{t,\tau}}(\phi)\right);$$
where $\phi$ all satisfy $\Mas(\phi)=0$, and 
$\phi'$ satisfy $\Mas(\phi')=1$. 
Counting these ends with sign, we obtain the relation
$$\Gamma^\infty_{\Psi_t}\circ \Gamma^\infty_{\Psi_{1-t}} = \Id + \partial^\infty\circ H^\infty - H^\infty\circ \partial^\infty.$$ 
Switching the roles of $\Psi_t$ and
$\Psi_{1-t}$, it follows that $\Gamma^\infty$ induces an isomorphism in
homology.

This 
same technique applies to prove the result for $\HFa$, $\HFp$, $\HFm$,
and the $U$ equivariance proves the result for $\HFred$.
\qed 
\vskip.3cm

The case where $b_1(Y)>0$ proceeds much as before. However, special
care must be taken to see that the sum appearing in the definition of
the chain map induced by an isotopy of the $\alpha$-curves, as given
in Equation~\eqref{eq:DefGamma} is, in fact, a finite sum for each
given $[\x,i]$. Again, this is done with the help of admissibility hypotheses.

In a single pair creation, the domains for the diagrams
$(\Sigma,\alphas,\betas)$, and
$(\Sigma,\{\psi(\alpha_1),\alpha_2,...,\alpha_g\},\betas)$ do not
coincide: the latter has a new domain. Correspondingly a homotopy class
$\pi_2^{\Psi}(\x,\y)$ does not have a well-defined multiplicity at
this new domain, since the $\Psi_t(\Ta)$ crosses the subvariety
$\{w\}\times \Sym^{g-1}(\Sigma)\subset \Sym^g(\Sigma)$, where $w$ is
any point in this new domain. 

However, the multiplicities at the other domains are still
well-defined; i.e. if $\cald_i$ is any domain which exists before the
pair-creation, and $w_i\in\cald_i$ is a point in the interior of this
domain, then the intersection number $u\cap
\left(\{w_i\}\times \Sym^{g-1}(\Sigma)\right)$ (where $u$ is any map representing
$\phi\in\pi_2^{\Psi_t}(\x,\y)$) is independent of the choice of
representative $u$ and the point $w_i$ (we choose the isotopy $\Psi_t$ to
be constant near $\{w_i\}\times\Sym^{g-1}(\Sigma)$). We call this collection of
multiplicities the domain of $\phi$. In fact, in our
fixed one-parameter family of nerly-symmetric complex structures, we 
choose our basepoints $\{w_i\}_{i=1}^m\subset \CurveComp$ with one in each 
of these domains $\cald_i$.

\begin{lemma}
\label{lemma:StrongDynamicIsotopies}
Fix  $(\Sigma,\alphas,\betas,z)$ be a strongly $\spinc$-admissible pointed
Heegaard diagram, and an isotopy $\Psi_t$ as above. 
Then, for each pair of integers $j$, and
for each $\x\in\Ta\cap \Tb$, $\y\in\Ta'\cap\Tb$, there are only finitely many homotopy classes
$\psi\in\pi_2^{\Psi_t}(\x,\y)$ with $\Mas(\psi)=j$ which support
$J_s$-holomorphic representatives.
\end{lemma}

\begin{proof}
Let $w_1,...,w_m$ be points contained in the interiors of the domains
before the pair-creation, and $w_{m+1}$ be a point in the new domain. 
Let $\Ta'$ be the torus $\psi_1(\alpha)\times...\times\alpha_g$.
As before, if $\x,\y\in\Ta\cap\Tb$, we let $\pi_2(\x,\y)$ denote the
space of homotopy classes of Whitney disks for $\Ta$, $\Tb$; if
$\y'\in\Ta'\cap\Tb$  we let
$\pi_2^{\Psi_t}(\x,\y)$ denote the homotopy classes with moving
boundary conditions defined above, and we let $\pi_2'(\x,\y)$ denote
the homotopy classes of Whitney disks for the pair $\Ta'$ and $\Tb$,
now thinking of $\x$ and $\y$ as intersections between those tori.

Fix $\x\in\Ta\cap\Tb$, and $\y\in\Ta'\cap\Tb$.
It is easy to see that each homotopy class $\pi_2^{\Psi_t}(\x,\y)$ has
a representative $u(s,t)$ which is constant for $t\leq 1$. As such,
$u$ can be thought of as representing a class $\pi_2'(\x,\y)$. Indeed,
this induces a one-to-one correspondence $\pi_2^{\Psi_t}(\x,\y)\cong
\pi_2'(\x,\y)$. In a similar manner, if $\x\in \Ta\cap\Tb$, we have
identifications $\pi_2'(\x,\x)\cong \pi_2^{\Psi_t}(\x,\x)\cong
\pi_2(\x,\x)$, which preserve all the local multiplicities $n_{w_j}$
for all $j=1,...,m$.

Let $\{\psi_i\}$ be a sequence of homotopy classes in
$\pi_2^{\Psi_t}(\x,\y)$ which support holomorphic representatives, and
have a fixed Maslov index. Next, fix $\phi_0\in\pi_2(\y,\x)$. Since
the $\psi_i$ all support holomorphic representatives, the local
multiplicities at the $w_j$ for $j=1,...,m$ are non-negative; it
follows that for $j=1,...,m$, $n_{w_j}(\psi_i*\phi_0)\geq
n_{w_j}(\phi_0)$. But $\psi_i*\phi_0$ is a homotopy class connecting
$\x$ with $\x$, which are intersection points which existed before the
pair creation, so we can consider the corresponding element of
$\pi_2(\x,\x)$. From the above observations, the multiplicities at all
$w_i$ for $i=1,...,m$ are are bounded below, and the Maslov index is
fixed, so there can be only finitely many such homotopy classes,
according to Lemma~\ref{lemma:StrongFiniteness}. The lemma now
follows.
\end{proof}

\vskip.2cm
\noindent{\bf{Proof of Theorem~\ref{thm:Isotopies} when $b_1(Y)>0$.}}
In view of Lemma~\ref{lemma:StrongDynamicIsotopies}, 
it follows that $\Gamma^{\infty}_{\Psi_t}$ as defined above is
a finite sum for each fixed $[\x,i]$, so the earlier proof in the 
case where $b_1(Y)=0$ applies. Note that the
Lemma~\ref{lemma:StrongDynamicIsotopies} also holds for isotopies
obtained by juxtaposing $\Psi_t$ with $\Psi_{1-t}$. 

Establishing its $H_1(Y;\Z)/\Tors$-equivariance of the map
$\Phi^\infty$ follows as in the proof of
Lemma~\ref{lemma:ActionCommutes} above.
\qed
\vskip.2cm

\subsection{Weakly admissible Heegaard diagrams}
\label{subsec:Weakly}

\begin{theorem}
\label{thm:WeakIsotopies}
If $Y$ is a three-manifold equipped with a $\SpinC$ structure
$\spinc$, and endowed with a weakly $\spinc$-admissible Heegaard diagram.
Then, there is an isotopic strongly $\spinc$-admissible Heegaard diagram,
a identification between equivalence
classes of orientation systems for the two diagrams and corresponding
identification between the homology groups 
(thought of as $\Z[U]\otimes_\Z\Wedge^* H_1(Y;\Z)/\Tors$-modules)
\begin{eqnarray*}
\HFa(\Sigma,\alphas,\betas,\orient)&{\cong}&
\HFa(\Sigma',\alphas',\betas',\orient') \\
\HFp(\Sigma,\alphas,\betas,\orient)&{\cong}&
\HFp(\Sigma',\alphas',\betas',\orient')
\end{eqnarray*} 
\end{theorem}

\begin{proof}
According to Lemma~\ref{lemma:WeaklyAdmissibleIsotopy}, we can find a
weakly $\spinc$-isotopic, strongly $\spinc$-admissible Heegaard
diagram. (Indeed, it is easy give exact Hamiltonian isotopies
of $\Sigma$ which are suitable for the purposes of
Lemma~\ref{lemma:StronglyAdmissible}.)
Note that the analogue of Lemma~\ref{lemma:StrongDynamicIsotopies}
also holds in the weakly admissible context (where now both
$\Mas(\psi)$ and $n_z(\psi)$ are fixed), so we can construct chain
homotopy equivalences $\Gamma^+_{\Psi_t}$ and ${\widehat
\Gamma}_{\Psi_t}$ by modifying Equation~\eqref{eq:DefGamma} (for
example, in the definition of $\Gamma^+_{\Psi_t}$ we drop terms
involving $[\y,j]$ with $j<0$) showing that the corresponding groups
are isomorphic.
\end{proof}

\section{Holomorphic triangles}
\label{sec:HolTriangles}

Maps between Floer homologies can be constructed by counting
pseudo-holomorphic triangles in a given equivalence class. This
construction is fundamental to establishing the handleslide invariance
of the Floer homologies considered here. As we shall see in the sequel
(c.f. Section~\ref{HolDiskTwo:sec:Surgeries} of~\cite{HolDiskTwo}),
they are useful also when comparing the Floer homology groups of
three-manifolds which differ by surgeries on a knot. More applications
are also given in~\cite{HolDiskThree}.  Thus, we allow ourselves now a
lengthy digression into the properties of these maps.

Since holomorphic triangles fit naturally into a four-dimensional
framework, we begin the section by setting up the relevant
(four-dimensional) topological preliminaries, including the map from
homotopy classes of triangles to $\SpinC$ structures over an
associated four-manifold. In Subsection~\ref{subsec:Orientations}, we
discuss issues of orientations. In
Subsection~\ref{subsec:HolTriangleMaps}, we discuss admissibility
issues, and then set up the maps induced by holomorphic triangles,
discussing various invariance properties of the maps.  The maps enjoy
a certain associativity property, which we will make use of in the
proof of handleslide invariance, and in the sequel. This associativity
is treated in Subsection~\ref{subsec:Assoc}, see also~\cite{FOOO},
\cite{Silva}.  

\subsection{Topological preliminaries on triangles}
\label{subsec:TriangleTopPrelim}

A {\em Heegaard triple-diagram of genus $g$} is an oriented two-manifold and three $g$-tuples
$\alphas$, $\betas$, and $\gammas$ which are complete sets of attaching
circles for handlebodies $U_{\alpha}$, $U_\beta$, and $U_{\gamma}$
respectively. 
Let 
$Y_{\alpha,\beta}=U_{\alpha}\cup U_{\beta}$, 
$Y_{\beta,\gamma}=U_{\beta} \cup U_{\gamma}$, and
$Y_{\alpha,\gamma}=U_{\alpha} \cup U_{\gamma}$ denote the three
induced three-manifolds.
A Heegaard triple-diagram  naturally
specifies a cobordism $X_{\alpha,\beta,\gamma}$ between these
three-manifolds. The cobordism is constructed as follows. 

Let $\Delta$ denote the two-simplex, with vertices $v_{\alpha},
v_{\beta}, v_{\gamma}$ labeled clockwise, and let $e_{i}$ denote the
edge $v_{j}$ to $v_{k}$, where $\{i,j,k\}=\{\alpha,\beta,\gamma\}$. 
Then, we form the identification space
$$X_{\alpha,\beta,\gamma}=\frac{\left(\Delta\times \Sigma\right)
\coprod \left(e_\alpha \times U_\alpha\right) \coprod
\left(e_\beta\times U_\beta \right)\coprod
\left(e_\gamma\times U_\gamma \right)}
{ \left(e_{\alpha}\times \Sigma\right) \sim \left(e_\alpha \times
\partial U_\alpha\right), \left(e_{\beta}\times \Sigma\right) \sim
\left(e_\beta \times \partial U_\beta\right), \left(e_{\gamma}\times
\Sigma\right) \sim \left(e_\gamma \times \partial U_\gamma\right) }.
$$ Over the vertices of $\Delta$, this space has corners, which can be
naturally smoothed out to obtain a smooth, oriented, four-dimensional
cobordism between the three-manifolds $Y_{\alpha,\beta}$,
$Y_{\beta,\gamma}$, and $Y_{\alpha,\gamma}$ as claimed.

We will call the cobordism $X_{\alpha,\beta,\gamma}$ described above a
{\em pair of pants connecting $Y_{\alpha,\beta}$, $Y_{\beta,\gamma}$,
and $Y_{\alpha,\gamma}$}.  Note that $$\partial
X_{\alpha,\beta,\gamma}=-Y_{\alpha,\beta}-Y_{\beta,\gamma}+Y_{\alpha,\gamma},$$
with the obvious orientation.

\begin{example}
\label{ex:Cobordism}
Let $(\Sigma,\alphas,\betas)$ be a Heegaard diagram for $Y$, and let $\gammas$ be a
$g$-tuple of curves which are isotopic to $\betas$. Then the
triple-diagram 
$$(\Sigma,\alphas,\betas,\gammas)$$ is a diagram for the
cobordism between $-Y$, $Y$, and $\#^g(S^1\times S^2)$ obtained from
$Y\times [0,1]$ by deleting a regular neighborhood of $U_{\beta}\times
\OneHalf$. 
\end{example}

\subsubsection{Two-dimensional homology}

We can think of the two-dimensional homology of
$X=X_{\alpha,\beta,\gamma}$ in terms of the $\alphas$, $\betas$, and
$\gammas$ as follows:

\begin{prop}
\label{prop:HomologyOfX}
Let $\SpanA\subset H_1(\Sigma;\Z)$ denote the
lattice spanned by the one-dimensional homology classes induced by the
$\alphas$. Then, there are natural identifications
\begin{equation}
H_2(X;\Z)\cong \Ker\Big(\SpanA\oplus\SpanB\oplus\SpanC
\longrightarrow H_1(\Sigma;\Z)\Big);
\end{equation}
or, equivalently,
\begin{equation}
H_2(X;\Z)\cong \Ker\Big(H_1(\Ta;\Z)\oplus H_1(\Tb;\Z)\oplus H_1(\Tc;\Z)
\longrightarrow H_1(\Sym^g(\Sigma);\Z)\Big).
\end{equation}
Similarly, we have 
\begin{equation}
H_1(X;\Z)\cong \Coker\Big(\SpanA\oplus\SpanB\oplus\SpanC 
\longrightarrow H_1(\Sigma;\Z)\Big);
\end{equation}
\end{prop}

\begin{proof}
First, note that the boundary homomorphism $\partial \colon
H_2(U_\alpha,\Sigma;\Z)\longrightarrow H_1(\Sigma;\Z)$ is injective, and its image
is $\SpanA$. 
The first isomorphism then
follows from the long exact sequence in homology for the pair
$(X,\Delta\times \Sigma)$, bearing in mind that
$$H_2(X,\Delta\times \Sigma) \cong H_2(U_\alpha,\Sigma)\oplus
H_2(U_\beta,\Sigma)\oplus H_2(U_\gamma,\Sigma)$$
(by excision), and
that the map $H_2(\Sigma)\longrightarrow H_2(X)$ is trivial: the
Heegaard surface is obviously null-homologous in $X$.

The second isomorphism follows from the fact that 
under the natural identification 
$H_1(\Sym^g(\Sigma);\Z)\cong H_1(\Sigma;\Z)$,
the image of $H_1(\Ta;\Z)$ is identified with $\SpanA$.

The final isomorphism follows from the fact that
$$H_1(X,\Delta\times \Sigma) \cong H_1(U_\alpha,\Sigma)\oplus
H_1(U_\beta,\Sigma)\oplus H_1(U_\beta,\Sigma)\cong H^2(U_\alpha)\oplus
H^2(U_\beta)\oplus H^2(U_\gamma)=0.$$
\end{proof}

Suppose $(a,b,c)\in\SpanA\oplus\SpanB\oplus\SpanC$ satisfies
$a+b+c=0$. Then, of course, $a+b+c$ spans some two-chain in
$\Sigma$. Two-chains of this type which also vanish at a given
base-point $z$ (lying outside the collection of attaching circles) are
natural analogues of the periodic domains considered earlier. We call
such two-chains {\em triply-periodic domains}. In keeping with earlier
terminology, the data 
$(\Sigma,\alphas,\betas,\gammas,z)$ where we choose a reference point
$z\in\ThreeCurveComp$ is called a {\em pointed Heegaard
triple-diagram}.

\subsubsection{Homotopy classes of triangles}

Let $\x\in\Ta\cap \Tb$, 
$\y\in\Tb\cap \Tc$, $\w\in\Ta\cap\Tc$. Consider the map 
$$u \colon \Delta \longrightarrow \Sym^{g}(\Sigma) $$
with the boundary conditions that 
$u(v_{\gamma})=\x$, $u(v_{\alpha})=\y$, and $u(v_{\beta})=\w$, and 
$u(e_{\alpha})\subset \Ta$, $u(e_{\beta})\subset \Tb$, 
$u(e_{\gamma})\subset 
\Tc$. Such a map is called a {\em Whitney triangle connecting $\x$, 
$\y$, and $\w$}. Two Whitney triangles are homotopic if the maps are 
homotopic through maps which are all Whitney triangles. We let 
$\pi_{2}(\x,\y,\w)$ denote the space of homotopy classes of Whitney 
triangles connecting $\x$, $\y$, and $\w$. In the case where $g=2$,
this will actually denote equivalence classes of homotopy classes of Whitney
triangles, corresponding to the action of
$\pi_1(\Sym^g(\Sigma))$
on $\pi_2(\Sym^g(\Sigma))$. (Note that one
can alternately stabilize, as in Section~\ref{sec:Stabilization} to 
reduce to the case where $g>2$.)

As in the definition of Whitney disks, we have an obstruction
$$\epsilon \colon (\Ta\cap\Tb)\times (\Tb\cap\Tc) \times
(\Ta\cap \Tc)\longrightarrow
\frac{H_1(\Sym^g(\Sigma))}{H_1(\Ta)+H_1(\Tb)+H_1(\Tc)}\cong H_1(X;\Z)$$
which vanishes if $\pi_2(\x,\y,\w)$ is non-empty. The obstruction is
defined as follows. Choose an arc $a\subset \Tb$ from $\x$ to $\y$,
$b\subset \Tc$ from $\y$ to $\w$, and an arc $c\subset \Ta$ from $\w$ to $\x$. 
Then, $\epsilon(\x,\y,\w)$ is the equivalence class of the closed path
$a+b+c$.

Using a base-point 
$z\in\Sigma-\alpha_1-...-\alpha_g-\beta_1-...-\beta_g-\gamma_1-...-\gamma_g$,
we obtain an intersection number
$$n_z\colon \pi_2(\x,\y,\w)\longrightarrow \Z.$$

\begin{prop}
\label{prop:CalcPiTwo}
Given $\x\in\Ta\cap\Tb$, $\y\in\Tb\cap\Tc$, $\w\in\Ta\cap\Tc$, then
$\pi_2(\x,\y,\w)$ is non-empty if and only if $\epsilon(\x,\y,\w)=0$. 
Moreover, if $g>1$ and $\epsilon(\x,\y,\w)=0$ then 
$$\pi_2(\x,\y,\w)\cong \Z \oplus H_2(X;\Z).$$
\end{prop}

\begin{proof}
The first statement is clearly true.
Let $\Map^W(\Delta,\Sym^g(\Sigma))$ denote the space of Whitney triangles
connecting $\x,\y,\w$. Then, evaluation along the boundary gives a
fibration
$$\Map^W(\Delta,\Sym^g(\Sigma))\longrightarrow \Omega_{\Ta}(\w,\x)\times
\Omega_{\Tb}(\x,\y)\times\Omega_{\Tc}(\y,\w),$$
whose fiber is homotopy equivalent to the space of pointed maps from
the sphere to $\Sym^g(\Sigma)$ (the base space here is a product of
path spaces). In the case where $g>2$, this gives us an exact sequence
$$\begin{CD} 0@>>> \Z @>>> \pi_0(\Map^W(\Delta,\Sym^g(\Sigma))) @>>>
H_1(\Ta)\oplus H_1(\Tb)\oplus H_1(\Tc).
\end{CD}$$
By definition, $\pi_0(\Map^W(\Delta,\Sym^g(\Sigma)))\cong
\pi_2(\x,\y,\w)$. The evaluation $n_z$ provides a splitting for the
first inclusion, so that 
$$\pi_2(\x,\y,\w)\cong \Z \oplus \Image
\Big(\pi_2(\x,\y,\w)\longrightarrow H_1(\Ta)\oplus H_1(\Tb)\oplus
H_1(\Tc)\Big).$$ That image, in turn, is clearly identified with the
kernel of the natural map $H_1(\Ta)\oplus H_1(\Tb) \oplus
H_1(\Tc)\longrightarrow H_1(\Sym^g(\Sigma))$ (we are using here the
fact that $\pi_1(\Sym^g(\Sigma))$ is Abelian). The proposition when $g>2$ then
follows from Proposition~\ref{prop:HomologyOfX}. The case where $g=2$
follows as above, after dividing by the action of $\pi_1(\Sym^g(\Sigma))$.
\end{proof}

Note that the identification $\pi_2(\x,\y,\w)\cong \Z \oplus
H_2(X;\Z)$ is not canonical, but rather it is affine. Specifically, if
we fix a homotopy class: $\psi_0\in\pi_2(\x,\y,\w)$, then any other
homotopy class $\psi\in \pi_2(\x,\y,\w)$ differs from $\psi_0$ by an
integer $n_z(\psi)-n_z(\psi_0)$, and a triply-periodic domain
$\cald(\psi)-\cald(\psi_0)-(n_z(\psi)-n_z(\psi_0))[\Sigma]$ (which in
turn can be thought of as a two-dimensional homology class in $X$).

\subsubsection{$\SpinC$ structures}

There is a geometric interpretation of $\SpinC$ structures in four
dimensions, analogous to Turaev's interpretation of $\SpinC$
structures in three-dimensions, compare~\cite{KMcontact} and 
\cite{GompfStipsicz}.

Let $X$ be a four-manifold. We consider pairs $(J,P)$, where $P\subset
X$ is a collection of finitely many points in $X$, and $J$ is an
almost-complex structure defined over $X-P$. We say that two pairs
$(J_1,P_1)$ and $(J_2,P_2)$ are {\em homologous} if there is a compact
one-manifold with boundary $C\subset X$ containing $P_1$ and $P_2$,
with the property that $J_1|X-C$ is isotopic to $J_2|X-C$. We can
think of a $\SpinC$ structure on $X$ as a homology class of such pairs
$(J,P)$.

The identification with a more traditional definition is as
follows. Note that an almost-complex structure over $X-P$ has a
canonical $\SpinC$ structure, and that can be uniquely extended
over the points $P$ (the obstruction to extending lies in
$H^3(X,X-P)=0$, and the indeterminacy in extending lies in
$H^2(X,X-P)=0$). Conversely, given a $\SpinC$ structure with spinor
bundle $W^+$, a generic section $\Phi\in \Sections(X,W^+)$ vanishes at
finitely many points, away from which Clifford multiplication on
$\Phi$ sets up an isomorphism between $TX$ and $W^-$, hence endowing
$TX$ with a complex structure.  

Given a pair $(J,P)$, the first Chern class of the induced complex
tangent bundle of $X-P$ canonically extends to give a two-dimensional
cohomology class $c_1(J,P)\in H^2(X;\Z)$. In fact, this agrees with
the first Chern class $c_1(\spinc)$ of the spinor bundle $W^+$.

\subsubsection{Triangles and $\SpinC$ structures}
\label{subsec:TrianglesAndSpinCStructures}

The base-point $z\in\ThreeCurveComp$ gives rise to a relationship
between $\SpinC$ structures on $X$ and holomorphic triangles,
analogous to the construction of the $\SpinC$ structure on a
three-manifold belonging to intersection point between $\Ta\cap\Tb$
together with the basepoint $z$.

To describe this, fix ``height functions'' over the handlebodies
$f_i\colon U_i\longrightarrow [0,1]$ where $i=\alpha$, $\beta$, or
$\gamma$ with only $g$ index one critical points and one index zero
critical point, with $f_i(\partial U_i)=1$. 

Now, given a generic map $u\colon \Delta\longrightarrow \Sym^g(\Sigma)$
representing $\phi\in\pi_2(\x,\y,\w)$, there
is an immersed surface-with-boundary $F=F_0\cup F_1\subset
X_{\alpha,\beta,\gamma}$ constructed as follows. The intersection of
the component $F_0$ with $U_\xi\times e_\xi$, is the product of
$e_\xi$ with the upward gradient connecting the index zero critical
point with the point $z\in\Sigma$; its intersection with $\Delta\times
\Sigma$ is simply $\Delta\times\{z\}$. The intersection of $F_1$ with
$U_{\xi}\times e_\xi$ is given by the $g$-tuple of $f_\xi$ gradient
flow-lines connecting the various index one critical points with the
$g$ points over $(x,u(x))$ (where $x\in e_\xi$). Finally, in the
inside region $\Delta\times \Sigma$, the subset $F_1$ consists of
points $(x,\sigma)$, where $\sigma\in u(x)$. Note that in the
complement $X-(F_0\cup F_1)$, there is a well-defined oriented
two-plane field ${\mathcal L}$ which is tangent to $\Sigma$ inside
$\Delta\times \Sigma$, and agrees with the kernel of $df_\xi$ in
$TU_\xi\subset T(U_\xi\times e_\xi)$.

In fact, we extend the two-plane field further. Fix a central point
$x\in\Delta$, and three straight paths $a$, $b$, and $c$ from $x$ to the edges
$e_\alpha$, $e_\beta$, and $e_\gamma$ respectively. In the complement
$\Delta-a\cup b \cup c$, there is a foliation by line segments
which connect pairs of edges. For example, there is a family
$\ell_{\alpha,\beta}(t)$ of paths connecting $e_\alpha$ to $e_\beta$
which degenerates as $t\goesto 0$ to the vertex $v_\gamma$, and as
$t\goesto 1$ it degenerates to $a\cup b$. There are analogous families
of leaves $\ell_{\beta,\gamma}(t)$ and $\ell_{\alpha,\gamma}(t)$

There is a natural map $\pi\colon X \longrightarrow \Delta$. The
preimage under $\pi$ of $\ell_{\alpha,\beta}(t)$ for $t\in [0,1)$,
which we denote ${\widetilde \ell}_{\alpha,\beta}(t)$, is identified
with $Y_{\alpha,\beta}$. For all but finitely many $t$ in the 
open interval, the
intersection of $F$ with ${\widetilde \ell}_{\alpha,\beta}(t)$
consists of $g+1$ disjoint paths which connect the critical points of
$f_{\alpha}$ in $U_{\alpha}$ to critical points of $f_{\beta}$ in
$U_{\beta}$. For these generic 
$t$, we extend the oriented two-plane field in over a
neighborhood of these $g+1$ paths (as in
Subsection~\ref{subsec:SpinCStructures}) in a continuous manner. In this way,
we have extended ${\mathcal L}$ across the intersection of $F$ with
${\widetilde \ell}_{\alpha,\beta}(t)$ for all but finitely many $t$.

\begin{figure}
	\mbox{\vbox{\epsfbox{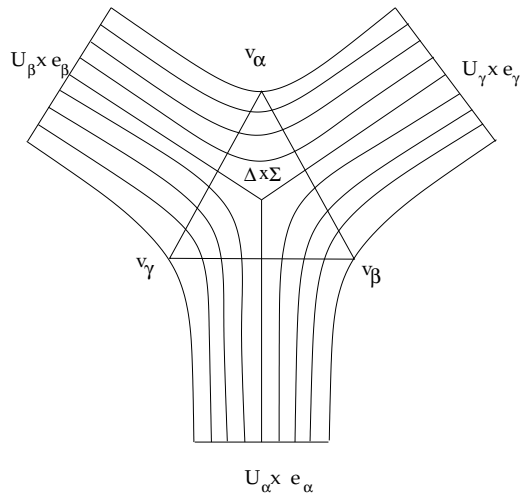}}}
	\caption{\label{fig:Cobord}
	Schematic for the cobordism $X$. We have illustrated the foliation of
	the triangle by segments whose preimages are the three-manifolds
	${\widetilde \ell}_{\xi,\eta}(t)$.}
\end{figure}

We proceed in the analogous manner to extend
over the ${\widetilde \ell}_{\beta,\gamma}(t)$ and ${\widetilde
\ell}_{\alpha,\gamma}(t)$. 

We have now extended ${\mathcal L}$ over $X$, except for the
intersection of $F$ with certain excluded leaves in the foliation of
$\Delta$. These excluded leaves fall into two categories. First, there
is the singular leaf $a\cup b \cup c$; and then there are those leaves
in $\Delta$ which contain a point $x$ for which $\sigma(x)$
has either a repeated entry, or $\sigma(x)$ contains the basepoint
$z\in\Sigma$.  
These are the points where the paths of $F$
cross. One can see that generically the intersection of $F$ with the
preimages of these special leaves is a collection of contractible
one-complexes; so its tubular neighborhood consists of a finite
collection of disjoint four-balls embedded in $X$.

The two-plane field ${\mathcal L}$ and the orientation on $X$ 
determine a complex structure over
the complement of finitely many balls in $X$, and hence a $\SpinC$
structure over $X$.

\begin{prop}
\label{prop:AssocSpinC}
The above construction induces a map
$$\spinc_z\colon \pi_2(\x,\y,\w)\longrightarrow
\SpinC(X).$$
\end{prop}

\begin{proof}
Recall that $\x$, $\y$, and $\w$ determine the two-plane field on the
boundary minus finitely many three-balls. Fix $u$ and choose
extensions over the three-balls
(c.f. Section~\ref{subsec:SpinCStructures}). This data specifies the two-plane
field over $X_{\alpha,\beta,\gamma}-
(\interior(\Delta)\times
\Sigma)-\interior F_0-\interior F_1$. The above discussion shows that the $\SpinC$
structure extends over this region, and, indeed, since the deleted region is
topologically a $\Delta\times
\Sigma$, it follows from a cohomology long exact sequence that the
extension is unique. It is easy to see also that the induced $\SpinC$
structure does not depend on the extension of the two-plane fields to
the three-balls in the boundary.

Changing $u$ by a homotopy moves $F_1$ by an isotopy, so it is easy to
see that the induced $\SpinC$ structure depends only on the homotopy
class of $u$.

Indeed, it is easy to see from the above discussion that the
restriction $u|\partial \Delta$ determines the $\SpinC$ structure. It
follows that adding spheres to the homotopy class in $\pi_2(\x,\y,\w)$
does not change the induced $\SpinC$ structure.
\end{proof}

Homotopy classes of Whitney triangles can be collected into
$\SpinC$-equivalence classes, as follows. 
Let $\x,\x'\in \Ta\cap\Tb$, $\y,\y'\in\Tb\cap\Tc$, and ${\mathbf v},{\mathbf
v}'\in\Ta\cap\Tc$.  We say that two homotopy classes $\psi\in
\pi_2(\x,\y,{\mathbf v})$ and $\psi'\in\pi_2(\x',\y',{\mathbf v}')$
are {\em $\SpinC$-equivalent}, or simply equivalent, if there are
classes $\phi_1\in\pi_2(\x,\x')$, $\phi_2\in\pi_2(\y,\y')$, and
$\phi_3\in\pi_2({\mathbf v},{\mathbf v}')$ with
$$\psi'=\psi+\phi_1+\phi_2+\phi_3.$$ Let $S_{\alpha,\beta,\gamma}$
denote the space of homotopy classes of such triangles. 

To justify the terminology, we claim that the $\SpinC$ structure
constructed above depends only on its $\SpinC$-equivalence class as
follows:

\begin{prop}
\label{prop:SpinCForTriangles}
The map from Proposition~\ref{prop:AssocSpinC} descends to a map
$$\spinc_z \colon S_{\alpha,\beta,\gamma} \longrightarrow \SpinC(X_{\alpha,\beta,\gamma})$$
which is one-to-one, with image consisting of those
$\SpinC$-structures whose restrictions to the boundary are realized by
intersection points. 
\end{prop}

\begin{proof}
Since addition of Whitney disks to a given homotopy class
$\psi\in\pi_2(\x,\y,\w)$ changes the induced two-plane field in a
neighborhood of the boundary, it follows that the above constructions
descend to give a map $$\spinc_z \colon S_{\alpha,\beta,\gamma}
\longrightarrow \SpinC(X_{\alpha,\beta,\gamma});$$ what remains to
verify that we have identified the image, and also that $\spinc_z$ is
one-to-one.

To see that we have characterized the image.  Recall that
for an arbitrary four-manifold-with-boundary $(X,Y)$ there is a
canonical map $\epsilon' \colon \SpinC(Y) \longrightarrow
H^3(X,Y;\Z)$  which is defined as follows.  Choose a $\SpinC$
structure $\spinc_0$ over $X$, and let
$$\epsilon'(\spinct)=\delta(\spinct-\spinc_0|_{Y}),$$ where
$\delta\colon H^2(Y;\Z)\longrightarrow H^3(X,Y;\Z)$ is the coboundary
map. It is easy to see that $\epsilon'$ is independent of the choice
of $\spinc_0$, and that it vanishes if and only if $\spinct$ extends
over $X$. Next, we argue that
$\PD[\epsilon(\x,\y,\w)]=\pm \epsilon'(\x,\y,\w)$. To see this, isotope $\Ta$,
$\Tb$, and $\Tc$ so that there are intersection points $\x'$, $\y'$,
and $\w'$ for which $\epsilon(\x',\y',\w')=0$, so that there is a
triangle connecting them. We have explicitly constructed the
corresponding $\SpinC$ structure, thus $\epsilon'(\x',\y',\w')=0$, as
well.  It is easy to see that
$$\epsilon(\x,\y,\w)-\epsilon(\x',\y',\w')=
\pm \PD \circ \delta\Big(\PD(\epsilon(\x,\x'))\oplus 
\PD(\epsilon(\y,\y'))\oplus \PD(\epsilon(\w,\w'))
\Big).$$
Similarly, 
$$\epsilon'(\x,\y,\w)-\epsilon'(\x',\y',\w')=
\delta\Big(\big(\spinc_z(\x)-\spinc_z(\x')\big)\oplus 
\big(\spinc_z(\y)-\spinc_z(\y')\big)
\oplus
\big(\spinc_z(\w)-\spinc_z(\w')\big) \Big).$$
It follows that $\PD(\epsilon(\x,\y,\w))=\pm\epsilon'(\x,\y,\w)$: the
obstructions to extending a $\SpinC$ structure are the same as the
obstruction to finding a Whitney triangle.

Suppose that $u$ and $v$ are a pair of triangles in $\pi_2(\x,\y,\w)$
with $n_z(u)=n_z(v)$, so that their difference from
Proposition~\ref{prop:CalcPiTwo} can be interpreted as the
triply-periodic domain $\cald(u)-\cald(v)$.  We claim that this
triply-periodic domain corresponds to a relative cohomology class in
$H^2(X,\partial X)$ whose image in $H^2(X)$ is the difference
$\spinc_z(u)-\spinc_z(v)$, which in turn verifies that the map
$\spinc_z$ is injective as claimed.  But this now is a local
calculation since, as is easy to verify, the restriction map
$$H^2(X,\partial X)
\longrightarrow H^2(U_\alpha \times (e_\alpha,\partial e_\alpha))\oplus H^2(U_\beta \times
(e_\beta,\partial e_\beta))
\oplus H^2(U_\gamma \times (e_\gamma,\partial e_\gamma))$$ is injective, and each of the
latter groups is generated by the Poincar\'e duals to the cylinders
$[\xi_i^*]\times e_\xi$ (where $\xi=\alpha$, $\beta$, or $\gamma$, and
$i=1,...,g$). On the one hand, the evaluation of a triply-periodic
domain on, say, $\alpha_1^*\times[0,1]$ is easily seen to be simply
the multiplicity of $\alpha_1$ in the boundary of the triply-periodic
domain. On the other hand, the pair of two-plane fields
representing $\spinc_z(u)$ and $\spinc_z(v)$
differ over $\alpha_1^*\times e_\alpha$ only at those
points where one of $u(e_\alpha)$ or $v(e_\alpha)$ contains
$\alpha_1\cap\alpha_1^*$. The fact that the constant appearing here is
one could be determined by calculating a model case (see~\cite{HolDiskThree}).
\end{proof}

\subsubsection{Higher polygons}
\label{subsec:HigherPolygons}

The above results for triangles admit straightforward generalizations
to arbitrarily large collections of $g$-tuples, which call 
{\em Heegaard multi-diagrams} (or {\em pointed Heegaard
multi-diagrams}, when they are equipped with a basepoint $z$ in the
complement of all the attaching circles). In fact, the only
other case we will require in the present work is the case of
squares. Specifically, an oriented two-manifold $\Sigma$ and four
$g$-tuples of attaching circles $\alphas$, $\betas$, $\gammas$, and
$\deltas$ specify a four-manifold $X_{\alpha,\beta,\gamma,\delta}$
which provides a cobordism between $Y_{\alpha,\beta}$,
$Y_{\beta,\gamma}$, $Y_{\gamma,\delta}$ and $Y_{\alpha,\delta}$. It
admits two obvious decompositions
$$X_{\alpha,\beta,\gamma,\delta}=X_{\alpha,\beta,\gamma}\cup_{Y_{\alpha,\gamma}}
X_{\alpha,\gamma,\delta}=X_{\alpha,\beta,\delta}\cup_{Y_{\beta,\delta}}
X_{\beta,\gamma,\delta}.$$

We can define homotopy classes of squares $\pi_2(\x,\y,{\mathbf v},\w)$
in $\Sym^g(\Sigma)$, and equivalence classes of homotopy classes
$S_{\alpha,\beta,\gamma,\delta}$ -- i.e. two squares
$\varphi\in\pi_2(\x,\y,{\mathbf v},\w)$ 
and $\varphi'\in\pi_2(\x',\y',{\mathbf v}',\w')$ 
are equivalent if there are $\phi_1\in\pi_2(\x,\x')$, 
$\phi_2\in\pi_2(\y,\y')$, $\phi_3\in\pi_2({\mathbf v},{\mathbf v}')$, 
and $\phi_4\in\pi_2(\w,\w')$ with
$$\varphi+\phi_1+\phi_2+\phi_3+\phi_4=\varphi'.$$
Proposition~\ref{prop:SpinCForTriangles} admits a straightforward
generalization, giving a map from $S_{\alpha,\beta,\gamma,\delta}$ to
the space of $\SpinC$ structures over
$X_{\alpha,\beta,\gamma,\delta}$. 

\subsection{Orienting spaces of pseudo-holomorphic triangles}
\label{subsec:Orientations}

We will be counting pseudo-holomorphic triangles. To achieve the
required transversality, we allow $J$ to be a function from $\Delta$
to the space of (nearly-symmetric) almost-complex structures over
$\Sym^g(\Sigma)$ chosen to be compatible near the corners with the
paths $J_s$ used to define the notion of pseudo-holomorphic disk.
Moreover, we will use a class of perturbations of the constant complex
structure for which the analogue of Lemma~\ref{lemma:NonNegativity}
still holds: if $u$ is a $J$-holomorphic triangle, the domain
associated to $u$ is non-negative.

Now, we can collect the space of $J$-holomorphic Whitney triangles
representing a fixed homotopy class into a moduli space, which we
denote $\ModFlow(\psi)$. This moduli space has an expected dimension,
which we will denote $\Mas(\psi)$.

With the transversality in place, the modulo two count of
$\ModFlow_J(\psi)$ is straightforward to define. When we wish to work
over $\Z$, however, we must use a refined count. Again this can be
done since the determinant line bundle of the tangent space admits an
extension $\det(D_u)$ as a trivial line bundle over each component
$\psi\in\pi_2(\x,\y,\w)$. 

Let $\spinc$ be a $\SpinC$ structure over the four-manifold $X$
specified by a pointed Heegaard triple
$(\Sigma,\alphas,\betas,\gammas,z)$, and let $\orient_{\alpha,\beta}$,
$\orient_{\beta,\gamma}$ and $\orient_{\alpha,\gamma}$ be coherent
systems of orientations for the three bounding three-manifolds.

\begin{defn}
A {\em coherent system of orientations for $\spinc$ $\orient_{\alpha,\beta,\gamma}$, 
compatible with
$\orient_{\alpha,\beta}$, $\orient_{\beta,\gamma}$, and
$\orient_{\alpha,\gamma}$} is a
collection of sections $\orient_{\alpha,\beta,\gamma}$
of 
the determinant line bundle $\det(D_u)$ for each homotopy class of triangle
$\psi$ representing the $\SpinC$ structure
$\spinc$, which is compatible with splicing in the
sense that if $\psi\in\pi_2(\x,\y,\w)$, $\psi_1\in\pi_2(\x,\x')$,
$\psi_2\in\pi_2(\y,\y')$, $\psi_3\in\pi_2(\w,\w')$ are any three
Whitney disks, then:
$$\orient_{\alpha,\beta,\gamma}(\psi+\phi_1+\phi_2+\psi_3)=
\orient_{\alpha,\beta,\gamma}(\psi)\wedge \orient_{\alpha,\beta}(\phi_1)
\wedge \orient_{\beta,\gamma}(\phi_2)
\wedge \orient_{\alpha,\gamma}(\phi_3),$$
under the identification coming from splicing.
\end{defn}

Existence is ensured by the following:

\begin{lemma}
\label{lemma:OrientTriangles}
Let $(\Sigma,\alphas,\betas,\gammas,z)$ be a pointed Heegaard triple,
and fix a $\SpinC$ structure $\spinc$ over $X_{\alpha,\beta,\gamma}$
whose restrictions $\spinct_{\alpha,\beta}$, $\spinct_{\beta,\gamma}$
and $\spinct_{\alpha,\gamma}$ are all realized by intersection points.
For coherent  systems ${\mathfrak o}_{\alpha,\beta}$ and
${\mathfrak o}_{\beta,\gamma}$ for two of the boundary components,
there always exists at least one coherent orientation system
$\orient_{\alpha,\gamma}$ for the remaining boundary component, and a
coherent system ${\mathfrak o}_{\alpha,\beta,\gamma}$ which is
compatible with the $\orient_{\alpha,\beta}$,
$\orient_{\beta,\gamma}$, and $\orient_{\alpha,\gamma}$. 
\end{lemma}

\begin{proof}
Let $\psi_0\in\pi_2(\x_0,\y_0,\w_0)$ be a fixed homotopy class representing $\spinc$. 
Fix an arbitrary orientation ${\mathfrak o}_{\alpha,\beta,\gamma}(\psi_0)$. 

Next, we construct ${\mathfrak o}_{\alpha,\gamma}(\phi_{\alpha,\gamma})$, where
$\phi_{\alpha,\gamma}\in\PerClasses{\w_0}\subset \pi_2(\w_0,\w_0)$ are periodic
classes. To this end, let $K$ denote the subgroup of $\PerClasses{\w_0}$ 
consisting of elements $\phi_3\in\PerClasses{\w_0}$ 
which satisfy the property that
$$\psi_0+\phi_3=\psi_0+\phi_1+\phi_2$$ for some periodic classes
$\phi_1$ and $\phi_2$ in $\PerClasses{\x_0}$ and $\PerClasses{\y_0}$
respectively. We claim that (1) given $\phi_3\in K$, the corresponding
$\phi_1$ and $\phi_2$ are uniquely determined, and also that
(2) the quotient $Q$ of $\PerClasses{\w_0}$ is a free Abelian group.
According to Claim (2), we obtain a splitting
$$\PerClasses{\w_0}\cong K\oplus Q.$$ Now, for 
$\phi_3\in K$, we define $\orient_{\alpha,\gamma}(\phi_3)$ so that
$$\orient_{\alpha,\beta,\gamma}(\psi_0)\wedge
\orient_{\alpha,\gamma}(\phi_3)=
\orient_{\alpha,\beta,\gamma}(\psi_0)\wedge \orient_{\alpha,\beta}(\phi_1)\wedge 
\orient_{\beta,\gamma}(\phi_2).$$
According to Claim (1), this is a consistent definition.
We then define $\orient_{\alpha,\gamma}(\phi)$ arbitrarily on a
basis of generators for $Q$, and allow that to induce the orientation on all
$\psi\in\pi_2(\x_0,\y_0,\w_0)$.
Both Claim (1) and (2) follow from the following diagram:
$$
\begin{CD}
&	&	&	&	H_2(Y_{\alpha,\gamma}) \\
&	&	&	&	@V{i}VV	\\
H_3(X,Y_{\alpha,\beta}\coprod Y_{\beta,\gamma})
@>>>
H_2(Y_{\alpha,\beta}\coprod Y_{\beta,\gamma}) @>{j}>>H_2(X)
@>{q}>> H_2(X,Y_{\alpha,\beta}\coprod Y_{\beta,\gamma}),
\end{CD}$$
where $X=X_{\alpha,\beta,\gamma}$.  Interpreting the second homology
as periodic domains, we see that the quotient group $Q$ is the image
of $q\circ i$.  Now, by excision, 
$$H_*(X,Y_{\alpha,\beta}\coprod
Y_{\beta,\gamma})\cong H_*(U_{\beta}\times (e_\beta, \partial e_\beta))$$
(where here $U_\beta$ is a handlebody and $e_\beta$ is an interval); in particular, $j$ 
is injective, establishing Claim (1), and also the image of $q$ is
non-torsion, establishing Claim (2).

We have thus consistently oriented all triangles in $\pi_2(\x_0,\y_0,\w_0)$.
As a final step, we choose a complete set of paths
$\{\theta_i\}_{i=1}^m$ for $Y_{\alpha,\gamma}$ over which we choose
our orientations (for $\orient_{\alpha,\gamma}$) arbitrarily, and use them to define the orientation
for all the remaining $\psi\in\pi_2(\x,\y,\w)$ in the given $\SpinC$-equivalence class.
\end{proof}

\subsection{Holomorphic triangles and maps between Floer homologies}
\label{subsec:HolTriangleMaps}

Our aim is to use these counts to define maps between Floer
homologies. To do this, we will need our triple-diagram to satisfy
some admissibility hypotheses, which are direct generalizations of the
admissibility conditions from Subsection~\ref{subsec:Admissibility}.

\begin{defn}
\label{def:AdmissibleTriple}
A pointed Heegaard triple-diagram is
called {\em weakly admissible} if each non-trivial triply-periodic domain
which can be written as a sum of doubly-periodic domains has both
positive and negative coefficients. A pointed triple-diagram
is called {\em strongly admissible} for the $\SpinC$ structure $\spinc$
if for each triply-periodic domain ${\mathcal D}$ which can be written
as a sum of doubly-periodic domains 
$${\mathcal D}={\mathcal D}_{\alpha,\beta}+{\mathcal D}_{\beta,\gamma}+{\mathcal
D}_{\alpha,\gamma}$$
with the property that 
$$\langle c_1(\spinc_{\alpha,\beta}), H({\mathcal D}_{\alpha,\beta})
\rangle + \langle c_1(\spinc_{\beta,\gamma}), H({\mathcal D}_{\beta,\gamma})
\rangle + \langle c_1(\spinc_{\alpha,\gamma}), H({\mathcal D}_{\alpha,\gamma})
\rangle =2n\geq 0,$$
there is some coefficient of ${\mathcal D}>n$.
(In the above expression, of course, $\spinc_{\xi,\eta}$ is the restriction of $\spinc$
to the boundary component $Y_{\xi,\eta}$).
\end{defn}

Note that the above notion of weak admissibility is independent of
$\SpinC$ structures -- it corresponds to the notion of weak
admissibility for any torsion $\SpinC$ structure, for an ordinary
pointed Heegaard diagram. (We could, of course, have given a slightly
weaker formulation depending on the $\SpinC$ structure, more parallel
to the definition of weakly admissible for pointed Heegaard diagrams
given earlier, but we have no particular use for this presently.)

The following are analogues of Lemmas~\ref{lemma:WeakFiniteness} and
\ref{lemma:StrongFiniteness}:

\begin{lemma}
\label{lemma:WeakFiniteTriple}
Let $(\Sigma,\alphas,\betas,\gammas,z)$ be weakly admissible Heegaard
triple, with underlying four-manifold $X$. Fix
intersection points $\x$, $\y$, and $\w$ and a $\SpinC$ structure
$\spinc$ over $X$. Then, for each integer $k$, 
there are only finitely many
homotopy classes $\psi\in\pi_2(\x,\y,\w)$ with $n_z(\psi)=k$ with
$\spinc_z(\psi)=\spinc$, and which support holomorphic representatives.
\end{lemma}

\begin{proof}
Given $\psi,\psi'\in\pi_2(\x,\y,\w)$ with $n_z(\psi)=n_z(\psi')$ and
$\spinc_z(\psi)=\spinc_z(\psi')$, the difference
$\cald(\psi)-\cald(\psi')$ is a triply-periodic domain which, in view
of Proposition~\ref{prop:SpinCForTriangles}, can be written as a sum
of doubly-periodic domains. Given this, finiteness follows as in the
proof of Lemma~\ref{lemma:WeakFiniteness}.
\end{proof}

\begin{lemma}
\label{lemma:StrongFiniteTriple}
For a strongly admissible pointed Heegaard triple-diagram for a given
$\SpinC$ structure $\spinc$, and an integer $j$, there are only
finitely many $\psi\in\pi_2(\x,\y,\w)$ representing $\spinc$
with $\Mas(\psi)=j$ and which
support holomorphic representatives.
\end{lemma}

\begin{proof}
Suppose that $\psi, \psi'\in\pi_2(\x,\y,\w)$ satisfy
$\spinc_z(\psi)=\spinc_z(\psi')$, and $\Mas(\psi)=\Mas(\psi')$. Then
we can write
$\psi'=\psi+\phi_1+\phi_2+\phi_2$; so by the additivity of the index, 
it follows that $\Mas(\phi_1)+\Mas(\phi_2)+\Mas(\phi_3)=0$ (which is
identified with the first Chern class evaluation). 
The proof then follows from the proof of Lemma~\ref{lemma:StrongFiniteness}. 
\end{proof}

Existence of admissible triples follows along the lines of
Section~\ref{sec:Special}.

\begin{lemma}
Given a Heegaard triple-diagram $(\Sigma,\alphas,\betas,\gammas,z)$,
there is an isotopic weakly admissible Heegaard triple
diagram. Moreover, given a $\SpinC$ structure $\spinc$ over $X$, there
is an isotopic strongly $\spinc$-admissible Heegaard triple diagram.
\end{lemma}

\begin{proof}
This follows as in Lemma~\ref{lemma:StronglyAdmissible}: we wind transverse 
to all of the $\alphas$, $\betas$, and $\gammas$
simultaneously.
\end{proof}

A $\SpinC$ structure over $X$ gives rise to a map
$$f^{\infty}(~\cdot~;\spinc)\colon \CFinfty(Y_{\alpha,\beta},\spinc_{\alpha,\beta})
\otimes \CFinfty(Y_{\beta,\gamma},\spinc_{\beta,\gamma}) \longrightarrow
\CFinfty(Y_{\alpha,\gamma},\spinc_{\alpha,\gamma}) $$
by the formula:
\begin{equation}
\label{eq:DefTriangle}
f^{\infty}_{\alpha,\beta,\gamma}([\x,i]\otimes[\y,j];\spinc)
=\sum_{\w\in\Ta\cap\Tc}
\sum_{\{\psi\in\pi_2(\x,\y,\w)\big|\spinc_z(\psi)=\spinc,\Mas(\psi)=0\}}
\Big(\#\ModFlow(\psi)\Big)\cm {[\w,i+j-n_z(\psi)]}.
\end{equation}
For each fixed $[\x,i]$ and $[\y,j]$ the above is a finite sum when
the triple is strongly admissible for $\spinc$.

In fact, for each fixed $[\x,i]$ and $[\y,j]$, the $[\w,k]$ coefficient is a
sum of $\#\ModFlow(\psi)$, where $\psi$ ranges over
$\psi\in\pi_2(\x,\y,\w)$ with $\spinc_z(\psi)=\spinc$ and
$n_z(\psi)=i+j-k$. Thus (according to Lemma~\ref{lemma:WeakFiniteTriple}), the $[\w,k]$ coefficient is given by a finite
sum under the weak admissibility hypothesis. 

Hence, if the triple is weakly admissible, the above sum induces a map
$$f_{\alpha,\beta,\gamma}^+(\cdot~;\spinc)\colon 
\CFp(Y_{\alpha,\beta},\spinc_{\alpha,\beta})
\otimes {CF}^{\leq 0}(Y_{\beta,\gamma},\spinc_{\beta,\gamma})
\longrightarrow \CFp(Y_{\alpha,\gamma},\spinc_{\alpha,\gamma}),$$
where, $\CF^{\leq 0}(Y,\spinc)\subset \CFinf(Y,\spinc)$ is the
subcomplex generated by $[\x,i]$ with $i\leq 0$. Of course, $\CFleq(Y,\spinc)$ is
isomorphic to $\CFm(Y,\spinc)$ as a chain complex (but the latter is
generated by $[\x,i]$ with $i<0$). 

Similarly, we can define a map 
$$
{\widehat f}_{\alpha,\beta,\gamma}(\x\otimes\y;\spinc)
=\sum_{\w\in\Ta\cap\Tc}\sum_{\{\psi\in\pi_2(\x,\y,\w)\big|\spinc_z(\psi)=\spinc,\Mas(\psi)=0,n_z(\psi)=0\}}
\left(\#\ModFlow(\psi)\right)\w.
$$
Again, this is a finite sum under the weak admissibility hypothesis.

\begin{theorem}
\label{thm:HolTriangles} 
Let $(\Sigma,\alphas,\betas,\gammas,z)$ be a pointed Heegaard
triple-diagram, which is strongly $\spinc$-admissible for some $\SpinC$
structure $\spinc$ over the underlying  four-manifold $X$. 
Then the sum on
the right-hand-side of Equation~\eqref{eq:DefTriangle} is finite,
giving rise to a $U$-equivariant chain map which also induces maps on homology:
\begin{eqnarray*}
F^{\infty}_{\alpha,\beta,\gamma}(\cdot,\spinc_{\alpha,\beta,\gamma})
\colon \HF^{\infty}(Y_{\alpha,\beta},\spinc_{\alpha,\beta})
\otimes \HF^{\infty}(Y_{\beta,\gamma},\spinc_{\beta,\gamma})
&\longrightarrow & \HF^{\infty}(Y_{\alpha,\gamma},\spinc_{\alpha,\gamma})
\\
F^{\leq 0}_{\alpha,\beta,\gamma}(\cdot,\spinc_{\alpha,\beta,\gamma})
\colon \HF^{\leq 0}(Y_{\alpha,\beta},\spinc_{\alpha,\beta})
\otimes \HF^{\leq 0}(Y_{\beta,\gamma},\spinc_{\beta,\gamma})
&\longrightarrow & \HF^{\leq 0}(Y_{\alpha,\gamma},\spinc_{\alpha,\gamma}).
\end{eqnarray*}
The induced ($U$-equivariant) chain map
$$
f^+_{\alpha,\beta,\gamma}(\cdot,\spinc_{\alpha,\beta,\gamma})\colon
\CFp(Y_{\alpha,\beta},\spinc_{\alpha,\beta})
\otimes CF^{\leq 0}(Y_{\beta,\gamma},\spinc_{\beta,\gamma})
\longrightarrow \CFp(Y_{\alpha,\gamma},\spinc_{\alpha,\gamma}) 
$$
gives a well-defined  chain map when the triple diagram is only weakly
admissible, and the Heegaard diagram $(\Sigma,\betas,\gammas,z)$ is
strongly admissible for $\spinc_{\beta,\gamma}$. 
In fact, the induced map
$$
{\widehat f}_{\alpha,\beta,\gamma}(\cdot,\spinc_{\alpha,\beta,\gamma})\colon
\CFa(Y_{\alpha,\beta},\spinc_{\alpha,\beta})
\otimes \CFa(Y_{\beta,\gamma},\spinc_{\beta,\gamma}) \longrightarrow
\CFa(Y_{\alpha,\gamma},\spinc_{\alpha,\gamma})
$$
gives a well-defined chain map when the diagram is weakly admissible.
There are induced maps on homology:
\begin{eqnarray*}
{\widehat F}_{\alpha,\beta,\gamma}(\cdot,\spinc_{\alpha,\beta,\gamma})
\colon \HFa(Y_{\alpha,\beta},\spinc_{\alpha,\beta})
\otimes \HFa(Y_{\beta,\gamma},\spinc_{\beta,\gamma})
&\longrightarrow & \HFa(Y_{\alpha,\gamma},\spinc_{\alpha,\gamma}), \\
F^+_{\alpha,\beta,\gamma}(\cdot,\spinc_{\alpha,\beta,\gamma})
\colon \HFp(Y_{\alpha,\beta},\spinc_{\alpha,\beta})
\otimes \HF^{\leq 0}(Y_{\beta,\gamma},\spinc_{\beta,\gamma})
&\longrightarrow& \HFp(Y_{\alpha,\gamma},\spinc_{\alpha,\gamma}),
\end{eqnarray*}
the latter of which is also $U$-equivariant.
\end{theorem}

\begin{proof}
The fact that $f^\infty_{\alpha,\beta,\gamma}$ is a chain map follows
by counting ends of one-dimensional moduli spaces of
holomorphic triangles (compare ~\cite{McDuffSalamon}). Fix $\x\in\Ta\cap\Tb$,
$\y\in\Tb\cap\Tc$, $\w\in\Ta\cap\Tc$,
and consider moduli spaces of holomorphic triangles
$\ModFlow(\psi)$ where $\psi\in\pi_2(\x,\y,\w)$,
$\spinc_z(\psi)=\spinc$, and
$\Mas(\psi)=1$. The ends of this moduli space are modeled on:
\begin{eqnarray*}
&\left(\coprod_{\x'\in\Ta\cap\Tb}
\coprod_{\phi_{\alpha,\beta}\csum\psi_{\alpha,\beta}=\psi}
\UnparModFlow(\phi_{\alpha,\beta})\times
\ModFlow(\psi_{\alpha,\beta}) \right) & \\
&\coprod& \\
&\left( \coprod_{\y'\in\Tb\cap\Tc}
\coprod_{\phi_{\beta,\gamma}\csum\psi_{\beta,\gamma}=\psi}
\UnparModFlow(\phi_{\beta,\gamma})\times \ModFlow(\psi_{\beta,\gamma}) \right)& \\
&\coprod& \\
&\left( \coprod_{\w'\in\Ta\cap\Tc}
\coprod_{\phi_{\alpha,\gamma}\csum\psi_{\alpha,\gamma}=\psi}
\UnparModFlow(\phi_{\alpha,\gamma})\times
\ModFlow(\psi_{\alpha,\gamma}) \right). & \\
\end{eqnarray*}
In the above expression, the pairs of
homotopy classes $\phi_{\alpha,\beta}$ and $\psi_{\alpha,\beta}$ range
over  $\phi_{\alpha,\beta}\in\pi_2(\x,\x')$ and
$\psi_{\alpha,\beta}\in\pi_2(\x',\y,\w)$
with
$\Mas(\phi_{\alpha,\beta})=1$, 
$\Mas(\psi_{\alpha,\beta})=0$,
$\phi_{\alpha,\beta}\csum\psi_{\alpha,\beta}=\psi$ (with analogous
conditions for the $\phi_{\beta,\gamma}\in\pi_2(\y,\y')$ and 
$\phi_{\alpha,\gamma}\in\pi_2(\w',\w)$).
Counted with signs, the first two unions
give the
$[\w,i+j-n_z(\psi)]$-coefficient of
$f^{\infty}_{\alpha,\beta,\gamma}\circ \partial([\x,i]\otimes [\y,j])$
(using the natural differential on the tensor product), while the last
gives the $[\w,i+j-n_{z}(\psi)]$-coefficient of $\partial\circ
f^{\infty}([\x,i]\otimes[\y,j])$. 

Recall that if $\psi$ has a holomorphic representative, then
$n_z(\psi)\geq 0$. Thus, $f^{\infty}$ maps the subcomplex 
$$\CFme(Y_{\alpha,\beta},\spinc_{\alpha,\beta})\otimes
\CFme(Y_{\beta,\gamma},\spinc_{\beta,\gamma})\subset 
\CFinf(Y_{\alpha,\beta},\spinc_{\alpha,\beta})\otimes
\CFinf(Y_{\beta,\gamma},\spinc_{\beta,\gamma})$$ 
into $\CFme(Y_{\alpha,\gamma},\spinc_{\alpha,\gamma})$. Similarly, 
$f^+_{\alpha,\beta,\gamma}$ as above also gives a chain map.

The $U$-equivariance 
$$f^\infty_{\alpha,\beta,\gamma}\left(U\left([\x,i]\right)\otimes [\y,j]\right)= U\cm
f^\infty_{\alpha,\beta,\gamma}\left([\x,i]\otimes [\y,j]\right)$$ (and indeed for
the other induced maps, where stated) follows immediately from the
definitions.
\end{proof}

Now familiar arguments can be used to establish invariance properties
of these maps, as in the following:

\begin{prop}
\label{prop:TrianglesJIndep}
The maps on homology listed in Theorem~\ref{thm:HolTriangles} are
independent of the choice of family $J$ 
(and underlying complex structure $\sj$ over $\Sigma$)
used in its definition.
\end{prop}

\begin{proof}
Fix first the complex structure $\sj$ over $\Sigma$.
Consider a one-parameter variation family of maps
$J_\tau$ from $\Delta$
into the space of almost-complex structures over $\Sym^{g}(\Sigma)$,
where $\tau$ is a real parameter $\tau\in[0,1]$ (which are perturbations
of the symmetrized complex structure $\Sym^g(\sj)$ over $\Sym^g(\Sigma)$). 
We 
write down the case
of $\CFinfty$; the other homology theories work the same way, with
only notational changes.
Consider the map
$$H^\infty \colon \CFinfty(Y_{\alpha,\beta},\spinc_{\alpha,\beta})
\otimes \CFinfty(Y_{\beta,\gamma},\spinc_{\beta,\gamma})
\longrightarrow \CFinfty(Y_{\alpha,\gamma},\spinc_{\alpha,\gamma})$$
defined by
$$H^\infty([\x,i]\otimes [\y,j],\spinc) = \sum_{\w\in\Ta\cap\Tc}
\sum_{\{\psi\in\pi_2(\x,\y,\w)\big|
\spinc_z(\psi)=\spinc,\Mas(\psi)=-1\}}
\#\left(\bigcup_{\tau\in[0,1]}\ModFlow_{J_\tau}(\psi)\right)[\w,i+j-n_z(\psi)].$$
Now, the ends of 
$$\coprod_{\{\psi\in\pi_2(\x,\y,\w)\big|\spinc_z(\psi)=\spinc, \Mas(\psi)=0\}}
\left(\bigcup_{\tau\in[0,1]}\ModFlow_{J_\tau}(\psi)\right)$$
count 
$$f^\infty_{J_0}([\x,i]\otimes[\y,j];\spinc)-f^\infty_{J_1}([\x,i]\otimes[\y,j];\spinc)+
\partial\circ H^{\infty}([\x,i]\otimes[\y,j]) - H^{\infty}\circ
\partial([\x,i]\otimes[\y,j]);$$
i.e. $f^{\infty}_{J_0}$ and $f^{\infty}_{J_1}$ are chain homotopic.

The above argument shows that if we fix a map $J$ into the space of
almost-complex structures over $\Sym^g(\Sigma)$, with specified
behaviour near the vertices of the triangle, then the induced map is
uniquely specified up to chain homotopy. Now, we investigate
invariance as one changes the one-parameter families of paths used at
the corners.  To this end, we claim that the following diagram
$$\begin{CD}
(\CFinf(Y_{\alpha,\beta},\spinc_{\alpha,\beta}),\delta^\infty_{J_s(0)})
\otimes \CFinf(Y_{\beta,\gamma},\spinc_{\alpha,\beta}) 
@>f_{\alpha,\beta,\gamma}>> (\CFinf(Y_{\alpha,\gamma}),\delta^\infty_{J_s(0)}) \\
@V{\Phi^{\infty}_{\alpha,\beta}\otimes \Id} VV	@VV{\Phi^{\infty}_{\alpha,\gamma}}V \\
(\CFinf(Y_{\alpha,\beta},\spinc_{\alpha,\beta}),\delta^\infty_{J_s(1)})
\otimes \CFinf(Y_{\beta,\gamma},\spinc_{\alpha,\beta})
@>f_{\alpha,\beta,\gamma}>> (\CFinf(Y_{\alpha,\gamma}),\delta^\infty_{J_s(1)}),
\end{CD}$$
commutes up to homotopy. Here, the maps maps $\Phi^{\infty}$ are the
homomorphisms induced by variations of these paths (which we
encountered in the proof of Theorem~\ref{thm:IndepCxStruct}; note that
we have now dropped the variation of path from the notation, and
replaced it with the an indication of the Heegaard diagram we are
working with). The homotopy now is obtained by counting triangles with
in a formally $-1$-dimensional space, but which are
$J_\tau$-holomorphic with respect to any $J_\tau$ for $\tau$ in a
one-parameter family. We do not spell out the details here: a very
similar argument is given in the next proposition. Similar remarks hold
as we vary the boundary behaviour at the other vertices.

As a final remark, since the induced maps are invariant under small
perturbations of the family $J$, it follows also that the induced map
is independent under variations of the complex structure $\sj$ over
$\Sigma$.
\end{proof}

\begin{prop}
\label{prop:TriangleIsotopyInvariance}
The maps on homology listed in Theorem~\ref{thm:HolTriangles} 
are invariant under isotopies of the $\alphas$,
$\betas$, and $\gammas$ preserving all the admissibility hypotheses.
\end{prop}

\begin{proof}
We begin with isotopies of the $\alphas$.
As in the proof of isotopy invariance of Floer homologies, we let
$\Psi_\tau$ be an isotopy (induced from an exact Hamiltonian isotopy of the 
$\alphas$ in $\Sigma$), and we consider moduli spaces with dynamic
boundary conditions.
Specifically, let
$E_{\alpha}\colon \R \longrightarrow \Delta$
be a parameterization of the edge $e_{\alpha}$, with 
\begin{eqnarray*}
\lim_{t\goesto -\infty}E_{\alpha}(t)=v_{\gamma} &{\text{and}}&
\lim_{t\goesto +\infty}E_{\alpha}(t)=v_{\beta} 
\end{eqnarray*}
Consider moduli spaces indexed by a real parameter $\tau\in\R$:
$$
\ModFlow_\tau=
\left\{
u\colon\Delta\longrightarrow \Sym^{g}(\Sigma) \Bigg|
\begin{array}{c}
u\circ E_{\alpha}(t)\in \Psi_{t+\tau}(\Ta) \\
u(e_{\beta})\subset \Tb,
u(e_{\gamma})\subset \Tc
\end{array}
\right\},$$
and divide them into homotopy classes $\pi_2^{\Psi_t}(\x,\y,\w)$, with
$\x\in\Ta\cap\Tb$, $\y\in\Tb\cap\Tc$, $\w\in\Tc\cap \Psi_1(\Ta)$.

Note that if $\Mas(\psi)=-1$, then
$\bigcup_{\tau\in\R}\ModFlow_\tau(\psi)$ is generically a compact
zero-dimensional manifold, so we can define 
$$H^\infty([\x,i]\otimes
[\y,j];\spinc)=\sum_{\w}\sum_{\{\psi\in\pi_2^{\Psi_t}(\x,\y,\w) |
\Mas(\psi)=-1,\spinc_z(\psi)=\spinc\}}\left(\#\bigcup_{\tau\in \R}\ModFlow_\tau(\psi)\right)
[\w,i+j-n_z(\psi)].$$

Fix, now, 
any homotopy class $\psi\in\pi_2^{\Psi_t}(\x,\y,\w)$ with 
$\spinc_z(\psi)=\spinc$ and $\Mas(\psi)=0$, and consider the one-manifold
$$\bigcup_{\tau\in\R}\ModFlow_\tau(\psi).$$
This has ends as $\tau\goesto\pm \infty$, which are modeled on
$$\left(\coprod_{\phi_{\alpha,\beta}\csum\psi_{\alpha,\beta}=\psi}\Mod_{\Psi_t}(\phi_{\alpha,\beta})\times 
					\ModFlow(\psi_{\alpha,\beta})\right)
\coprod
\left(\coprod_{\phi_{\alpha,\gamma}\csum\psi_{\alpha,\gamma}=\psi}
				\Mod_{\Psi_t}(\phi_{\alpha,\gamma})\times \ModFlow(\psi_{\alpha,\gamma})\right),$$
where the first union is over all 
$\x'\in \Psi_1(\Ta)\cap\Tb$ 
with $\spinc_z(\x')=\spinc_{\alpha,\beta}$,
$\phi_{\alpha,\beta}\in\pi_2^{\Psi_t}(\x,\x')$ (in the sense of Subsection~\ref{subsec:Isotopies}),
$\psi_{\alpha,\beta}\in\pi_2(\x',\y,\w)$, and $\Mas(\phi_{\alpha,\beta})=\Mas(\psi_{\alpha,\beta})=0$
(with analogous conditions on the second union). There are also ends of the form
$$\begin{array}{c}
\left(\coprod_{\phi_{\alpha,\beta}\csum\psi_{\alpha,\beta}=\psi}\UnparModFlow(\phi_{\alpha,\beta})\times 
\left(\bigcup_{\tau\in\R} \ModFlow_\tau(\psi_{\alpha,\beta})\right)\right)\\
\coprod \\
\left(\coprod_{\phi_{\beta,\gamma}\csum\psi_{\beta,\gamma}=\psi}\UnparModFlow(\phi_{\beta,\gamma})\times 
\left(\bigcup_{\tau\in\R} \ModFlow_\tau(\psi_{\beta,\gamma})\right)\right) \\
\coprod \\
\left(\coprod_{\phi_{\alpha,\gamma}\csum\psi_{\gamma,\alpha=\psi}}\UnparModFlow(\phi_{\alpha,\gamma})\times 
\left(\bigcup_{\tau\in\R} \ModFlow_\tau(\psi_{\alpha,\gamma})\right)\right),
\end{array}$$
where the first union is over all $\x'\in \Ta\cap\Tb$ in the same
equivalence class as $\x$, $\phi_{\alpha,\beta}\in\pi_2(\x,\x')$ (in
the sense of Subsection~\ref{subsec:Isotopies}),
$\psi_{\alpha,\beta}\in\pi_2^{\Psi_t}(\x',\y,\w)$, and
$\Mas(\phi_{\alpha,\beta})=1$ and $\Mas(\phi_{\alpha,\beta})=1$ while
$\Mas(\psi_{\alpha,\beta})=-1$ (with analogous conditions over the
other two unions). Counting ends with sign, we get that $$
\Gamma_{\alpha,\alpha',\gamma}\circ f_{\alpha,\beta,\gamma} -
f_{\alpha',\beta,\gamma}\circ \Gamma_{\alpha,\alpha'\beta} = \partial\circ H
- H\circ
\partial,$$
where 
\begin{eqnarray*}
\Gamma_{\alpha,\alpha',\beta}\colon \CFinf(Y_{\alpha,\beta},\spinc_{\alpha,\beta}) \longrightarrow \CF(Y_{\alpha',\beta},\spinc_{\alpha',\beta})
&{\text{and}}&
\Gamma_{\alpha,\alpha',\gamma}\colon \CF(Y_{\alpha,\gamma},\spinc_{\alpha,\gamma}) \longrightarrow \CF(Y_{\alpha',\gamma},\spinc_{\alpha',\gamma})
\end{eqnarray*}
are the chain maps induced by the isotopy $\Psi_t$, as constructed in
Subsection~\ref{subsec:Isotopies} (note that here we have suppressed
the isotopy $\Psi_t$ from the notation.

Isotopies of the $\gammas$ work the same way; 
we now set up isotopies of 
the $\betas$. Consider moduli spaces indexed by a
real $\tau\in [0,\infty)$ 
$$
\ModFlow_\tau=
\left\{
u\colon\Delta\longrightarrow \Sym^{g}(\Sigma) \Bigg|
\begin{array}{c}
u\circ E_{\gamma}(t)\in \Psi_{t+\tau}^{-1}\circ 
\Phi_{t-\tau}(\Tc) \\
u(e_\alpha)\subset \Ta,
u(e_{\gamma})\subset \Tc
\end{array}
\right\}.$$
These moduli spaces partition according to 
homotopy
classes $\psi\in\pi_2(\x,\y,\w)$ (with $\x\in\Ta\cap\Tb$,
$\y\in\Tb\cap\Tc$, $\w\in\Ta\cap\Tc$).
Note that for $\tau=0$, this is the usual moduli space
for holomorphic triangles for $\alphas$, $\betas$, and
$\gammas$.
Again, when $\Mas(\psi)=0$, the union $\bigcup_{\tau\in[0,\infty)}
\ModFlow_\tau(\psi)$ is generically a compact, zero-dimensional manifold, and we can define
Note that if $\Mas(\psi)=-1$, then
$\bigcup_{\tau\in[0,\infty)}\ModFlow_\tau(\psi)$ is generically a compact
zero-dimensional manifold, so we can define 
$$H^\infty([\x,i]\otimes
[\y,j];\spinc)=\sum_{\w}\sum_{\{\psi\in\Psi_t |
\Mas(\psi)=-1,\spinc_z(\psi)=\spinc\}}\left(\#\bigcup_{\tau\in[0,\infty)}\ModFlow_\tau(\psi)\right)
[\w,i+j-n_z(\psi)].$$

Fix a homotopy class $\psi\in\pi_2(\x,\y,\w)$ with
$\spinc_z(\psi)=\spinc$ and $\Mas(\psi)=0$, and consider the ends of
$$\bigcup_{\tau\in[0,\infty)}\ModFlow_\tau(\psi).$$
The ends as $\tau\goesto\infty$ are modeled on
$$\bigcup_{\phi_{\alpha,\beta}\csum\psi\csum\phi_{\beta,\gamma}}
\ModFlow_{\Psi_t}(\phi_{\alpha,\beta})
\times\ModFlow_{\Psi_t}(\phi_{\beta,\gamma})\times 
\ModFlow(\psi_{\alpha,\beta',\gamma}),$$
where the union is over all $\x'\in\Ta\cap\Psi_1(\Tb)$,
$\y'\in\Psi_1(\Tb)\cap\Tc$ with 
$\spinc_z(\y')\in\spinc_{\beta',\gamma}$,
$\phi_{\alpha,\beta}\in\pi_2^{\Psi_t}(\x,\x')$,
$\phi_{\beta,\gamma}\in\pi_2^{\Psi_t}(\y,\y')$, 
and $\psi_{\alpha,\beta',\gamma}\in\pi_2(\x,\y',\w)$ with $\Mas(\phi_{\alpha,\beta})=\Mas(\phi_{\beta,\gamma})
=\Mas(\psi_{\alpha,\beta,\gamma})=0$.
Counting these ends with sign, we get a contribution of 
$$f_{\alpha,\beta',\gamma}\circ (\Gamma_{\alpha,\beta,\beta'}([\x,i])\otimes 
\Gamma_{\beta,\beta',\gamma}([\y,j])),$$
while the end as $\tau\goesto\infty$ corresponds simply to
$f_{\alpha,\beta,\gamma}([\x,i]\otimes [\y,j])$. There are other ends
as before, whose contribution is $$\partial \circ H^\infty -
H^\infty\circ \partial.$$ Thus, we have exhibited a chain homotopy
from $f_{\alpha,\beta,\gamma}$ with $f_{\alpha,\beta',\gamma}\circ 
(\Gamma_{\alpha,\beta,\beta'}\otimes 
\Gamma_{\beta,\beta',\gamma})$.
\end{proof}

\subsection{Associativity of holomorphic triangles}
\label{subsec:Assoc}

The map induced by triangles satisfies an associativity property, on
the level of homology. As is familiar in Lagrangian Floer homology,
the chain homotopies required for the associativity is provided by a
count of holomorphic squares, see~\cite{QAssoc}, \cite{FOOO},
\cite{Silva}.

\begin{figure}
\mbox{\vbox{\epsfbox{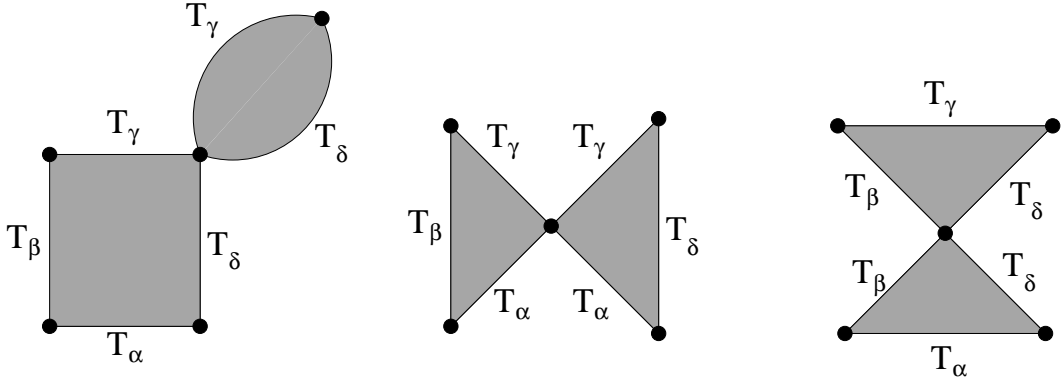}}}
\caption{\label{fig:deg} 
{\bf{Ends of the moduli spaces holomorphic maps of rectangles.}}
The two degenerations on the right correspond to degenerations
of the conformal type of the underlying rectangle. (There are
three additional degenerations like the one on the left, corresponding
to the other three vertices.)}
\end{figure}

Loosely speaking, this count is done as follows. Let $\Rect$ denote
the ``rectangle'': the unit disk with four marked boundary points
which we denote $a_0$, $a_1$, $a_2$, and $a_3$ (labeled in clockwise
order). Observe that the space of conformal structures on the
rectangle, $\Mod(\Rect)$, is identified with $\R$ under the map which
is the logarithmic difference between the length and the width. We
would like to count pseudo-holomorphic maps of the rectangle in
$\Sym^g(\Sigma)$ (without fixing the conformal structure of the
domain). We wish to consider moduli spaces of pseudo-holomorphic
Whitney rectangles with formal dimension one. Some ends of these
moduli spaces are modeled on flowlines breaking off at the corners,
but there is another type of end not encountered before in the counts
of trianges, arising from the non-compactness of $\Mod(\Rect)\cong
\R$. As this parameter goes to $\pm
\infty$, the corresponding rectangle breaks up conformally into a pair
of triangles meeting at a vertex (in two different ways, depending on
which end we are considering), as illustrated in Figure~\ref{fig:deg}.
This is how a count of holomorphic squares induces a chain homotopy
between two different compositions of holomorphic triangle counts..

The ingredients required in this construction are: the role of
$\SpinC$ structures in the construction, admissibility hypotheses
required to make the holomorphic rectangle counts to be finite,
compatibility with holomorphic triangle counts. With these components in
place, the proof of associativity proceeds in the same way as it does in 
the usual Lagrangian Floer homology.

\subsubsection{$\SpinC$ structures on rectangles}

Fix a pointed Heegaard quadruple
$(\Sigma,\alphas,\betas,\gammas,\deltas,z)$, and let
$X_{\alpha,\beta,\gamma,\delta}$ be the corresponding cobordism. We
have, of course, restriction maps:
$$\SpinC(X_{\alpha,\beta,\gamma,\delta})\longrightarrow
\SpinC(X_{\alpha,\beta,\gamma})\times
\SpinC(X_{\alpha,\gamma,\delta}),$$
which correspond to splitting the cobordism along an embedded copy of
$Y_{\alpha,\gamma}$. 
There is a subgroup 
$\delta H^1(Y_{\alpha,\gamma})\subset
H^2(X_{\alpha,\beta,\gamma,\delta})$, whose orbits on
$\SpinC(X_{\alpha,\beta,\gamma,\delta})$ are the fibers of this
restriction map.
Similarly, we have a restriction map
$$\SpinC(X_{\alpha,\beta,\gamma,\delta})\longrightarrow
\SpinC(X_{\alpha,\beta,\delta})\times
\SpinC(X_{\beta,\gamma,\delta}),$$
which corresponds to splitting along $Y_{\beta,\delta}$. 
In view of this, we will find it convenient to fix not  a single $\SpinC$ structure
over $X_{\alpha,\beta,\gamma,\delta}$, but rather a
$\delta H^1(Y_{\beta,\delta})+\delta H^1(Y_{\alpha,\gamma})$ orbit.

\subsubsection{Admissibility for rectangles}
\label{subsec:AdmissRect}

There are notions of admissibility for Heegaard quadruples (and, in
general, multi-diagrams), which generalize the corresponding notions
for triangles. For instance, a Heegaard quadruple
$(\Sigma,\alphas,\betas,\gammas,\deltas,z)$ is called {\em weakly
admissible} if every periodic domain which can be written as sums of
doubly-periodic domains for $Y_{\alpha,\beta}$, $Y_{\beta,\gamma}$,
$Y_{\gamma,\delta}$, and $Y_{\alpha,\delta}$ has both positive and
negative coefficients. Existence is achieved by winding, as in
Section~\ref{sec:Special}.

Strong admissibility requires fixing a $\delta
H^1(Y_{\beta,\delta})+\delta H^1(Y_{\alpha,\gamma})$-orbit ${\mathfrak S}$
of a fixed $\SpinC$ structure over $X_{\alpha,\beta,\gamma,\delta}$.
We say that a Heegaard quadruple is {\em strongly admissible} for the
orbit ${\mathfrak S}$ if for each $\SpinC$ structure $\spinc\in{\mathfrak S}$
and each quadruply-periodic domain which can be written as a sum of
doubly-periodic domains: 
\begin{equation}
\label{eq:DecomposePerDom}
\PerDom=\sum_{\{\xi,\eta\}\subset
\{\alpha,\beta,\gamma,\delta\}}
\cald_{\xi,\eta}
\end{equation}
with the property that $$\sum \langle
c_1(\spinc_{\xi,\eta}),H(\cald_{\xi,\eta})\rangle = 2n\geq 0,$$
(i.e. where here $\spinc_{\xi,\eta}$ denotes the corresponding
restriction of $\spinc$), it follows that some local multiplicity of
$\PerDom$ is strictly greater than $n$.  Note that this notion
involves the orbit ${\mathfrak S}$ only through its restrictions to
the six three-manifolds $Y_{\xi,\eta}$, for all subsets
$\{\xi,\eta\}\subset\{\alpha,\beta,\gamma,\delta\}$. Note also that if
a Heegaard quadruple is strongly  ${\mathfrak S}$-admissible, then the
associated Heegaard diagrams for $Y_{\xi,\eta}$ (for all
$\{\xi,\eta\}\subset\{\alpha,\beta,\gamma,\delta\}$) are automatically
strongly $\spinc_{\xi,\eta}$-admissible.

Since the orbit ${\mathfrak S}$ {\em a priori} may
contain infinitely many $\SpinC$ structures,
it might be impossible to achieve strong admissibility. However, suppose that
the Heegaard quadruple satisfies the hypothesis that:
\begin{eqnarray}
\label{eq:GoodDecomposition}
\delta H^1(Y_{\beta,\delta})|Y_{\alpha,\gamma}=0
&{\text{and}}&
\delta H^1(Y_{\alpha,\gamma})|Y_{\beta,\delta}=0.
\end{eqnarray}
In this case, if we fix any
$\{\xi,\eta\}\subset\{\alpha,\beta,\gamma,\delta\}$, the restriction
$\spinc_{\xi,\eta}$ to $Y_{\xi,\eta}$ of any $\spinc\in{\mathfrak S}$
is independent of the choice of $\spinc$. Equivalently, if we choose
any quadruply-periodic domain which can be written as a
sum of doubly-periodic domains as in Equation~\eqref{eq:DecomposePerDom}
we have that 
$$\sum_{\{\xi,\eta\}\subset\{\alpha,\beta,\gamma,\delta\}}
\langle c_1(\spinc_{\xi,\eta}), H(\cald_{\xi,\eta}) \rangle $$
is a function of the periodic domain ${\mathfrak Q}$ and the orbit
${\mathfrak S}$ (i.e. it is independent of the choice of
$\SpinC$ structure
$\spinc\in{\mathfrak S}$). Thus, the proof of
Lemma~\ref{lemma:spincRealized} adapts immediately to show that
${\mathfrak S}$-strong admissibility can be achieved.

\subsubsection{Compatibility with counts of squares}

Having established the necessary admissibility requirements for
defining counts of pseudo-holomorphic rectangles, we need two more 
ingredients: transversality and orientations.

We must set up the transverse perturbations for
$J$-holomorphic rectangles in a manner which is
compatible with the $J_s$ paths used over strips and the
``nearly-symmetric families'' used on the
triangles. To do this, we must fix a map $$J\colon \ModRect \times
\Rect \longrightarrow {\mathcal U},$$ 
with certain properties. First, we arrange for ${\mathcal U}$ to
consist of $(\sj,\eta,V)$-nearly symmetric almost-complex structures,
with $V$ chosen to contain $\bigcup_{z_i}\{z_i\}\times
\Sym^{g-1}(\Sigma)$, where the $\{z_i\}$ are a finite collection of
points, one from each domain in
$\Sigma-\alpha_1-...-\alpha_g-\beta_1-...-\beta_g-\gamma_1-...-\gamma_g-\delta_1-...-\delta_g$. Moreover,
${\overline V}$ is chosen to be disjoint from all four tori $\Ta$,
$\Tb$, $\Tc$, and $\Td$.  To arrange for compatibility with strips, we
fix path $J_s^{(i)}$ for $i=0,...,3$, and assume that for each
conformal structure in $\ModRect$, the map $J$ agrees with the
$J_s^{(i)}$ path near the $i^{th}$ vertex (under a conformal
identification between a neighborhood of each corner and strip).  To
arrange for compatibility with triangles, we restrict the behaviour of
$J$ as the parameter in $\ModRect$ goes to infinity.  As the parameter
in $\ModRect$ approaches one extreme, the complex structure over
$\Rect$ is identified with $[0,1]\times [-T,T]$ for some large $T$;
when the parameter approaches the other extreme, the complex structure
is identified with $[-T,T]\times [0,1]$.  For the first degeneration
into triangles, we assume that $J|[0,1]\times [-T,0]$ agrees, after
translation by $T$ on the second factor, with the restriction of a
some fixed admissible family on $[0,1]\times[0,+\infty)\cong
\Delta$, and similarly, that $J|[0,1]\times [0,T]$ agrees, after
translation by $-T$ on the second factor, with the restriction
of an admissible family on $[0,1]\times(-\infty,0]\cong \Delta$. We
also assume analogous properties for the other degeneration.

Adapting the transversality proof, one can see that for generic
$J\colon \ModFlow(\Rect)\times\Rect\longrightarrow {\mathcal U}$
satisfying the compatibilities near infinity as above, the
corresponding moduli spaces with formal dimension $\leq 1$
are smooth. 

When working with quadruples, and $\Z$ coefficients, we need yet
another generalization of the notion of coherent systems of
orientations. We now fix a $\delta H^1(Y_{\beta,\delta})+\delta
H^1(Y_{\alpha,\gamma})$ orbit in
$\SpinC(X_{\alpha,\beta,\gamma,\delta})$, which we denote ${\mathfrak
S}$. A coherent system of orientations for ${\mathfrak S}$, then, is a
collection of non-vanishing sections indexed by subsets
$\{\xi_1,...,\xi_\ell\}\subset \{\alpha,\beta,\gamma,\delta\}$ with
$\ell=2,3,4$,
$\orient_{\xi_1,...,\xi_\ell}(\phi_{\xi_1,...,\xi_\ell})$, for the
determinant line bundle defined over the homotopy class of polygons
$\phi_{\xi_1,...,\xi_{\ell}}$ (i.e. this can be a Whitney disk,
triangle, or square) representing the restriction of some $\spinc\in{\mathfrak S}$ to 
$Y_{\xi_1,\xi_2}$ when $\ell=2$ or $X_{\xi_1,...,\xi_{\ell}}$ if
$\ell=3,4$. These are required to be compatible with the gluings in the sense that
$$\orient_{\xi_1,...,\xi_\ell}(\phi_{\xi_1,...,\xi_\ell})\wedge 
\orient_{\eta_1,...,\eta_m}(\phi_{\eta_1,...,\eta_m})
= \orient_{\xi_1,...,\xi_\ell,\eta_1,...,\eta_m}
(\phi_{\xi_1,...,\xi_\ell}\csum\phi_{\eta_1,...,\eta_m}),
$$ under gluing
maps which are defined whenever we have subsets
$\{\xi_1,...,\xi_\ell\}$ and $\{\eta_1,...,\eta_m\}$ with two
elements, say $\xi_1$ and $\xi_2$, in common, and for which the
polygons $\phi_{\xi_1,...,\xi_\ell}$ and $\phi_{\eta_1,...,\eta_m}$
meet in a single intersection point for ${\mathbb T}_{\eta_1}\cap {\mathbb T}_{\eta_2}$. 

The following construction of coherent orientation systems will suffice for our
purposes:

\begin{prop}
\label{prop:BuildOrientationSystem}
Let $X_{\alpha,\beta,\gamma,\delta}$ belong to a Heegaard quadruple
for which 
\begin{equation}
\label{eq:TrivialCoboundary}
\Big(\delta H^1(Y_{\beta,\delta})\subseteq
H^2(X_{\alpha,\beta,\gamma,\delta},\partial
X_{\alpha,\beta\gamma,\delta})\Big)=0.
\end{equation}
Fix a $\delta H^1(Y_{\alpha,\gamma})$-orbit ${\mathfrak S}$ in the set
of $\SpinC$ structures over $X_{\alpha,\beta,\gamma,\delta}$.  Fix
also orientation systems $\orient_{\alpha,\beta,\gamma}$ and
$\orient_{\alpha,\gamma,\delta}$ for the two corresponding Heegaard
subtriples, which induce the same orientation over
$Y_{\alpha,\gamma}$.  Then, the given orientation systems can be
completed compatibly to give an orientation system over the Heegaard
quadruple $X_{\alpha,\beta,\gamma,\delta}$.
\end{prop}
 
\begin{proof}
The given data already allow us to orient all squares representing
${\mathfrak S}$. Specifically, any such square can be decomposed as a
sum of Whitney triangles
$\varphi_{\alpha,\beta,\gamma,\delta}=\psi_{\alpha,\beta,\gamma}
+\psi_{\alpha,\gamma,\delta}$, where here the $\psi$ denote homotopy
classes of Whitney triangles for the subtriples indicated in the
notation, and which represent the $\SpinC$-equivalence class induced
from ${\mathfrak S}$. We then define the orientation for the given
square by
$$\orient(\varphi_{\alpha,\beta,\gamma,\delta})
=\orient(\psi_{\alpha,\beta,\gamma})\wedge
\orient(\psi_{\alpha,\gamma,\delta})$$
(note that we have suppressed the subscripts from $\orient$ indicating
the diagram used for the polygon, as our
notation for the homotopy classes already carries this
information). Since the induced orientations over $Y_{\alpha,\gamma}$
are compatible, this definition is independent of the decomposition.
In particular, if it follows that for any four Whitney disks
$\phi_1,...,\phi_4$, for the Heegaard pairs on the boundary
($Y_{\alpha,\beta}$, $Y_{\beta,\gamma}$, $Y_{\gamma,\delta}$, and
$Y_{\alpha,\delta}$ respectively), we already have that
\begin{equation}
\label{eq:WhatWeHave}
\orient(\phi_{\alpha,\beta,\gamma,\delta}+\phi_1+...+\phi_4)=
\orient(\phi_{\alpha,\beta,\gamma,\delta})\wedge\orient(\phi_1)\wedge...\wedge\orient(\phi_4).
\end{equation}

Next, fix a pair of triangles
$\psi_{\alpha,\beta,\delta}\in\pi_2(\x_0,{\mathbf u}_0,\w_0)$,
$\psi_{\beta,\gamma,\delta}\in\pi_2(\y_0,{\mathbf v}_0,{\mathbf u}_0)$ for the
corresponding Heegaard triples, and fix an arbitrary orientation
$\orient(\psi_{\alpha,\beta,\delta})$, and then orient
$\psi_{\beta,\gamma,\delta}$ so that
$$\orient(\psi_{\alpha,\beta,\delta})\wedge\orient(\psi_{\beta,\gamma,\delta})
=\orient(\psi_{\alpha,\beta,\delta}+\psi_{\beta,\gamma,\delta})$$
Given this, we can now orient all $\phi_{\beta,\delta}
\in\pi_2({\mathbf u}_0,{\mathbf u}_0)$: given any such $\phi_{\beta,\delta}$, in view of
Equation~\eqref{eq:TrivialCoboundary} we can write
$$\psi_{\alpha,\beta,\delta}+\phi_{\beta,\delta}+\psi_{\alpha,\gamma,\delta}
=
\psi_{\alpha,\beta,\delta}+\psi_{\alpha,\gamma,\delta} +
\phi_1+\phi_2+\phi_3+\phi_4,$$ where $\phi_1\in\pi_2(\x_0,\x_0)$,
$\phi_2\in\pi_2(\y_0,\y_0)$, $\phi_3\in\pi_2({\mathbf v}_0,{\mathbf v}_0)$, and
$\phi_4\in\pi_2(\w_0,\w_0)$. We now define $\orient(\phi_{\beta,\gamma})$ to be
compatible with this decomposition, and observe that, according to
Equation~\eqref{eq:WhatWeHave}, the answer is
independent of the choices of $\phi_1,...,\phi_4$. We extend these
orientations on $\pi_2({\mathbf u}_0, {\mathbf u}_0)$ to an orientation
system for $Y_{\beta,\delta}$ in the usual manner. (Note that the isomorphism
class of this orientation system is uniquely determined.)

This orientation system of $Y_{\beta,\delta}$ and the 
fixed orientation $\orient(\psi_{\alpha,\beta,\delta})$ on the 
initial triangle $\psi_{\alpha,\beta,\delta}$ suffices to orient the
all the remaining triangles for the triples $X_{\alpha,\beta,\delta}$
in the $\SpinC$-equivalence class corresponding to the restriction of
${\mathfrak S}$: any other such triangle
$\psi'_{\alpha,\beta,\delta}$ can be written in the form
$$\psi'_{\alpha,\beta,\delta}=\psi_{\alpha,\beta,\delta}+\phi_{\alpha,\beta}+\phi_{\beta,\delta}+\phi_{\alpha,\delta},$$
and can be oriented correspondingly. The orientation obtained
in this manner is well-defined, owing to Equation~\eqref{eq:WhatWeHave}.
The same remarks hold for the remaining subtriple $X_{\beta,\gamma,\delta}$,
and since the induced orientation systems on the $Y_{\beta,\delta}$ 
agree, it is straightforward now to verify that the induced orientation
system on the $X_{\alpha,\beta,\delta}$ and $X_{\beta,\gamma,\delta}$ are
consistent with the orientation system on the squares in $X_{\alpha,\beta,\gamma,\delta}$, completing the construction.
\end{proof}

\subsubsection{The associativity theorem}

\begin{theorem}
\label{thm:Associativity}
Let $(\Sigma,\alphas,\betas,\gammas,\deltas,z)$ be a
pointed Heegaard quadruple which is strongly
${\mathfrak S}$-admissible, where ${\mathfrak S}$ is
a $\delta H^1(Y_{\beta,\delta})+ \delta
H^1(Y_{\alpha,\gamma})$-orbit in
$\SpinC(X_{\alpha,\beta,\gamma,\delta})$.
Then, we have
\begin{eqnarray*}
\lefteqn{\sum_{\spinc\in{\mathfrak S}}
\Fstar{\alpha,\gamma,\delta}
(\Fstar{\alpha,\beta,\gamma}
(\xi_{\alpha,\beta}\otimes
\theta_{\beta,\gamma}; 
\spinc_{\alpha,\beta,\gamma})\otimes 
\theta_{\gamma,\delta};\spinc_{\alpha,\gamma,\delta})} \\
&=&\sum_{\spinc\in{\mathfrak S}}
\Fstar{\alpha,\beta,\delta}(\xi_{\alpha,\beta}\otimes 
\Fleq{\beta,\gamma,\delta}
(\theta_{\beta,\gamma}\otimes\theta_{\gamma,\delta};\spinc_{\beta,\gamma,\delta});
\spinc_{\alpha,\beta,\delta}),
\end{eqnarray*}
where $\Fstar{}=\Finf{}$, $\Fp{}$ or $\Fm{}$,
$\xi_{\alpha,\beta}\in \HFstar(Y_{\alpha,\beta})$, 
$\theta_{\beta,\gamma}$ and $\theta_{\gamma,\delta}$ lie in
$\HFleq(Y_{\beta,\gamma})$ and $\HFleq(Y_{\gamma,\delta})$ respectively; also,
\begin{eqnarray*}
\lefteqn{\sum_{\spinc\in{\mathfrak S}}
\Fa{\alpha,\gamma,\delta}(\Fa{\alpha,\beta,\gamma}
(\xi_{\alpha,\beta}\otimes
\xi_{\beta,\gamma}; 
\spinc_{\alpha,\beta,\gamma})\otimes 
\xi_{\gamma,\delta};\spinc_{\alpha,\gamma,\delta})} \\
&=&\sum_{\spinc\in{\mathfrak S}}
\Fa{\alpha,\beta,\delta}(\xi_{\alpha,\beta}\otimes 
\Fa{\beta,\gamma,\delta}
(\xi_{\beta,\gamma}\otimes\xi_{\gamma,\delta};\spinc_{\beta,\gamma,\delta});
\spinc_{\alpha,\beta,\delta}), \\
\end{eqnarray*}
where now $\xi_{\alpha,\beta}$, $\xi_{\beta,\gamma}$, and
$\xi_{\alpha,\gamma}$ lie in $\HFa$ for the corresponding three-manifolds.
When working over $\Z$, we assume a consistent family of orientations for all the 
$\SpinC$ structures in ${\mathfrak S}$, used in the definitions of the maps on triangles.
\end{theorem}

\begin{proof}

We define a map
$$H^\infty(~\cdot~,{\mathfrak S})\colon \bigoplus_{\spinc\in{\mathfrak S}}\CFinf(Y_{\alpha,\beta},\spinc_{\alpha,\beta})\otimes 
\CFinf(Y_{\beta,\gamma},\spinc_{\beta,\gamma})\otimes 
\CFinf(Y_{\gamma,\delta},\spinc_{\gamma,\delta})\longrightarrow
\bigoplus_{\spinc\in{\mathfrak S}}\CFinf(Y_{\alpha,\delta})$$
by 
$$H^\infty([\x,i]\otimes[\y,j]\otimes[\w,k],{\mathfrak S})=\sum_{\p\in\Ta\cap\Td}
\sum_{\left\{\varphi\in\pi_2(\x,\y,\w,\p)\big|
\begin{tiny}\begin{array}{c}
\spinc_z(\varphi)\in{\mathfrak S}\\
\Mas(\varphi)=0
\end{array}
\end{tiny}
\right\}}\Big(\#\Mod(\varphi)\Big)[\p,i+j+k-n_z(\varphi)].$$

Note that above map is a finite sum by the strong admissibility requirement on
the Heegaard quadruple; indeed, it also implies that the 
sums appearing in the statement of the theorem are finite sums.

Counting ends of one-dimensional moduli spaces $\ModFlow(\varphi)$
with $\Mas(\varphi)=1$, we see that $H$ induces a chain homotopy
between the maps
$$\xi_{\alpha,\beta}\otimes\theta_{\beta,\gamma}\otimes\theta_{\gamma,\delta}\mapsto
\sum_{\spinc\in{\mathfrak S}}
f_{\alpha,\gamma,\delta}(f_{\alpha,\beta,\gamma}
(\xi_{\alpha,\beta}\otimes
\theta_{\beta,\gamma},\spinc_{\alpha,\beta,\gamma})\otimes \theta_{\gamma,\delta},\spinc_{\alpha,\gamma,\delta})$$
and 
$$\xi_{\alpha,\beta}\otimes\theta_{\beta,\gamma}\otimes\theta_{\gamma,\delta}\mapsto
\sum_{\spinc\in{\mathfrak S}}
f_{\alpha,\beta,\delta}(\xi_{\alpha,\beta}\otimes
f_{\beta,\gamma,\delta}(\theta_{\beta,\gamma}\otimes
\theta_{\gamma,\delta},\spinc_{\beta,\gamma,\delta}),\spinc_{\alpha,\beta,\delta}).$$
Again, the other cases are established in the same manner.
\end{proof}

\section{Handleslide invariance}
\label{sec:HandleSlides}

Our aim is to prove handleslide invariance of the Floer homology
groups.  In Subsection~\ref{subsec:STwoTimesSOne} we establish such a
result for the three-manifold $\#^g(S^2\times S^1)$ (equipped with a
particular $\SpinC$ structure) and a specific handleslide, by
explicitly calculating the Floer homologies for both Heegaard
diagrams. In Subsection~\ref{subsec:Naturality}, we use the
holomorphic triangle construction to transfer the result for this
specific three-manifold to a general three-manifold.

\subsection{The Floer homologies of $\#^g (S^2\times S^1)$.}
\label{subsec:STwoTimesSOne}
We give now a model calculation of a handleslide.

Let $\Sigma$ be an oriented surface of genus $g$ equipped with a
basepoint $z$.  Fix a $g$-tuple $\betas$ of attaching circles for a
handlebody bounding a genus $g$ surface $\Sigma$ (disjoint from $z$),
and let $\gammas$ be another $g$-tuple obtained by handlesliding
$\beta_1$ over $\beta_2$ (in the complement of $z$). Fix also another
$g$-tuple $\deltas$ which is isotopic to the $\betas$. More precisely,
let $\beta_1'$ be a curve obtained by handle-sliding $\beta_1$ over
$\beta_2$. Move $\beta_1'$ by a small isotopy to $\gamma_1$, so that
it meets $\beta_1$ in a pair of transverse intersection points with
opposite signs. Let $\gamma_i$ for $i>1$ be obtained by small
isotopies of $\beta_i$. Furthermore let $\delta_i$ for $i=1,...,g$ be
obtained by small isotopies of $\beta_i$. We also require that for
each $i=1,...,g$ the pairwise intersections $\beta_i
\cap \gamma _i$, $\beta_i \cap\delta _i$, $\gamma _i\cap \delta _i$
consist of two transverse points with opposite signs.  We denote these
points $y_i ^\pm$, $w_i ^\pm$, $v_i ^\pm$ respectively. All the above
isotopies are taken to be small enough to be disjoint from the initial
basepoint $z$. Moreover, we arrange that $z$ is also outside of the
pair of pants domain which bounds $\beta_ 1,
\beta_2$ and $\beta_1 '$.  See the corresponding
Figure~\ref{fig:bcd} for $g=2$, where the small circles are identified
by their vertical pairs to give the genus two surface.

It is easy to see that the three Heegaard diagrams
$(\Sigma,\betas,\gammas,z)$, $(\Sigma,\gammas,\deltas,z)$ all
$(\Sigma,\betas,\gammas,z)$ represent the three-manifold
$\#^g(S^2\times S^1)$. Let $\spinc_0$ denote the $\SpinC$ structure
on $\#^g(S^2\times S^1)$ with $c_1(\spinc_0)=0$. 

We will calculate the Floer homologies of $\#^g(S^2\times S^1)$ in the
$\SpinC$ structure $\spinc_0$ with respect to all three Heegaard
diagrams. We start with the easiest case, $(\Sigma,\betas,\deltas,z)$.
Note that since we can choose the curves $\beta_i$ to be exact
Hamiltonian isotopic to the $\delta _i$, so the following can be
thought of as a natural analogue of a basic result of Floer.

\begin{lemma} 
    \label{lemma:Tori1} 
        The pointed Heegaard diagram $(\Sigma,\betas,\deltas,z)$ is 
        admissible for the $\SpinC$ structure $\spinc_0$. Moreover, there is a choice
        of orientation conventions $\orient_{\beta,\gamma}$ with the property that
        \begin{eqnarray*}
                \HFa(\Sigma,\betas,\deltas,\orient_{\beta,\delta})&\cong& H_*(T^g;\Z) \\
                \HFm(\Sigma,\betas,\deltas,\orient_{\beta,\delta})&\cong& \Z[U]\otimes H_*(T^g;\Z), \\
        \end{eqnarray*}
        where $T^g$ denotes the $g$-dimensional torus.
\end{lemma}

\begin{proof}
There are altogether $2^g$ intersection points between $\Tb\cap\Tc$.
Fix some $i=1,...,g$, and let 
$\x=\{x_1,...,x_g\}$ and $\y=\{y_1,...,y_g\}$ be intersection points
with $x_i= w_i^+$, 
$y_i=w_i ^-$ and for all $j\neq i$ we have $x_j=y_j$.
Note that
$\beta_i$, $\delta_i$ bound two pairs of disks, which we denote by
$D_1$ and $D_2$, as in Figure~\ref{fig:Disks}. Let $\phi_1, \phi _2 \in \pi _2(\x,\y)$ 
be homotopy classes given by $\cald (\phi_1)= D_1$ and $\cald (\phi _2)= D_2$ respectively.
It is easy to see that $\mu(\phi _1)=\mu (\phi _2)=1$. 
By juxtaposing the disks we get a periodic class $\phi _1-\phi _2$ 
in $\pi _2(\x,\x)$, with $\mu(\phi _1-\phi _2)=0$.
It easy to see that $\pi _2(\x,\x)$ and similarly $\pi _2(\y,\y)$ are generated by
$g$ classes, all of which can be realized in this way, so it follows that all $2^g$
intersection points represent the trivial $\SpinC$ structure. Moreover, since
each non-trivial combination of these periodic classes has both positive and negative 
coefficients, it follows that $(\Sigma,\betas,\deltas,z)$ is admissible for
$\spinc_0$ (since $c_1(\spinc_0)=0$, the strong and weak admissibility notions coincide).
It also follows that if $d(\w)$ denotes the number of positive coordinates
for $\w= \{w_1 ^\pm ,...,w _g ^\pm\}$, then the relative grading is given by
$\gr(\w,\w')= d(\w)-d(\w')$. We still have to show that all boundary maps
are zero.

    \begin{figure}
\mbox{\vbox{\epsfbox{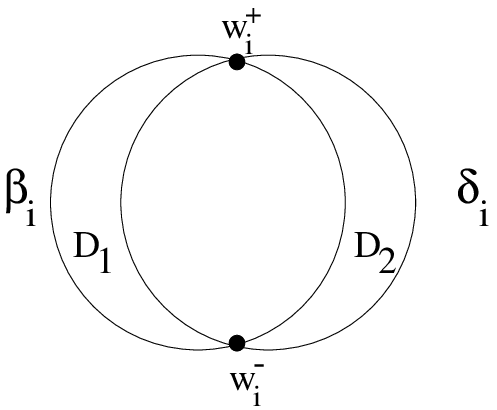}}}
\caption{\label{fig:Disks}}
\end{figure}

For the purpose of these calculations, we will use the path $J_t\equiv
\Sym^g(\sj)$, where $\sj$ is some fixed complex structure over
$\Sigma$. In these calculations, transversality for the flow-lines
considered
is either immediate or
it is achieved by moving the curves, as in
Proposition~\ref{prop:MoveToriTransversality}.

First let $\x$ and $\y$ be as above. Then for all $\phi\in
\pi_2(\x,\y)$, with $n_z(\phi)=0$ and $\phi \neq
\phi _1, \phi _2$, the moduli space $\ModFlow(\phi)$ is empty, since $\cald (\phi)$ has some negative
coefficients. Also $\UnparModFlow(\phi _1)$ contains a unique solution that maps the trivial $g$-fold cover of the disk
to $\Sigma $, so that one of the sheets is mapped to $D_1$, and other sheets map to  $x_j$
for $j\neq i$ respectively. Similarly $\UnparModFlow(\phi _2)$ also contains a 
unique smooth solution. Since the two domains differ by a periodic domain, 
we are free to choose an orientation system for which
the signs corresponding to these solutions
are different. For this system, the 
the component of $\y$ in
${\widehat \partial}\x$ is zero.

Now suppose that $\x$ and $\y$ differ in at least two coordinates, and
$\gr(\x,\y)=1$. Then it is easy to see that for all $\phi \in \pi
_2(\x,\y)$, with $n_z(\phi)=0$ the moduli space $\ModFlow(\phi)$ is
empty, since $\cald(\phi)$ has some negative coefficients. This shows
that ${\widehat \partial}\equiv 0$, and consequently
$\HFa(\betas,\deltas)\cong H_{*}(T^{g};\Z)$. In fact, since the two disks
$D_1$ and $D_2$ separately have a unique holomorphic representative,
it follows that the $H_1(\#^g(S^2\times S^1))$-module structure is
given by the identification $H_*(T^g)\cong \Wedge^* H_1(\#^g(S^2\times
S^1))$.

Let $\y+=w_1^+\times...\times w_g^+$ be the intersection point
representing a top-dimensional homology class in
$\CFa(\betas,\deltas,z)$. We have shown that for all $\y\neq \y^+\in
\Tb\cap \Tc$, $\gr (\y^+,\y)\geq 1$. It follows that for all $\phi \in
\pi_2(\y^+,\y)$ with  $\mu(\phi)=1$ we have $n_z(\phi)\leq 0$. By Lemma~\ref{lemma:NonNegativity},
when $n_z(\phi)<0$, the moduli space $\Mod(\phi)$ is empty. The case
$n_z(\phi)=0$ corresponds to the boundary map in
$\CFa(\betas,\gammas,\spinc_0)$ which was just shown to be trivial. It
follows that $\partial^\infty ([\y^+,i])=0$. The algebra action on
$[\y^+,-1]$ gives then the map $$\Z[U]\otimes_\Z \Wedge^*
H^1(\#^g(S^2\times S^1);\Z)\longrightarrow
\CFm(\betas,\deltas,\spinc_0),$$ which is easily seen to be an
isomorphism by properties of the short exact sequence $$
\begin{CD}
0@>>> \CFm(\betas,\deltas,\spinc_0)@>{U}>>\CFm(\betas,\deltas,\spinc_0)@>>>
\CFa(\betas,\deltas,\spinc_0)@>>> 0.
\end{CD}
$$
\end{proof}

\begin{remark}
\label{rmk:CoeffSystems}
In fact, the above proof shows that for any orientation system
$\orient$, the chain complex for
$\CFa(\betas,\deltas,\spinc_0,\orient)$ is a $g$-fold tensor product
of two-step complexes $\Z\longrightarrow \Z$, where the boundary map
is multiplication by $0$ or $\pm 2$. Of the $2^g$ possible
(isomorphism classes of) complexes, the orientation of
Lemma~\ref{lemma:Tori1} is characterized as the only one with a
non-zero, $g$-dimensional cycle.
\end{remark}

In order to calculate the complexes associated to
$(\Sigma,\betas,\gammas,z)$,  it is useful to have the following:

\begin{figure}
\mbox{\vbox{\epsfbox{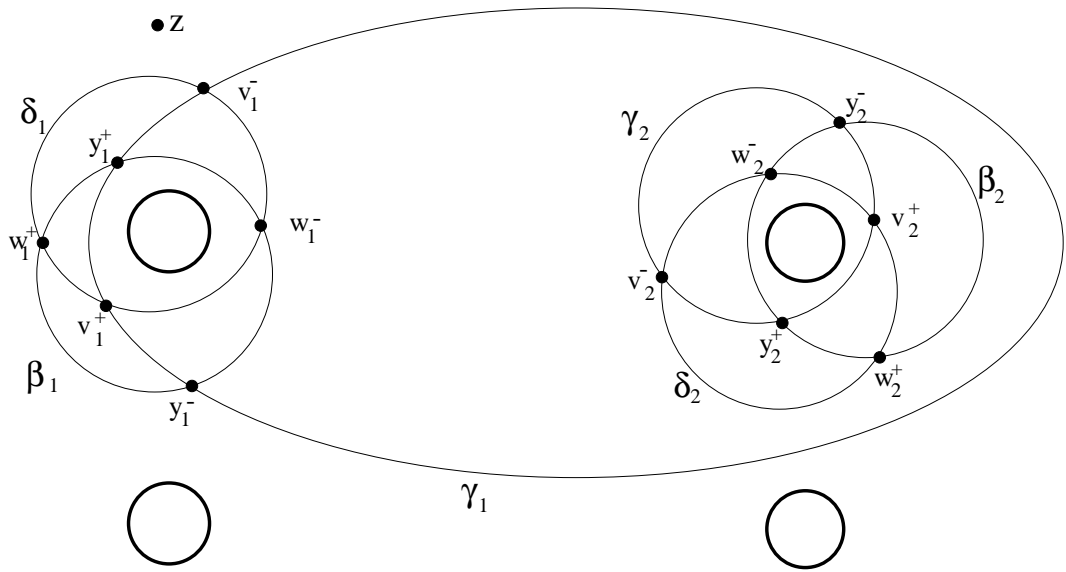}}}
\caption{\label{fig:bcd} $\beta_i$, $\gamma_i$, $\delta _i$ in
the genus 2 surface.}
\end{figure} 

\begin{lemma}
\label{lemma:Annuli}
Let $A$ be an annulus, and fix arcs $\xi_1$ and $\xi_2$ on the outer
and inner boundaries respectively. Then there is a holomorphic
involution of $A$
switching $\xi_1$ and $\xi_2$ if and only if the angles swept out by
$\xi_1$ and $\xi_2$ agree. Moreover, if
this condition is satisfied, the involution is unique. 
\end{lemma}

\begin{proof}
We can think of $A$ as the set $\{z\in \C \big| r<|z|<1/r\}$.  By the
Schwartz reflection principle, we can extend any involution of $A$ to
the Riemann sphere, so that it switches $0$ and $\infty$. Thus, this
involution has the form $z\mapsto c/z$ for some $c\in \C^*$. Since the 
involution takes the set of complex numbers with modulus $r$ to those 
with modulus $1/r$, 
it follows that $|c|=1$, and the lemma follows.
\end{proof}

\begin{lemma} 
    \label{lemma:Tori2} The pointed Heegaard diagram
    $(\Sigma,\betas,\gammas,z)$ is admissible for the $\SpinC$
    structure $\spinc_0$. Moreover, there is a choice of orientation
    conventions $\orient_{\beta,\gamma}$ with the property that \begin{eqnarray*}
    \HFa(\Sigma,\betas,\gammas,\orient_{\beta,\gamma})&\cong&
    H_*(T^g;\Z) \\
    \HFm(\Sigma,\betas,\gammas,\orient_{\beta,\gamma})&\cong&
    \Z[U]\otimes H_*(T^g;\Z), \\ \end{eqnarray*} where $T^g$ denotes
    the $g$-dimensional torus.
\end{lemma}

\begin{proof}
We work out the $g=2$ case, for notational convenience (the general
case follows easily).  

Let $D_i$ for $1\leq i\leq 5$ denote the
connected components of $\Sigma_2-\beta_1-\beta_2-\gamma_1-\gamma_2$,
see Figure~\ref{fig:regions}.  Let $\x=\{y_1^+, y_2 ^+\}$, and
$\y=\{y_1^-, y_2^+\}$.  Then there are exactly three classes, $\phi \in
\pi_2(\x,\y)$, with $n_z(\phi)=0$ and $\cald (\phi)\geq 0$, and these
are given by $\cald(\phi_1)= D_1$, $\cald (\phi_2)= D_2+D_3$, and
$\cald (\phi_3)= D_2+D_4$.  As before, $\mu (\phi_1)=1$. We claim that
$\mu(\phi_2)=\mu(\phi_3)=1$.  Since these cases are symmetric we deal
with $\phi_2$.  

According to Lemma~\ref{lemma:Correspondence}, a holomorphic disk $u:
D \rightarrow
\Sym^2 (\Sigma)$ representing $\phi_2$ gives rise to a branched
double cover $\pi: F\rightarrow D $ and a map ${\widehat u}:
F\rightarrow \Sigma$.  In our case ${\widehat u}$ has degree 1 in
$D_2$ and $D_3$ and has degree 0 on the other regions.  Here, $F$ is
an annulus, and the image of $\partial F$ lies in the union of
$\beta_1$, $\beta_2$, $\gamma_1$, $\gamma_2$.  In fact the part of the
image that lies in $\gamma_2$ is an arc starting at $y^+_2$, see
Figure~\ref{fig:cut}.  More generally, for each $t\in [0,1)$, we can
consider the subset $B_{t}\subset \Sigma$, which is the interior of
the region in $D_{2}\cup D_{3}$ obtained by removing a length $t$
subarc of $\gamma_{2}$ starting at $y_{2}^{+}$, with arc-length
normalized so that $t=1$ corresponds to the endpoint $y_{2}^{-}$.

The region $B_{t}$ is topologically an (open) annulus, and hence can
be identified conformally with a standard (open) annulus
$A_{t}^{\circ}=\{z\in\C \big| r_{t}<|z|<1\}$, where, $r_{t}$ is a
non-zero real number depending on $t$. The identification is given by
a map $\Phi_{t}\colon A_{t}\longrightarrow B_{t}$ which extends to a
continuous map from the closure $A_{t}$ of $A_{t}^{\circ}$ to the
closure of $B_{t}$, which maps the boundary of $A_{t}$ into the union
of $\beta_{1}$, $\beta_{2}$, $\gamma_{1}$, and $\gamma_{2}$, so as to
map onto the length $t$ sub-arc of $\gamma_{2}$. Let $ \xi_1(t)$ and
$\xi_2(t)$ denote the subsets of $\partial A_{t}$ which map to
$\beta_{1}$ and $\beta_{2}$ respectively, and let $f_{1}(t)$ and
$f_{2}(t)$ denote the angles swept out by $\xi_1(t)$ and
$\xi_2(t)$. 
When the map $\Phi_{t}\colon A_{t}\longrightarrow
B_{t}\subset
\Sigma$ is induced from a holomorphic disk in the second symmetric 
product, the corresponding branched double-cover (of the disk, by
$A_{t}$) induces an involution on $A_t$ which switches $\xi_1(t)$ to
$\xi_2(t)$. These assertions follow from simple modifications of the
Riemann mapping theorem, see for example~\cite{Ahlfors}. Indeed, it follows from
these classical considerations that if the curves are smoothly embedded,
then the objects $r_t$ and $\Phi_t$ also depend smoothly on the parameter $t$.

In light of Lemmas~\ref{lemma:Correspondence} and \ref{lemma:Annuli}
points in the moduli space $\UnparModFlow(\phi_2)$ are in one-to-one
correspondence with $t\in (0,1)$ with $f_1(t)=f_2(t)$.  For certain
(generic) choices of the curve $\beta_1$, we can arrange that the
graph of $f_1$ and $f_2$ have non-empty, transversal intersection;
i.e. that there is a finite, non-empty collection of holomorphic
disks.  By a slight perturbation of the curves, it then follows that
the formal dimension of $\UnparModFlow(\phi_2)$ is 0, so
$\mu(\phi_2)=1$.

\begin{figure}
\mbox{\vbox{\epsfbox{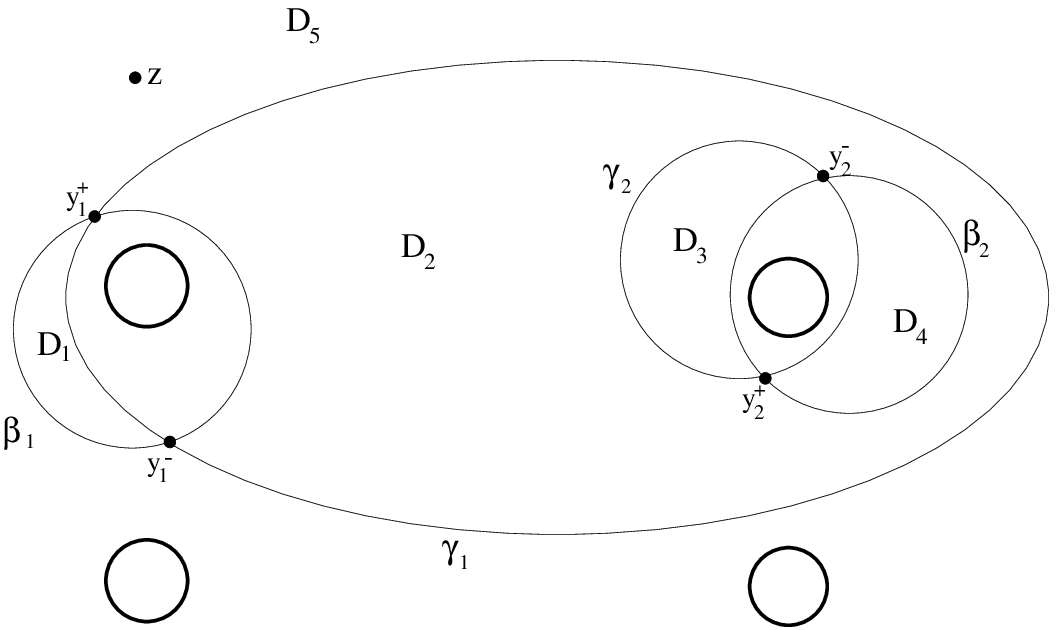}}}
\caption{\label{fig:regions}}
\end{figure}

Now let $\x'=\{ y_1^+, y_2^- \}$, and let $\phi_4, \phi_5 \in
\pi_2(\x,\x')$ be given by $\cald (\phi_4)= D_3$. $\cald (\phi_5)=
D_4$. Clearly $\phi_4-\phi_5$ and $\phi_1-\phi_2$ generate the
periodic classes of $\pi_2(\x,\x)$, and both of them has zero Maslov
index.  Clearly any non-trivial linear combination has positive and
negative coefficients as well, so our diagram $(\Sigma,\betas,\gammas,\spinc_0)$
is $\spinc_0$-admissible.

In order to compute the boundary maps, note that
$\#\UnparModFlow(\phi_1)=\pm 1$, where the unique solution is again a
product of a holomorphic disk and a constant map.  We also claim that
$$\#\UnparModFlow(\phi_2)+\#\UnparModFlow(\phi_3)=\pm 1.$$ This
computation follows easily from the previous discussion and some
complex analysis. Let $A$ denote the annulus given by the domain $D_2$
in $\Sigma$, and let $\nu_1$ and $\nu_2$ denote the conformal angles
of the parts of $\partial A$ that correspond to $\beta_1$ and
$\beta_2$ respectively. By general position we can assume that
$\nu_1\neq \nu_2$.

Let $f_1$, $f_2$ be defined as above. Clearly $f_2(0)=2\pi $, so for small enough
$t$ we have $f_1(t)< f_2(t)$. 
We claim that $\lim_{t\goesto 1} A_{t}=A$,  $\lim_{t\goesto 1}f_{1}(t)=\nu_{1}$, 
$\lim_{t\goesto 1}f_{2}(t)=\nu_{2}$. 
This follows readily from Gromov's compactness
(although in this case, it could also be proved using classical 
conformal analysis, see for example~\cite{Ahlfors}).
It follows that if $\nu_1< \nu_2$, then
$\#(\UnparModFlow(\phi_2))=0 $, and if  $\nu_1 >\nu_2$, then
$\#(\UnparModFlow(\phi_2))=\pm 1$.

\begin{figure}
\mbox{\vbox{\epsfbox{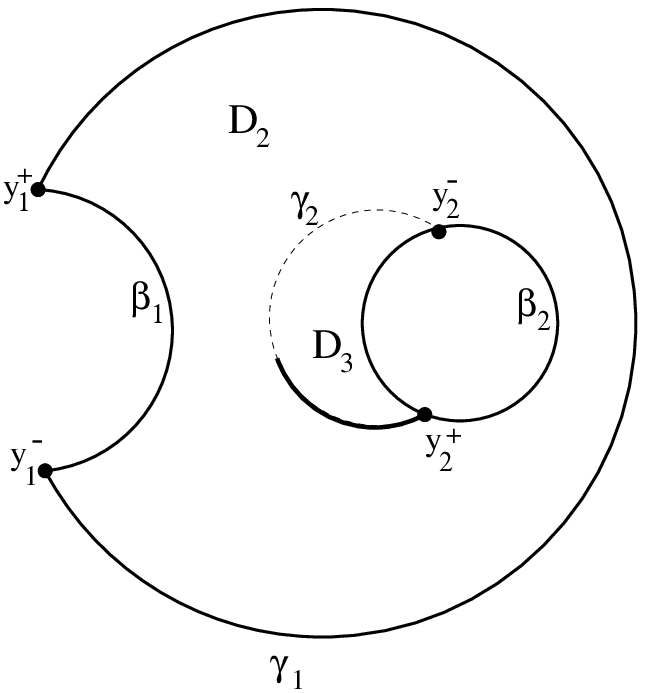}}}
\caption{\label{fig:cut}}
\end{figure} 

Thus, we have calculated $\#\ModFlow(\phi_{2})$.
The computation for  $\phi_3$ is similar. Again we have an annulus with a
cut of length $t\in [0,1)$. Let $g(t)$ denote the conformal angle corresponding
to $\beta_1$, $\beta_2$ respectively. In this case $g_2(0)=0$, so for small
enough $t$ we have $g_1(t)>g_2(t)$. As above, as $t\rightarrow 1$, the conformal angle  $g_i(t)$ converges to $\nu_i$.
It follows that if $\nu_1 <\nu_2$, then  we have 
$\#(\UnparModFlow(\phi_3))=\pm 1 $, and if  $\nu_1 >\nu_2$ then
$\#(\UnparModFlow(\phi_3))=0$. This proves  
$\#\UnparModFlow(\phi_2)+\#\UnparModFlow(\phi_3)=\pm 1$.

Again, we are free to choose an orientation system 
so that the sign of the flow obtained here cancels the sign in 
$\#\ModFlow(\phi_{1})$ (since the domains differ by a periodic domain), 
so that the $\y$ component of $\partial \x$ vanishes.

The same argument shows also that the $\y'$-component of ${\widehat\partial}\x'=0$, 
where $\y'= \{y_1^-,y_2^-\}$. The remaining components
of ${\widehat \partial}$ can be shown to vanish, using the arguments
from Lemma~\ref{lemma:Tori1}.  This shows that
$\HFa(\betas,\gammas,z)\equiv H_*(T^2,\Z)$. The corresponding
statement about $\HFm(\betas,\gammas,z)$ also follows as before.

The case where $g>2$ follows from the $g=2$ case (where we have the same
holomorphic maps multiplied with $g-2$ constant maps).
\end{proof}

\subsection{Naturality}
\label{subsec:Naturality}

Continuing notation from the previous section, we 
let 
$$
\begin{array}{lll}
{\widehat\Theta}_{\beta,\gamma}\in \HFa(\betas,\gammas,
\spinc_0,\orient_{\beta,\gamma}),&
{\widehat\Theta}_{\gamma,\delta}\in\HFa(\gammas,\deltas,\spinc_0,\orient_{\gamma,\delta}),
&{\widehat\Theta}_{\beta,\delta}\in\HFa(\betas,\deltas,\spinc_0,\orient_{\beta,\delta})
\end{array}$$
and ${\Theta}^{\leq}_{\beta,\gamma}$,
${\Theta}^{\leq}_{\gamma,\delta}$, and
${\Theta}^{\leq}_{\beta,\delta}$ denote the top-dimensional generators
of $\HFleq$ coming from Lemmas~\ref{lemma:Tori1} and \ref{lemma:Tori2};
i.e. ${\widehat \Theta}_{\beta,\gamma}$ is represented by the intersection point
$\y^+=\{y_1^+,...,y_g^+\}$, ${\widehat \Theta}_{\gamma,\delta}$ 
is represented by ${\mathbf v}^+=\{v_1^+,...,v_g^+\}$, and 
${\widehat \Theta}_{\beta,\delta}$ is represented by
$\w^+=\{w_1^+,...,w_g^+\}$.

Let $Y$ be a three-manifold equipped with a $\SpinC$ structure
$\spinc$, and $(\Sigma,\alphas,\betas,z)$ be a strongly
$\spinc$-admissible Heegaard diagram, and let $\gammas$ be obtained
from $\betas$ by a handleslide as chosen in the isotopy class
described in Subsection~\ref{subsec:STwoTimesSOne}. Note that the
Heegaard triple $(\Sigma,\alphas,\betas,\gammas,z)$ represents the
cobordism $Y\times [0,1]$, with a bouquet of $g$ circles removed
(c.f. Example~\ref{ex:Cobordism}). It follows that for each $\SpinC$
structure $\spinc$ over $Y$, there is a uniquely induced $\SpinC$
structure on the Heegaard triple $X_{\alpha,\beta,\gamma}$ whose
restriction to $Y=Y_{\alpha,\beta}$ is $Y$ and $\#^g(S^2\times
S^1)=Y_{\beta,\gamma}$ has trivial first Chern class.

It is also easy to see that if we
specify some arbitrary isomorphism classes of orientation system
$\orient_{\alpha,\beta}$ for $(\Sigma,\alphas,\betas,z)$, the
orientation system $\orient_{\beta,\gamma}$ given to us by
Lemma~\ref{lemma:Tori2}, then these uniquely induce an orientation
system $\orient_{\alpha,\gamma}$ on $(\Sigma,\alphas,\gammas,z)$.

\begin{theorem}
\label{thm:HandleslideInvariance}
Let $Y$ be a closed three-manifold equipped with a $\SpinC$ structure
$\spinc\in\SpinC(Y)$, and let $(\Sigma,\alphas,\betas,z)$ be a
strongly $\spinc$-admissible Heegaard diagram. Fix an arbitrary
isomorphism class of orientation system $\orient_{\alpha,\beta}$ for
$(\Sigma,\alphas,\betas,z)$, and let $\gammas$ be the $g$-tuple of
circles obtained from $\betas$ by a handleslide as above. Then, for
the induced orientation system $\orient_{\alpha,\gamma}$, the chain map
$\xi\mapsto f^{\infty}(\xi\otimes\Theta^{\leq}_{\beta,\gamma})$ given
by counting holomorphic triangles over the Heegaard triple
$(\Sigma,\alphas,\betas,\gammas,z)$ induces an isomorphism
of $\Z[U]\otimes_{\Z}\Wedge^* H_1(Y;\Z)/\Tors$-modules
$$\HFinf(\Sigma,\alphas,\betas,\spinc)\longrightarrow
\HFinf(\Sigma,\alphas,\gammas,\spinc);$$
and indeed the triangle count induces isomorphisms on all the other
homology groups $\HFm$, $\HFp$, and $\HFa$.
\end{theorem}

\begin{lemma} 
\label{lemma:Assoc}
Let $(\Sigma,\alphas,\betas,\gammas,\deltas,z)$ as before.
Then, for any $\xi\in\HFa(\alphas,\betas, \spinc)$, we have that
\begin{eqnarray*}
{\widehat F}_{\alpha ,\gamma, \delta} 
({\widehat F}_{\alpha, \beta, \gamma} (\xi\otimes{\widehat\Theta}_{\beta,\gamma})\otimes {\widehat\Theta}_{\gamma,\delta})
={\widehat F}_{\alpha, \beta, \delta}(\xi\otimes {\widehat F}_{\beta, \gamma,\delta} 
({\widehat\Theta}_{\beta,\gamma}\otimes{\widehat\Theta}_{\gamma,\delta})), 
\end{eqnarray*}
and for $\xi\in\HFinf(\alphas,\betas,\spinc)$,
\begin{eqnarray*}
F^{\infty}_{\alpha ,\gamma, \delta} 
(F^{\infty}_{\alpha, \beta, \gamma} (\xi\otimes\Theta^{\leq}_{\beta,\gamma})\otimes \Theta^{\leq}_{\gamma,\delta})
=F^{\infty}_{\alpha, \beta , \delta}(\xi \otimes F^{\leq}_{\beta, \gamma,\delta} (\Theta^{\leq}_{\beta,\gamma}\otimes\Theta^{\leq}_{\gamma,\delta})).
\end{eqnarray*}
\end{lemma}

\begin{proof}
This follows from associativity of the triangle construction. 

First, we begin with some remarks on admissibility for the pointed
Heegaard quadruple $(\Sigma,\alphas,\betas,\gammas,\deltas,z)$ Note
that $$\delta H^1(Y_{\beta,\delta})+\delta H^1(Y_{\alpha,\gamma})=0.$$
To see this, note that all quadruply-periodic domains for our Heegaard
quadruple can be written as sums of doubly-periodic domains for the
four bounding three-manifolds (this is an easy consequence of the fact
that the spans in $H_1$ of the three $g$-tuples $\betas$, $\gammas$,
and $\deltas$ all coincide). Since the map from doubly-periodic
domains to quadruply-periodic domains models the map $H^1(\partial
X_{\alpha,\beta,\gamma,\delta})$ to
$H^2(X_{\alpha,\beta,\gamma,\delta}, \partial
X_{\alpha,\beta,\gamma,\delta})$, it follows that the map
$H^2(X_{\alpha,\beta,\gamma,\delta})\longrightarrow H^2(\partial
X_{\alpha,\beta,\gamma,\delta})$ is injective. Since the restriction
of $\delta H^1(Y_{\beta,\delta})+\delta H^1 (Y_{\alpha,\gamma})$ to
the boundary is trivial, it follows that the subgroup itself is
trivial. Thus, the orbits ${\mathfrak S}\subset
\SpinC(X_{\alpha,\beta,\gamma,\delta})$ needed for
admissibility (in the sense of Subsection~\ref{subsec:AdmissRect})
are trivial.  Moreover, it is now easy to see that if we
choose $\betas$, $\gammas$ and $\deltas$ all sufficiently close to one
another, so the six Heegaard diagrams are strongly admissible, then
the quadruple $X_{\alpha,\beta,\gamma,\delta}$ is strongly admissible.

The result is now a direct application of Theorem~\ref{thm:Associativity}.
\end{proof}

\begin{lemma}
\label{lemma:CalcProd}
The triple $(\Sigma,\betas,\gammas,\deltas,z)$ is admissible in
respect to $z$, and (for the orientation conventions from
Lemmas~\ref{lemma:Tori1} and \ref{lemma:Tori2}), we have that
\begin{eqnarray*}
{\widehat F}_{\beta,\gamma,\delta}(\Theta_{\beta,\gamma}\otimes{\Theta_{\gamma,\delta}})
&=&{\widehat\Theta}_{\beta,\delta}, \\
{F}^\leq_{\beta,\gamma,\delta}(\Theta^\leq_{\beta,\gamma}\otimes
{\Theta^\leq_{\gamma,\delta}})&=&
{\Theta}^{\leq}_{\beta,\delta}, 
\end{eqnarray*}
\end{lemma}

\begin{proof}
The admissibility is easy to check. Now let $\Delta_i$ denote the
triangular region in $\Sigma -\{\beta_i\}-\{\gamma_i\}-\{\delta_i\}$
with vertices $y_i^+, v_i^+, w_i^+$, c.f. Figure~\ref{fig:bcd}. Let
$\psi \in \pi_2(\y^+, {\mathbf v}^+, \w^+)$ be given by $\cald(\psi)=
\sum_{i=1}^g \Delta_i$. Then all $\psi '\neq \psi$, with $n_z(\psi)=0$
has some negative coefficients so $\ModFlow (\psi')=\emptyset$.  To
study $\ModFlow(\psi)$, note that for each $i$ there is a unique map
from the two-simplex to $\Delta_i$ that satisfies the boundary
conditions. The corresponding $g$-tuple of maps to $\Sigma$ gives the
unique solution in $\ModFlow(\psi)$, which is easily seen to be
smooth.  It follows that the coefficient of $\w^+$ in 
${\widehat f} (\y^+ \otimes {\mathbf v}^+)$ 
is $\pm 1$.  For all $\w\neq \w^+$ we
have $\gr(\y^+,{\mathbf v}^+,\w)= \gr(\y^+,{\mathbf
v}^+,\w^+)+\gr(\w^+,\w) >0$. Thus, ${\widehat
F}(\Theta_{\beta,\gamma}\otimes\Theta_{\gamma,\delta})=\Theta_{\beta,\delta}$. More
precisely, choosing orientation conventions on
$(\Sigma,\betas,\gammas,z)$ and $(\Sigma,\gammas,\deltas,z)$ as in
Lemmas~\ref{lemma:Tori1} and \ref{lemma:Tori2}, i.e. for which $\y^+$
and ${\mathbf v}^+$ are cycles, then for the induced orientation
convention on $(\Sigma,\betas,\deltas,z)$, $\w^+$ is a cycle so there,
too, we are using the orientation convention of
Lemma~\ref{lemma:Tori1}.
The corresponding result for $f^{\leq}$ follows similarly.
\end{proof}

\begin{prop}
\label{prop:Isomorphism}
If the $\deltas$ are sufficiently close to the $\betas$, then the maps
$f_{\alpha,\beta,\delta}(\cdot\otimes\Theta_{\beta,\delta})$ induce
isomorphisms $$\HFinf(\Sigma,\alphas,\betas,\spinc)\cong
\HFinf(\Sigma,\alphas,\deltas,\spinc);$$
and indeed isomorphisms of all the other homology groups $\HFa$, $\HFm$, and $\HFp$.
\end{prop}

In fact, the above proof is slightly simpler in the case where
$c_1(\spinc)$ is torsion, so we start with a proof in this case. Indeed,
we will find it useful to introduce here an ``energy filtration'' on $\CFinf$.
To this end, we introduce some terminology. 

\begin{defn}
\label{def:FilteredGroup}
A {\em filtered group} is
a free Abelian group $C$ which is freely generated by a distinguished set
of generators $\SpinCz$, and equipped with a map $$\Filt\colon
\SpinCz\longrightarrow \R.$$ If $\xi,\eta$ are any two elements in
$C$, we say that $\xi<\eta$, if, writing $\xi=\sum_{\x\in\SpinCz}
\xi_\x\cm \x$ and $\eta=\sum_{\x\in\SpinCz}\eta_\x \cm \x$
(where $\xi_\x,\eta_x\in\Lambda$), we have
that $$\max\{\Filt(\x)\big|\xi_\x\neq 0\}<
\min\{\Filt(\x)\big|\eta_\x\neq 0\}.$$
We say that a filtration is {\em bounded below} if for each real number $r$,
$(\Filt(\SpinCz))\cap (-\infty,r]$ is a finite set.
\end{defn}

There are, of course, more general notions where $\R$ is replaced by an 
arbitrary partially ordered set, but we have no need for those presently.

The basic property is the following:

\begin{lemma}
\label{lemma:IsoFilter}
Let $F\colon A\longrightarrow B$ be a map of filtered groups, which
can be decomposed as as a sum $F=F_0 + \ell$ where $F_0$ is a
filtration-preserving isomorphism, and $\ell(\x)<F_0(\x)$ for all
$\x\in{\mathfrak S}$.  Then, if the filtration on $B$ is bounded
below, we can conclude that $F$ is an isomorphism of groups.
\end{lemma}

\begin{proof}
Straightforward.
\end{proof}

When $c_1(\spinc)$ is torsion and $(\Sigma,\alphas,\betas,z)$ is
strongly $\spinc$-admissible, then that Heegaard diagram is, of
course, weakly admissible for all $\SpinC$ structures.  Fix a
reference point $\x_0\in\Ta\cap\Tb$ which represent $\spinc$, and
equip $\Sigma$ with a volume form for which the signed area of each
periodic domain is zero.  The {\em energy filtration} on
$\CFinf(\Sigma,\alphas,\betas,\spinc)$ is the filtration which
defined by $\Filt[\x,i]=-\Area(\phi)$ where $\phi\in\pi_2(\x_0,\x)$ is
any homotopy class with $n_z(\phi)=-i$, and $\Area(\phi)$ refers to
its signed area.

\vskip.2cm
\noindent{\bf{Proof of Proposition~\ref{prop:Isomorphism} when $c_1(\spinc)$ is torsion.}}
We endow $(\Sigma,\alphas,\betas,\spinc)$ with a volume form as
above, and arrange for the $\delta_i$ to be exact Hamiltonian
translates of the $\beta_i$, with respect to this volume form (so that
all the $\betas$-$\deltas$ periodic domains also have total signed
area equal to zero).  When all the $\delta_i$ are sufficiently close
to their corresponding $\beta_i$, and all the $\alpha_j$ intersect the
$\beta_i$ transversally (as we always assume), then to each
$\x\in\Ta\cap\Tb$, there is a uniquely associated point $\x'\in
\Ta\cap\Td$ which is closest to $\x$. We let $${\widehat f}_0\colon
\CFa(\alphas,\betas,\spinc)\longrightarrow\CFa(\alphas,\deltas,\spinc)$$
denote the group homomorphism induced by this closest point
map. Clearly, this is an isomorphism of groups.

Let $\w^+\in\Tb\cap\Td$ denote the intersection point representing the
class ${\widehat\Theta}_{\beta,\delta}$. Note that there is a unique
``small triangle'' homotopy class $\psi_0^\x\in\pi_2(\x,\w^+,\x')$,
which is the only triangle satisfying the two properties that
$\cald(\psi_0^\x)\geq 0$, and $\cald(\psi_0^\x)$ is supported inside
the support of the isotopy between $\betas$ and $\deltas$. Indeed,
$\cald(\psi_0^\x)$ has multiplicity one on the small triangles
$\Delta(x_i,w_i^+,x'_i)$ connecting $x_i$, $w_i^+$, and $x'_i$ (for
$i=1,...,g$), and zero everywhere else.  We claim that if $J$ is a
family which is a sufficiently small perturbation of the constant
complex structure $\sj$, then $\#\ModFlow_J(\psi_0^\x)=1$. To see
this, we first consider complex structures of the form $J\equiv
\Sym^g(\sj)$.  By elementary complex analysis we see that there is a
unique holomorphic map $f_i$ from the 2-simplex $\Delta$ to the
triangle $\Delta(x_i, w_i^+, x_i')$ that lines up the corresponding
boundary segments. It follows that the moduli space
$\ModFlow_{\Sym^g(\sj)}(\psi)$ has a unique element, which is given as
the product of $f_1,...,f_g$ from the trivial $g$-fold cover of
$\Delta$ to $\Sigma$ (with respect to the constant almost-complex
structure $\Sym^g(\sj)$). It is easy to see that this is a smooth
solution.  Counted with the appropriate orientation we have
$\#(\ModFlow_{\Sym^g(\sj)}(\psi))=1$. Thus, this assertion persists
for all sufficiently small perturbations of the constant family.

It follows that we can decompose ${\widehat
f}_{\alpha,\beta,\delta}(\cdot\otimes{\widehat\Theta}_{\beta,\delta})$
as $${\widehat f}_{\alpha,\beta,\delta}(\cdot\otimes
\Theta_{\beta,\delta})={\widehat f}_0 + {\widehat\ell},$$ where
${\widehat\ell}$ counts holomorphic triangles $\psi$ with
$n_z(\psi)=0$ and whose domains with are not supported inside the
small region bounded by $\betas$ and $\deltas$. In fact, let
$\epsilon$ denote the total (unsigned) area swept out by the isotopy
between the $\betas$ and the $\deltas$ which, as its notation
suggests, can be made arbitrarily small by varying $\deltas$. Let
$M_0$ denote the minimal area of any (non-zero) domain in
$\Sigma-\alphas-\betas$ (note that this is independent of $\deltas$).
It is easy to see that the triangles counted by $\ell$ all have total
area bounded below by $M_0-\epsilon$. Choose a reference point
$\x_0\in\Ta\cap\Tb$, and consider the induced area filtration on
$(\Sigma,\alphas,\deltas,\spinc)$ induced from $\x_0'$. Since for
each $\x\in\Ta\cap\Tc$, the area of the canonical triangle $\psi_0^\x$
is bounded by $\epsilon$, if we arrange that $\epsilon<M_0/2$, then we
have that $${\widehat\ell}(\x)<{\widehat f}_0(\x),$$ with respect to
the energy filtration induced on $\CFa(\alphas,\deltas,\spinc)$.
Applying Lemma~\ref{lemma:IsoFilter}, it follows that ${\widehat f}$,
which we already know is a chain map, is actually isomorphism of chain
complexes.

By modifying the above constructions, we see that
$$f^\infty_{\alpha,\beta,\delta}(\cdot\otimes{\Theta^\leq_{\beta,\delta}})
=f_0^{\infty}+ {\text{lower order}},$$ where
$$f_0^\infty[\x,i]=[\x',i].$$ Although now the filtration on
$\CFinf(\alphas,\betas,z)$ is no longer bounded below, it comes with
a relative $\Z$-grading, and the map $f^\infty$ preserves this
relative grading, in the sense that if we have
$\xi,\xi'\in\CFinf(\alphas,\betas,\spinc)$ with $\gr(\xi,\xi')=0$,
then
$$\gr(f^\infty_{\alpha,\beta,\delta}(\xi\otimes{\Theta^\leq_{\beta,\delta}}),
f^\infty_{\alpha,\beta,\delta}(\xi'\otimes{\Theta^\leq_{\beta,\delta}}))=0.$$
Now, the filtration induced on the subset of
$\CFinfty(\alphas,\betas,\spinc)$ with fixed relative grading is
easily seen be bounded below (since, in fact, there are only finitely
many generators representing any given degree.
\qed
\vskip.2cm

When $c_1(\spinc)$ is non-torsion, we must use a refinement, since the
relative grading on $\CFinf$ is no longer a relative $\Z$-grading, but
only a $\Zmod{\delta}$-grading (where $\delta=\delta({\spinc})$ is the
indeterminacy from Equation~\eqref{eq:Indeterminacy}), and the complex
for each fixed degree might be an infinitely generated $\Z$-module.

So, we equip $\Sigma$ with a volume form for which each ${\mathfrak
s}$-renormalized periodic domain has total signed area zero. This can
be arranged as in the proof of Lemma~\ref{lemma:EnergyZero}.  
Given any $[\x,i]$ and $[\y,j]$ with the same relative grading,
we can find some disk $\phi\in \pi_2(\x,\y)$
with $n_z(\phi)=i-j$ and $\Mas(\phi)=0$. 
In each relative degree, then, we have a filtration (uniquely
defined up to an additive constant) defined by
$$\Filt([\x,i])-\Filt([\y,j])=-\Area(\cald(\phi)).$$ Since any two possible
choices of such disk $\phi$, $\phi'$ differ by a renormalized periodic
domain, it follows that the difference in filtration defined above is independent of
the the choice of disk.

\vskip.2cm
\noindent{\bf{Proof of Proposition~\ref{prop:Isomorphism} 
when $c_1(\spinc)$ is non-torsion.}}  We consider $\CFm$, and use the
notation from the earlier proof. We write
$$f^\infty_{\alpha,\beta,\delta}=f_0^\infty+\ell^\infty,$$ where
$f_0^\infty[\x,i]=[\x',i]$ as before. Clearly, for given $[\x,i]$,
$f_0^\infty([\x,i])$ and $\ell^\infty([\x,i])$ have the same relative degree. 

We claim that, just as before, if the $\deltas$ are sufficiently close
to the $\betas$, then $\ell^\infty$ has lower order (than $f_0^\infty$) with
respect to the above refined energy filtration.  To this see this, let
$[\y,j]$ be any element appearing with non-zero multiplicity in
$\ell^\infty([\x,i])=[\y,j]$ (so in particular $\y\neq
\x'$).  Then, there is a pseudo-holomorphic triangle
$\psi\in\pi_2(\x,\w^+,\y)$ with $\cald(\psi)\geq 0$ and
$\Mas(\psi)=0$. It is easy to see that in this case, there is also a
$\phi\in\pi_2(\x',\y)$ with $\psi=\psi^\x_0 * \phi$, so that
$\Mas(\phi)=0$ and $n_z(\phi)=i-j$. Since $\cald(\psi)$ is not
supported inside the region between $\betas$ and $\gammas$, we know
that $\Area(\psi)>M_0-\epsilon$. Thus, it follows that
$$\Area(\phi)=\Area(\psi)-\Area(\psi_0)>M_0-2\epsilon,$$ so, since
$-\Area(\phi)$ determines the filtration difference between
$f_0^\infty[\x,i]$ and $\ell^\infty[\x,i]$, if we choose the exact Hamiltonian
translates to be sufficiently close, $\ell^\infty([\x,i])<f_0^\infty[\x,i]$.

We claim that the refined filtration is bounded below. To see this,
observe that if $\delta$ denotes the divisibility of $c_1(\spinc)$,
then for each $[\x,i]$ representing a fixed relative, we can find a
$\SpinC$-renormalized periodic domain connecting $[\x,i]$ to
$[\x,i+\delta]$, so these elements all have the same filtration.
Thus, Lemma~\ref{lemma:IsoFilter} applies,
to prove that $f^\infty_{\alpha,\beta,\delta}$ induces an
isomorphism of chain complexes.
\qed

\vskip.2cm
\noindent{\bf Proof of Theorem~\ref{thm:HandleslideInvariance}.}
This follows easily from Lemma~\ref{lemma:Assoc} and
Proposition~\ref{prop:Isomorphism}.  A direct application of these
results shows that that the map induced from $\otimes
\Theta_{\beta,\gamma}$ is injective, while the map induced from
${\gamma,\delta}$ is surjective. But the roles of these $g$-tuples is
symmetric: we can introduce a fifth tuple of circles $\etas$ which are
small isotopic translates of the $\gammas$, and apply the above
reasoning to see that the map induced from $\otimes
\Theta_{\beta,\gamma}$ is surjective, as well. 

To see that the map commutes with the $H_1(Y;\Z)$-module structure,
we represent the action by a codimension-one constraint
$V\in \Ta$ as in Remark~\ref{rmk:GeoRep}, and
consider the moduli space
$$\ModFlow_{V}(\psi)=\bigcup_{\tau\in\R}
\left\{u\in \ModFlow(\psi)\big| u\circ E_{\alpha}(\tau)\in V\right\},$$
where $E_{\alpha}\colon \R$ is a parameterization of the $\alpha$-edge
as in the proof of Proposition~\ref{prop:TrianglesJIndep}.  As usual,
when $\Mas(\psi)=-1$, this space is compact, and can be used to
construct a chain homotopy $$H([\x,i]\otimes[\y,j])=
\sum_{\w}\sum_{\psi}\#\Big(\ModFlow_{V}(\psi)
\Big)[\w,i+j-n_z(\psi)]$$
Consider homotopy classes 
of triangles $\pi_2(\x,\Theta_{\beta,\gamma},\y)$
with $\Mas(\psi)=0$. The ends as
$\tau\goesto\infty$ correspond to the commutator of $F^\infty(\cdot\otimes
\Theta_{\beta,\gamma})$ with the
action of $V$; the other ends correspond to the commutator of $H$ with
the boundary maps.
\qed

\subsection{The Maslov index of a periodic domain}
\label{subsec:MasIndex}

We can now prove Theorem~\ref{thm:Grading}, which was used in the
definition of the relative gradings. 

\vskip.3cm
\noindent{\bf{Proof of Theorem~\ref{thm:Grading}.}}
Fix $\x\in\Ta\cap \Tb$, and
consider the map $\Hom(H_2(Y;\Z),2\Z)$ which, given $c\in H_2(Y;\Z)$,
calculates $\mu(\psi(c))$, where $\psi(c)\in\pi_2(\x,\x)$ is the periodic
class associated to $c\in H_2(Y;\Z)$. Note that this is a
homomorphism, since the Maslov index is additive.  Indeed, this
assignment depends on the point $\x\in\Ta\cap\Tb$
only through its induced $\SpinC$ structure $\spinc=s_z(\x)$, by
the additivity of the Maslov index. We denote the map by 
$m_\spinc\in\Hom(H_{2}(Y;\Z),\Z)$. 

We argue that $m_\spinc$ depends on $Y$ alone, i.e. it is invariant
under pointed isotopies, pointed handle-slides, and stabilization
(c.f. Proposition~\ref{prop:PointedHeegaardMoves}). To see
stabilization invariance, it suffices to see how the Maslov index
changes by adding $S\in\pi_2(\Sym^g(\Sigma))$, and thereby reducing to
the case where the coefficient of the domain is zero on the
two-torus. Handle-slide invariance follows from 
index calculations parallel to (but much
simpler than) Theorem~\ref{thm:HandleslideInvariance}.
Specifically, let
$\alphas$, $\betas$, $\gammas$ be attaching circles, 
where $\gammas$
are obtained from $\betas$ by a handle slide and a small Hamiltonian
isotopy. Indeed, if the pair of pants used for
the handleslide is sufficiently small, we can associate
to $\x\in\Ta\cap\Tb$ a nearest intersection point
$\x'\in\Ta\cap\Tc$,
which is connected to $\x$ by a small triangle
$\Delta\in
\pi_2(\x,\theta,\x')$ 
with $n_z=0$ and $\Mas(\Delta)=0$ (here,
$\theta$ is an intersection point representing the 
homology class
${\widehat\Theta}_{\beta,\gamma}$ from Lemma~\ref{lemma:Tori2}).
Now, in view of the
affine identification $\pi_2(\x,\theta,\x')\cong \Z \oplus
H^1(Y;\Z)$
(c.f. Propositions~\ref{prop:CalcPiTwo} 
and \ref{prop:HomologyOfX}), 
if $p_\x$ is a periodic class for
$(\Sigma,\alphas,\betas)$, and $p_\x'$ is the corresponding periodic
class for $(\Sigma,\alphas,\gammas,z)$, 
then there is a periodic class
$\delta\in\pi_2(\theta,\theta)$ for $\Tb\cap\Tc$ with the property
that $p_\x+\Delta = p_\x'+\Delta +
\delta$. Moreover, since the Maslov index on any such
element $\delta$ vanishes (this is established in the course of the
proof of Lemma~\ref{lemma:Tori2}), it follows that
$\Mas(p_\x)=\Mas(p_{\x'})$. Isotopy invariance is straightforward,
except in the case where the isotopy cancels all intersection points
belonging to the given $\SpinC$ structure $\spinc$. To avoid this we
use only special isotopies, as in Lemma~\ref{lemma:spincRealized} and
\ref{lemma:StronglyAdmissibleIsotopy} (see Remark~\ref{rmk:SpinCRealized}). 

Now, we argue that if $\spinc, \spinc'\in\SpinC(Y)$ are represented by
intersection points, then we claim that $$ m_{\spinc}=m_{\spinc'}+2c ,
$$ where $c\in H^2(Y;\Z)$ is the class for which
$\spinc'=\spinc-c$. To see this, it suffices to consider the effect of
moving the base-point $z$ across some fixed circle, say,
$\alpha_1$. Note then that $s_{z'}=s_{z}+\alpha_1^*$, according to
Lemma~\ref{lemma:VarySpinC}.  If $\psi$ is the periodic class
corresponding under the basepoint $z$ to some $v\in H_2(Y;\Z)$ then
clearly $n_{z'}(\psi)= -\langle \alpha_1^*,v\rangle$.  Moreover, the
periodic class for $\psi(z',v)=\psi(z,v)-n_{z'}(\psi(z,v))[S]$.  It
follows that $m_{\spinc}=m_{\spinc'}+2c$.

It follows that $m_{\spinc}=c_1(\spinc)+K$, in $\Hom(H_2(Y;\Z),\Z)$
for some $K$ which is independent of $\spinc$. We wish to show that
$K=0$. To this end, we compare $m_\spinc$ and $m_{{\overline\spinc}}$.
Switching the roles of $\Ta$ and $\Tb$ and reversing the orientation
of $\Sigma$, we get a new Heegaard diagram describing $Y$, and an
obvious identification of intersection points; letting $s_z'(\x)$ be
the $\SpinC$ structure with respect to this new data, it is clear that
$s_z'(\x)={\overline{s_z(\x)}}$.  Note that switching the two tori and
the orientation of $\Sigma$ simultaneously leaves holomorphic data,
such as the Maslov index of a given periodic domains,
unchanged. In particular, if 
${\mathcal P}=\sum a_i\cald_i$ is a periodic domain, and ${\mathcal
P}'=\sum a_i\cald_i'$, where the $\cald_i'$ have the opposite
orientation to the $\cald_i$, then 
$\langle m_{\spinc},H({\mathcal P})\rangle=
\langle m_{\overline\spinc},H({\mathcal P}')\rangle$. 
However, $H({\mathcal P})=-H({\mathcal P}')$.  Thus,
$m_{\overline\spinc}=-m_\spinc$. Since it is also true that
$c_1(\spinc)=-c_1({\overline\spinc})$, it follows that $K=0$.
\qed
\vskip.3cm

\section{Stabilization}
\label{sec:Stabilization}

The final step in establishing topological invariance of the Floer
homologies is stabilization invariance. Fix a strongly
$\spinc$-admissible pointed Heegaard diagram $(\Sigma,\alphas,\betas,z)$ for
$Y$, where $\Sigma$ is a surface of genus $g$ and
$\alphas=\{\alpha_1,...,\alpha_g\}$, $\betas=\{\beta_1,...,\beta_g\}$.
The stabilized diagram $(\Sigma',\alphas',\betas',z)$ is obtained by
forming the connected sum $\Sigma'=\Sigma\#E$, which is the connected
sum of $\Sigma$ with $E$, a surface of genus $g$, and letting
$\alphas'=\{\alpha_1,...,\alpha_g,\alpha_{g+1}\}$ and
$\betas'=\{\beta_1,...,\beta_g,\beta_{g+1}\}$, where $\alpha_{g+1}$
and $\beta_{g+1}$ are a pair of circles in the $E$ summand which meet
transversally in a single positive point $c$.  For simplicity, we
choose the point in $\Sigma$ along which we perform the connected sum
to lie in the same path-component of $\CurveComp$ as $z$. (For this
choice, there is a one-to-one correspondence between periodic domains
for the two diagrams, and hence the notions of admissibility coincide.)

We begin with the much simpler case of $\CFa$.

\begin{theorem}
\label{thm:StabilizeHFa}
Let $(\Sigma,\alphas,\betas,z)$ be a strongly $\spinc$-admissible
pointed Heegaard diagram for $Y$, and let $(\Sigma',\alphas',\betas',z)$
denote its stabilization. Then, for each orientation system $\orient$
on the original Heegaard diagram, there is an induced orientation
system $\orient'$ on its stabilization, and a corresponding isomorphism
of $\Wedge^* H_1(Y;\Z)/\Tors$-modules.
$$\HFa(\alphas,\betas,\spinc,\orient)\cong
\HFa(\alphas',\betas',\spinc,\orient')$$
\end{theorem}

\begin{proof}
It is easy to see that the intersection points correspond:
$\Ta'\cap \Tb'=\left(\Ta\cap\Tb\right)\times \{c\}$. 
We choose a basepoint $z$ in the same
path-component of $\Sigma-\alpha_1-...-\alpha_g-\beta_1-...-\beta_g$ as
$\sigma_1\in \Sigma$, the point along which we perform the connected sum 
to obtain $\Sigma'$. Let $z'$ denote the corresponding basepoint in $\Sigma'$.
If $\x\in\Ta\cap\Tb$, and $\x'\in\Ta'\cap\Tb'$ 
is the corresponding point $\x\times\{c\}$, then the induced $\SpinC$ 
structures agree $\spinc_{z}(\x)=\spinc_{z'}(\x')$, since
the corresponding vector 
fields agree away from the three-ball where the stabilization occurs.

Let ${\mathbf x},{\mathbf y}\in\Ta\cap\Tb$,
and $\phi\in\pi_2({\mathbf x},{\mathbf y})$ be the class with
coefficient $n_z(\phi)=0$. Let ${\bf x}'={\bf x}\times \{c\}$, ${\bf
y}'={\bf y}\times \{c\}$, and $\phi'\in\pi_2({\mathbf x}',{\mathbf
y}')$ be the class with $n_{z'}(\phi')=0$. Then, we argue that for 
certain special paths of almost-complex structures, the moduli space
$\Mod(\phi')$ is identified with $\Mod(\phi)\times \{c\}$
(together with its deformation theory, 
including the determinant line bundles).  Hence, the chain complexes
are identical.

To set up the complex structure, we let $\sigma_1$ and $\sigma_2$
denote the connected sum points for $\Sigma$ and $E$ respectively.
Recall that there is a holomorphic map 
$$\Sym^{g}(\Sigma-B_{r_1}(\sigma_1))\times (E-B_{r_2}(\sigma_2))
\longrightarrow
\Sym^{g+1}(\Sigma'),$$
which is a diffeomorphism onto its image.
Suppose that $J_s$ is a (generic) 
family of almost-complex structures over $\Sym^g(\Sigma)$ which agrees
with $\Sym^g(\sj)$ in a tubular neighborhood of $\{\sigma_1\}\times
\Sym^{g-1}(\Sigma)$, a neighborhood which is holomorphically
identified with $B_\epsilon({\sigma_1})\times
\Sym^{g-1}(\Sigma)$.  
Using the above product map, we can transfer $J_s\times \sj_E$ to an
open subset of $\Sym^{g+1}(\Sigma')$, and extend it to all of
$\Sym^{g+1}(\Sigma')$.  It follows from our choice of
$J_s$ that if $n_z(\phi)=0$, then any $J_s$-holomorphic representative
for $\phi$ must have its image in $\Sym^{g}(\Sigma-B_{\epsilon})$, and
hence its product with the constant map is
$(J_s\times{\sj}_E)$-holomorphic in $\Sym^{g+1}(\Sigma')$. Conversely,
a $(J_s\times{\sj}_E)$-holomorphic curve in $\Sym^{g+1}(\Sigma')$ with
$n_{z'}(\phi')=0$ must be contained in
$\Sym^g(\Sigma)\times\{c\}$. The identification of deformation
theories is straightforward.

Since the map we describe identifies individual moduli spaces
it is clear (c.f. Remark~\ref{rmk:GeoRep}) that the map is equivariant
under the action of $H_1(Y;\Z)/\Tors$.
\end{proof}

The corresponding fact for $\HFpm$ and $\HFinf$ is more subtle, and depends on a
gluing theorem for holomorphic curves.

\begin{theorem}
\label{thm:StabilizeHFb}
Let $(\Sigma,\alphas,\betas,z)$ be a strongly $\spinc$-admissible
pointed Heegaard diagram for $Y$, and let $(\Sigma',\alphas',\betas',z)$
denote its stabilization. Then, for each orientation system $\orient$
on the original Heegaard diagram, there is an induced orientation
system $\orient'$ on its stabilization, and a corresponding 
$\Z[U]\otimes_\Z\Wedge^* H_1(Y;\Z)/\Tors$-module isomorphisms
\begin{eqnarray*}
\HFinf(\alphas,\betas,\spinc,\orient)\cong  
\HFinf(\alphas',\betas',\spinc,\orient')
&{\text{and}}&
\HFpm(\alphas,\betas,\spinc,\orient)\cong  
\HFpm(\alphas',\betas',\spinc,\orient')
\end{eqnarray*}
\end{theorem}

The proof of this theorem occupies the rest of this section. However,
we give a brief outline of the proof presently. 

\vskip.3cm
\noindent{\bf{Sketch of proof.}}
As in the proof of Theorem~\ref{thm:StabilizeHFa} above, the
generators for the complexes are identified.  We fix some complex
structure $\sj$ over $\Sigma$, and $\sj_E$ over $E$.  Let $\sj'(T)$
denote the complex structure on $\Sigma'$ obtained by inserting a
cylinder $[-T,T]\times S^1$ between the $\Sigma$ and
$E$. Correspondingly, the symmetric product $\Sym^{g+1}(\Sigma')$ is
endowed with a complex structure $\Sym^{g+1}(\sj'(T))$, which admits
an open subset which is holomorphically identified with
$$\Sym^{g}(\Sigma-B_{r_1}(\sigma_1))\times
\Sym^1(E-B_{r_2}(\sigma_2)),$$ given the complex structure obtained by
restricting $\Sym^{g}(\sj)\times \sj_E$.  Using this parameterization,
we can transfer some fixed initial path $J_s$ of nearly-symmetric
almost complex structures over $\Sym^g(\Sigma)$ to a path $J_s'(T)$ of
$\sj'(T)$-nearly symmetric almost-complex structures which extend over
over $\Sym^{g+1}(\Sigma')$.

For $T$ is sufficiently large, we show that for any pair $\x, 
\y\in\Ta\cap \Tb$, if 
$\phi\in\pi_{2}(\x,\y)$ 
is the homotopy class with $\Mas(\phi)=1$, 
then there is an identification
$\Mod_{J_s}(\phi)\cong \Mod_{J'_s(T)}(\phi')$, where $\phi'\in\pi_{2}(\x',\y')$ is 
the corresponding class with $\Mas(\phi')=1$, and 
$\x'=\x\times\{c\}$, $\y'=\y\times \{c\}$ (see 
Theorem~\ref{thm:Gluing} below). It is an easy consequence of this 
that all the relevant corresponding chain complexes are identified:
\begin{eqnarray*}
    (\CF^{\pm}(\alphas,\betas,\spinc),\partial^{\pm}_{J_s})&\cong&
    (\CF^{\pm}(\alphas',\betas',\spinc),\partial^{\pm}_{J_s(T)}) \\
    (\CF^{\infty}(\alphas,\betas,\spinc),\partial^{\infty}_{J_s})
    &\cong&
    (\CF^{\infty}(\alphas',\betas',\spinc),\partial^{\infty}_{J_s(T)}).
\end{eqnarray*}
Since the homology groups of these chain complexes are independent of
the paths of almost-complex structures, it follows that the groups are
isomorphic, c.f. the last paragraph in Subsection~\ref{sec:Gluing}. 

The identification between the moduli spaces is given by a now 
familiar gluing construction. Suppose that all flows in 
$\ModFlow_{J_s}(\phi)$ meet the subvariety 
$\{\sigma\}\times\Sym^{g-1}(\Sigma)$ in general position: i.e. suppose 
that each Maslov index one flow-line $u\in\ModFlow(\phi)$ in 
$\Sym^{g}(\Sigma)$ meets the above subvariety corresponding to the 
connected sum point $\sigma\in\Sigma$ transversally (this can 
be easily arranged), 
then, given such a flow $u\in\ModFlow(\phi)$,
we have that the pre-image 
$$u^{-1}(\{\sigma\}\times\Sym^{g-1}(\Sigma))=\{q_{1},\ldots,q_{n}\}, $$
where $n=n_{z}(\phi)$. Taking the product of $u$ with $\{c\}$, we 
obtain a map to $\Sym^{g}(\Sigma)\times E$, which we can view
as a subset of $\Sym^{g+1}(\Sigma\vee E)$, which in turn can be 
thought of as a degenerated version of $\Sym^{g+1}(\Sigma')$. 
Splicing in spheres 
in $\Sym^{g-1}(\Sigma)\times \Sym^{2}(E)$ (which can also be thought 
of as a subset of the $(g+1)$-fold symmetric product of the wedge), 
we obtain for any sufficiently large $T$, a nearly $J'_s$-holomorphic map 
into $\Sym^{g+1}(\Sigma\#_{T}E)$ for all sufficiently large $T$. As 
the parameter $T\goesto \infty$, this spliced map becomes closer to 
being $J'_s(T)$-holomorphic (in an appropriate norm, as discussed below). 
Using the inverse function theorem in the usual manner, we obtain 
a nearby pseudo-holomorphic map $u'\in\ModFlow_{J'_s(T)}(\phi')$. 

A Gromov compactness argument shows that for sufficiently large $T$,
all of the moduli spaces $\ModFlow_{J'_s(T)}(\phi')$ lie in the domain
of this gluing map.  (Note also that in the case where $g=1$ not every
homotopy class $\phi'\in\pi_2(\x',\y')$ is obtained by stabilizing a
class $\phi\in\pi_2(\x,\y)$, e.g. there may not be any homotopy class
$\phi\in\pi_2(\x,\y)$ with $\Mas(\phi)=1$. Nonetheless, the
compactness argument applies, showing that then the moduli space
$\ModFlow_{J'_s(T)}(\phi')$ are empty for sufficiently large $T$.)

Note that the $\Z[U]$-equivariance of the identifications is a formal
consequence of the fact that $n_z(\phi)=n_{z'}(\phi')$; the
$H_1(Y;\Z)/\Tors$-equivariance also follows similarly.
\qed
\vskip0.3cm

An important ingredient in the above story is a description of the 
spheres in $\Sym^{2}(E)$. In particular, it is crucial that 
through each pair of points in $E$, there is a unique 
holomorphic sphere in $\Sym^{2}(E)$ in the positive 
generator of $\pi_{2}(\Sym^{2}(E))$. This follows from the 
following:

\begin{lemma}
    \label{lemma:RSurf}
    Let $E$ be a genus one Riemann surface. 
    The second symmetric product $\Sym^{2}(E)$ is naturally a 
    ruled surface over $E$.
\end{lemma}
    
\begin{proof}
    The map from  $\Sym^{2}(E)$ to the base $E$
    is the map sending the pair $\{x,y\}$ to the sum $x+y$. 
    The fiber over a base point $a$ is represented by pairs of the 
    form $\{a+w,-w\}$ with $w\in E$. Two pairs $(a+w_{1},-w_{1})$ 
    and $(a+w_{2},-w_{2})$  are equivalent if
    $-a-w_{1}=w_{2}$. Now, the 
    map $w\mapsto -a-w$ is an involution whose quotient is a 
    projective line. 
\end{proof}

Ultimately, the proof of Theorem~\ref{thm:StabilizeHFb} outlined above
is a modification of the usual picture for dealing with
non-compactness in Gromov theory (\cite{ParkerWolfson},
\cite{McDuffSalamon}, \cite{Liu}).  The gluing here is analytically
closely related to gluing problems which arose in the study of the
Yang-Mills equations (see for example~\cite{Taubes}, \cite{Mrowka}).
Moreover, the degeneration of $\Sigma'$ above is closely related to
the degenerations studied by Ionel-Parker and Li-Ruan which concern
degenerations of symplectic manifolds along codimension one symplectic
submanifolds (see~\cite{IonelParker}, \cite{LiRuan}). The rest of this
section is devoted to proving Theorem~\ref{thm:StabilizeHFb}.

\subsection{Gluing: the statement}
\label{sec:Gluing}

In this subsection, we state the result which is the cornerstone of 
Theorem~\ref{thm:StabilizeHFb} above, allowing us to
use flows in $\Sym^{g}(\Sigma)$ to 
construct flows in $\Sym^{g+1}(\Sigma\#_{T}E)$ for 
sufficiently large $T$. 

For simplicity, we assume henceforth that the path $J_s$ of
almost-complex structures over $\Sym^g(\Sigma)$ agrees with
$\Sym^g(\sj)$ for some complex structure $\sj$ over $\Sigma$ in a
neighborhood of the subset $\Sym^{g-1}(\Sigma)\times \{\sigma\}\subset
\Sym^g(\Sigma)$ (i.e. we are using $(\sj,\eta,V)$-nearly symmetric
almost-complex structures for a choice of $V$ containing
$\{\sigma\}\times\Sym^{g-1}(\Sigma)$). 

Moreover, as indicated earlier,
the path $J_s$ can be used to construct a one-parameter family
$J'_s(T)$ of paths over $\Sym^{g+1}(\Sigma')$ compatibly with $J_s$ in
the following sense. Recall that 
$\Sym^{g+1}(\Sigma\#_T E)$ (endowed with the complex structure
$\Sym^{g+1}(\sj'(T))$, where $\sj'(T)=\sj\#_T\sj_E$) has an open subset
holomorphically identified with:
$$\Sym^{g}(\Sigma-B_{r_1}(\sigma_1))\times
\Sym^{1}(E-B_{r_2}(\sigma_2)),$$
where $r_1$, $r_2$ are non-negative real numbers. Fix another pair of
real numbers $R_1>r_1$, $R_2>r_2$.  We choose $J_s'$ with the
following properties:

\begin{itemize}
\item over $\Sym^g(\Sigma-B_{R_1}(\sigma_1))\times
\Sym^1(E-B_{R_2}(\sigma_2))$, the path $J_s'$ agrees with 
$J_s\times{\sj_E}$,
\item 
over $\Sym^g(\Sigma-B_{r_1}(\sigma_1))\times
\Sym^1(B_{R_2}(\sigma_2)-B_{r_2}(\sigma_2))$, as the normal parameter
to $\sigma_2$ goes from $R_2$ to $r_2$,
$J_s'$ connects from $J_s\times \sj_E$ to
$\Sym^g(\sj)\times \sj_E$ (for each $s$), always splitting as a product of some
almost-complex structure on the first factor with $\sj_E$ on the
second,
\item 
over the rest of $\Sym^{g+1}(\Sigma')$, 
$J_s'\equiv \Sym^{g+1}(\sj'(T))$.
\end{itemize}
Note that this description fits together continuously, 
in light of our hypothesis that $J_s\equiv \Sym^{g}(\sj)$ 
near $\Sym^{g-1}(\Sigma)\times\{\sigma_1\}$.

Fix a class $\phi\in \pi_{2}({\mathbf x},{\mathbf y})$, with
$n_{\sigma}(\phi)=n$ for
$\sigma\in\Sigma-\alpha_{1}-\ldots-\alpha_{g}-\beta_{1}-\ldots-\beta_{g}$.
There is a naturally induced map $$\rho_{1} \colon
\ModFlow_{J_s}(\phi)\longrightarrow
\Sym^{n}\left(\Sym^{g-1}\left(\Sigma)\times 
\Sym^{1}(E\right)\right)$$
given by 
$$\rho_{1}(u)= u\left(u^{-1}\left(\Sym^{g-1}(\Sigma)
\times \{\sigma\}\right)\right)\times 
\{c\}.$$ 
There is a complex codimension one subset in the range corresponding
to the subvariety $\Sym^{g-2}(\Sigma)\times\{\sigma\}\subset
\Sym^{g-1}(\Sigma)$.  For general choice of $\sigma$, all
one-dimensional moduli spaces of flows miss this special locus --
i.e. the images of $u$ in the symmetric product never contain the
prospective connected sum point $\sigma$ with multiplicity greater
than one.

Consider the moduli space of (unparameterized) holomorphic 
maps
$$\ModFlow\left(\bigcup_{i=1}^{n}{\mathbb S}\longrightarrow 
\Sym^{g-1}(\Sigma)\times \Sym^2(E)\right)
$$ from $n$ disjoint Riemann spheres, whose restriction to each
component represents the generator of $\pi_{2}(\Sym^{2}(E))$ (and
hence is constant on the first factor).

There is a map
$$
\rho_{2}\colon 
\ModFlow\left(\bigcup_{i=1}^{n}{\mathbb S}\longrightarrow 
\Sym^{g-1}(\Sigma)\times \Sym^2(E)\right)\longrightarrow 
\Sym^{n}\left(\Sym^{g-1}\left(\Sigma\right)\times 
\Sym^{1}\left(E\right)\right),$$
given by
$$ \rho_{2}(v)= v\left(v^{-1}\left(\Sym^{g-1}(\Sigma)\times \{\sigma\times 
\Sym^{1}(E)\}\right)\right)
$$

In the next section, we show how to splice the spheres coming from 
this moduli space to the disks in 
$\ModFlow(\phi)\times\{c\}$, giving rise to a map from the 
fibered product
of $\rho_{1}\times \rho_{2}$ to the space of maps connecting $\x'=\x\times\{c\}$ to 
$\y'=\y\times\{c\}$ which are nearly holomorphic. Indeed, this
fibered product description is particularly simple in view of the 
fact that $\rho_{2}$ is a diffeomorphism (see 
Lemmas~\ref{lemma:ModSphere} and ~\ref{lemma:IdentDefTheorySphere} 
below). In particular, each $u\in\ModFlow_{J_s}(\phi)$ has a unique 
corresponding $v\in\ModSpheres$ with $\rho_{1}(u)=\rho_{2}(v)$.

Our aim, then is to prove the following:

\begin{theorem}
        \label{thm:Gluing} Fix $\x,\y\in\Ta\cap\Tb$ and a class
        $\phi\in\pi_{2}({\mathbf x},{\mathbf y})$ with $\Mas(\phi)=1$.
        For $\x'=\x\times\{c\}$ and $\y'=\y\times\{c\}$, let
        $\phi'\in\pi_2(\x',\y')$ be the class with the property that
        $n_\sigma(\phi)=n_{\sigma'}(\phi')$ where $\sigma$ is any
        point in $\CurveComp$, and $\sigma'$ is the corresponding
        point in $\Sigma'$.  Then, for all sufficiently large $T$,
        there is a diffeomorphism $\ModFlow_{J_s}(\phi)\cong
        \ModFlow_{J_s'(T)}(\phi')$.
\end{theorem}

\begin{remark}
\label{rmk:IdentDefTheory}
    The diffeomorphism statement above is to be interpreted as an 
    identification between deformation theories. In particular, 
    $\Mas(\phi)=\Mas(\phi')$. 
\end{remark}
   
We argue that the $J_s'(T)$ are $\sj'(T)$-nearly symmetric, and
can be used to calculate the Floer
homologies in $\Sym^{g+1}(\Sigma')$. 
Fix a $T$ so that Theorem~\ref{thm:Gluing} holds for all $\phi$ with $\Mas(\phi)=1$, 
then, we fix a generic path $I_s'$
lying in the open set ${\mathcal U}'$ of
Theorem~\ref{thm:GromovInvariant}. 
Connecting $J_s'(T)$ and $I_s$ by a one-parameter family of
$\sj(T)$-nearly symmetric paths of almost-complex structures
(as in the proof of Theorem~\ref{thm:IndepCxStruct}), we see that the Floer
homology calculated using the path $J_s'(T)$ is the same as that
calculated using $I_s'$. 

\subsection{Approximate Gluing}
\label{subsec:Splicing}

We describe how to splice spheres to flows for $\Sym^{g}(\Sigma)$ to 
obtain flows in $\Sym^{g+1}(\Sigma')$. First, we must introduce some 
notation. 

Write $W=\Sym^{g-1}(\SigOne)\times\{\sigma_1\}$ for the subvariety of 
$V=\Sym^{g}(\SigOne)$. We consider flows $u\in\ModFlow({\bf x},{\bf 
y})$ which meet $W$ transversally in $n=n_{\sigma_{1}}[u]$ distinct 
points $\{q_{1},\ldots,q_{n}\}$. For such a flow $u$, 
it is easy to see that there 
are constants $0<r_1<R_1$ with the property that
the intersection of the image of $u$ with 
$B_{r_1}(\sigma_{1})\times 
\Sym^{g-1}(\SigOne-B_{r_1}(\sigma_{1}))$ is contained in the subset
\begin{equation}
\label{eq:ProdIdent}
B_{r_1}(\sigma_{1})\times 
\Sym^{g-1}(\SigOne-B_{R_1}(\sigma_{1}))\subset
\Sym^{g}(\SigOne).
\end{equation}
We introduce the corresponding product-like metric on this region. 
Moreover, we fix a conformal identification
\begin{equation}
    \label{isom:ConfIdent}
    B_{1}(\sigma_{1})-\sigma_{1}\cong [0,\infty)\times S^{1}
\end{equation}
(observe that all the $J_s$ agree with the complex structure in this region).
We can find an $n$-tuple of disjoint balls $B_{\epsilon}(q_{i})$ which
are mapped into the subspace from Equation~\eqref{eq:ProdIdent}
by $u$. By elementary complex
analysis, we can find an $n$-tuple
$(w_{1},t_{1}+i\theta_{1},\ldots,w_{n},t_{n}+i\theta_{n})$, where $w_i\in \Sym^{g-1}(\SigOne)$, $t_i\in \R$, 
$\theta_i\in [0,2\pi)$, 
with the
property that, with respect to a conformal identification with the
unit disk with the ball around $q_{i}$
\begin{equation}
    \label{eq:ConfIdentDisk}
    \CDisk\cong B_{\epsilon}(q_{i})
\end{equation}
and the identification of the range with the product from Equation~\eqref{eq:ProdIdent}, 
the restriction of $u$ to 
$B_{\epsilon}(q_{i})$ can be 
written as 
\begin{equation}
    \label{eq:AsymptoticsOne}
    z\mapsto (w_{i},e^{t_{i}+i\theta_{i}}z)+\BigO(|z|^{2}).
\end{equation}
It is convenient to recast this in terms of cylindrical coordinates, 
with respect to fixed conformal identifications 
$$B_{\epsilon}(q_{i})-q_{i}\cong [0,\infty)\times S^{1},$$
We use weighted Sobolev spaces 
with weight function $e^{\delta\tau_{1}}$, where
$$\tau_{1}\colon \SigOne-\{q_{1},\ldots,q_{n}\}\longrightarrow 
[0,\infty)$$ is a smooth function satisfying:
\begin{itemize}
    \item $\tau_{1}$ vanishes in the complement of 
    $B_{\epsilon}(q_{i})$,
    \item $\tau_{1}(s+i\varphi)\equiv s$ for $s\geq 1$ over each 
    cylinder.
\end{itemize}

In terms of the cylindrical coordinates, then, the Taylor expansion 
in Equation~\eqref{eq:AsymptoticsOne} gives the following.
For each $i$, 
there is $w_{i}\in \Sym^{g-1}(\SigOne)$ and $t_{i}+i\theta_{i}$,
for which the restriction of
$u$ to $B_{\epsilon}(q_{i})-q_i\cong [0,\infty)\times S^{1}$
differs by a $\Sobol{p}{1,\delta}$ map (for some $\delta>0$)
from the smooth map 
$$a_{t_{i}+i\theta_{i},w_{i}}\colon 
[0,\infty)\times S^{1}\longrightarrow 
\Sym^{g-1}(\SigOne)\times [0,\infty)\times S^{1}\subset 
\Sym^{g}(\SigOne)$$  
defined by
$$a_{t_{i}+i\theta_{i},w_{i}}(s+i\varphi)=(w_{i},(s+t_{i})+i(\varphi+\theta_{i})),$$
where we have used the conformal identification
$$\Sym^{g-1}(\SigOne-B_{R_1}(\sigma_{1}))\times \left([0,\infty)\times S^{1}\right)
\cong \Sym^{g-1}(\SigOne-B_{R_1}(\sigma_{1}))\times 
\left(B_{r_1}(\sigma_{1})-\sigma_1\right)\subset 
\Sym^{g}(\SigOne)$$
(strictly speaking, when $t_{i}<0$, we must cut off $s+t_{i}$ 
in the region where $s<-t_{i}$).

We can use these asymptotics to ``cut off'' the pseudo-holomorphic map $u$ to 
construct a nearly pseudo-holomorphic map into $\Sym^{g}(\SigOne-\sigma_{1})$, 
which agrees with a standard map between cylinders, as follows. Given 
a real number $T>0$, 
let $X_{1}(T)$ denote the subset of $\Strip-\{q_{1},\ldots,q_{n}\}$ 
consisting of points with cylindrical coordinate $\leq T$, i.e. 
$X_{1}(T)=\tau_{1}^{-1}((0,T])$, and $X_1(\infty)\cong 
\Strip-\{q_{1},\ldots,q_{n}\}$. 
Consider the map
$${\widetilde u}_{T}\colon X_1(\infty)\longrightarrow \Sym^{g}(\SigOne)$$
defined to agree with $u$ away from the balls $B_{\epsilon}(q_{i})$
and over $B_{\epsilon}(q_{i})\cong 
[0,\infty)\times S^{1}$, defined by
$${\widetilde u}_{T}(s+i\varphi)=
h(s-T)a_{t_{i}+i\theta_{i},w_{i}}(s+i\varphi)+(1-h(s-T))u(s+i\varphi),$$
where $h\colon \R\longrightarrow [0,1]$ is a smooth, 
increasing cut-off function which is identically $0$ for $t<0$ and 
$1$ for $t>1$. In the latter formula, the convex combination is to be 
interpreted using the exponential map on the range, with respect to 
the product metric.

When considering $\delta$-weighted spaces for functions defined over 
the entire $\Strip-\{q_{i}\}$, we use the weight 
function $e^{\delta \tau_{1}}$ in addition to the weights on $\Strip$ 
described in Section~\ref{subsec:FredholmTheory}. There, we chose some 
smooth function $\tau_{0}\colon [0,1]\times \R \cong \Strip 
\longrightarrow \R$ with $\tau_{0}(s+it)=|t|$ for all sufficiently 
large $t$ (and here we can assume that the support of $\tau_{0}$ 
misses the balls $B_{\epsilon}(q_{i})$), using the weight function 
$e^{\delta_{0}\tau_{0}}$, for some $\delta_{0}$ depending on the 
local geometry around the intersections $\Ta\cap \Tb$. Now, 
suppose 
$$f\colon X_{1}(\infty)\longrightarrow \Sym^{g}(\SigOne-\sigma_1)$$ 
is a map which sends the cylindrical ends of $X_{1}(\infty)$ to the locus
$\sigma_1\times\Sym^{g-1}(\SigOne)$. 
Then, we define the $\Sobol{p}{k;\delta_{0},\delta}$ norm of 
a section $\xi\in \Sections(X_{1}(\infty),f^{*}(T\Sym^{g}(\Sigma)))$ by
$$\|\xi\|_{\Sobol{p}{k;\delta_{0},\delta}(f)}
= \sum_{i=0}^{p}\int_{\Strip-\{q_{i}\}}|\nabla^{(i)}\xi|^{p}
e^{\delta_{0}\tau_{0}}e^{\delta\tau_{1}}d\mu,$$
where $\mu$ is the volume form for $\Strip-\{q_{i}\}$ with the 
cylindrical-end metric (obtained using the cylindrical identification 
of Equation~\eqref{eq:ConfIdentDisk}), and the norm of $\nabla^{(i)}$ denotes the 
metric induced from the pull-back of the cylindrical metric on 
$T\Sym^{g}(\SigOne-\sigma_{1})$
(obtained from the identification~\eqref{eq:ProdIdent} together with~\eqref{isom:ConfIdent}). 
The Banach space $\Sobol{p}{k;\delta_{0},\delta}(f)$ 
is the space of sections
$$\Sobol{p}{k;\delta_{0},\delta}(f)
:=\left\{\xi\in \Sobol{p}{k;\delta_{0},\delta}(\Strip-\{q_{i}\},
f^{*}(T\Sym^{g}(\Sigma)))\Big| 
\begin{array}{l}
\xi(0,t)\in T_{f(0+it)}(\Tb), \forall t\in\R,\\
\xi(1,t)\in T_{f(1+it)}(\Ta), \forall t\in\R
\end{array}\right\}.$$
Similarly, we let $\Sobol{p}{\delta_0,\delta, k-1}(\Lambda^{0,1} f)$
denote the Banach space of
$\Sobol{p}{k-1;\delta_{0},\delta}$ sections of
$f^{*}(T\Sym^{g}(\Sigma))$ 
with a doubly-weighted norm, without the boundary conditions (note
that, as in Section~\ref{sec:Analysis} the notation is justified by
the fact that the bundle $\Wedge^{0,1}$ is trivial over $\Strip$).
For our purposes, it suffices to use $k=1$, and
fix some $p>2$ once and for all.  In fact, in the interest of notational 
expediency, we typically drop the $\delta_{0}$ from the notation, 
simply writing $\Sobol{p}{1,\delta}(f)$ and 
$\Sobol{p}{\delta}(\Lambda^{0,1} f)$, with some $\delta_{0}$ being 
understood as fixed.
Moreover, the results below will hold for all $\delta>0$ sufficiently small (indeed,
any $\delta\in (0,p)$ will do).

\begin{lemma}
    \label{lemma:ErrorTermOne}
Fix $u$ as above. There are constants $t$, $k>0$, $S_{0}>0$,
and $C>0$ so that for each
$S>S_{0}$, the map ${\widetilde u}_{S}$ constructed above
satisfies
\begin{equation}
    \label{eq:ApproxHolomorphic} \|\DBar_{J_s} {\widetilde
    u}_{S}\|_{\Sobol{p}{\delta}(\Lambda^{0,1} {\widetilde u}_{S})}\leq
    Ce^{-kS}.
\end{equation}
Moreover, for each $S$, there is $T_{0}$ so that for all $T>T_{0}$,
${\widetilde u}_S$ maps $X_{1}(T)$ into 
$$\Sym^{g}(\SigOne-[T+t,\infty)\times S^{1}).$$
\end{lemma}

\begin{proof}
    Choose holomorphic coordinates 
    $(z_{1},\ldots,z_{g})$ centered at
    $$w_{i}\times \sigma_{1}\in  
    \Sym^{g-1}(\SigOne-B_{R_{1}}(\sigma_{1}))\times 
    B_{r_{1}}(\sigma_{1}),$$
    with respect to which $z_{g}=0$ corresponds to the intersection of 
    the coordinate patch with the subvariety 
    $\Sym^{g-1}(\Sigma)\times\{\sigma_{1}\}$. 
    (These are $\Sym^g(j)$-holomorphic coordinates,
    which are actually holomorphic for all the
    $J_s$, in the specified region.)
    With respect to these coordinates, we can write $u$ as 
    $(u_{1},\ldots,u_{g-1},f)$. Next, 
    fix a local coordinate function $z$ around the point $q_{i}$, so that $f(0)=0$. 
    We can factor $f(z)=c z g(z)$, 
    where $c=f'(0)$ and $g(0)=1$. Thus, in cylindrical coordinates, we can write 
    this as the map sending $s+i\varphi\in [0,\infty)\times S^{1}$ to 
    $$-\log f(e^{-(s+i\varphi)})=-\log c + s + i\varphi - \log 
    g(e^{-(s+i\varphi)}).$$
    Now $z\mapsto \log g(z)$ vanishes at the origin, so that there is 
    a constant $c_{1}$ 
    with the property that $|\log g(e^{-(s+i\varphi)})|\leq 
    c_{1}e^{-s}$. Applying the above cut-off construction forces us 
    to consider the new function
    $${\widetilde f}(e^{-(s+i\varphi)})=-\log c + s+i \varphi + (h(s-S)-1)\log 
    g(e^{-(s+i\varphi)}).$$
    Now,  $\DBar$ of this is given by
    $$\DBar {\widetilde f} = (\dbar h(s-S))\log g(e^{-(s+i\varphi)}),$$
    so we have a constant $c_{2}$ with the property that
    $$|\DBar \widetilde f(s+i\varphi)|\leq c_{3} e^{-s}$$
    pointwise; moreover, this tensor is supported in the strip where 
    $s$ is constrained to lie in the interval $[S,S+1]$. Thus
    there is a $c_{3}$ with 
    \begin{equation}
        \label{eq:NormalComponent}
        \int_{[0,\infty)\times S^{1}}|\DBar \widetilde f|^{p}e^{\delta s}\leq c_{3} 
        e^{(-p+\delta)S}.
    \end{equation}

    To handle the other components $u_{i}$, we proceed in a similar 
    manner. Noting that $u_{i}(0)=0$, the cutting off construction 
    gives
    ${\widetilde u}_{i}(e^{-(s+i\varphi)})=(1-h(s-S))u_{i}(e^{-(s+i\varphi)})$.
    Applying $\DBar$ to this, we obtain
    $$\DBar {\widetilde u}_{i}(e^{-(s+i\varphi)})=-(\DBar 
    h(s-S))u_i(e^{-(s+i\varphi)}).$$ We have a constant $c_{4}$ with 
    $|u_{i}(z)|\leq c_{4} z$, so that 
    $$|\DBar {\widetilde u}_{i}(s+i\varphi)|\leq c_{4}e^{-s}$$
    pointwise. Once again, since the tensor is supported in a bounded 
    strip, we have that 
    \begin{equation}
        \label{eq:OtherComponents}
        \int_{[0,\infty)\times S^{1}}
        |\DBar \widetilde u_{i}|^{p}e^{\delta s}\leq c_{3} 
        e^{(-p+\delta)S}.
    \end{equation}
    Together, Inequalities~\eqref{eq:NormalComponent} and 
    \eqref{eq:OtherComponents} give 
    Inequality~\eqref{eq:ApproxHolomorphic}, with 
    respect to a cylindrical-end metrics on the range 
    $\Sym^{g}(\SigOne-\sigma_{1})$ (which, 
    over the subset $\Sym^{g}(\SigOne-B_{r_{1}})$, can be made to 
    be isometric to the corresponding subset of
    $\Sym^{g}(\SigOne)$).

    The second part of the lemma is straightforward.
\end{proof}
    
Next, we turn our attention to spheres in $\Sym^2(E)$. More 
precisely, consider  holomorphic maps
$$ v \colon {\mathbb S}\longrightarrow \Sym^{g-1}(\SigOne)
\times \Sym^2(E) $$ which are
constant on the first factor, and which represent the positive generator of
$\pi_2(\Sym^2(E))$ on the second (i.e. they satisfy 
$n_{\sigma_{2}}[v]=1$).  We denote the space of such maps, modulo 
holomorphic reparameterization, by
$$\ModSphere.$$

\begin{lemma}
    \label{lemma:ModSphere}
    The map 
    $$\ModSphere\longrightarrow \Sym^{g-1}(\SigOne)\times 
    E$$ which takes $[v]\in\ModSphere$ to the unique pair $(w,y)$, with
    the property that $(w,\{y,\sigma_{2}\})$ is in the image of $v$, 
    induces a one-to-one correspondence.
\end{lemma}

\begin{proof}
    This is true because $\pi\colon \Sym^{2}(E)\longrightarrow 
    E$ is a ruled surface (see Lemma~\ref{lemma:RSurf}), 
    the generator corresponds to 
    holomorphic spheres representing 
    the fiber class (according to the proof of that lemma),
    and the map $E\longrightarrow 
    \Sym^{2}(E)$ given by $\sigma\mapsto 
    \{\sigma,\sigma_{2}\}$ is a section. 
    Thus, the composite $\pi\circ v$ is a holomorphic 
    map from the sphere to the torus, which must be constant. Hence, 
    $v$ maps to a fiber. Indeed, from our assumption on the homotopy 
    class, it follows that $v$ must map isomorphically to a fiber, and 
    each fiber is determined by its intersection with the section.
\end{proof}
    
Let $v$ be as above. Thinking of ${\mathbb S}$ as the Riemann sphere 
${\mathbb S}\cong \C\cup\{\infty\}$, we 
normalize $v$ so that $v(0)=w\times \{y,\sigma_2\}$ (this can be achieved by
precomposing $v$ with a M\"obius transformation if necessary), and
suppose moreover that $y\neq \sigma_2$. A neighborhood of the image
$w\times \{y,\sigma_{2}\}$ can be identified with a product
$$\Sym^{g-1}(\SigOne)\times B_{r_{2}}(\sigma_{2})\times 
\left(E-B_{R_{2}}\left(\sigma_{2}\right)\right)\subset 
\Sym^{g-1}(\SigOne)\times\Sym^{2}(E),$$
for some $r_{2}<R_{2}$. 
With respect to these coordinates, near $0$, $v$ then takes the form
$$z\mapsto (w,\sigma_2+\BigO(|z|),e^{t+i\theta}z)+\BigO(|z|^2)$$
for some $t+i\theta$ (when the third coordinate $=0$ corresponds to
the point $y\in E$).
In terms of cylindrical coordinates in $[0,\infty)\times S^1\cong 
B_{\epsilon}(0)-0\subset {\mathbb S}$ 
and 
$[0,\infty)\times S^{1}\cong B_{r_{2}}(\sigma_{2})$, this can
be phrased as follows:
$$v\colon [0,\infty)\times S^1\longrightarrow 
\Sym^{g-1}(\SigOne)\times \left([0,\infty)\times S^{1}\right)\times E$$
satisfies the property that there is a $w\in W$ and $t+i\theta$ so
that $v$ differs from the map
$$b_{(t+i\theta,w,y)}(s+i\varphi)=(w,(s+t)+i(\theta+\varphi),y)$$ by a
map which is in $\Sobol{p}{1,\delta}$ for any $\delta>0$. (Here, 
$\Sobol{p}{1,\delta}$ is defined with respect to a weight function 
$e^{\delta\tau_{2}}$, for $\tau_{2}\colon {\mathbb 
S}-0\longrightarrow \R^{+}$ defined in a manner analogous to 
$\tau_{1}$.)

In view of these asymptotics, we can form a ``cutting off'' construction as before.
Given a real number $S>0$, let ${\mathbb S}(S)$ denote the subspace 
of $x\in{\mathbb S}$ where $\tau_{2}(x)\leq S$. 
Next, consider the map
$$v_{S}\colon {\mathbb S}-\infty\longrightarrow 
\Sym^{g-1}(\SigOne)\times \Sym^{2}(E-\sigma_{2})$$ 
defined to agree with $v$ over ${\mathbb S}(S)$, and given by the 
formula
$$
v_{S}(s+i\varphi)= 
h(s-S) b_{t+i\theta}+(1-h(s-S))v(s+i\varphi)
$$
over the cylinder $[0,\infty)\times S^{1}\cong B_{r_{2}}(0)-0$.

The six-dimensional group ${\mathbb P}{\mathrm{Sl}}_{2}(\C)$ acts on the space of 
parameterized holomorphic maps from ${\mathbb S}$. 
The normalization condition 
on $v$ above cuts out two dimensions from this automorphism group. An 
additional two dimensions could be cut out by the following conditions. 
The holomorphic map $v$ has an ``energy measure'' obtained by pulling 
back the symplectic form from $\Sym^{2}(E)$. This measure in 
turn has a center of mass in the three-ball. Specifying the center of
mass lies on the $z$-axis (the axis connecting the points
corresponding to $0$ and $\infty$ in ${\mathbb S}$) 
cuts out an additional 
two-dimensions from the automorphism group, but instead we find it 
convenient to formulate this condition as follows.
Let $v_{0}$ be a normalized holomorphic sphere 
with $v_{0}(0)=w\times\{c,\sigma_2\}$, 
whose center of mass
lies on the $z$ axis, and choose  a point $\sigma_{2}'\in E$ with the 
property that $v_{0}(\infty)$ lies in the submanifold 
$W'=\Sym^{g-1}(\SigOne)\times \sigma_{2}'\times E$ (there
are two possible choices for $\sigma_{2}'$).

\begin{defn}
    \label{def:Centered}
A holomorphic map 
$$ v \colon {\mathbb S}
\longrightarrow \Sym^{g-1}(\SigOne)\times \Sym^2(E) $$ is called
{\em centered} if the following two conditions are satisfied:
\begin{itemize}
    \item $v(0)=w\times\{y,\sigma_{2}\}$ for some $w\in\Sym^{g-1}(\SigOne)$
     and $y\in E-\sigma_2$.
    \item $v(\infty)\in W'$.
\end{itemize}
\end{defn}

We can collect the set of centered holomorphic maps into a moduli
space, which we denote by
$$\ModCent.$$ Since the space of holomorphic automorphisms which respect
the centered condition is identified with $\C^*$, this moduli space
has a $\C^*$ action, and fibers over the region in $\ModSphere$ 
corresponding to $\Sym^{g-1}(\SigOne)\times 
\left(E-B_{R_2}(\sigma_{2})\right)$.
Now, we can view the assignment $v\mapsto (w, t+i\theta,y)$ as a
map from the moduli space of centered maps
$$\rho_2\colon \ModCent\longrightarrow
\Sym^{g-1}(\SigOne)\times \R\times S^1\times E.$$
It is an easy consequence of Lemma~\ref{lemma:ModSphere} that 
$\rho_{2}$ is a $\C^{*}$-equivariant diffeomorphism.
Given 
$(w,t,i\theta)\in \Sym^{g-1}(\SigOne)\times \R\times S^{1}$, let 
$$v_{(w,t,i\theta)}\colon{\mathbb S}\longrightarrow 
\Sym^{g-1}(\SigOne)\times\Sym^{2}(E)$$
denote the centered map with the property 
$\rho_{2}v_{(w,t,i\theta)}=(w,t+i\theta,c)$.

Given the $J_s$-holomorphic disk $u$, 
we would like to splice in a number of centered spheres, to construct a 
nearly $J_s'(T)$-holomorphic curve in $\Sigma'(T)$. 
For instance, at the marked 
point $q_{i}$, with value $w_{i}$ we splice in a copy of 
$v_{(w_{i},-t_{i},\theta_{i})}$.

More precisely, using the given conformal identifications 
of Equation~\eqref{isom:ConfIdent}
of the neighborhoods of the puncture points,
the connected sum $\Sigma'(T)$ can be thought of as the space obtained 
from $\SigOne(2T)$ and $ E(2T)$, 
by identifying the cylinders 
\begin{eqnarray*}
[0,2T]\times S^{1}\subset \SigOne(2T)
&{\text{and}}&[0,2T]\times S^{1}\subset  E(2T)
\end{eqnarray*}
using the involution
$(t,i\theta)\sim (2T-t,i\theta)$. 
Let 
$$ X_{2}(T)=\bigcup_{i=1}^{n}{\mathbb S}(T)_{i},$$
and $X_{1}\cup_{T}X_{2}$ denote the union of $X_{1}(T)$ with 
$X_{2}(T)$ glued along their common boundary.
Now, given a holomorphic disk $u$ 
and a pair $S$ and $T$ of real numbers with $0<S<T-t$, 
let 
$${\widetilde v}_{S}\colon X_{2}(\infty)\longrightarrow 
\Sym^{g-1}(\SigOne)\times \Sym^{2}(E-\sigma_{2})$$
be the map
whose restriction to ${\mathbb 
S}(T)_{i}\subset X_{2}(\infty)$ is 
the cut-off of $v_{(w_{i},t_{i}+i\theta_{i})}$. 
We can form a spliced map 
$${\widehat \gamma}={\widehat \gamma}(u,S,T) 
\colon X_{1}\cup_{T} X_{2}\longrightarrow
\Sym^{g+1}(\SigOne\#_{T} E),$$
which is defined to agree with ${\widetilde u}_{S}\times \{c\}$ 
over $X_{1}(T)$, and to agree with
the map ${\widetilde v}_{S}$ over $X_{2}(T)$.
Note that $X_{1}\cup_{T}X_{2}$ 
is conformally equivalent to $\Strip$, though 
its geometry is different. Reflecting these geometries, we find it 
convenient to introduce weighted Sobolev spaces over $X_{1}\cup_{T}X_2$. 
Fix some $p>2$, and let 
${\mathcal Z}(T)$ denote the space of sections 
$\Gamma\left(X_{1}\cup_{T} X_{2},
\Wedge^{0,1}\otimes {\widehat \gamma}^{*}(T\Sym^{g+1}(\Sigma'))\right)$, 
endowed with the 
norm
$$\|\phi\|_{{\mathcal Z}_{\delta}(T)} = \left(\int_{X_{1}\cup_{T}X_{2}}|\phi|^{p}e^{\delta_{0}\tau_{0}}e^{\delta \tau}\right)^{1/p},$$
where $\tau_{0}$ is the weight function used near the two ends of the 
strip $$[0,1]\times \infty \cong X_{1}\cup_{T}X_{2},$$
and 
$\tau\colon X_{1}\cup_{T}X_{2}\longrightarrow \R$ is a smooth 
function (depending, of course, on $T$) which 
agrees with $\tau_{1}$ and $\tau_{2}$ over $X_{1}(T-1)$ and 
$X_{2}(T-1)$, and satisfies a uniform $\Cinfty$ (i.e. independent of 
$T$) bound on $d\tau$.

\begin{lemma}
\label{lemma:ErrorTerm}
    Given any sufficiently large $S$, there is a $T_{0}$ so that for 
    all $T>T_{0}$, the map 
    ${\widehat \gamma}(u,S,T)$ determines a smooth map
    $${\widehat \gamma}(u,S,T)\colon 
    X_{1}\cup_{T}X_{2}\longrightarrow \Sym^{g+1}(\SigOne\#_{T} E).$$
    Moreover,  there are positive constants $k$ and $C$, so that for all 
    $S>S_{0}$ and $T>T_{0}$, 
    ${\widehat \gamma}(u,S,T)$ satisfies
    $$\|\DBar_{J_s'} {\widehat \gamma}(u,S,T)\|_{{\mathcal Z}_{\delta}(T)}\leq C 
    e^{-k S}.$$
\end{lemma}

\begin{proof}
    This follows immediately from Lemma~\ref{lemma:ErrorTermOne}, and the
    corresponding estimate for ${\widetilde v}_{S}$ (which follows in
    the same manner). Note that the image of ${\widehat\gamma}$ is
    contained in the region where the divisor on the $E$ side is
    bounded away from $\sigma_2$, so that the almost-complex
    structures $J_s'$ split as $J_s\times {\sj_E}$.
\end{proof} 
    
\subsection{Paramatrices}~\newline
\label{subsec:FredTheory}

To perform the gluing, we must construct an inverse for the linearization of the 
$\DBar$-operator for the spliced map from  $X_{1}\cup_{T}X_{2}$, 
whose norm is bounded independent of $T$ (see 
Proposition~\ref{prop:UniformParamatrix} below). To set this up, we 
describe suitable function spaces for the operators over the
cylindrical-end versions:
$$D_{\widetilde u}\colon 
\Sections(X_{1}(\infty),E_{1})\longrightarrow 
\Sections(X_{1}(\infty),F_{1})$$
and 
$$D_{\widetilde v}\colon
\Sections(X_{2}(\infty),E_{2})\longrightarrow 
\Sections(X_{2}(\infty),F_{2}),$$
where:
\begin{eqnarray*}
    E_{1}&=&{\widetilde u}^{*}\left(T\Sym^{g}
    \left(\SigOne-\sigma_{1}\right)\right)\\
    F_{1}&=&\Hom(TX_{1}(\infty),{\widetilde u}^{*}\left(T\Sym^{g}
    \left(\SigOne-\sigma_{1}\right)\right)), \\
    E_{2}&=&{\widetilde v}^{*}\left(T\left(\Sym^{g-1}
    \left(\SigOne-\sigma_{1}\right)\times\Sym^{2}
    \left(E-\sigma_{2}\right) \right)\right) \\
    F_{2}&=&\Hom(TX_{2}(\infty),{\widetilde v}^{*}\left(T\left(\Sym^{g-1}
    \left(\SigOne-\sigma_{1}\right)\times\Sym^{2}
    \left(E-\sigma_{2}\right) \right)\right)).
\end{eqnarray*}
Note that, given the complex structures on $X_{1}$ and $X_{2}$, there 
are canonical identifications
\begin{eqnarray*}
    F_{1}\cong \Wedge^{(0,1)}\otimes E_{1} &{\text{and}}&
    F_{2}\cong \Wedge^{(0,1)}\otimes E_{2},
\end{eqnarray*}
but we 
have chosen to write $F_{1}$ and $F_{2}$ as above, as we
will allow the complex structure on $X_{1}(\infty)$ to vary.

Indeed, to allow for this variation, we enlarge the domain of 
$D_{\widetilde u}$, by $\C^{n}\cong \bigoplus_{i=1}^{n}T_{q_{i}}\Strip$.
Let $\phi_{\xi}\in \Diff(\Strip)$ be a family of 
diffeomorphisms of $\Strip$ indexed by $\xi\in\C^n$,
which satisfies the following properties:
\begin{itemize}
    \item $\phi_{0}=\Id$
    \item $\phi_{\xi}|\Strip-\bigcup B_{\epsilon}(q_{i})\equiv \Id$,
    \item If $\xi$ is sufficiently small, then 
    $\phi_{\xi}|B_{\epsilon/2}(q_{i})\longrightarrow \Strip$ is 
    translation by $\xi_{i}$.
\end{itemize}
A family of diffeomorphisms $\phi_{\xi}$ as above
is constructed, for example, by fixing a cut-off 
function $h$ supported in $\bigcup B_{\epsilon}(q_{i})$ which is 
identically one on $B_{3\epsilon/4}(q_{i})$. Then, if $\xi\in 
\bigoplus T_{q_{i}}(\Strip)$, we let ${\widetilde \xi}$ be its 
locally constant extension to $\bigcup_{i=1}^{n}B_{\epsilon}(q_{i})$. 
Then, $\phi_{\xi}$ is obtained by exponentiating the compactly 
supported vector field $h {\widetilde \xi}$.
We then define $j_{\xi}=\phi_{\xi}^{*}j_{0}$, where $j_{0}$ is the 
standard complex structure on $\Strip$. 

We then consider a modified Cauchy-Riemann operator
$$
        \DBar^{\ext}_{J_s}\colon 
        \Map(X_{1}(\infty),\Sym^{g}(\Sigma))
        \times 
        \C^{n}
        \longrightarrow
        \Sections(X_{1}(\infty), F_{1}),
$$
which assigns to $(u,\xi)\in 
\Maps(X_{1}(\infty),\Sym^{g}(\Sigma))
\times \C^{n}$ the bundle map from $T\Strip$ 
to $F_{1}$ given by 
$$
        \DBar^{\ext}_{J_s}(u,\xi)=(u^{*}
        J_{\phi_\xi^{-1}})\circ du - du\circ j_{\xi},
$$ 
where
$J_{\phi_\xi^{-1}}$ is the family of almost-complex structures over
$\Strip$ given by
$$
        J_{\phi_\xi^{-1}}(s,t) = J(\phi_\xi^{-1}(s,t))
$$ 
(where we view our original path $J$ as a map on $\Strip$ which is
constant in the $t$ directions).  Clearly,
$\DBar^{\ext}_{J_s}(u,\xi)=0$ if and only if the map $u\circ
\phi_{\xi}^{-1}$ is $J_s$-holomorphic.
The linearization of $\DBar^{\ext}_{J_s}$ is the map
$$
        D_{(u,\xi)}\colon \Sections(X_{1}(\infty),E_{1})\oplus \C^{n}
        \longrightarrow \Sections(X_{1}(\infty), F_{1})
$$
given by
\begin{eqnarray}
        D_{(u,\xi)}(\nu,x)
        &=&
        (u^*\nabla_\nu J_{\phi_{\xi}^{-1}})\circ du + 
        (u^* J_{\phi_{\xi}^{-1}})\circ (\nabla\nu)
                \label{eq:Linearization} \\
        &&
        - (\nabla\nu)\circ j_{\xi}
        -(du)\circ({\mathcal L}_x j_\xi)
        + (d_x(J_{\phi_{\xi}^{-1}}))\circ du. \nonumber
\end{eqnarray}
The first term here involves a covariant derivative
of the tensor $J_{\phi_\xi^{-1}(s,t)}$ (defined over
$T\Sym^{g}(\Sigma)$);
the fourth term involves the Lie derivative of the complex structure
$j_\xi$ in the direction specified by the vector field corresponding
to $x\in\C^n$, the final term involves a directional derivative of the map 
$J_{\phi_\xi^{-1}}$, thought of as a map from $\Strip$ to the space of
almost-complex structures, in the direction specified by $x$, thought
of as a vector field along $\Strip$.

With these remarks in place, we turn attention to the Sobolev 
topologies to be used on the domain and range. 
For the operator 
$D_{\widetilde u}$, we consider the Banach spaces
$$
        {\sW}_{1}=
        \left(\Sobol{p}{1;\delta_{0},\delta}({\widetilde u})
        +\Harm\right)\oplus \C^n,
$$
and
$$
        {\sZ}_{1}=\Sobol{p}{\delta_{0},\delta}(\Lambda^{1}\otimes
        {\widetilde u}),
$$
where $\Harm$ is the finite-dimensional vector space
$$
        \Harm=\bigoplus_{i=1}^{n}
        \left(T_{w_{i}}\Sym^{g-1}(\SigOne)\oplus \C \right).
$$
Recall that, by definition, 
$\Sobol{p}{1;\delta_{0},\delta}({\widetilde 
u})$ and 
$\Sobol{p}{\delta_{0},\delta}(\Lambda^{1}\otimes 
{\widetilde u})
\cong 
\Sobol{p}{\delta_{0},\delta}(\Lambda^{(0,1)}\otimes 
{\widetilde u})$ are spaces of sections of $E_{1}$ and $F_{1}$ 
respectively (satisfying the appropriate boundary and decay 
conditions, specified earlier). 

To make sense of ${\sW}_{1}$, we 
view $\Harm$ as a space of sections of $E_{1}$, as follows. 
First,
consider
the tangent vectors $T_{w_{i}}\Sym^{g-1}(\SigOne)$ 
as vector fields in ${\widetilde u}^{*}T\Sym^{g-1}(\SigOne)$
which are locally constant over the cylindrical ends. More precisely,
a tangent vector $T_{w_{i}}(\Sym^{g-1}(\SigOne))$
gives rise to a section  of $E_{1}$ over the $i^{th}$
cylindrical end. These vector fields are then
to  be multiplied by a fixed cut-off function (supported in the 
$i^{th}$ end),  to give rise to global
sections of $E_{1}$.
The $\C$ summand in $\Harm$ corresponds to
tangent vector fields 
to the cylinders which are constant, as follows. We identify
the constant vector fields over $S^{1}\times \R$ with $\C$. 
Again, using a smooth cut-off function, we can transfer the 
vector fields over $S^{1}\times \R$ to $X_{1}(\infty)$ (to get 
a vector space of tangent vector fields, which are constant 
at the cylindrical ends).  We can then use the derivative of 
${\widetilde u}$ to push these forward to get sections of $E_{1}$
(which are holomophic over the ends).

The vector space $\C^{n}$ corresponds to variations in the complex 
structure over $X_{1}(\infty)$, corresponding to variations 
of the puncture points on the disks, as described above. 
Then, the differential operator $D_{({\widetilde u},\xi)}$ extended as above
induces a continuous linear map
$$ 
        \Lu\colon \sW_{1}\longrightarrow \sZ_{1}. 
$$
(Perhaps a more suggestive notation for $\sW_{1}$ is to write it as:
$$
        \Sobol{p}{1;\delta_{0},\delta}({\widetilde u})
        +\bigoplus_{i=1}^{n}\left(T_{w_{i}}\Sym^{g-1}(\SigOne)\oplus 
        T_{q_{i}}\Strip\oplus N_{w_{i}}\Sym^{g-1}(\SigOne)\right),
$$
where $N_{w}\Sym^{g-1}(\SigOne)$ is the normal vector space to
$\Sym^{g-1}(\SigOne)\times\{\sigma_{1}\}$ at the point $w\times 
\sigma_{1}$.)
There is also a natural linear projection
$$
        \rho_{1}\colon \sW_{1}\longrightarrow \Harm.
$$
Note that these operators (and indeed the Banach spaces on which they 
are
defined) depend on the real number $S>0$, which we have suppressed 
from the notation. 

\begin{lemma} 
    \label{lemma:IdentDefTheoryStrip}
    For all sufficiently large $S>0$, the operator 
    $\Lu$ has a one-dimensional kernel and no cokernel.
    Moreover, it admits a right inverse, with an inverse whose 
    norm is bounded above independent of $S$. 
\end{lemma}

\begin{proof}
    Consider the operator 
    $$
        D_{(u,\xi)}^{\cyl}\colon 
            \left(\Sobol{p}{1;\delta_{0},\delta}(u)+\Harm
            \right)\oplus \C^{n}\longrightarrow 
            \Sobol{p}{\delta_{0},\delta}(\Lambda^{1}\otimes u).
    $$
    Here, $u$ is viewed as a map between $X_{1}(\infty)$ and 
    $\Sym^{g}(\SigOne-\sigma_{1})$, endowed with cylindrical metrics.
    By the removable singularities theorem,
    the kernel and cokernel are identified with kernel and cokernel of 
    the flow-line in $\Sym^{g}(\SigOne)$, 
    $$
        D_{u} \colon \Sobol{p}{1,\delta_{0}}(u)
        \longrightarrow 
        \Sobol{p}{\delta_{0}}(\Lambda^{1}\otimes u).
    $$
    (See the corresponding discussion in~\cite{APSI}.) By assumption,
    the cokernel vanishes and its kernel is one-dimensional, since we
    assumed that $u$ was a generic flow-line. Moreover, the operators
    $D_{({\widetilde u}_{S},\xi)}$ converge to 
    $D_{(u,\xi)}^{\cyl}$ as $S\goesto
    \infty$. The lemma then follows.
\end{proof}
    
An analogous discussion applies on the other side, as well.
However, to compensate for the $\C^{n}$ increase of the domain on 
the $X_{1}$ side, we must decrease the domain on the $X_{2}$ side.
We define Banach spaces
$$
        {\mathcal W}_{2}^{\ext}=
        \Sobol{p}{1,\delta}(X_{2}(\infty),E_{2})+\Harm
$$
and 
$$
        {\mathcal Z}_{2}=
        \Sobol{p}{\delta}(X_{2}(\infty),F_{2})
$$
(where $\Harm$ is the finite dimensional vector space defined
earlier), 
and once again we have a map:
$$
        \Lv^{\ext}\colon {\mathcal W}_{2}^{\ext}\longrightarrow {\mathcal 
        Z}_{2}
$$
(inherited by the differential operator $D_{\widetilde v}$),
and the linear projection
$$ 
        \rho_{2}\colon \sW_{2}^{\ext}\longrightarrow \Harm. 
$$

As it turns out, $\Lv $ on 
${\sW}_{2}^{\ext}$ has kernel.
To get rid of the kernel, we consider a linearization of the centered 
condition described in the previous subsection. Specifically,
if $x_{1},\ldots x_{n}\in X_{2}(\infty)$ are the points corresponding 
to the origins under the identification $\bigcup_{i=1}^{n}\C\cong X_{2}(\infty)$, 
we assume that ${\widetilde v}(x_{i})\in W'$. Then, we can define
$$\sW_{2}=\{\xi\in\sW^{\ext}_{2}\big| \xi(x_{i})\in T_{v(x_{i})}W' 
~~~\forall i=1,\ldots,n\}.$$

\begin{lemma}
    \label{lemma:IdentDefTheorySphere}
    For all sufficiently large $S>0$, the map $\Lv$ has 
    no cokernel, and its kernel is identified under $\rho_{2}$ with 
    $\Harm$. 
\end{lemma}

\begin{proof}
    This follows as in the proof of 
    Lemma~\ref{lemma:IdentDefTheoryStrip}.  Note that (by a removable singularities theorem) the kernel and 
    cokernel of 
    $$D_{v}^{\ext}\colon \sW_{2}^{\ext}\longrightarrow \sZ_{2}$$
    can be identified with the corresponding spaces for the operator
    $$D_{v}\colon \Sobol{p}{1}(v)\longrightarrow 
    \Sobol{p}{}(\Lambda^{1}\otimes v),$$
    where 
    $$\Sobol{p}{1}(v)=
    \Sobol{p}{1}\left(\bigcup_{i=1}^{n}
    {\mathbb S},v^{*}(T\Sym^{g}(\SigOne))\right)$$
    and 
    $$\Sobol{p}{1}(\Lambda^{1}\otimes v)=
    \Sobol{p}{1}
    \left(\bigcup_{i=1}^{n}{\mathbb S},\Lambda^{1}\otimes 
    v^{*}(T\Sym^{g}(\SigOne))\right).$$
    In turn, these spaces are given by the deformation theory of 
    (parameterized) holomorphic 
    spheres, up to reparameterizations. 
    The fact that the cokernel vanishes follows from the 
    fact that the fiber class of a ruled surface form a smooth moduli 
    space, which in turn follows from an easy Leray-Serre spectral 
    sequence argument, which we defer to the next paragraph.
    For each sphere, then, there is a complex three-dimensional family of 
    holomorphic vector fields which fix the image curve. These generate 
    a three-dimensional space of kernel elements (for each sphere). 
    However, if we restrict the domain
    to ${\sW}_{2}$, we are considering spheres 
    which are specified to lie in the 
    subvariety $W\times W'$ at a fixed pair of points: hence, there 
    is only a one-dimension space $\C^{*}$ of automorphisms left.
    This kernel  is easily seen to be captured by $\rho_{2}$.
    
    To see that the $n$-fold fiber class of $\Sym^{2}(E)$ over $E$
    considered above has a smooth deformation theory, note 
    that each $n$-tuple of fibers is given as the zero set of a section $\sigma$ of 
    a bundle ${\mathcal L}$, which is the pull-back of a Chern class 
    $n$ line bundle ${\mathcal L}_{0}$ over the base $E$. The cokernel of the 
    deformation complex is identified with the one-dimensional 
    cohomology of the quotient sheaf of ${\mathcal L}$ by $\sigma$,
    $$H^{1}(\Sym^{2}(E), {\mathcal L}/\sigma{\mathcal O}).$$
    First, we prove the vanishing
    \begin{equation}
        \label{eq:VanishHOne}
        H^{1}(\Sym^{2}(E),{\mathcal L})=0.
    \end{equation}
    The Leray-Serre 
    spectral sequence for the fibration has an $E_{2}^{p,q}$ term with
    $$H^{p}(E, R^{q}\pi_{*}{\mathcal L})\Rightarrow
    H^{p+q}(\Sym^{2}(E),\pi^{*}{\mathcal 
    L}), $$
    where 
    $R^{q}\pi_{*}$ is the derived functor of the push-forward map 
    $\pi_{*}$. It suffices to prove that both 
    $H^{1}(E,\pi_{*}{\mathcal L})=0$ and 
    $H^{0}(E,R^{1}\pi_{*}{\mathcal L})=0$. The first vanishing 
    follows from Serre duality, since $\pi_{*}{\mathcal L}={\mathcal 
    L}_{0}$, a positive line bundle over an elliptic curve. The second 
    vanishing statement follows from the projection formula 
    $R^{1}\pi_{*}(\pi^{*}{\mathcal L}_{0})\cong {\mathcal 
    L}_{0}\otimes R^{1}\pi_{*}{\mathcal O}=0$ (see for 
    example Chapter III.8 of~\cite{Hartshorne}).
    This proves the vanishing in Equation~\eqref{eq:VanishHOne}. The 
    vanishing of $H^{1}(\Sym^{2}(E), {\mathcal L}/\sigma{\mathcal 
    O})$, then, follows from Equation~\eqref{eq:VanishHOne},
    the long exact sequence in cohomology 
    associated to the short exact sequence of sheaves:
    $$\begin{CD}
    0@>>> {\mathcal O} @>{\sigma}>> {\mathcal L} @>>> {\mathcal 
    L}/\sigma {\mathcal O} @>>> 0,
    \end{CD},$$
    and the fact that $H^{2}(\Sym^{2}(E),{\mathcal O})=0$. 
\end{proof}

We can consider norms on spaces of sections 
$${\sX}(T)^{\ext}=\Sobol{p}{1}(X_{1}\cup_{T} X_{2}, {\widetilde 
u}\#_{T}{\widetilde v}^{*}(T(\Sym^{g+1}(\Sigma')))\oplus \C^{n},$$
which 
reflect the neck stretching. These norms are obtained by using a partition 
of unity
to transfer sections to $X_{1}(\infty)$ and 
$X_{2}(\infty)$ and use the norms of ${\sW}_{1}$ and ${\sW}_{2}$ on those 
spaces. More precisely, choose a partition of unity  $\{\phi_{1},\phi_{2}\}$ on $X_{1}\cup_{T} X_{2}$ 
subordinate to its cover 
$\{X_{1}(T+1), X_{2}(T+1)\}$, constructed so that $d\phi_{1}$ is 
uniformly $\Cinfty$ bounded, independent of $T$. There is an 
injection $$\iota\colon {\sX}(T)^{\ext}\longrightarrow {\sW}_{1}\oplus {\sW}_{2}^{\ext}$$
given by
$$\iota(a,\xi)=
\left((\phi_{1} a + (1-\phi_{1})h, \xi), (\phi_{2} a + 
(1-\phi_{2}) h)\right),$$
where $h=\Proj_{\mathcal H}\left(a|_{S^{1}\times \{0\}}\right)$
over each cylinder,
and 
where $\phi_{1}a + (1-\phi_{1})h$ and $\phi_{2} a + (1-\phi_{2})h$ are 
to be thought of as sections over $X_{1}(\infty)$ and $X_{2}(\infty)$ 
respectively in the obvious manner. 
The norm on ${\sX}(T)^{\ext}$ is defined to make $\iota$ an isometry onto 
its image; i.e. the norm in ${\sX}(T)^{\ext}$ is defined by
$$\|(a,\xi)\|_{{\mathcal X}(T)}=\|\phi_{1}a+(1-\phi_{1})h\|_{\Sobol{p}{1,\delta}(X_{1})+\mathcal H}
+ \|\phi_{2}a+
(1-\phi_{2})h\|_{\Sobol{p}{1,\delta}(X_{2})+\mathcal H}
+\|\xi\|_{\C^{n}},$$

Let $\sX(T)\subset \sX(T)^{\ext}$ denote the subset whose 
image under $\iota$ maps to $\sW_{1}\oplus \sW_{2}\subset \sW_{1}\oplus\sW_{2}^{\ext}$;
i.e. where the vectors at the centers of $X_{2}$ are tangent to $W'$, 
$a(x_{i})\in T_{v(x_{i})}W'$ for $i=1,\ldots,n$.
Similarly, we consider the Banach space
$$\sZ(T)=
\Sobol{p}{}\left(X_{1}\cup_{T}X_{2},\Wedge^{0,1}\otimes {\widetilde u}
\#_{T}{\widetilde v}^{*}(T\Sym^{n}(\SigOne))\right)$$
with norm given by
$$\|b\|_{{\sZ}(T)}
=\|\phi_{1}b\|_{{\sZ}_{1}}+\|\phi_{2}b\|_{{\sZ}_{2}}.$$
It is easy to see that for any fixed $T$, 
the induced topology on $\sZ(T)$ is the same 
as the ordinary unweighted Sobolev topologies
on the corresponding bundle over $X_{1}\#_{T}X_{2}$.
These norms were chosen instead because they satisfy the following property:

\begin{prop}
    \label{prop:UniformParamatrix}
    There is a constant  $C$  with the property that for all $T>T_{0}$,
    $D_{T}=D_{{\widetilde u}\#_T{\widetilde v}}$ has a continuous right inverse 
    $$P_{T}\colon \sZ(T)\longrightarrow {\mathcal X}(T),$$
    whose operator norm is uniformly bounded above (independent of $T$).
\end{prop}

The inverse is constructed by patching together right inverses 
over the cylindrical end versions $X_{1}(\infty)$ and $X_{2}(\infty)$. 
We describe these inverses, and prove 
Proposition~\ref{prop:UniformParamatrix} at the end of the subsection. 

Let ${\mathcal X}_{\infty}$ denote the fibered product of $\rho_{1}$ 
and $\rho_{2}$; i.e. it is the Banach space which fits into the short 
exact sequence
$$\begin{CD}
0@>>>{\mathcal X}_{\infty}@>>> {\mathcal W}_{1}\oplus{\mathcal 
W}_{2}
@>{\rho_{1}-\rho_{2}}>> {\mathcal H} @>>>0
\end{CD}.$$
Letting $\sZ_{\infty}=\sZ_{1}\oplus \sZ_{2}$, the maps $\Lu$ and $\Lv$ induce a map
$$\Dinf \colon \sX_{\infty} \longrightarrow \sZ_{\infty},$$
by 
$$\Dinf(a,b)=(\Lu a,\Lv b).$$

\vskip.5cm
\noindent{\bf{Proof of Proposition~\ref{prop:UniformParamatrix}.}}
The restriction of $\rho_{2}$ to the 
kernel of $D_{\widetilde v}$ in $\sW_{2}$ is an isomorphism onto $\Harm$  
(according to Lemma~\ref{lemma:IdentDefTheorySphere}). It follows 
from this, and Lemma~\ref{lemma:IdentDefTheoryStrip} that 
$\Dinf\colon \sX_{\infty}\longrightarrow \sZ_{\infty}$ is 
invertible.
Let $P_{\infty}=(P_{\infty}^{(1)},P_{\infty}^{(2)})$ 
denote the inverse of $\Dinf$.
A paramatrix $Q_{T}$ for $D_{T}$ can be 
defined by
$$Q_{T}(y)=\psi_{1;T}P_{\infty}^{(1)}(\psi_{1;T}y,\psi_{2;T}y)+
\psi_{2;T}P_{\infty}^{(2)}(\psi_{1;T}y,\psi_{2;T}y),$$
where $\psi_{1;T}$, $\psi_{2;T}$ are functions defined on
$X_{1}\cup_{T}X_{2}$,  which satisfy the 
property that $\psi_{1;T}^{2}+\psi_{2;T}^{2}\equiv 1$, and 
$$\sup_{X_{1}\cup_{T}X_{2}}|d\psi_{1;T}|\leq \frac{c}{T},$$
for some fixed constant $c$. Such a family can be obtained from some 
fixed $\{\psi_{1;1},\psi_{2;1}\}$ supported in the tube by
rescaling the tube by $T$.

The operator norms of the $Q_{T}\colon \sZ(T)\longrightarrow \sX(T)$ 
is uniformly bounded independent of $T$. Indeed,
\begin{eqnarray*}
    \|Q_{T} y\|_{\sX(T)}
    &\leq & \|P_{\infty}(\psi_{1} y, \psi_{2} y)\|_{\sX_{\infty}} \\
    &\leq &\|P_{\infty}\| \|(\psi_{1}y,\psi_{2}y)\|_{\sZ_{\infty}} \\
    & \leq & \|P_{\infty}\|\|y\|_{\sZ(T)},
\end{eqnarray*}
where $\|P_{\infty}\|$ denotes the operator norm of 
$$P_{\infty}\colon \sZ_{\infty}\longrightarrow \sX_{\infty}.$$

Moreover, $Q_{T}$ is nearly an inverse for $D_{T}$, as
\begin{eqnarray*}
    D_{T}\circ Q_{T}(y) &= &
    D_{T}\psi_{1;T}P_{\infty}^{(1)}(\psi_{1;T}y,\psi_{2;T}y)+
    \psi_{2;T}P_{\infty}^{(2)}(\psi_{1;T}y,\psi_{2;T}y) \\ &=&
    (d\psi_{1;T})P_{\infty}^{(1)}(\psi_{1;T}y,\psi_{2;T}y)+
    (d\psi_{2;T})P_{\infty}^{(2)}(\psi_{1;T}y,\psi_{2;T}y) \\ && + \psi_{1;T}\Lu
    P_{\infty}^{(1)}(\psi_{1;T}y,\psi_{2;T}y)+ \psi_{2;T}\Lv
    P_{\infty}^{(2)}(\psi_{1;T}y,\psi_{2;T}y) \\ &=&
    (d\psi_{1;T})P_{\infty}^{(1)}(\psi_{1;T}y,\psi_{2;T}y)+
    (d\psi_{2;T})P_{\infty}^{(2)}(\psi_{1;T}y,\psi_{2;T}y) \\ && +
    \psi_{1;T}\psi_{1;T}y+ \psi_{2;T}\psi_{2;T}y;
\end{eqnarray*}
i.e. the operator $D_{T}\circ Q_{T}$ differs from the identity map 
(on $\sZ(T)$) by the operator 
$$y\mapsto     (d\psi_{1;T})P_{\infty}^{(1)}(\psi_{1;T}y,\psi_{2;T}y)+
    (d\psi_{2;T})P_{\infty}^{(2)}(\psi_{1;T}y,\psi_{2;T}y).$$ 
But
\begin{eqnarray*}
    \|(d\psi_{1;T}P_{\infty}^{(1)}(\psi_{1;T}y,\psi_{2;T}y)\|_{\sX(T)}
    & \leq & 
    \frac{c}{T}\|P_{\infty}^{(1)}(\psi_{1;T}y,\psi_{2;T}y)\|_{\sW_{1}} \\
    & \leq &
    \frac{c}{T}\|P_{\infty}(\psi_{1;T}y,\psi_{2;T}y)\|_{\sX_{\infty}} \\
    &\leq &
    \frac{c}{T}\|P_{\infty}\| \|(\psi_{1;T}y,\psi_{2;T}y)\|_{\sZ_{\infty}} \\
    &\leq & 
    \frac{c}{T}\|P_{\infty}^{(1)}\| \|y\|_{\sZ(T)},
\end{eqnarray*}
where $\|P_{\infty}\|$ denotes the operator norm of 
$$P_{\infty}\colon \sZ_{\infty}\longrightarrow \sX_{\infty}.$$
A similar estimate holds for 
$\|(d\psi_{2;T}P_{\infty}^{(2)}(\psi_{1;T}y,\psi_{2;T}y)\|_{\sX(T)}$.
Thus, if $T$ is sufficiently large relative to $\|P_{\infty}\|$, then 
we can arrange for, say $$\|D_{T}\circ Q_{T}-\Id\|\leq \OneHalf,$$
so that the von Neumann expansion gives the inverse of $D_{T}\circ 
Q_{T}$, with $\|D_{T}\circ Q_{T}\|\leq 2$. 
Letting $P_{T} = Q_{T}\circ (D_{T}\circ Q_{T})^{-1}$, it follows that 
$P_{T}$ is a right inverse for $D_{T}$, with norm bounded by $2 
\|P_{\infty}\|$ (again, for sufficiently large $T$).
\qed

\subsection{Proof of the gluing result}

Given the approximately holomorphic disks ${\widehat \gamma}={\widehat 
\gamma}(u)$
constructed in Subsection~\ref{subsec:Splicing} and the deformation
theory from Subsection~\ref{subsec:FredTheory}, the construction of
the holomorphic disks in $\Sym^{g+1}(\Sigma')$, as stated in 
Theorem~\ref{thm:Gluing},
follows an application
of the inverse function theorem.

We give a suitable modification of Floer's set-up 
(see~\cite{FloerUnregularized}).
There is an ``exponential mapping'' 
$$\Exp_{\widehat \gamma}\colon {\mathcal X}(T) \longrightarrow \Map(X_{1}\cup_{T} 
X_{2}, \Sym^{g+1}(\Sigma'))\times \Diff(\Strip),$$
which is obtained as follows.
Let $(a,\xi)\in {\mathcal X}(T)$; i.e. 
$a$ is a section of ${\widehat \gamma}_{T}^{*}(T\Sym^{g+1}(\Sigma'))$, 
and $\xi\in \C^{n}\subset 
\Sections(X_{1}\cup_{T}X_{2},T(X_{1}\cup_{T}X_{2})$. Fix a path 
$\met_{s}$ of metrics for which $\Tb'$ is $\met_{0}$-totally geodesic, and 
$\Ta'$ is $\met_{1}$-totally geodesic; this family of metrics is also 
required to by cylindrical in the connected sum region (for the 
appropriate uniformity in $T$). 
In fact, we find it convenient to use an $s$-independent
product K\"ahler metric on the
region $\Sym^{g+1}(\Sigma-B_{R_1})\times 
\Sym^1(E-B_{R_2})\subset \Sym^{g+1}(\Sigma')$. 
Then, we define
$$\Exp(a,\xi)(s+it) = \exp_{{\widehat 
\gamma}_{T}(s+it)}^{\met_{s}}(a(s+it))\times \exp_{s+it}(\xi),$$
where in the first factor
$$\exp^{\met}_{x}\colon T\Sym^{g+1}(\Sigma')\longrightarrow 
\Sym^{g+1}(\Sigma')$$
denotes the exponential map for the metric $\met$ over 
$\Sym^{g+1}(\Sigma')$ at the point $x\in\Sym^{g+1}(\Sigma')$; while 
the second denotes simply the exponential map for $\Strip$. The 
range of this exponential map consists of holomorphic disks with the 
appropriate boundary conditions.

We can think of the $\DBar$-operator extended to
$\Map(\Strip,\Sym^{g+1}(\Sigma'))\times \Diff(\Strip)$, which sends the 
pair $(u,\phi)$ to the section $\Sections(\Strip, \End(T\Strip, 
u^{*}(T\Sym^{g+1}(\Sigma'))))$, defined by
$$
        \DBar^{\ext}_{J_s'}(u,\xi)=(u^{*}
        J'_{\phi_\xi^{-1}})\circ du - du\circ j_{\xi},
$$ 
Clearly, a 
section $(u,\phi)$ lies in the kernel of this operator iff $u\circ 
\phi$ is a $J_s'$-holomorphic map from $\Strip$ to $\Sym^{g+1}(\Sigma')$.

As in~\cite{FloerUnregularized}, this section has a second-order 
expansion 
$$\DBar \Exp_{\widehat \gamma}(a,\xi) = \DBar {\widehat\gamma} + 
D_{\widehat\gamma}(a,\xi) + N_{\widehat \gamma_{T}}(a,\xi), $$
for some operator
$$N_{\widehat\gamma_{T}}\colon \sX(T)\longrightarrow \sZ(T).$$

\begin{lemma}
    \label{lemma:QuadraticTerm}
    There are constants $\epsilon$ and $c$  depending on $u$ with 
    the property that for all $T$ large enough for ${\widehat\gamma}_{T}$ 
    to be defined, if $\|\alpha_{1}\|_{\sX(T)}, 
    \|\alpha_{2}\|_{\sX(T)} < \epsilon$, then
    $$
         \|N_{{\widehat \gamma}_{T}}(\alpha_{1}) -
         N_{{\widehat \gamma}_{T}}(\alpha_{2})\|_{\sZ(T)}
         \leq c 
         \|\alpha_{1}-\alpha_{2}\|_{\sX(T)}
         \left(\|\alpha_{1}\|_{\sX(T)}+\|\alpha_{2}\|_{\sX(T)}\right).
    $$
\end{lemma}

\begin{proof}
        We are free to use, metrics near
        $\{w_i\}\times \{c\}\in 
        \Sym^{g}(\Sigma-B_{r_1})\times\Sym^{1}(E-B_{r_2})\subset
        \Sym^{g+1}(\Sigma')$
        which are Euclidean. In these neighborhoods, then, the
        non-linear part of the second-order expansion vanishes. But these
        are the regions where the map depends on the neck-length $T$. 
\end{proof}

We can now apply    
Newton's iteration scheme (see 
for instance~\cite{DonKron}, or~\cite{Liu} in a more closely related 
context) to find the unique holomorphic map
$\gamma=\gamma(u)$ in a 
sufficiently small neighborhood of ${\widehat \gamma}$. 

\begin{prop}
    \label{prop:TaubesGlue}
    There is an $\epsilon>0$ with the property that
    for all sufficiently large $T>0$, there is a
    unique holomorphic curve $\gamma(u)$ which lies in an 
    $\epsilon$-neighborhood of ${\widehat\gamma}(u)$, measured in 
    the ${\mathcal X}(T)$ norm (for $T$ sufficiently large relative 
    to $S$). 
\end{prop}

\begin{proof}
Given ${\widehat \gamma}$, we find, for all sufficiently large $T$,
a holomorphic 
$$\gamma={\widehat \gamma} + P_{T}\eta,$$
by finding
$\eta\in{\mathcal Z}(T)$ with the property that
$$\eta + N_{\widehat \gamma}(P_{T}\eta)=\DBar_{J_s'(T)}{\widehat \gamma}.$$
This can be done, since the map
$$\eta\mapsto N_{\widehat \gamma}\circ P_{T}\eta $$
is uniformly quadratically  contracting (this is a consequence of 
Lemmas~\ref{prop:UniformParamatrix} and \ref{lemma:QuadraticTerm}), and 
the error term $\DBar_{J'_s(T)}{\widehat \gamma}$ goes to zero in ${\mathcal 
Z}(T)$, according to Lemma~\ref{lemma:ErrorTerm}. 
\end{proof}

In fact, from the construction of the right inverses, it is clear that
the kernel of $D_{\widehat \gamma}$ is identified with the kernel of
the original $D_u$ (see Lemma~\ref{lemma:IdentDefTheoryStrip} and
Proposition~\ref{prop:UniformParamatrix}). Moreover, since $D_\gamma$
is a small perturbation of $D_{\widehat\gamma}$, their kernels are
identified as well. Hence the deformation theory of $u$ is identified with
the deformation theory of $\gamma(u)$ (this was the claim made in 
Remark~\ref{rmk:IdentDefTheory}). This can also be used to 
identify the corresponding determinant line bundles, so that the
signs appearing in the signed counts agree. 

We wish to show that for sufficiently large $T>0$, all the holomorphic 
curves in $\ModFlow(\x',\y')$ are contained in the domain of the
gluing map constructed in Proposition~\ref{prop:TaubesGlue}. 
This can be seen after an application of
the Gromov compactness theorem. To apply Gromov's theorem, it is 
important to set up a uniform version of the energy bounds from 
Subsection~\ref{subsec:EnergyBounds} which hold as $T\goesto \infty$. 
To this end, we find it convenient to think of the degeneration (as 
$T\goesto \infty$) as the formation of a node. Specifically, we 
consider a holomorphic one-parameter family of holomorphic curves 
parameterized by the complex disk
$$f \colon X \longrightarrow \CDisk,$$ 
whose fiber is singular
at $t=0$, and which contains the complex structures $\sj'(T)$ 
for all large $T$ 
(i.e. $X$ is an algebraic variety, together with a map to the disk, whose 
fiber $X_{t}$ at over a point $t\in\CDisk$ is 
a curve of genus $g+1$ whenever $t\neq 0$; 
and the limit $t\goesto 0$ corresponds to $T\goesto \infty$).
Forming the fiberwise 
symmetric product, we obtain an algebraic variety
$$f \colon S(X)\longrightarrow \CDisk, $$
whose fiber over $t$ is identified as $S(X)_{t}\cong
\Sym^{g+1}(X_{t})$ (see also~\cite{DonaldsonSmith} for a 
symplectic construction
of this object).

\begin{prop}
\label{prop:GromovCompactness}
    Fix a class $\phi\in\pi_{2}(\x',\y')$. Then any sequence of
    holomorphic disks
    $u_{t}\in \ModFlow_{J_{t}}(\phi)$ with $t\goesto 0$
    has a Gromov limit 
    which is a holomorphic disk in $\Sym^{g+1}(\SigOne\vee E)$.
\end{prop}
    
\begin{proof}
We can work in an ambient manifold by embedding the algebraic variety
$S(X)$ into a complex projective space, from which it inherits the
K\"ahler form $\omega$. As in Subsection~\ref{subsec:EnergyBounds}, we
must obtain an {\em a priori} bound on the $\omega$-energy of any disk
$u\in\ModFlow(\phi)$ (thought of as a disk in $S(X)$).  We imitate the
discussion from Subsection~\ref{subsec:EnergyBounds}, only now in
families.
    
The product form $\omega_{0}$ gives a K\"ahler form on the product
$X^{\times(g+1)}$.  Restricting this to the subvariety $P(X)\subset
X^{\times(g+1)}$, defined by $$P(X)=\{(x_{1},\ldots,x_{g+1})\in
X^{\times(g+1)}\big| f(x_{1})=\ldots=f(x_{g+1})\},$$ we obtain a
quadratic form $\met_{0}(v)=\omega_{0}(v,Jv)$ on the tangent cone of
$P(X)$. There is a holomorphic map $\pi\colon P(X)\longrightarrow
S(X)$, under which we can pull back $\omega$ to the subvariety, to
obtain another quadratic form $g(v)=\omega(v,Jv)$ on the tangent
cone. Note that $\met_{0}$ is nowhere vanishing. Thus, we can form the
ratio $\met/\met_{0}$, to obtain a continuous function on the
projectivized tangent cone of $P(X)$. By compactness, this function is
bounded above by some constant which we will denote by $C_3$ (since it
is the constant $C_3$ appearing in
Inequality~\eqref{ineq:BoundEnergyDiag}). It follows that if we have a
$\Sym^{g+1}(\sj'(T))$-holomorphic map $u\colon \Omega\longrightarrow
P(X)$, then $$ \int_{\Omega}u^*(\omega)\leq
C_3\int_{\widehat\Omega}{\widehat u}^{*}(\omega_{0}).$$ But we must
consider perturbations of these.

To control the energy integrand in the region where the complex
structures are varying, consider first
the (holomorphic) commutative diagram
$$
\begin{CD}
\Big(\Sigma-B_{R_1}(\sigma_1)\Big)^{\times g}\times
(E-B_{R_2}(\sigma_2))\times\CDisk @>{\widetilde \iota}>>P(X) \\
@V{\pi}VV  @V{\pi'}VV \\
\Sym^{g}(\Sigma-B_{R_1}(\sigma_1))\times (E-B_{R_2}(\sigma_2)) \times\CDisk
@>\iota>> S(X),
\end{CD},
$$
where the vertical maps are induced by the quotients, and
the maps all commute with projection to the base $\CDisk$.
Letting $V$ be some open subset of the diagonal in
$\Sym^{g}(\Sigma-B_{R_1}(\sigma_1))$, and ${\widetilde V}$ be the pre-image of 
its closure in $(\Sigma-B_{R_1}(\sigma_1))^{\times g}$, then the restriction 
of $\pi$ to ${\widetilde V}\times\CDisk$ is a covering space. Thus,
we can consider the differential
form $\pi_*{\widetilde \iota}^*(\omega_0)$ over 
$V \times
(E-B_{R_2}(\sigma_2))\times\CDisk\subset\Sym^{g}(\Sigma-B_{R_1}(\sigma_1))\times
(E-B_{R_2}(\sigma_2))\times\CDisk$. Since the maps are all
holomorphic local diffeomorphisms, 
observe that the complex structure $\Sym^{g}(\sj)\times
\sj_E \times\sj_\CDisk$ tames the form $\pi_*{\widetilde \iota}^*(\omega_0)$.
It follows that if $J_s'$ is a sufficiently small perturbation of 
the constant path $\Sym^{g}(\sj)\times \sj_E$ (which is the
restriction of $\Sym^{g+1}(\sj)$), then the taming condition is
preserved. With respect to such a path $J_s'$, then, the energy of
a $J_s'$-holomorphic 
map $u$ from $\Strip$ to any given $t$-fiber of $P(X)$ can be
controlled as in Lemma~\ref{lemma:EnergyBound} by the
$\omega_0$-integral of its associated branched cover ${\widehat u}$. 
On the other hand, for a given $\phi$,
this integral is controlled by some multiplicity (depending on the
maximal multiplicity in $\cald(\phi)$, which is, of course, a 
topological quantity) 
of the the $\omega_0$-area a fiber of $X$. But the latter area is bounded
in the family, since it is obtained from a smooth symplectic form
which extends over the family $X$. 

This gives us the uniform energy bound on the $u_{t}$ required by
Gromov's compactness theorem, and hence gives rise to a limiting
cusp-curve which maps into $S(X)$. Indeed, by continuity of the
limiting process, it follows that the cusp-curve maps to the $t=0$
fiber, $\Sym^{g+1}(\SigOne\vee E)$.
\end{proof}
    
With respect to the cylindrical geometries on the domain and the 
range, this gives us the following:

\begin{prop}
    \label{prop:Compactness}
        Fix $\x,\y\in\Ta\cap\Tb$, and let $\x',\y'$ be their
        stabilizations. Fix a homotopy class $\phi'\in\pi_2(\x',\y')$
        with $\Mas(\phi')=1$. Given any sequence of holomorphic curves
        $\{u_{i}\}_{m=1}^{\infty}$ with 
        $u_{m}\in\Mod_{J'_s(T_{m})}(\phi')$ for some unbounded sequence 
        of real numbers $\{T_{m}\}$, with $\Mas(u_{m})=1$. Then, for each 
        $\epsilon>0$, after passing to a suitable 
        subsequence of the $\{u_{m}\}$ if necessary, we obtain 
        the following data:

    \begin{itemize}

        \item a collection of points $\{w_{1},\ldots,w_{n}\}\in 
        \Sym^{g}(\SigOne)$
        \item a pair $u$ and $v$, where 
        $u\in\ModFlow(\phi)$, and $v\in\ModCent$, 

        \item
        a decomposition of the 
        strip $\Strip=X_{1}\cup_{T'_{m}} X_{2}$
        (and, of course, a corresponding subset
        $\{q_{1},\ldots,q_{n}\}\subset \Strip$),

        \item a sequence $\{\xi_{m}\}\in \C^{n}$ converging to zero, with the 
        property that $u_{m}\circ \phi_{\xi_{m}}(q_m)\in W'$ 
        (recall that $\phi_{\xi_{m}}$ is 
        some ``translation'' of the strip $\Strip=X_{1}\cup_{T'_{m}} X_{2}$ 
        which carries the center $q_{i}$ of the $i^{th}$ component of $X_{2}$ 
        to some nearby point).
        
    \end{itemize}
        
    which satisfy the following properties:
    
    \begin{list}
        {(\arabic{bean})}{\usecounter{bean}\setlength{\rightmargin}{\leftmargin}}
        \item
        \label{prop:MapFirst}
        each $u_{m}\circ \phi_{\xi_{m}}$ maps $X_{1}$ into 
        $$\Sym^{g}(\SigOne-B_{r_{1}}(\sigma_{1}))\times 
        \Sym^{1}(E-B_{r_{2}}(\sigma_2))\subset \Sym^{g+1}(\Sigma')$$

        \item
        \label{prop:MapSecond}
        each $u_{m}\circ \phi_{\xi_m}$ maps the the $i^{th}$ component of
        $X_{2}$ into 
        $$B_{\epsilon}(w_{i})\times B_{\epsilon}(c)
        \subset \Sym^{g-1}(\SigOne-B_{r_{1}}(\sigma_1))\times 
        \Sym^{1}(E-B_{r_{2}}(\sigma_2))\subset \Sym^{g+1}(\Sigma').$$
        
        \item
        \label{prop:CInftyConv}
        for any $S>0$,
        the restrictions $\{u_{m}\circ \phi_{\xi_{m}}|X_{1}(S)\}_{m=1}^{\infty}$ and 
        $\{u_{m}\circ \phi_{\xi_{m}}|X_{2}(S)\}_{m=1}^{\infty}$ converge in the 
        $\Cinfty$ topology to the restrictions $u\times\{c\}|X_{1}(S)$ and 
        $v|X_{2}(S)$ 
        respectively, where we think of $u\times\{c\}$ and $v$ as maps
        $$u\times\{c\}\colon X_{1}(\infty)\longrightarrow 
        \Sym^{g}(\SigOne-\sigma_{1})\times \{c\}$$
        and 
        $$v\colon X_{2}(\infty)\longrightarrow 
        \Sym^{g-1}(\SigOne-B_{r_{1}}(\sigma_{1}))\times
        \Sym^{2}(E-\sigma_{2})\subset \Sym^{g+1}(\Sigma')$$
        respectively.
        
        \item 
        \label{prop:BubblePoint}
        under $u_{m}\circ \phi_{\xi_{m}}$, each cylinder $[-T'_{m},T'_{m}]\times S^{1}_{i}\subset 
        X_{1}\cup_{T'_{m}}X_{2}$ maps to the product
        $$B_{\epsilon}(w_{i})\times \left([-T_{m},T_{m}]\times S^{1}\right)
        \times B_{\epsilon}(c)
        \subset \Sym^{g+1}(\Sigma'(T_{m})).$$
        
    \end{list}
\end{prop}

\begin{proof}
    Gromov's compactness theorem adapted to our situation
    (Proposition~\ref{prop:GromovCompactness}) gives a map
    $$u_{\infty}\colon {\mathbb B}
    \longrightarrow \Sym^{g+1}(\SigOne\vee E),$$ 
    which is the limit of
    the $\{u_{m}\}$, where ${\mathbb B}$ here is a bubbletree
    obtained by attaching $n$ spheres ${\mathbb S}_i$ to the disk $\CDisk$
   . For generic choice of $\sigma_1$,
    its image lies in the union
    of $\Sym^{g}(\SigOne)\times\Sym^{1}(E)$ and
    $\Sym^{g-1}(\SigOne)\times \Sym^{2}(E)$. Since $\Ta'\cup\Tb'$ is
    contained in the first set, it follows that the main component of
    the bubble-tree is mapped to the first set.  Indeed, this gives
    rise to the $u\times\{c\}$. Moreover, by dimension counts, it
    follows that the spheres are mapped into the second subset, giving
    rise to the map $v$. (Note also that the limiting almost-complex
    structures $J_s'$ all agree with the constant complex structure,
    over this latter subset.)
    
    To view the domains of all the $\{u_{m}\}$ as
    $X_{1}\cup_{T'_{m}}X_{2}$ for fixed $X_{1}$, we use the
    diffeomorphisms generated by vector fields $\xi_{i}$
    ($i=1,\ldots,n$). Specifically, in the proof of the compactness
    theorem (see, for instance~\cite{ParkerWolfson}) one rescales
    around a point where energy is accumulating, to construct the
    bubbles. These points where the energy concentrates form
    subsequences of $n$-tuples
    $\q_{m}=\{q_{1}^{(m)},\ldots,q_{n}^{(m)}\}$, which converge to
    $\q$.  Now, after rescaling, we obtain a $\CInftyLoc$
    convergent subsequence to spheres 
    whose energy-measures are centered at the origin.  This
    means that there is a sequence of points
    $\q'_{m}=\{{q'}_{1}^{(m)},\ldots {q'}_{n}^{(m)}\}$ also converging
    to $\q$, with the property that $u_{m}(q'_{i})\in W'$. We let
    $\xi_{m}$ be the corresponding translation taking $\q$ to
    ${\q'}_{m}$, so that $\phi_{\xi_{m}}(\q)={\q'}_{m}$. Thus, the
    pairs $(u_{m}\circ \phi_{\xi_{m}}, \phi_{\xi_{m}}^{-1})$ satisfy
    the extended holomorphic condition on
    $\Map(X_{1}\cup_{T_{m}}X_{2},\Sym^{g+1}(\Sigma'))\times
    \Diff(\CDisk)$ considered earlier. Moreover, the
    $\CInftyLoc$ convergence of the rescaled $u_m$ ensures
    $\CInftyLoc$ convergence of the restrictions of
    $u_m\circ\phi_{\xi_m}|X_2(S)$. The limit is a 
    holomorphic map which is centered
    (in the sense of Definition~\ref{def:Centered}).

    The usual $\Cinfty_{\loc}$ on the main 
    component of the bubble-tree
    gives us Property~(\ref{prop:MapFirst}) and the 
    corresponding part of Property~(\ref{prop:CInftyConv}), while the 
    $\Cinfty_{\loc}$ convergence on the spheres gives us 
    Property~(\ref{prop:MapSecond}), and the rest of 
    Property~(\ref{prop:CInftyConv}). 
    
    The $C^{0}$ convergence at 
    the bubble points is equivalent to 
    Property~(\ref{prop:BubblePoint}), in light of the following.
    Given $w\in\Sym^{g-1}(\SigOne-B_{r_{1}}(\sigma_1))$, 
    $\sigma\in\Sym^{1}(E-B_{r_{2}}(\sigma_2))$, and letting $s$ be the 
    wedge-point in $\SigOne\vee E$, 
    a neighborhood of 
    $w\times s \times \sigma\in\Sym^{g+1}(\SigOne\vee E)$ 
    meets the nearby fiber $\Sigma'(t)$ of the one-parameter family 
    in a subset identified with the product
    \begin{eqnarray*}
        \lefteqn{B_{\epsilon}(w) \times \left([-T,T]\times S^{1}\right) \times
    B_{\epsilon}(\sigma_{2})\subset} && \\
    &&  
    \Sym^{g-1}(\SigOne-B_{r_{1}}(\sigma_1))\times \left([-T,T]\times 
    S^{1}\right) \times \Sym^{1}(E-B_{r_{2}}(\sigma_2)
    \subset 
    \Sym^{g+1}(\Sigma').
    \end{eqnarray*}
\end{proof}

We would like to prove that the properties arising from the 
compactness above, specifically the behaviour of the maps on the 
cylinders $u_{m}|[-T'_{m},T'_{m}]\times S^{1}\subset 
X_{1}\cup_{T'_{m}}X_{2}$ above,
imply stronger decay conditions. According to the proposition, 
we can view these as maps (which are now holomorphic in the usual
sense) 
to the product space $B_{\epsilon}(w_{i})\times 
\left([-T_{m},T_{m}]\times S^{1}\right)\times B_{\epsilon}(c)$. There 
is an a priori bound on the energy of these maps, where the energy is 
calculated with respect to the symplectic form on the target with 
collapsing connected sum region. But we wish to use the cylindrical 
geometry on the target instead, to stay in the geometric framework of 
Section~\ref{sec:Gluing}. The $\Sobol{p}{1,\delta}$ convergence we are 
aiming for can be divided into two parts, then, corresponding to the two 
natural projections on the target, which we denote
$$
      \Pi_{2}\colon \Sym^{g-1}(\SigOne-B_{r_{1}}(\sigma_1))\times \left([-T,T]\times 
      S^{1}\right) \times \Sym^{1}(E-B_{r_{2}}(\sigma_2))
      \longrightarrow \R\times S^{1}
$$
and 
\begin{eqnarray*}
    \lefteqn{\Pi_{1,3}\colon \Sym^{g-1}(\SigOne-B_{r_{1}}(\sigma_1))\times \left([-T,T]\times 
      S^{1}\right) \times \Sym^{1}(E-B_{r_{2}}(\sigma_2))} && \\
      &&
      \longrightarrow 
      \Sym^{g-1}(\SigOne-B_{r_{1}}(\sigma_1))\times 
      \Sym^{1}(E-B_{r_{2}}(\sigma_2)).
\end{eqnarray*}
Our aim is to show that the projections
$\Pi_{1,3}\circ u_{m}|[-T'_{m},T'_{m}]\times S^1$ and $\Pi_{2}\circ 
u_{m}|[-T'_{m},T'_{m}]\times S^1$ all lie in a uniform $\Sobol{p}{1,\delta}$ ball.
Note that, thanks to the $\Cinfty$ convergence on the compact pieces
(Property~(\ref{prop:CInftyConv})), we 
have that the restrictions $\{u_{m}|_{\{-T'_{m}\}\times S^{1}}\}$ 
and $\{u_{m}|_{\{T'_{m}\}\times S^{1}}\}$ lie in $\Cinfty$-bounded 
subsets, up to translations on the cylindrical factor on the range.
In handling the $\Pi_{1,3}$ projection, we do not use this 
translation, and we can appeal to the following 
elementary result:

\begin{lemma}
    \label{lemma:NeckSqueezeOne}
    Let $T_{i}$ be an increasing sequence of real numbers, and 
    suppose that $u_{i}\colon [-T_{i},T_{i}]\times S^{1}
    \longrightarrow \C$ is a sequence of holomorphic maps of the cylinders 
    into the complex plane, with the property that
    $u_{i}|_{([-T_{i},-T_{i}+1] \cup [T_{i}-1,T_{i}])\times 
    S^{1}}$ is uniformly $\Cinfty$ bounded, then there is 
    a subsequence $\{w_{m}\}_{m=1}^{\infty}\in\CDisk$ 
    converging to some $w\in\CDisk$, 
    with the property that 
    $u_{m}-w_{m}$ is uniformly 
    $\Sobol{p}{1,\delta}$-bounded.
\end{lemma}

\begin{proof}
    The point here is that the 
    norm of a holomorphic function $u$ on a cylinder is controlled by 
    its behaviour near the boundary, provided that
    $\int_{\{T\}\times S^{1}}u=0$. More precisely, 
    for each integer $k\geq 0$, we have constants $c_{k}$ with the property that if 
    $$u\colon [-T,T]\times S^{1}\longrightarrow \C$$ is a holomorphic 
    map with $\int_{\{-T\}\times S^{1}}u=0$, then we have for each sufficiently small $\delta>0$,
    $$\|u\|_{\Sobol{2}{k,\delta}}\leq 
    c_{k}\left(\|u\|_{\Sobol{2}{k,\delta}([-T,-T+1]\times S^{1})}
    +\|u\|_{\Sobol{2}{k,\delta}([T-1,T]\times S^{1})}\right).$$
    This follows from looking at the 
    Fourier coefficients of $u$. We content ourselves here with the 
    case where $k=0$.
     
    Write $u(t,\theta)=\sum_{n\neq 0}a_{n}e^{n(t+i\theta)}$. Then, we 
    have
    \begin{eqnarray*}
        \|u\|_{\Sobol{2}{\delta}}^{2}
        &=&
        \sum_{n>0}|a_{n}|^{2}\int_{-T}^{T}e^{(2n+\delta)t}
        + 
        \sum_{n<0}|a_{n}|^{2}\int_{-T}^{T}e^{(2n+\delta)t} \\
        &=&
        \sum_{n>0}|a_{n}|^{2}\int_{-T}^{T}e^{(2n+\delta)t}
        + 
        \sum_{n<0}|a_{n}|^{2}\int_{-T}^{T}e^{(2n+\delta)t} \\
        &\leq & 
        C\left(\sum_{n>0}|a_{n}|^{2}\int_{T-1}^{T}e^{(2n+\delta)t}
        + 
        \sum_{n<0}|a_{n}|^{2}\int_{-T}^{-T+1}e^{(2n+\delta)t}\right) \\
        &=&
        C
        \left(\|u\|_{\Sobol{2}{\delta}([-T,-T+1]\times S^{1})}^{2}
        +\|u\|_{\Sobol{2}{\delta}([T-1,T]\times S^{1})}^{2}\right).
    \end{eqnarray*}
    The estimate for larger $k$ follows in a similar manner. The 
    corresponding estimate for $\Sobol{p}{1,\delta}$ follows from the 
    uniform, $T$-independent inclusion (provided $T$ is bounded below)
    $$\Sobol{2}{2,\delta}([-T,T]\times S^{1}) \longrightarrow 
    \Sobol{p}{1,\delta'}([-T,T]\times S^{1}),$$
    where $\delta'=\delta p/2$, which in turn follows from the fact that the 
    norm of the inclusion $\Sobol{2}{2}\longrightarrow \Sobol{p}{1}$ 
    is uniformly bounded, provided $T$ is bounded below.

    To apply this, 
    let $w_{i}=\int_{\{-T_{i}\}\times S^{1}}u_{i}$. Clearly, the 
    $w_{i}$ have a convergent subsequence. 
    The above argument gives a uniform $\Sobol{p}{k,\delta}$-bound on
    $u_{i}-w_{i}$.
\end{proof}

To handle the $\Pi_{2}$ projection, we must study a similar problem, 
only where now the target is also a cylinder (rather than a disk). 
To set up notation, let 
$$j_{1}\colon [0,1]\times S^{1}\longrightarrow \R\times S^{1}$$
denote the standard inclusion of an annulus into the 
cylinder. 

\begin{lemma}
    \label{lemma:NeckSqueezeTwo}
    There is an $\epsilon>0$ with the property that
    if $u_{m}\colon [-T_{m},T_{m}]\times S^{1}\longrightarrow 
    \R\times S^{1}$ is a sequence of maps with the property that, up to 
    translations and rotations on the cylinder, both sequences
    $\{u_{m}|[-T_{m},-T_{m}+1]\times S^{1}\}$ and 
    $\{u_{m}|[T_{m}-1,T_{m}]\times S^{1}\}$ line in an $\epsilon$
    neighborhood of $j_1$. Then, up to 
    translations and rotations, the $u_{m}$ lie in a uniform
    $\Sobol{2}{\delta}$-neighborhood of the standard inclusion 
    $[-T_{m},T_{m}]\times S^{1}\subset \R\times S^{1}$.
\end{lemma}
    
\begin{proof}
    Let 
    $$j_{T_{m}}\colon [-T_{m},T_{m}]\times S^{1}\longrightarrow \R\times 
    S^{1}$$
    denote the natural inclusion. 
    Consider next the difference $j_{T_m}-u_{m}$. If the $C^{0}$ norm of 
    the difference between $u_{m}|\{T_{m}\}\times S^{1}$ and the 
    standard inclusion of the circle is sufficiently small (up to 
    translation), then the two maps are homotopic. Hence, we can lift 
    the difference $u_{m}-j_{T_{m}}$ to a sequence of holomorphic 
    maps 
    $$u_{m}-j_{T_{m}}\colon [-T_{m},T_{m}]\times S^1\longrightarrow \C,$$
    where $\C=\R\times\R$ is thought of as the universal covering 
    space of the annulus. We then apply the argument from 
    Lemma~\ref{lemma:NeckSqueezeOne}.
\end{proof}
    
\vskip.3cm
\noindent
{\bf Proof of Theorem~\ref{thm:Gluing}.}
    Given $\epsilon>0$, we can find a $T_{0}$ so that for all 
    $T>T_{0}$, all holomorphic curves $\gamma\in\ModFlow(\phi')$ lie 
    in the image of the map constructed in 
    Proposition~\ref{prop:TaubesGlue}. This follows from 
    Lemma~\ref{lemma:NeckSqueezeOne} and 
    Lemma~\ref{lemma:NeckSqueezeTwo}.
\qed

\section{Conclusion: topological invariance}
\label{sec:Conclusion}

We have established all the pieces now to conclude topological
invariance of the homology groups, first stated in
Theorem~\ref{intro:Invariance}, which justifies our dropping the
Heegaard diagram from the notation for the Floer homology groups. As a
preliminary remark, recall that Theorem~\ref{thm:IndepCxStruct} shows
that the homology groups are independent of the complex structure.  We
can now give the following more precise statement of
Theorem~\ref{intro:Invariance}:

\begin{theorem}
\label{thm:Invariance}
If $(Y,\spinc)$ is a three-manifold equipped with a $\SpinC$ structure
$\spinc\in\SpinC(Y)$, then the Floer homology groups are topological
invariants of $(Y,\spinc)$ in the following sense. There is a strongly
$\spinc$-admissible Heegaard diagram, and if two different strongly
$\spinc$-admissible Heegaard diagrams
$(\Sigma_1,\alphas_1,\betas_1,z_1)$ and
$(\Sigma_2,\alphas_2,\betas_2,z_2)$ represent the same three-manifold,
then there is a one-to-one correspondence between isomorphism classes
of orientation conventions for the first and the second Heegaard
diagram and corresponding
$\Z[U]\otimes_\Z\Wedge^*H_1(Y;\Z)/\Tors$-module isomorphisms
(represented by the vertical arrows): 
\begin{equation}
\label{eq:LongExactSeq}
\begin{CD}
...@>>>\HFm(\Sigma_1,\alphas_1,\betas_1,\spinc)@>>>
\HFinf(\Sigma_1,\alphas_1,\betas_1,\spinc)@>>>
\HFp(\Sigma_1,\alphas_1,\betas_1,\spinc)@>>>... \\
&& @V{\Phi^-}VV @V{\Phi^\infty}VV @V{\Phi^+}VV  \\
...@>>>\HFm(\Sigma_2,\alphas_2,\betas_2,\spinc)@>>>
\HFinf(\Sigma_2,\alphas_2,\betas_2,\spinc)@>>>
\HFp(\Sigma_2,\alphas_2,\betas_2,\spinc)@>>>... 
\end{CD}
\end{equation}
(where the groups in the first row can be calculated using any
orientation system $\orient_1$, and the second are calculating using
the induced orientation system $\orient_2$). Indeed, if
$(\Sigma_1,\alphas_1,\betas_1,z_1)$ and
$(\Sigma_2,\alphas_2,\betas_2,z_2)$ are weakly
$\spinc$-admissible, we have isomorphisms:
$$
\begin{CD}
...@>>>\HFa(\Sigma_1,\alphas_1,\betas_1,\spinc)@>>>
\HFp(\Sigma_1,\alphas_1,\betas_1,\spinc)@>>>
\HFp(\Sigma_1,\alphas_1,\betas_1,\spinc)@>{U}>>... \\
&& @V{\widehat\Phi}VV @V{\Phi^+}VV @V{\Phi^+}VV  \\
...@>>>\HFa(\Sigma_2,\alphas_2,\betas_2,\spinc)@>>>
\HFp(\Sigma_2,\alphas_2,\betas_2,\spinc)@>{U}>>
\HFp(\Sigma_2,\alphas_2,\betas_2,\spinc)@>{U}>>... 
\end{CD}
$$
\end{theorem}

\begin{proof}
The strongly $\spinc$-admissible Heegaard diagrams exist according to
Proposition~\ref{prop:VerySpecialMoves}, which also shows that they
can be connected by a sequence of strongly $\spinc$-admissible
isotopies, handleslides, and stabilizations. Invariance under those
isotopies in the above sense was established in
Theorem~\ref{thm:Isotopies}. Invariance under handleslides was
established in Theorem~\ref{thm:HandleslideInvariance} for a special
choice of isotopy type of the handleslide. Of course, any handleslide
(which does not cross the basepoint) can be brought to this form after
an isotopy of the attaching circles (which does not cross the
basepoint), so general handleslide invariance follows. Finally,
stabilization invariance was established in
Theorem~\ref{thm:StabilizeHFb} (and Theorem~\ref{thm:StabilizeHFa} for
$\HFa$). When the Heegaard diagram is only weakly $\spinc$-admissible,
Lemma~\ref{lemma:WeaklyAdmissibleIsotopy} gives us a strongly
$\spinc$-admissible diagram which is isotopic to the given diagram,
and Theorem~\ref{thm:WeakIsotopies} gives the necessary
identifications between $\HFa$ and $\HFp$. We then reduce to to
strongly $\spinc$-admissible case considered before.
\end{proof}

Recall that $\HFred(Y,\spinct)$ can be determined from the long exact
sequence in Equation~\eqref{eq:LongExactSeq}, so it, too is a
topological invariant.

When $b_1(Y)>0$, we stress that there is still the auxilliary choice
of an (isomorphism class of) orientation system, giving rise to
$2^{b_1(Y)}$ different candidates for the ``Floer homologies.'' In
fact, in Theorem~\ref{HolDiskTwo:thm:HFinfTwist} of~\cite{HolDiskTwo},
we show how to identify a canonical orientation system.

It is not difficult to establish naturality of the maps appearing in
the above theorem; we return to this point in~\cite{HolDiskThree}.

\subsubsection{Further remarks}

There are other variants of the Floer homologies defined in the
present paper. For example, there are variants 
twisted with a $\Z[H^1(Y;\Z)]$-module $M$. We will describe this
construction in Section~\ref{HolDiskTwo:sec:TwistedCoeffs}
of~\cite{HolDiskTwo}.

There is another construction, which works even in the absence of
admissibility hypotheses, which we sketch now.
For this construction, we will work over the Novikov ring $\Ring$
consisting of formal power series $\sum_{r\geq 0} a_r e^r$, for which
the support of the $a_r$ (in $r$) is discrete, endowed with the
multiplication law: $$\left( \sum_{r\geq 0} a_r e^r \right) \cm
\left( \sum_{r\geq 0} b_r e^r \right) = \sum_{r\geq 0} \left(\sum_{s\geq 0} a_{s}
b_{r-s}\right) e^r$$
(we emphasize that the symbol $e$ here is a formal variable).

For a pointed Heegaard diagram $(\Sigma,\alphas,\betas,z)$ for $Y$, we
consider a the chain complex $\CFp_{{\mathrm Nov}}$ which is freely
generated (over $\Ring$) by pairs $[\x,i]\in \left(\Ta\cap
\Tb\right)\times \Z^{\geq 0}$, endowed with the boundary map by
$$\partial^+_{\mathrm Nov}[\x,i]=\sum_{\{\y\in {\mathcal S}\}}
\sum_{\{\phi\in\pi_2(\x,\y)\big| n_z(\phi)\leq i\}} e^{{\mathcal A}(\phi)}\left(\#\ModFlow(\phi)\right)
\cdot [\y,i-n_z(\phi)],$$ where ${\mathcal A}(\phi)$ denotes the area of
the domain $\cald(\phi)$. Observe that this construction depends on
the choice of volume form for $\Sigma$ through the induced areas of
each periodic domain -- a real valued function on $H_2(Y;\R)$. That
datum, in turn, can be thought of as a real two-dimensional cohomology
class $\eta\in H^2(Y;\R)$.

Taking homologies, we obtain homology groups
$\HFp_{\mathrm{Nov}}(Y,\spinc,\eta)$ which are invariants of the
underlying topological data (and an orientation system, which we
suppress), and which require no admissibility hypotheses to define.
We will have no further use for this construction
in~\cite{HolDiskTwo}, though it may turn out to be useful in other
applications.  In particular, this construction is analogous to the
Seiberg-Witten-Floer homology perturbed by a real two-dimensional
cohomology class.

\commentable{
\bibliographystyle{plain}
\bibliography{biblio1}

\begin{thebibliography}{10}

\bibitem{Ahlfors}
L.~V. Ahlfors.
\newblock {\em Conformal invariants: topics in geometric function theory}.
\newblock McGraw-Hill, 1973.

\bibitem{KTheory}
M.~F. Atiyah.
\newblock {\em $K$-theory}.
\newblock Advanced Book Classics. Addison-Wesley Publishing Company, 1989.
\newblock Notes by {D}. {W}. {A}nderson.

\bibitem{APSI}
M.~F. Atiyah, V.~K. Patodi, and I.~M. Singer.
\newblock Spectral asymmetry and {R}iemannian geometry, {I}.
\newblock {\em Math. Proc. Camb. Phil. Soc.}, 77:43--69, 1975.

\bibitem{Silva}
V.~de~Silva.
\newblock {\em Products in the symplectic {F}loer homology of {L}agrangian
  intersections}.
\newblock PhD thesis, Oxford University, 1999.

\bibitem{DonKron}
S.~K. Donaldson and P.~B. Kronheimer.
\newblock {\em The Geometry of Four-Manifolds}.
\newblock Oxford Mathematical Monographs. Oxford University Press, 1990.

\bibitem{DonaldsonSmith}
S.~K. Donaldson and I.~Smith.
\newblock Lefschetz pencils and the canonical class for symplectic
  {$4$}-manifolds.
\newblock math.SG/0012067, 2000.

\bibitem{FloerLag}
A.~Floer.
\newblock Morse theory for {L}agrangian intersections.
\newblock {\em J. Differential Geometry}, 28:513--547, 1988.

\bibitem{FloerMaslov}
A.~Floer.
\newblock A relative {M}orse index for the symplectic action.
\newblock {\em Comm. Pure Appl. Math.}, 41(4):393--407, 1988.

\bibitem{FloerUnregularized}
A.~Floer.
\newblock The unregularized gradient flow of the symplectic action.
\newblock {\em Comm. Pure Appl. Math.}, 41(6):775--813, 1988.

\bibitem{FloerHofer}
A.~Floer and H.~Hofer.
\newblock Coherent orientations for periodic orbit problems in symplectic
  geometry.
\newblock {\em Math. Z.}, 212(1):13--38, 1993.

\bibitem{FloerHoferSalamon}
A.~Floer, H.~Hofer, and D.~Salamon.
\newblock Transversality in elliptic {M}orse theory for the symplectic action.
\newblock {\em Duke Math. J}, 80(1):251--29, 1995.

\bibitem{FOOO}
K.~Fukaya, Y-G. Oh, K.~Ono, and H.~Ohta.
\newblock {\em Lagrangian intersection Floer theory---anomaly and obstruction}.
\newblock Kyoto University, 2000.

\bibitem{GompfStipsicz}
R.~E. Gompf and A.~I. Stipsicz.
\newblock {\em {$4$}-manifolds and Kirby calculus}, volume~20 of {\em Graduate
  Studies in Mathematics}.
\newblock American Mathematical Society, 1999.

\bibitem{Gromov}
M.~Gromov.
\newblock Pseudo holomorphic curves in symplectic manifolds.
\newblock {\em Invent. Math.}, 82:307--347, 1985.

\bibitem{Harris}
J.~Harris.
\newblock {\em Algebraic Geometry: a first course}.
\newblock Number 133 in Graduate Texts in Mathematics. Springer-Verlag, 1995.

\bibitem{Hartshorne}
R.~Hartshorne.
\newblock {\em Algebraic Geometry}.
\newblock Number~52 in Graduate Texts in Mathematics. Springer-Verlag, 1977.

\bibitem{Hattori}
A.~Hattori.
\newblock Topology of {$\C^n$} minus a finite number of affine hyperplanes in
  general position.
\newblock {\em J. Fac. Sci. Univ. Tokyo}, 22(2):205--219, 1975.

\bibitem{IonelParker}
E-N. Ionel and T.~H. Parker.
\newblock Relative {G}romov-{W}itten invariants.
\newblock math.SG/9907155.

\bibitem{KMcontact}
P.~B. Kronheimer and T.~S. Mrowka.
\newblock Monopoles and contact structures.
\newblock {\em Invent. Math.}, 130(2):209--255, 1997.

\bibitem{LiRuan}
A-M. Li and Y.~Ruan.
\newblock Symplectic surgery and {G}romov-{W}itten invariants of {C}alabi-{Y}au
  3-folds.
\newblock {\em Invent. Math.}, 145(1):151--218, 2001.

\bibitem{Liu}
G.~Liu.
\newblock Associativity of quantum multiplication.
\newblock {\em Comm. Math. Phys.}, 191(2):265--282, 1998.

\bibitem{MacDonald}
I.~G. MacDonald.
\newblock Symmetric products of an algebraic curve.
\newblock {\em Topology}, 1:319--343, 1962.

\bibitem{McDuffSalamon}
D.~McDuff and D.~Salamon.
\newblock {\em {$J$}-holomorphic curves and quantum cohomology}.
\newblock Number~6 in University Lecture Series. American Mathematical Society,
  1994.

\bibitem{Mrowka}
T.~S. Mrowka.
\newblock A local {M}ayer-{V}ietoris principle for the {Y}ang-{M}ills moduli
  space.
\newblock Berkeley Ph.D. Thesis, 1989.

\bibitem{OhTransversality}
Y-G. Oh.
\newblock Fredholm theory of holomorphic discs under the perturbation of
  boundary conditions.
\newblock {\em Math. Z}, 222(3):505--520, 1996.

\bibitem{OhFloer}
Y-G. Oh.
\newblock On the structure of pseudo-holomorphic discs with totally real
  boundary conditions.
\newblock {\em J. Geom. Anal.}, 7(2):305--327, 1997.

\bibitem{HolDiskGraded}
P.~S. Ozsv{\'a}th and Z.~Szab{\'o}.
\newblock Absolutely graded {F}loer homologies and intersection forms for
  four-manifolds with boundary.
\newblock To appear in {\em Adv. Math.}, math/0110170, 2001.

\bibitem{HolDiskTwo}
P.~S. Ozsv{\'a}th and Z.~Szab{\'o}.
\newblock Holomorphic disks and three-manifold invariants: properties and
  applications.
\newblock math.SG/0105202, 2001.

\bibitem{HolDiskThree}
P.~S. Ozsv{\'a}th and Z.~Szab{\'o}.
\newblock Holomorphic triangles and invariants for smooth four-manifolds.
\newblock math/0110169, 2001.

\bibitem{ParkerWolfson}
T.~H. Parker and J.~G. Wolfson.
\newblock Pseudo-holomorphic maps and bubble trees.
\newblock {\em J. Geom. Anal.}, 3(1):63--98, 1993.

\bibitem{Reidemeister}
K.~Reidemeister.
\newblock Zur dreidimensionalen {T}opologie.
\newblock {\em Abh. Math. Sem. Univ. Hamburg}, (9):189--194, 1933.

\bibitem{SpecFlow}
J.~Robbin and D.~Salamon.
\newblock The spectral flow and the {M}aslov index.
\newblock {\em Bull. London Math. Soc.}, 27(1):1--33, 1995.

\bibitem{SalamonZehnder}
D.~Salamon and E.~Zehnder.
\newblock Morse theory for periodic solutions of {H}amiltonian systems and the
  {M}aslov index.
\newblock {\em Comm. Pure Appl. Math.}, 45(10):1303--1360, 1992.

\bibitem{Scharlemann}
M.~Scharlemann.
\newblock Heegaard splittings of compact {$3$}-manifolds.
\newblock math.GT/0007144.

\bibitem{Singer}
J.~Singer.
\newblock Three-dimensional manifolds and their {H}eegaard diagrams.
\newblock {\em Trans. Amer. Math. Soc.}, 35(1):88--111, 1933.

\bibitem{Taubes}
C.~H. Taubes.
\newblock {\em Metrics, connections and gluing theorems}.
\newblock Number~89 in CBMS Regional Conference Series in Mathematics. AMS,
  1996.

\bibitem{QAssoc}
G.~Tian.
\newblock Quantum cohomology and its associativity.
\newblock In R.~Bott, M.~Hopkins, A.~Jaffe, I.~Singer, D.~Stroock, and S-T.
  Yau, editors, {\em Current Developments in Mathematics, 1995}, pages
  361--401. Internat. Press, 1994.

\bibitem{Turaev}
V.~Turaev.
\newblock Torsion invariants of {S}pin{$^c$}-structures on $3$-manifolds.
\newblock {\em Math. Research Letters}, 4:679--695, 1997.

\bibitem{Viterbo}
C.~Viterbo.
\newblock Intersection de sous-vari{\'e}t{\'e}s lagrangiennes, fonctionnelles
  d'action et indice des syst{\`e}mes hamiltoniens.
\newblock {\em Bull. Soc. Math. France}, 115(3):361--390, 1987.

\bibitem{Ye}
R.~Ye.
\newblock Gromov's compactness theorem for pseudo holomorphic curves.
\newblock {\em Trans. Amer. Math. Soc.}, 342(2):671--694, 1994.

\bibitem{Zieschang}
H.~Zieschang.
\newblock On {H}eegaard diagrams of {$3$}-manifolds.
\newblock In {\em On the Geometry of Differentiable Manifolds}, volume 163-164,
  pages 247--280. Ast{\'e}risque, 1988.

\end{thebibliography}
}

\end{document}